\documentclass[11pt]{amsart}
\usepackage{amssymb, amsmath,latexsym,amsfonts,amsbsy, amsthm,mathtools,graphicx,color}
\usepackage{float}
\usepackage{hyperref}

\numberwithin{equation}{section}
\allowdisplaybreaks

\usepackage{mathrsfs,enumerate,extarrows,bm,tikz}
\usetikzlibrary{decorations.markings,arrows.meta}

\let\pa=\partial

\def\cE{{\mathcal E}}

\def\cM{{\mathcal M}}

\def\cP{{\mathcal P}}

\def\C{\mathbb C}
\def\R{\mathbb R}
\def\T{\mathbb T}
\def\Z{\mathbb Z}
\def\N{\mathbb N}

\def\e{\mathrm e}

\def\dive{\mathop{\rm div}\nolimits}

\newcommand{\beq}{\begin{equation}}
	\newcommand{\eeq}{\end{equation}}
\newcommand{\ben}{\begin{eqnarray}}
	\newcommand{\een}{\end{eqnarray}}
\newcommand{\beno}{\begin{eqnarray*}}
	\newcommand{\eeno}{\end{eqnarray*}}

\renewcommand{\Re}{\operatorname{Re}}

\newtheorem{theorem}{Theorem}[section]

\newtheorem{lemma}[theorem]{Lemma}
\newtheorem{proposition}[theorem]{Proposition}

\theoremstyle{remark}
\newtheorem{remark}[theorem]{Remark}
\newtheorem*{theorem*}{Theorem}

\newcommand{\eps}{\varepsilon}
\newcommand{\dd}{\,\mathrm d}
\newcommand{\Div}{\operatorname{div}}
\newcommand{\supp}{\operatorname{supp}}
\newcommand{\tr}{\operatorname{tr}}
\newcommand{\Id}{\mathrm I}

\newcommand{\norm}[1]{\left\|#1\right\|}

\newcommand{\HH}{\mathcal H}

\newcommand{\Pu}{\mathbb P_{\rm u}}

\begin{document}

\title[Blow-up for the cut-off Boltzmann equation] {On blow-up for the cut-off Boltzmann equation with very hard potentials}

\author[N. Jiang]{Ning Jiang}
\address[N. Jiang]{School of Mathematics and Statistics, Wuhan University, Wuhan 430070, China}
\email{njiang@whu.edu.cn}

\author[F. Shao]{Feng Shao}
\address[F. Shao]{Academy of Mathematics \& Systems Science and Hua Loo-Keng Center for Mathematical Sciences, Chinese Academy of Sciences, Beijing 100190, China}
\email{fshao@amss.ac.cn}

\author[D. Wei]{Dongyi Wei}
\address[D. Wei]{School of Mathematical Sciences, Peking University, Beijing 100871,  China}
\email{jnwdyi@pku.edu.cn}

\author[Z. Zhang]{Zhifei Zhang}
\address[Z. Zhang]{School of Mathematical Sciences, Peking University, Beijing 100871, China}
\email{zfzhang@math.pku.edu.cn}

\author[K. Zhou]{Kai Zhou}
\address[K. Zhou]{School of Mathematical Sciences, Nanjing Normal University, Nanjing 210023, China}
\email{kzyd@njnu.edu.cn}

\date{\today}

\begin{abstract}
	We construct smooth, strictly positive solutions to the three-dimensional cut-off Boltzmann equation with very hard potentials whose densities become unbounded in finite time. A key ingredient is a new velocity weight, which is Gaussian on a ball and has a polynomial tail of fixed order. We establish coercive estimates for the linearized collision operator and bilinear estimates for the collision term, with constants independent of the radius of the Gaussian region. Combining these weighted estimates with an imploding solution to the compressible Euler equations, we develop an energy framework that controls the Boltzmann solution up to the blow-up time. The initial data agree with a constant Maxwellian outside a bounded spatial region and have finite relative entropy.
\end{abstract}

\maketitle

\tableofcontents

\section{Introduction}\label{Sec:intro}
This paper concerns the inhomogeneous Boltzmann equation
\begin{equation}\label{Eq.Boltzmann}
	\partial_t f+v\cdot\nabla_x f=\frac1\varepsilon Q(f,f), \qquad (t,x,v)\in[0,T)\times\mathbb R^3\times\mathbb R^3,
\end{equation}
where $T>0$ is the blow-up time under consideration, $\varepsilon>0$ is the Knudsen number, and the nonnegative function $f=f(t,x,v)$ denotes the phase-space density of particles at position $x$ with velocity $v$. The collision term on the right-hand side of \eqref{Eq.Boltzmann} is defined by
\begin{equation}\label{Eq.collision-operator}
	Q(f,g)(v):=\int_{\R^3}\int_{\mathbb S^2}B_\gamma(v-v_*)\bigl(f_*'g'-f_*g\bigr)\dd\omega\dd v_*, \qquad B_\gamma(v-v_*):=\frac{|v-v_*|^\gamma}{4\pi}.
\end{equation}
Here $\mathbb S^2$ is the unit sphere in $\mathbb R^3$, and $\dd\omega$ is its standard surface measure. All functions in the collision integral are evaluated at the same $(t,x)$. For any phase-space function $h=h(t,x,v)$, we use the abbreviations
\[h:=h(t,x,v),\qquad h_*:=h(t,x,v_*),\qquad h':=h(t,x,v'),\qquad h_*':=h(t,x,v_*').\]
The pre-collisional velocities $(v,v_*)$ and the post-collisional velocities $(v',v_*')$ are related by
\begin{equation}\label{eq:postcollisional-velocities}
	v':=\frac{v+v_*}{2}+\frac{|v-v_*|}{2}\omega, \qquad v_*':=\frac{v+v_*}{2}-\frac{|v-v_*|}{2}\omega, \qquad \omega\in\mathbb S^2.
\end{equation}
These relations preserve momentum and kinetic energy, i.e.,
\[v+v_*=v'+v_*',\qquad |v|^2+|v_*|^2=|v'|^2+|v_*'|^2.\]
The angular measure in \eqref{Eq.collision-operator} has total mass one after multiplication by $1/(4\pi)$, and hence the collision kernel satisfies Grad's angular cut-off condition. The cases $\gamma>0$, $\gamma=0$, and $\gamma<0$ are called hard potentials, Maxwell molecules, and soft potentials, respectively. The range $\gamma>1$ is commonly referred to as very hard potentials. The hard-sphere kernel corresponds to $\gamma=1$. In this paper, we consider
\begin{equation}\label{Eq.gamma-range}
	\sqrt3<\gamma\leq2.
\end{equation}
For repulsive inverse-power forces $F(r)\sim r^{-s}$ with $s>1$, the corresponding exponent is $\gamma=(s-5)/(s-1)<1$. Thus, \eqref{Eq.gamma-range} lies outside the classical inverse-power range. Nevertheless, as discussed in \cite[Section~4]{Villani2025}, super hard spheres with $\gamma=2$ have been used as a phenomenological model in kinetic theory. Toscani and Villani \cite{ToscaniVillani1999} also studied the super-hard-sphere kernel and obtained an entropy-entropy production inequality for the spatially homogeneous Boltzmann equation. Consequently, the whole range $-3\leq \gamma\leq2$ is physically meaningful.

The global existence theory for the spatially inhomogeneous Boltzmann equation was initiated by Ukai \cite{Ukai1974,Ukai1976}, who constructed global solutions near Maxwellians for hard potentials with angular cutoff in the periodic box and the whole space. For soft potentials, Caflisch \cite{CaflischSoft1980} treated the case $-1<\gamma<0$, and Guo \cite{Guo2003} covered the full range $-3<\gamma<0$ by a nonlinear energy method. See also \cite{Guo2004,LYY2004} for the development of energy methods near Maxwellians. For the hard-sphere Boltzmann equation, Gualdani, Mischler and Mouhot \cite{GMM2017} established nonlinear stability for perturbations with polynomial decay in velocity. Recently, Chen, Shen and Zhang \cite{CSZ2026} proved the sharp global stability.

Another class of global solutions arises from perturbations of the vacuum state $f=0$. Illner and Shinbrot \cite{IllnerShinbrot1984} proved global existence for a rare gas of hard spheres in the whole space. For $-\frac12\leq\gamma\leq0$, Chen, Shen and Zhang \cite{CSZ2023} obtained global existence and scattering for sufficiently regular initial data, by only assuming the smallness of a critical norm.

For the Boltzmann equation without angular cutoff, global solutions near Maxwellians were constructed by Gressman and Strain \cite{GressmanStrain2011} and Alexandre, Morimoto, Ukai, Xu and Yang \cite{AMUXY2011,AMUXY2011II,AMUXY2012}. Duan, Liu, Sakamoto and Strain \cite{DLSS2021} established global existence of mild solutions with low regularity in the spatial variable. For the stability of vacuum without angular cutoff, we refer to Chaturvedi \cite{Chaturvedi2021}.

For general (large) initial data, in their groundbreaking work, DiPerna and Lions \cite{DiPernaLions1989} proved the global existence of renormalized solutions to the Boltzmann equation with angular cutoff. The initial data are only required to satisfy
\[\int_{\R^3\times\R^3}f_0(x,v)(1+|x|^2+|v|^2+|\log f_0(x,v)|)\dd x\dd v<\infty,\]
under suitable assumptions on the collision kernel. For the Boltzmann equation without angular cutoff, Alexandre and Villani \cite{AlexandreVillani2002} constructed global renormalized solutions with a defect measure. For the periodic hard-sphere Boltzmann equation, the recent OpenAI preprint \cite{OpenAI2026Boltzmann} claims the non-uniqueness of global entropy solutions, satisfying the local conservation laws of mass, momentum and energy.

For the Boltzmann equation with angular cutoff and $-3<\gamma\leq1$, Duan, Huang, Wang and Yang \cite{DHWY2017} constructed global mild solutions in $\Omega\in\{\R^3,\T^3\}$ for a class of initial data with large amplitude. Duan and Wang \cite{DW2019} also constructed global solutions with large amplitude initial data in bounded domains with isothermal diffuse reflection boundary conditions.

For the long-time behavior of solutions to the spatially periodic Boltzmann equation, Desvillettes and Villani \cite{DesvillettesVillani2005} proved that
\[\|f(t)-\mu\|_{L^1_{x,v}}\leq C_k(1+t)^{-k},\qquad \forall\,k>0,\]
by assuming uniform Sobolev bounds, moment bounds and a suitable lower bound on $f$. For a class of non-cutoff collision kernels, Imbert and Silvestre \cite{ImbertSilvestre2022} established $C^\infty$ a priori estimates under the following assumptions on the macroscopic quantities:
\[0<m_0\leq\int_{\R^3}f\dd v\leq M_0,\qquad \int_{\R^3}|v|^2f\dd v\leq E_0,\qquad \int_{\R^3}f\log f\dd v\leq H_0.\]
See also Fern\'andez-Real, Ros-Oton and Weidner \cite{FRW2025} for conditional regularity estimates in terms of the mass density, the pressure tensor and higher-order velocity moments.

The derivation of the equations of fluid mechanics from kinetic theory goes back to the work of Maxwell \cite{Maxwell1867} and Boltzmann \cite{Boltzmann1894}, and the rigorous justification of such derivations is part of Hilbert's sixth problem \cite{Hilbert1902}. For the compressible Euler limit of \eqref{Eq.Boltzmann}, given a smooth solution $(\rho,u,\mathcal T)$ to the compressible Euler equations on a time interval $[0,\tau]$, one seeks to construct solutions $f^\varepsilon$ to \eqref{Eq.Boltzmann} such that
\[
f^\varepsilon(t,x,v)\longrightarrow\frac{\rho(t,x)}{(2\pi\mathcal T(t,x))^{3/2}}\exp\left(-\frac{|v-u(t,x)|^2}{2\mathcal T(t,x)}\right)
\]
as $\varepsilon\to0$ on $[0,\tau]$. Nishida \cite{Nishida1978} justified the compressible Euler limit in the framework of analytic functions using the spectral properties of the linearized Boltzmann operator. Caflisch \cite{Caflisch1980} justified the limit for smooth Euler solutions using the Hilbert expansion. For the hard-sphere Boltzmann equation, Guo, Jang and Jiang \cite{GJJ2009} established an $L^2-L^\infty$ energy framework. The range of collision kernels was extended to $-3<\gamma\leq1$ with angular cutoff in \cite{GJJ2010}. More precisely, in \cite{GJJ2009,GJJ2010}, the Euler solution is assumed to satisfy
\[\sup_{[0,\tau]\times\Omega}\mathcal T<2\inf_{[0,\tau]\times\Omega}\mathcal T,
\qquad \Omega=\R^3\ \text{or}\ \mathbb T^3.\]
For cutoff kernels with $-3<\gamma\leq1$, Jiang, Luo and Tang \cite{JLT2025} established Grad-Caflisch estimates for the inverse of the linearized Boltzmann operator.

Can smooth large data for the spatially inhomogeneous Boltzmann equation develop a singularity? Villani conjectured global smooth solvability for smooth, rapidly decaying data and cutoff hard potentials in certain periodic or bounded convex spatial domains \cite[Conjecture 1]{Villani2008}. See also the regularity and well-posedness questions of Villani and Mouhot \cite{VillaniMouhot2015}. There are several positive results concerning the global regularity problem for large classical solutions. In one dimension, Biryuk, Craig and Panferov \cite{BCP2006} and Wynter \cite{Wynter2023} constructed global strong solutions for large initial data. For bounded convex domains with specular reflection boundary conditions, Duan, Ko and Lee \cite{DKL2023} constructed global solutions for initial data with arbitrarily large amplitude but sufficiently small relative entropy. Recently, He, Jiang, Kim, Lee and Park \cite{HJKLP2026} proved the uniqueness and exponential convergence to equilibrium of global mild solutions for cutoff hard potentials under a uniform-in-time bound on a velocity-weighted $L^1_vL^\infty_x$ norm.

For both cutoff and non-cutoff collision kernels, Golding, Henderson and Silvestre \cite{GHS} proved that the mass and energy densities remain bounded on $[0,T)$ if
\begin{equation}\label{Eq.GHS-integrable}
	\int_0^T\|\langle v\rangle^{3+\gamma}f(t)\|_{L^\infty_{x,v}}\dd t<\infty,
\end{equation}
and the initial data have sufficiently fast decay in $v$. As a consequence, for $-3<\gamma<\sqrt3$, Golding, Henderson and Silvestre \cite[Section 2.5]{GHS} ruled out the construction of blow-up solutions to the Boltzmann equation as perturbations of local Maxwellians associated with the known smooth self-similar Euler implosions.

Nevertheless, in the present paper, we show that the construction of blow-up solutions from smooth self-similar Euler implosions is possible for $\sqrt3<\gamma\leq2$. More precisely, by making use of the smooth imploding solutions to the compressible Euler equations constructed in \cite{SWWZ2025}, we construct smooth solutions to the cut-off Boltzmann equation \eqref{Eq.Boltzmann} whose mass densities become unbounded in finite time. To our knowledge, this is the first construction of finite-time blow-up solutions to the spatially inhomogeneous Boltzmann equation from smooth initial data with finite relative entropy.

A closely related collisional kinetic model is the Landau equation
\begin{align*}
	&\qquad\partial_t f+v\cdot\nabla_x f=\frac1\varepsilon Q_{\rm L}(f,f),\\
	Q_{\rm L}(f,g)(v)&:=\nabla_v\cdot\int_{\R^3}\mathsf a_\gamma(v-v_*)\bigl(f(v_*)\nabla_vg(v)-g(v)\nabla_{v_*}f(v_*)\bigr)\dd v_*,\notag\\
	&\mathsf a_\gamma(z):=\frac1{8\pi}|z|^{\gamma+2}\left(\Id-\frac{z\otimes z}{|z|^2}\right),\qquad z\ne0.
\end{align*}
The Landau operator describes the grazing-collision limit of the Boltzmann operator. The Coulomb case corresponds to $\gamma=-3$ \cite{Villani2002}. For the spatially homogeneous Landau equation with hard potentials, Desvillettes and Villani \cite{DesvillettesVillani2000} proved global existence, uniqueness, and immediate smoothing for a large class of initial data. Guo \cite{Guo2002Landau} constructed global classical solutions near Maxwellians in a periodic box. Guillen and Silvestre \cite{GuillenSilvestre2025} proved the monotonicity of the Fisher information and obtained global classical solutions to the homogeneous equation for $\gamma\in[-3,1]$, including Coulomb interactions, for initial data bounded above by a Maxwellian. For the Landau equation with Coulomb potential, the compressible Euler limit was justified by Lei, Liu, Xiao and Zhao \cite{LLXZ2022} via the Hilbert expansion and by Duan, Yang and Yu \cite{DYY2025} via the macro-micro decomposition. 

A finite-time singularity with a bounded distribution function was constructed for the inhomogeneous Landau equation in \cite{BCGJVY}. The construction uses the monatomic Euler profiles of \cite{SWWZ2025} and gives, for $\sqrt3<\gamma\leq2$,
\[\sup_{0\leq t<T}\|f(t)\|_{L^\infty_{x,v}}<\infty,\qquad\sup_{0\leq t<T}\left\|\int_{\R^3}f(t,\cdot,v)\dd v\right\|_{L^\infty_x}=\infty.\]
The solution is smooth and positive for $t<T$, its spatial H\"older seminorms of positive order diverge, and in similarity variables it tends to the Euler Maxwellian at $T$. To our knowledge, \cite{BCGJVY} is the first construction of a singularity of a collisional kinetic equation by Euler implosion.

The smooth self-similar Euler profiles used in both the present paper and \cite{BCGJVY} were constructed by Shao, Wang, Wei and Zhang \cite{SWWZ2025} for the adiabatic exponent $5/3$. The study of self-similar implosion for compressible Euler equations dates back to the classical work of Guderley \cite{Guderley1942}, who constructed self-similar solutions with a converging shock. Motivated by Guderley's work, Merle, Rapha\"el, Rodnianski and Szeftel \cite{MRRS2022I,MRRS2022II} made a breakthrough in the study of singularity formation for compressible fluids. For a class of adiabatic exponents, they proved the existence of a sequence of smooth, radially symmetric, self-similar imploding solutions to the compressible Euler equations. Unlike the Guderley solutions, MRRS's profiles are smooth across the sonic point. By using these self-similar profiles, they also constructed finite-time blow-up solutions to the 3-D compressible Navier-Stokes equations from smooth initial data with finite energy. Later, using computer-assisted proofs, Buckmaster, Cao-Labora and G\'omez-Serrano \cite{BCGS2025} constructed smooth self-similar Euler profiles for all adiabatic exponents greater than one and proved the blow-up of the compressible Navier-Stokes equations for diatomic gases (with adiabatic exponent $7/5$). Nevertheless, the Euler profiles in \cite{BCGS2025} can not be used to construct blow-up solutions for Navier-Stokes, Boltzmann or Landau. See \cite{CLGSSS} for the construction of non-radial imploding solutions to the compressible Euler and Navier-Stokes equations.

The physically important case of monatomic gases, corresponding to the adiabatic exponent $5/3$, is not covered by the analysis in \cite{MRRS2022I,MRRS2022II}. For this particular exponent, three singular points in the phase portrait of the ODE system for the self-similar profiles merge into a single point as the self-similarity exponent $r$ approaches $3-\sqrt3$. The triple-point degeneracy changes the asymptotic behavior of the Taylor coefficients of the local solutions around the sonic point, and hence the analysis in \cite{MRRS2022I} needs to be substantially modified. In this case, Shao, Wang, Wei and Zhang \cite{SWWZ2025} established the existence of smooth self-similar Euler profiles for a discrete sequence of exponents $\{r_n\}_{n\geq1}$ satisfying
\[1<r_n<3-\sqrt3,\qquad \lim_{n\to\infty}r_n=3-\sqrt3.\]
Using these self-similar profiles and the abstract theorem in \cite{CLGSSS}, Shao, Wang, Wei and Zhang \cite{SWWZ2025} also constructed finite-time blow-up solutions to the 3-D compressible Navier-Stokes equations for monatomic gases. 

\medskip

As in \cite{BCGJVY}, our starting point is a smooth self-similar imploding solution to the compressible Euler equations for monatomic gases, constructed in \cite{SWWZ2025}. By truncating the self-similar solution and solving a correction equation, we obtain a smooth Euler solution $(\rho,u,\mathcal T)$ on $[0,T)$ which agrees with a constant state outside a bounded set (see Proposition~\ref{Thm.Euler-solution}). Since the Euler flow is isentropic, we normalize the entropy so that $\rho=\theta^3$ and $\mathcal T=\theta^2$ for some positive function $\theta$. For some self-similarity exponent $r>1$ and smooth profiles $U:\R^3\to\R^3$ and $\vartheta:\R^3\to(0,\infty)$, we have
\[(u,\theta)(t,x)=(T-t)^{\frac1r-1}(U,\vartheta)\left(\frac{x}{(T-t)^{1/r}}\right)\]
in a fixed neighborhood of the origin.

For $\rho,\mathcal T>0$ and $u\in\R^3$, we define the Maxwellian
\begin{equation}\label{eq:physical-Maxwellian}
	\mathcal M(\rho,u,\mathcal T)(v):=\frac{\rho}{(2\pi\mathcal T)^{3/2}}\exp\left(-\frac{|v-u|^2}{2\mathcal T}\right).
\end{equation}
In particular, the standard Maxwellian is given by
\[\mu(v):=\mathcal M(1,0,1)(v)=(2\pi)^{-3/2}\mathrm e^{-|v|^2/2}.\]
Consider the local Maxwellian associated with the above Euler solution
\[\bar M(t,x,v):=\mathcal M(\theta^3,u,\theta^2)(v)=\mu\left(\frac{v-u(t,x)}{\theta(t,x)}\right).\]
For $0\leq t<T$, we have
\[\|\bar M(t)\|_{L^\infty_{x,v}}=(2\pi)^{-3/2},\qquad\int_{\R^3}\bar M(t,0,v)\dd v=\vartheta(0)^3(T-t)^{-\frac{3(r-1)}r}.\]
Hence $\bar M$ remains uniformly bounded although its mass density at the origin blows up as $t\uparrow T$. Our goal is to construct a solution $f$ to \eqref{Eq.Boltzmann} which stays close to $\bar M$ up to time $T$ and exhibits the same blow-up behavior.

Before stating our main result, we introduce some notation. For a distribution function $f$, let $\rho_f(t,x):=\int_{\R^3}f(t,x,v)\dd v$. Given a positive global Maxwellian $M_\infty$, we define the relative entropy of $f$ with respect to $M_\infty$ by
\begin{equation}\label{eq:entropy}
	\mathcal H(f\mid M_\infty):=\iint_{\R^3\times\R^3}\left(f\log\frac{f}{M_\infty}-f+M_\infty\right)\dd v\dd x,
\end{equation}
with the convention $0\log(0/M_\infty)=0$. For $0<\delta_x\leq1$, we define
\begin{equation}\label{eq:holder-definition}
	[f]_{C_x^{\delta_x}L_v^\infty}:=\sup_{\substack{x_1,x_2\in\R^3,\ x_1\ne x_2\\v\in\R^3}}\frac{|f(x_1,v)-f(x_2,v)|}{|x_1-x_2|^{\delta_x}}.
\end{equation}
We say that $f$ is simultaneous-rotation-invariant if $f(t,Ox,Ov)=f(t,x,v)$ for every $O\in SO(3)$.

\begin{theorem}\label{thm:main}
	Let $\gamma$ satisfy \eqref{Eq.gamma-range}. There exist a similarity exponent $r$ from Theorem~\ref{Thm.profiles} and constants $T,\eps_0,\theta_\infty>0$ satisfying
	\[\frac{3+\gamma}{2+\gamma}<r<3-\sqrt3\]
	with the following property. We define
	\begin{equation}\label{eq:exponents}
		\alpha_{\rm ss}:=\frac{r-1}{r},\qquad p:=\frac{r}{r-1},\qquad \delta_{\rm gap}:=3+\gamma-p>0,\qquad \eta:=\alpha_{\rm ss}\delta_{\rm gap}>0.
	\end{equation}
	For each $0<\eps\le\eps_0$, there exist smooth positive simultaneous-rotation-invariant initial data $f_0$ and a positive solution
	\[f\in C^\infty([0,T)\times\R^3_x\times\R^3_v)\]
	to \eqref{Eq.Boltzmann} with $f(0)=f_0$ such that:
	\begin{enumerate}
		\item $f_0=M_\infty:=\mathcal M(\theta_\infty^3,0,\theta_\infty^2)$ outside a fixed spatial ball, and
		\[\mathcal H(f(t)\mid M_\infty) \le \mathcal H(f_0\mid M_\infty)<\infty,\qquad 0\le t<T.\]
		For every $t\in[0,T)$, the relative mass and kinetic-energy moments are absolutely finite, i.e.
		\[\iint_{\R^3_x\times\R^3_v}(1+|v|^2)|f(t,x,v)-M_\infty(v)|\dd v\dd x<\infty.\]
		\item The distribution function remains uniformly bounded, i.e.
		\begin{equation}\label{eq:main-uniform-bound}
			\sup_{0\leq t<T}\|f(t)\|_{L^\infty_{x,v}}\leq C,
		\end{equation}
		where $C$ is independent of $\eps\in(0,\eps_0]$.
		\item There exist smooth radial profiles $U:\R^3\to\R^3$ and $\vartheta:\R^3\to(0,\infty)$, specified in \eqref{eq:normalized-profile}, such that for every $L<\infty$,
		\begin{equation}\label{eq:main-rescaled-limit}
			\sup_{|y|\le L} \left\| f\left(t,(T-t)^{1/r}y, (T-t)^{-\alpha_{\rm ss}}(U(y)+\vartheta(y)\widehat v)\right) -\mu(\widehat v)\right\|_{L_{\widehat v}^2} \longrightarrow0 \quad\text{as }t\uparrow T,
		\end{equation}
		Moreover,
		\begin{equation}\label{eq:main-density}
			\rho_f(t,0) =\vartheta(0)^3(T-t)^{-3\alpha_{\rm ss}}(1+o(1)) \longrightarrow\infty.
		\end{equation}
		\item Given $\delta_x\in(0,1]$, one can take $c_{\delta_x}>0$ and $t_{\delta_x}<T$ such that
		\begin{equation}\label{eq:main-holder}
			[f(t)]_{C_x^{\delta_x}L_v^\infty} \ge c_{\delta_x}(T-t)^{-\delta_x/r}, \qquad t_{\delta_x}<t<T.
		\end{equation}
	\end{enumerate}
\end{theorem}

Several remarks are in order.

\begin{remark}[The range of $\gamma$]\label{rem:intro-euler-mechanism}
	The restriction $\gamma>\sqrt3$ comes from the comparison between the transport and collision terms near the Euler implosion. Let $(u,\theta)$ be the Euler solution given by Proposition~\ref{Thm.Euler-solution}. Under the change of velocity variables
	\begin{equation}\label{Eq.hat-v-def}
		\widehat v=\frac{v-u(t,x)}{\theta(t,x)},\qquad F(t,x,\widehat v)=f(t,x,u(t,x)+\theta(t,x)\widehat v),
	\end{equation}
	the Boltzmann equation \eqref{Eq.Boltzmann} becomes
	\[\partial_tF+(u+\theta\widehat v)\cdot\nabla_xF+\mathsf H_{\rm tr}\cdot\nabla_{\widehat v}F=\frac{\theta^{3+\gamma}}{\eps}Q(F,F),\]
	where
	\[\mathsf H_{\rm tr}=\left(\frac13\Div u-\widehat v\cdot\nabla_x\theta\right)\widehat v+5\nabla_x\theta-\widehat v\cdot\nabla_xu.\]
	By the estimates for the Euler solution in Proposition~\ref{Thm.Euler-solution}, we have $|\nabla_xu|+|\nabla_x\theta|\lesssim \theta^p$, where $p=\frac{r}{r-1}$. At the origin, we have
	\[\theta(t,0)^p=\frac{\vartheta(0)^p}{T-t},\qquad\frac1\eps\theta(t,0)^{3+\gamma}=\frac{\vartheta(0)^{3+\gamma}}{\eps(T-t)^{1+\eta}},\]
	where $\eta=\frac{r-1}{r}(3+\gamma)-1$. For $\eta>0$, there holds
	\[\left.\frac{\theta^p}{\eps^{-1}\theta^{3+\gamma}}\right|_{x=0}
	=\eps\vartheta(0)^{p-3-\gamma}(T-t)^\eta\to0\qquad\text{as }t\uparrow T.\]
	The condition $\eta>0$ is equivalent to
	\[\frac{r}{r-1}<3+\gamma\quad\Longleftrightarrow\quad r>\frac{3+\gamma}{2+\gamma}.\]
	Shao, Wang, Wei and Zhang \cite{SWWZ2025} constructed smooth self-similar profiles for a sequence of exponents $\{r_n\}_{n\geq1}$ satisfying $1<r_n<3-\sqrt3$ and $r_n\to3-\sqrt3$. Hence we can choose a self-similarity exponent $r=r_n$ satisfying $\frac{r}{r-1}<3+\gamma$ if and only if $\gamma>\sqrt3$.
	
	The restriction $\gamma>\sqrt3$ is related to the fact that the Euler profiles used in the present paper are isentropic. Recently, Chen, Shkoller and Vicol \cite{CSV2026} constructed smooth non-isentropic self-similar imploding solutions to the compressible Euler equations. As discussed in \cite[Section~1.6.2]{CSV2026}, for any $\gamma\in(-3,2]$, one can choose a non-isentropic profile such that the collision term is formally dominant over the transport terms in self-similar variables. In particular, the profiles in \cite{CSV2026} provide candidates for the construction of blow-up solutions to the Boltzmann equation for Maxwell molecules ($\gamma=0$) and hard spheres ($\gamma=1$). However, for the Euler solutions constructed in \cite{CSV2026}, the temperature vanishes at the origin for all $t<T$. For each fixed $t<T$, the associated local Maxwellian satisfies
	\[\mathcal M(\rho(t,x),u(t,x),\mathcal T(t,x))\rightharpoonup\rho(t,0)\delta_{u(t,0)}\qquad\text{as }x\to0\]
	in the sense of distributions. In particular, the local Maxwellian is already singular at the origin before the blow-up time. Our analysis relies on the strict positivity of the temperature of the background Euler solution. Thus, the construction of Boltzmann singularities from the non-isentropic profiles in \cite{CSV2026} requires new ideas. See also \cite[Section~2.5]{GHS}.
\end{remark}

\begin{remark}[The kinetic estimates]\label{rem:intro-kinetic-estimates}
	A main difficulty in the kinetic estimates is the loss of velocity weights caused by the transport operator. In the velocity variable $\widehat v$ defined by \eqref{Eq.hat-v-def}, consider the usual decomposition $F=\mu+\sqrt\mu g$. By \eqref{Eq.H-expression}, we have
	\begin{align*}
		\mathsf H_{\rm tr}\cdot\nabla_{\widehat v}\log\sqrt\mu=\frac12\left(\widehat v\otimes\widehat v-\frac13|\widehat v|^2\Id\right):\nabla_xu+\frac12(|\widehat v|^2-5)\widehat v\cdot\nabla_x\theta.
	\end{align*}
	In the $L^2$ energy estimate for $g$, we thus need to control $\theta^p\|\langle\widehat v\rangle^{3/2}g\|_{L^2_{\widehat v}}^2$,	while the dissipation from the linearized Boltzmann operator is of order $\frac{\theta^{3+\gamma}}{\eps}\|\langle\widehat v\rangle^{\gamma/2}g\|_{L^2_{\widehat v}}^2$. Since $\gamma<3$, there is a loss of velocity weights in the energy estimate. Although $p<3+\gamma$, we have
	\begin{equation}\label{eq:intro-Gaussian-transport-ratio}
		\frac{\theta^p\langle\widehat v\rangle^3}{\eps^{-1}\theta^{3+\gamma}\langle\widehat v\rangle^\gamma}
		=\eps\theta^{p-3-\gamma}\langle\widehat v\rangle^{3-\gamma}
		\to\infty\qquad\text{as }|\widehat v|\to\infty.
	\end{equation}
	
	To overcome the loss of velocity weights, one may consider a polynomially decaying weight instead of $\sqrt\mu$. However, the Gaussian weight is also used in the coupling of the macroscopic and microscopic energy estimates. We introduce a family of weights $w_{m,\ell}$ (see \eqref{eq:Phi-fixed-tail} and \eqref{eq:fixed-tail-weight}) satisfying
	\[w_{m,\ell}(v)=\begin{cases}
		\sqrt{\mu(v)},& |v|^2\leq m,\\
		c_{m,\ell}|v|^{-\ell/2},& |v|^2\geq2m,
	\end{cases}\]
	where $\ell$ is fixed and sufficiently large, $m\geq\ell$ and $c_{m,\ell}>0$. We prove
	\[|\mathsf H_{\rm tr}\cdot\nabla_{\widehat v}\log w_{m,\ell}|\leq\begin{cases}
		C\theta^p\langle\widehat v\rangle^3,&|\widehat v|^2\leq2m,\\
		C_\ell\theta^p\langle\widehat v\rangle,&|\widehat v|^2\geq2m.
	\end{cases}\]
	Since $\gamma>1$, it follows that
	\begin{equation}\label{eq:intro-tail-transport-ratio}
		\frac{|\mathsf H_{\rm tr}\cdot\nabla_{\widehat v}\log w_{m,\ell}|}
		{\eps^{-1}\theta^{3+\gamma}\langle\widehat v\rangle^\gamma}
		\leq C_\ell\eps\theta^{-\delta_{\rm gap}}(1+m^{(3-\gamma)/2}).
	\end{equation}
	
	One of the main ingredients of our proof is to establish collision estimates with constants independent of $m$. For $w=w_{m,\ell}$, let
	\[
	L_wg:=-\frac{Q(\mu,wg)+Q(wg,\mu)}{w},\qquad
	Q_w(g_1,g_2):=\frac{Q(wg_1,wg_2)}{w}.
	\]
	In Proposition~\ref{thm:new-collision}, we prove
	\begin{equation}\label{Eq.Lw-lower}
		\operatorname{Re}\langle L_wg,g\rangle\geq c_0\|\nu^{1/2}g\|_{L^2_v}^2
	\end{equation}
	for $g$ satisfying $\int_{\R^3}(1,v,|v|^2)wg\dd v=0$, and that
	\begin{equation}\label{Eq.bilinear-est}
		\|\nu^{-1/2}Q_w(g_1,g_2)\|_{L^2_v}\leq C_\ell\left(\|\nu^{1/2}g_1\|_{L^2_v}\|g_2\|_{L^2_v}+\|g_1\|_{L^2_v}\|\nu^{1/2}g_2\|_{L^2_v}\right),
	\end{equation}
	where $\nu(v)=\int_{\R^3}|v-v_*|^\gamma\mu(v_*)\dd v_*\asymp\langle v\rangle^\gamma$. The constants $c_0$ and $C_\ell$ are independent of $m$ for all sufficiently large $m$. Moreover, by \eqref{eq:new-tail-small}, we have
	\[
	\|w_{m,\ell}-\sqrt\mu\|_{L^2_v}
	+\|\mu/w_{m,\ell}-w_{m,\ell}\|_{L^2_v}
	\leq C_\ell\mathrm e^{-cm}.
	\]
	
	To reduce the errors arising from the replacement of $\sqrt\mu$ by $w_{m,\ell}$, we need to take $m$ large. On the other hand, by \eqref{eq:intro-tail-transport-ratio}, we need to keep $\eps\theta^{-\delta_{\rm gap}}m^{(3-\gamma)/2}$ small to absorb the transport term by the collision dissipation. Since the constants in \eqref{Eq.Lw-lower} and \eqref{Eq.bilinear-est} are independent of $m$, we choose $m$ as a function of $\theta(t,x)$ and $\eps$, i.e.
	\[m_\eps(t,x):=m_0+\kappa_{\rm w}\log\left(1+\frac{\theta(t,x)^{\delta_{\rm gap}}}{\eps}\right),\]
	where $m_0$ and $\kappa_{\rm w}$ are sufficiently large constants. By Proposition~\ref{thm:adaptive}, for every fixed $L\geq0$,
	\[
	\lim_{\eps\downarrow0}\sup_{0\leq t<T,\,x\in\R^3}
	\eps\theta^{-\delta_{\rm gap}}m_\eps^L=0.
	\]
	Hence the transport term in \eqref{eq:intro-tail-transport-ratio} can be absorbed by the collision dissipation for sufficiently small $\eps>0$. 
\end{remark}

\begin{remark}\label{rem:intro-conditional-regularity}
	Our blow-up result is consistent with \cite{GHS}. At $x=0$, we have $u(t,0)=0$. Take $v=\theta(t,0)e_1$ with $e_1=(1,0,0)$ in $\bar M=\mathcal M(\theta^3,u,\theta^2)$, then
	\[\bar M(t,0,\theta(t,0)e_1)=\mu(e_1),\qquad \|\langle v\rangle^{3+\gamma}\bar M(t)\|_{L^\infty_{x,v}}\geq\mu(e_1)\theta(t,0)^{3+\gamma}\gtrsim(T-t)^{-1-\eta}.\]
	Since $\eta>0$ in \eqref{eq:exponents}, we have
	\[\int_0^T\|\langle v\rangle^{3+\gamma}\bar M(t)\|_{L^\infty_{x,v}}\dd t=\infty,\]
	so $\bar M$ fails the time-integrability condition \eqref{Eq.GHS-integrable}.
\end{remark}

The use of self-similar solutions in the construction of blow-up solutions is not restricted to compressible fluids or kinetic equations. In \cite{MRRS2022NLS}, the self-similar Euler profiles were used to construct blow-up solutions to the defocusing nonlinear Schr\"odinger equation. More recently, the smooth self-similar imploding profiles of the relativistic Euler equations constructed in \cite{SWZ2026RelEuler} were used to construct blow-up solutions to the complex-valued and real-valued defocusing nonlinear wave equations in \cite{SWZ2025ComplexNLW,SWZ2026RealNLW}, respectively. See also \cite{BC2026NLW}.

\medskip

\noindent\textbf{Ideas of proof.} Let $(u,\theta)$ be the imploding Euler solution constructed in Proposition~\ref{Thm.Euler-solution} and let $\bar M=\mathcal M(\theta^3,u,\theta^2)$. Since $Q(\bar M,\bar M)=0$, we have
\[(\partial_t+v\cdot\nabla_x)(f-\bar M)=\eps^{-1}Q(f,f)-(\partial_t+v\cdot\nabla_x)\bar M.\]
To construct a solution $f$ close to $\bar M$ on $[0,T)$, we combine the stability estimates for the linearized Euler equations with the dissipation estimates for the Boltzmann collision operator.

\begin{enumerate}[(i)]
	\item \emph{Construction of the Euler background.} The densities of the self-similar Euler solutions constructed in \cite{SWWZ2025} tend to zero as $|x|\to\infty$. In our construction, we need an Euler solution with a strictly positive lower bound on the temperature and which agrees with a constant state outside a bounded set. Let
	\[(u_{\rm p},\theta_{\rm p})(t,x)=(T-t)^{-\alpha_{\rm ss}}(U,\vartheta)\left(\frac{x}{(T-t)^{1/r}}\right),\]
	and take a smooth cutoff function $\chi$ such that $\chi=1$ on $[0,1/2]$ and $\chi=0$ on $[1,\infty)$. For some large constant $R_0>0$, we seek an Euler solution of the form
	\begin{align*}
		u(t,x)&=u_{\rm p}(t,x)\chi\left(\frac{|x|}{R_0}\right)+u_{\rm r}(t,x),\\
		\theta(t,x)&=\theta_{\rm p}(t,x)\chi\left(\frac{|x|}{R_0}\right)+R_0^{1-r}\left(1-\chi\left(\frac{|x|}{R_0}\right)\right)+\theta_{\rm r}(t,x).
	\end{align*}
	The errors arising from the truncation are supported in $\{R_0/2\leq|x|\leq R_0\}$. We solve the equations for $(u_{\rm r},\theta_{\rm r})$ backward in time with zero terminal data at $t=T$. By the finite speed of propagation, for sufficiently small $T>0$, we have
	\[\supp u_{\rm r}(t)\cup\supp\theta_{\rm r}(t)\subset\{R_0/4<|x|<2R_0\}.\]
	Consequently, the Euler solution $(u,\theta)$ agrees with $(u_{\rm p},\theta_{\rm p})$ on $|x|\leq R_0/4$ and agrees with $(0,R_0^{1-r})$ on $|x|\geq2R_0$.
	Moreover, $\theta$ has a strictly positive lower bound on $[0,T)\times\R^3$. See Proposition~\ref{Thm.Euler-solution} for details.
	
	\item \emph{Macroscopic and microscopic perturbations.} Let $\rho_f$, $u_f$ and $T_f$ be the mass density, mean velocity and temperature of $f$, respectively, i.e.
	\[\rho_f=\int f\dd v,\qquad u_f=\frac1{\rho_f}\int vf\dd v,\qquad T_f=\frac1{3\rho_f}\int|v-u_f|^2f\dd v.\]
	We introduce the normalized macroscopic variables
	\[n=\frac{\rho_f}{\theta^3},\qquad V=\frac{u_f-u}{\theta},\qquad \Theta=\frac{T_f}{\theta^2}.\]
	Instead of estimating $(n-1,V,\Theta-1)$ directly, we use the pressure and entropy variables
	\[\theta_f=(\rho_fT_f)^{1/5},\qquad B=\frac{\theta_f}{\theta}-1,\qquad S=\log\frac{T_f}{\rho_f^{2/3}}.\]
	Denote $\mathbf z=(V,B,S)$. Since $n=(1+B)^3\mathrm e^{-3S/5}$ and $\Theta=(1+B)^2\mathrm e^{3S/5}$, estimating $\mathbf z$ is equivalent to estimating $(n-1,V,\Theta-1)$ near the Euler background.
	
	Let $M=\mathcal M(\rho_f,u_f,T_f)$ and $H=f-M$. Then
	\begin{equation}\label{eq:intro-macro-micro}
		f-\bar M=(M-\bar M)+H,\qquad \int(1,v,|v|^2)H\dd v=0.
	\end{equation}
	Let $D_f=\partial_t+u_f\cdot\nabla_x$ and $\mathcal C_f=\Pi:\nabla_xu_f+\Div q$, then the macroscopic system reads
	\begin{equation}\label{eq:intro-pressure-system}
		\begin{cases}
			D_fu_f+5\theta_f\mathrm e^{3S/5}\nabla_x\theta_f=-\rho_f^{-1}\Div\Pi,\\
			D_f\theta_f+\frac13\theta_f\Div u_f=-\frac2{15}\theta_f^{-4}\mathcal C_f,\\
			D_fS=-\frac23\theta_f^{-5}\mathcal C_f,
		\end{cases}
	\end{equation}
	where the stress tensor $\Pi$ and the heat flux $q$ are
	\[\Pi=\int(v-u_f)\otimes(v-u_f)H\dd v,\qquad q=\frac12\int|v-u_f|^2(v-u_f)H\dd v.\]
	Thanks to the pressure variable $\theta_f$, we note that no derivatives of $S$ occur in the acoustic part of \eqref{eq:intro-pressure-system}. We estimate $\mathbf z$ using the Euler linearization in \cite{BCGJVY}.
	
	By \eqref{eq:intro-macro-micro} and $Q(M,M)=0$, we have $Q(f,f)=Q(M,H)+Q(H,M)+Q(H,H)$. For the estimates of $H$, we apply the velocity weight constructed in Remark~\ref{rem:intro-kinetic-estimates} to the Maxwellian $M$ rather than $\bar M$. More precisely, we use
	\[\widetilde v=\frac{v-u_f}{\sqrt{T_f}},\qquad W=n^{1/2}\Theta^{-3/4}w_{m_\eps,\ell}(\widetilde v),\]
	where
	\[m_\eps(t,x)=m_0+\kappa_{\rm w}\log\left(1+\frac{\theta(t,x)^{\delta_{\rm gap}}}{\eps}\right).\]
	Since $M=n\Theta^{-3/2}\mu(\widetilde v)$, we have
	\begin{equation}\label{eq:intro-W-M}
		W^2=M\qquad\text{for }|\widetilde v|^2\leq m_\eps.
	\end{equation}
	The collision estimates \eqref{Eq.Lw-lower} and \eqref{Eq.bilinear-est} can be transferred to the linearization around $M$ by the change of velocity variables $v=u_f+\sqrt{T_f}\widetilde v$. See Section~\ref{sec:moments} for the equations of the macroscopic and microscopic perturbations.
	
	\item \emph{Construction of the coupled energy.} Let $K=2k$ be a sufficiently large even integer. In the $K$-th order energy estimate for \eqref{eq:intro-pressure-system}, the terms $\Delta_x^k\Div\Pi$ and $\Delta_x^k\Div q$ contain $K+1$ spatial derivatives of $H$. On the other hand, the equation for $H$ contains the source term $S_M=(I-\mathcal P_M)(v\cdot\nabla_xM)$, where $\mathcal P_M$ is the projection onto the tangent space of Maxwellians which preserves the mass, momentum and energy moments. Thus, in the $K$-th order energy estimate for $H$, the term $\Delta_x^kS_M$ contains $K+1$ spatial derivatives of $\mathbf z$. The separate estimates of $\mathbf z$ and $H$ would lead to the loss of spatial derivatives.
	
	To avoid the loss of derivatives, we estimate the sum of $\Delta_x^kH$ and the highest-order term in $\Delta_x^k(M-\bar M)$. By the chain rule, we have
	\[\Delta_x^k(f-\bar M)=\Delta_x^kH+D_{\mathbf z}M[\Delta_x^k\mathbf z]+\mathcal R_K,\]
	where $D_{\mathbf z}M$ is the differential of $M$ with respect to $\mathbf z$ with $(u,\theta)$ fixed, and $\mathcal R_K$ contains at most $K-1$ spatial derivatives of $\mathbf z$. All spatial derivatives are taken with the physical velocity $v$ fixed. Writing $D=\partial_t+v\cdot\nabla_x$, we have
	\[[D,\Delta_x^k]=0,\qquad D\Delta_x^k(f-\bar M)=\eps^{-1}\Delta_x^kQ(f,f)-\Delta_x^kD\bar M.\]
	By substituting the moment equations into $D\mathcal R_K$, we only need to estimate at most $K$ spatial derivatives of $\mathbf z$ and $H$. Hence the terms containing $K+1$ spatial derivatives of the perturbations cancel when we combine the macroscopic and microscopic equations.
	
	We next show that the energy of the combined perturbation controls both the macroscopic and microscopic perturbations. By the zero moments of $\Delta_x^kH$ and the fact that $D_{\mathbf z}M[\Delta_x^k\mathbf z]/M$ is a linear combination of $1,v_1,v_2,v_3,|v|^2$, we have
	\[\int\frac{\Delta_x^kH D_{\mathbf z}M[\Delta_x^k\mathbf z]}{M}\dd\widehat v=0.\]
	Since $W^2=M$ for $|\widetilde v|^2\leq m_\eps$, by the Gaussian decay of $D_{\mathbf z}M[\Delta_x^k\mathbf z]$, we obtain
	\begin{equation}\label{eq:intro-combined-energy-equivalence}
		\int\frac{|\Delta_x^kH+D_{\mathbf z}M[\Delta_x^k\mathbf z]|^2}{W^2}\dd\widehat v\asymp\int\frac{|\Delta_x^kH|^2}{W^2}\dd\widehat v+\int\frac{|D_{\mathbf z}M[\Delta_x^k\mathbf z]|^2}{M}\dd\widehat v.
	\end{equation}
	Moreover, by introducing suitable spatial weights in the combined energy, we can match the macroscopic part of the energy with the energy used in the stability estimates for the Euler equations. See Proposition~\ref{prop:global-matching} and Section~\ref{sec:energies} for details.
	
	\item \emph{The coupled energy estimates.} In addition to keeping the perturbations small in the whole space $\R^3_x$, we need to prove the convergence of the perturbations to zero in self-similar variables. Based on the analysis in (iii), we construct coupled energies with two different spatial weights $\Omega_0^*$ and $\Omega_a^*$, respectively, satisfying
	\[\Omega_a^*=\langle y\rangle^{-a}\Omega_0^*,\qquad y=\frac{x}{(T-t)^{1/r}},\qquad a=\frac r2.\]
	For each fixed $x\ne0$, we have $\Omega_a^*/\Omega_0^*\to0$ as $t\uparrow T$. 
	Let $\mathcal E$ be the coupled energy with the spatial weight $\Omega_a^*$ and a local quadratic correction arising from the stability estimates for the Euler equations. Let $\mathcal B_{\rm reg}$ be the coupled energy with the spatial weight $\Omega_0^*$, augmented by the energy of the mixed derivatives $\partial_x^\alpha\partial_v^\beta H$ with $|\beta|\geq1$ and $|\alpha|+|\beta|\leq K$. The mixed derivatives are needed to obtain the pointwise bounds on $H$ required for the continuation of the Boltzmann solution.
	
	The linearized Euler operator in self-similar variables has a finite-dimensional subspace $E_{\rm u}$ on which the required decay estimate does not hold. Let $\mathbf a(t)\in\R^J$, with $J=\dim E_{\rm u}$, be the coordinates of the localized macroscopic perturbation associated with $E_{\rm u}$. We write $|\mathbf a|_P^2=\mathbf a^{\mathsf T}P\mathbf a$ for some symmetric positive definite matrix $P$.
	
	Thanks to the cancellation of the highest-order derivatives in (iii), combining the Euler stability estimates with the collision estimates, for sufficiently small perturbations and $\eps>0$, we obtain
	\begin{align}
		\mathcal E'+\frac{\lambda}{T-t}\mathcal E+\frac c\eps\mathcal D_a&\leq\frac C{T-t}|\mathbf a|_P^2+C\eps(T-t)^{-1+\eta},\label{eq:intro-inner-energy}\\
		\mathcal B_{\rm reg}'+\frac{c_{\rm reg}}\eps\mathcal D_{\rm reg}&\leq\frac C{T-t}\mathcal E+C\mathcal B_{\rm reg}+C\eps(T-t)^{-1+\eta},\label{eq:intro-global-energy}
	\end{align}
	where $\lambda>0$ and $\mathcal D_a,\mathcal D_{\rm reg}$ are the corresponding microscopic dissipations. See Propositions~\ref{thm:new-closed-pair} and \ref{prop:mixed-velocity-estimate} for details.
	
	Since $\eta>0$, the term $\eps(T-t)^{-1+\eta}$ in \eqref{eq:intro-inner-energy} and \eqref{eq:intro-global-energy} is integrable on $[0,T)$. To close the energy estimates, it remains to choose the initial data such that
	\[\mathcal E(t)+|\mathbf a(t)|_P^2\leq C\eps T^{\eta-d_*}(T-t)^{d_*}\]
	for some $0<d_*<\min\{\lambda,\eta\}$. By \eqref{eq:intro-global-energy} and Gr\"onwall's inequality, the decay of $\mathcal E$ then gives a uniform bound for $\mathcal B_{\rm reg}$ on $[0,T)$.
	
	\item \emph{Selection of the initial data.} The equation for the unstable coordinates $\mathbf a$ takes the form
	\begin{equation}\label{eq:intro-unstable-equation}
		\mathbf a'=\frac1{T-t}A_{\rm u}\mathbf a+F_{\rm u},
	\end{equation}
	where $A_{\rm u}$ is the matrix of the restriction of the linearized Euler operator to $E_{\rm u}$. The source term $F_{\rm u}$ satisfies
	\begin{equation}\label{eq:intro-unstable-source}
		|F_{\rm u}|_P^2\leq C(T-t)^{-1+\eta}\mathcal D_a+C(T-t)^{-2}\mathcal E^2.
	\end{equation}
	Unlike the stability problem for the Euler equations, the equation \eqref{eq:intro-unstable-equation} contains a source term arising from the microscopic perturbation $H$. In particular, we need to use the collision dissipation $\mathcal D_a$ to control $F_{\rm u}$.
	
	For $\mathbf a_0\in\R^J$ sufficiently small, we construct a family of smooth positive initial data
	\[f_0(x,v;\mathbf a_0)=\mathcal M(\rho_{f,0}(x;\mathbf a_0),u_{f,0}(x;\mathbf a_0),T_{f,0}(x;\mathbf a_0))(v)\]
	satisfying
	\begin{equation}\label{eq:intro-initial-energy}
		H(0)=0,\qquad \mathbf a(0)=\mathbf a_0,\qquad \mathcal E(0)+\mathcal B_{\rm reg}(0)\leq C_I|\mathbf a_0|_P^2.
	\end{equation}
	The initial data agree with $\bar M(0)$ outside a bounded spatial region. See Proposition~\ref{prop:positive-realization} and Lemma~\ref{lem:selection-initial-energy} for details.
	
	Our goal is to choose $\mathbf a_0$ such that $|\mathbf a(t)|_P^2$ decays at the rate $(T-t)^{d_*}$ for some $0<d_*<\min\{\lambda,\eta\}$. By the spectral properties of $A_{\rm u}$, \eqref{eq:intro-unstable-equation} and \eqref{eq:intro-unstable-source}, we can choose $d_*$ such that
	\begin{equation}\label{eq:intro-unstable-growth}
		\frac\dd{\dd t}\big((T-t)^{-d_*}|\mathbf a|_P^2\big)\geq\frac{g_*}{2(T-t)^{1+d_*}}|\mathbf a|_P^2-C_g(T-t)^{\eta-d_*}\mathcal D_a-C_g(T-t)^{-1-d_*}\mathcal E^2
	\end{equation}
	for some $g_*>0$. To control the last two terms in \eqref{eq:intro-unstable-growth}, we use the damping term $\lambda\mathcal E/(T-t)$ and the collision dissipation $c\mathcal D_a/\eps$ in \eqref{eq:intro-inner-energy}. More precisely, for a sufficiently small constant $k_*>0$, we introduce
	\begin{equation}\label{eq:intro-shooting-functional}
		\mathcal Q_{\rm sh}(t):=(T-t)^{-d_*}(|\mathbf a(t)|_P^2-k_*\mathcal E(t)).
	\end{equation}
	By \eqref{eq:intro-unstable-growth} and \eqref{eq:intro-inner-energy}, for sufficiently small $\eps$ and $\mathcal E$, we obtain
	\begin{equation}\label{eq:intro-shooting-derivative}
		\mathcal Q_{\rm sh}'\geq\frac{g_*}{4(T-t)^{1+d_*}}|\mathbf a|_P^2+\frac{k_*(\lambda-d_*)}{2(T-t)^{1+d_*}}\mathcal E+\frac{k_*c}{2\eps}(T-t)^{-d_*}\mathcal D_a-Ck_*\eps(T-t)^{-1+\eta-d_*}.
	\end{equation}
	Thus, by subtracting $k_*\mathcal E$ in \eqref{eq:intro-shooting-functional}, we can use the collision dissipation to control the contribution of the microscopic perturbation to the equation for $\mathbf a$.
	
	For a sufficiently large constant $R_{\rm sh}>0$, let $L_{\rm sh}=R_{\rm sh}\eps T^{\eta-d_*}$. By \eqref{eq:intro-shooting-derivative}, we have
	\begin{equation}\label{eq:intro-outward-condition}
		\mathcal Q_{\rm sh}(t)=L_{\rm sh}\quad\Longrightarrow\quad\mathcal Q_{\rm sh}'(t)>0.
	\end{equation}
	On the other hand, as long as $\mathcal Q_{\rm sh}\leq L_{\rm sh}$, we have $|\mathbf a|_P^2\leq L_{\rm sh}(T-t)^{d_*}+k_*\mathcal E$. Substituting the bound for $|\mathbf a|_P^2$ into \eqref{eq:intro-inner-energy} and applying Gr\"onwall's inequality, for $|\mathbf a_0|_P^2\leq C\eps T^\eta$, we obtain
	\begin{equation}\label{eq:intro-selected-decay}
		\mathcal E(t)+|\mathbf a(t)|_P^2\leq C\eps T^{\eta-d_*}(T-t)^{d_*}.
	\end{equation}
	By \eqref{eq:intro-selected-decay}, we have
	\begin{equation}\label{eq:intro-time-integrals}
		\int_0^t\frac{\mathcal E(q)+|\mathbf a(q)|_P^2}{T-q}\dd q\leq\frac{C\eps T^\eta}{d_*},\qquad \eps\int_0^t(T-q)^{\eta-1}\dd q\leq\frac{\eps T^\eta}{\eta}.
	\end{equation}
	Applying Gr\"onwall's inequality to \eqref{eq:intro-global-energy} and using \eqref{eq:intro-time-integrals}, we obtain
	\begin{equation}\label{eq:intro-selected-global}
		\sup_{0\leq q\leq t}\mathcal B_{\rm reg}(q)+\frac{c_{\rm reg}}\eps\int_0^t\mathcal D_{\rm reg}(q)\dd q\leq C\eps T^\eta.
	\end{equation}
	For sufficiently small $\eps>0$, \eqref{eq:intro-selected-decay} and \eqref{eq:intro-selected-global} ensure the smallness of the perturbations. Hence the solution can be continued as long as $\mathcal Q_{\rm sh}\leq L_{\rm sh}$ and $t<T$.
	
	It remains to choose $\mathbf a_0$ such that $\mathcal Q_{\rm sh}(t)\leq L_{\rm sh}$ for all $t<T$. For a sufficiently small fixed constant $R_{\rm loc}>0$, we consider the set of initial parameters
	\[\mathcal A_\eps=\{\mathbf a_0\in\R^J:|\mathbf a_0|_P\leq R_{\rm loc},\ \mathcal Q_{\rm sh}(0;\mathbf a_0)\leq L_{\rm sh}\}.\]
	For sufficiently small $\eps$, the set $\mathcal A_\eps$ can be parametrized radially by the closed unit ball in $\R^J$, with $\mathcal Q_{\rm sh}(0;\mathbf a_0)=L_{\rm sh}$ on $\partial\mathcal A_\eps$. Suppose that for every $\mathbf a_0\in\mathcal A_\eps$, the corresponding solution reaches the level $\mathcal Q_{\rm sh}=L_{\rm sh}$ before some fixed time $S<T$. Let $t_{\rm exit}(\mathbf a_0)$ be the first time at which $\mathcal Q_{\rm sh}=L_{\rm sh}$. By \eqref{eq:intro-outward-condition}, the exit time depends continuously on $\mathbf a_0$. Since $t_{\rm exit}(\mathbf a_0)=0$ for $\mathbf a_0\in\partial\mathcal A_\eps$, the map
	\[\mathbf a_0\mapsto\frac{\mathbf a(t_{\rm exit}(\mathbf a_0);\mathbf a_0)}{|\mathbf a(t_{\rm exit}(\mathbf a_0);\mathbf a_0)|_P}\]
	would give, after the radial parametrization, a continuous retraction from the closed unit ball onto its boundary, contradicting Brouwer's fixed point theorem. Letting $S\uparrow T$ and using the compactness of $\mathcal A_\eps$ and the continuous dependence on the initial data, we obtain an initial parameter $\mathbf a_0$ with $|\mathbf a_0|_P\leq C\eps^{1/2}T^{\eta/2}$ such that \eqref{eq:intro-selected-decay} and \eqref{eq:intro-selected-global} hold for all $t<T$. In the case $J=0$, we simply take $f_0=\bar M(0)$. See Theorem~\ref{thm:coupled-shooting} for details.
	
	\item \emph{Formation of the singularity.} By \eqref{eq:intro-selected-global} and Sobolev embedding, we have
	\[\sup_{0\leq t<T}\|H/W\|_{L^\infty_{x,v}}\leq C\eps^{1/2}T^{\eta/2}.\]
	Since both $M$ and $W$ are uniformly bounded, we obtain $\sup_{0\leq t<T}\|f(t)\|_{L^\infty_{x,v}}\leq C$. On the other hand, by \eqref{eq:intro-selected-decay} and Sobolev embedding, for every fixed $R>0$, we have
	\begin{equation}\label{eq:intro-local-decay}
		\sup_{|x|\leq R(T-t)^{1/r}}\left(|n-1|+\|F(t,x)-\mu\|_{L^2_{\widehat v}}\right)\leq C_R\mathcal E(t)^{1/2}\to0
	\end{equation}
	as $t\uparrow T$. In particular, since $\rho_f=\theta^3n$, we obtain
	\[\rho_f(t,0)=\vartheta(0)^3(T-t)^{-3\alpha_{\rm ss}}(1+o(1))\to\infty\qquad\text{as }t\uparrow T.\]
	
	To prove the blow-up of the spatial H\"older seminorm, we choose $y_*\ne0$ such that $U(y_*)\ne0$. By \eqref{eq:intro-local-decay}, the limiting Maxwellians at $x=0$ and $x=(T-t)^{1/r}y_*$ in the velocity variable $(T-t)^{\alpha_{\rm ss}}v$ have different centers. Hence, for $t$ sufficiently close to $T$, we have
	\[\sup_{v\in\R^3}|f(t,(T-t)^{1/r}y_*,v)-f(t,0,v)|\geq c>0,\]
	which yields
	\[[f(t)]_{C_x^{\delta_x}L_v^\infty}\gtrsim(T-t)^{-\delta_x/r}.\]
	See Subsection~\ref{subsec:main-proof} for more details.
\end{enumerate}

\medskip

\noindent\textbf{Organization of the paper.} In Section~\ref{sec:roadmap}, we establish a road map, state the main propositions and prove Theorem~\ref{thm:main}. Section~\ref{Sec.backward-Euler-correction} is devoted to the construction of the imploding Euler solution. In Section~\ref{sec:velocity}, we construct the velocity weights and establish the corresponding collision estimates. In Section~\ref{sec:moments}, we derive the macroscopic and microscopic equations and study the linearized Euler operator. In Section~\ref{sec:energies}, we construct the coupled energies and establish the weighted norm estimates. Sections~\ref{sec:collision-hierarchy} and \ref{sec:expanded-critical} are devoted to the estimates of the collision and transport terms, respectively. In Section~\ref{sec:closed}, we derive the closed energy inequalities. In Section~\ref{sec:selection}, we use Brouwer's fixed point theorem to select the initial data. In Appendix~\ref{app:polynomial-local}, we establish the local existence theory and the continuation criterion.

\medskip

\noindent\textbf{Notations and conventions.} 
\begin{itemize}
	\item We denote $\N_0:=\Z_{\geq0}:=\{0,1,2,\ldots\}$, $\langle z\rangle:=(1+|z|^2)^{1/2}$ and $B_R:=\{z\in\R^3:|z|<R\}$.
	
	\item For multi-indices $\alpha,\beta\in\N_0^3$, we write
	\[|\alpha|:=\alpha_1+\alpha_2+\alpha_3,\qquad \partial_x^\alpha:=\partial_{x_1}^{\alpha_1}\partial_{x_2}^{\alpha_2}\partial_{x_3}^{\alpha_3},\qquad \partial_v^\beta:=\partial_{v_1}^{\beta_1}\partial_{v_2}^{\beta_2}\partial_{v_3}^{\beta_3}.\]
	
	\item The derivatives $\partial_t$ and $\partial_x^\alpha$ of functions in the physical coordinates $(t,x,v)$ are taken at fixed $v$. For functions in the normalized coordinates $(t,x,\widehat v)$, the derivatives $\partial_t$ and $\partial_x^\alpha$ are taken at fixed $\widehat v$. All velocity integrals without an indicated domain are over $\R^3$.
	
	\item For nonnegative quantities $A$ and $B$, we write $A\lesssim B$ if $A\leq CB$ for some constant $C>0$, and $A\asymp B$ if $A\lesssim B\lesssim A$. The positive constants $c,C$ may change from line to line and may depend on the fixed profiles, differentiation orders and weight parameters, but are independent of $t$ and $\eps$. We use $A\lesssim_\lambda B$ to indicate the dependence of the implicit constant on the parameter $\lambda$.
\end{itemize}

\section{A road map of the proof}\label{sec:roadmap}
The proof separates the imploding Euler background from the macroscopic and microscopic corrections to the Boltzmann distribution. In this section, we state the Euler background construction, the uniform collision estimates, the matched energy identities and the finite-dimensional correction of the initial data, and then prove Theorem~\ref{thm:main}. The detailed proofs of intermediate propositions are delayed to the subsequent sections.

\subsection{Imploding solutions of the Euler equations}
Formally, in the hydrodynamic limit $\varepsilon\to0$, equation \eqref{Eq.Boltzmann} is governed by the compressible Euler equations. We begin with the three-dimensional isentropic Euler system for a monatomic gas:
\begin{equation}\label{Eq.Euler-monatomic}
	\begin{cases}
		\partial_t\rho+\dive(\rho u)=0,\\
		\rho(\pa_tu+u\cdot\nabla u)+\nabla\left(\rho^{5/3}\right)=0.
	\end{cases}
\end{equation}
Assume that $(\rho,u)$ is a classical solution to \eqref{Eq.Euler-monatomic}, where
\[\rho:[0,T)\times\mathbb R^3\to\mathbb R_{>0}, \qquad u:[0,T)\times\mathbb R^3\to\mathbb R^3.\]
We define
\begin{equation}\label{Eq.theta-def}
	\theta(t,x):=\rho(t,x)^{1/3}.
\end{equation}
In terms of $(u,\theta)$, system \eqref{Eq.Euler-monatomic} is equivalently written as
\begin{equation}\label{Eq.Euler-u-theta}
	\begin{cases}
		\partial_tu+u\cdot\nabla u+5\theta\nabla\theta=0,\\
		\partial_t\theta+u\cdot\nabla\theta +\dfrac13\theta\,\dive u=0.
	\end{cases}
\end{equation}

Our analysis begins with the self-similar profiles constructed in \cite{SWWZ2025}. We fix a self-similarity exponent $r>1$. Substituting the self-similar ansatz
\begin{equation}\label{Eq.Euler-monatomic-ss-ansatz}
	\begin{aligned}
		u(t,x)&:=\frac{(T-t)^{\frac1r-1}}{r} u_{\rm s}\left(\frac{x}{(T-t)^{\frac1r}}\right),\\
		\theta(t,x)&:=\frac{(T-t)^{\frac1r-1}}{\sqrt{15}\,r} \theta_{\rm s}\left(\frac{x}{(T-t)^{\frac1r}}\right)
	\end{aligned}
\end{equation}
into \eqref{Eq.Euler-u-theta}, we obtain the following system for the profiles $(u_{\rm s},\theta_{\rm s})$:
\begin{equation}\label{Eq.Euler-monatomic-ss}
	\begin{cases}
		(r-1)u_{\rm s}+(x+u_{\rm s})\cdot\nabla u_{\rm s}+\frac13\theta_{\rm s}\nabla\theta_{\rm s}=0,\\
		(r-1)\theta_{\rm s}+(x+u_{\rm s})\cdot\nabla \theta_{\rm s}+\frac13\theta_{\rm s}\dive u_{\rm s}=0.
	\end{cases}
\end{equation}
The following result is proved in \cite[Theorem 1.1]{SWWZ2025}.

\begin{theorem}[Existence of self-similar Euler profiles]\label{Thm.profiles}
	There exists a discrete sequence
	\[\{r_n\}_{n\geq1}\subset(1,3-\sqrt3), \qquad \lim_{n\to\infty}r_n=3-\sqrt3,\]
	such that, for every $r\in\{r_n\}_{n\geq1}$, system \eqref{Eq.Euler-monatomic-ss} admits a globally defined, radially symmetric $C^\infty$ solution $(u_{\rm s},\theta_{\rm s})$ on $\mathbb R^3$. This solution satisfies
	\begin{equation}\label{eq:profile-decay}
		\begin{aligned}
			\theta_{\rm s}(x)&\gtrsim_r\langle x\rangle^{-(r-1)},\\
			|\nabla^j u_{\rm s}(x)|+|\nabla^j\theta_{\rm s}(x)| &\lesssim_{r,j}\langle x\rangle^{-(r-1)-j}, \qquad j\in\mathbb Z_{\geq0}.
		\end{aligned}
	\end{equation}
	Writing, for $R:=|x|>0$,
	\[u_{\rm s}(x)=U_{\rm s}(R)\frac{x}{R}, \qquad \theta_{\rm s}(x)=\Theta_{\rm s}(R),\]
	there exists a constant $c_{\rm rep}=c_{\rm rep}(r)>0$ such that, for every $R>0$,\footnote{Throughout the paper, the lowercase $r$ denotes the self-similar exponent and is a fixed parameter, while the capital $R$ denotes the spatial radial variable. These two notations should not be confused.}
	\[1+U_{\rm s}'(R)-\frac13|\Theta_{\rm s}'(R)|>c_{\rm rep}, \qquad 1+\frac{U_{\rm s}(R)}R-\frac13|\Theta_{\rm s}'(R)|>c_{\rm rep}.\]
	All radial expressions extend smoothly to $R=0$, and the inequalities continue to hold there.
\end{theorem}

We fix the exponent $r$ selected in \eqref{eq:exponents} from the sequence in Theorem~\ref{Thm.profiles}. The normalized Euler profiles and the sound-speed constant are given by
\begin{equation}\label{eq:normalized-profile}
	U:=\frac{u_{\rm s}}r,\qquad \vartheta:=\frac{\theta_{\rm s}}{\sqrt{15}\,r},\qquad c_{\rm s}:=\sqrt{\frac53}.
\end{equation}
Writing $U(y)=U_{\rm rad}(|y|)y/|y|$ and $\vartheta(y)=\vartheta_{\rm rad}(|y|)$ for $y\ne0$, the radial inequalities in Theorem~\ref{Thm.profiles} imply
\begin{equation}\label{eq:normalized-repulsivity}
	\frac1r+U_{\rm rad}'(\varrho)-c_{\rm s}|\vartheta_{\rm rad}'(\varrho)|\geq\frac{c_{\rm rep}(r)}r>0,\qquad \varrho\geq0.
\end{equation}
Here the radial expressions in \eqref{eq:normalized-repulsivity} are defined at $\varrho=0$ by smooth extension.

By \eqref{Eq.theta-def}, \eqref{Eq.Euler-monatomic-ss-ansatz}, and Theorem \ref{Thm.profiles}, every admissible profile generates a self-similar imploding solution to \eqref{Eq.Euler-monatomic}. More precisely, let the sequence $\{r_n\}$ be given by Theorem~\ref{Thm.profiles}, let $r\in\{r_n\}$, and let $(u_{\rm s}, \theta_{\rm s})$ be the smooth self-similar profile constructed in Theorem \ref{Thm.profiles}. Let $T>0$. We define
\begin{equation}\label{Eq.u_p-theta_p-def}
	\begin{aligned}
		u_{\rm p}(t,x)&:=\frac{(T-t)^{\frac1r-1}}{r} u_{\rm s}\left(\frac{x}{(T-t)^{\frac1r}}\right)=(T-t)^{\frac1r-1}U\left(\frac{x}{(T-t)^{\frac1r}}\right),\\
		\theta_{\rm p}(t,x)&:=\frac{(T-t)^{\frac1r-1}}{\sqrt{15}\,r} \theta_{\rm s}\left(\frac{x}{(T-t)^{\frac1r}}\right)=(T-t)^{\frac1r-1}\vartheta\left(\frac{x}{(T-t)^{\frac1r}}\right).
	\end{aligned}
\end{equation}
Setting $\rho_{\rm p}:=\theta_{\rm p}^3$, we obtain a solution $(\rho_{\rm p},u_{\rm p})$ to \eqref{Eq.Euler-monatomic}. The total Euler energy of this solution is infinite, and $\rho_{\rm p}(t,x)\to0$ as $|x|\to\infty$. The kinetic construction requires a compact perturbation of a strictly positive constant state. Therefore, we seek a truncated solution to \eqref{Eq.Euler-monatomic} of the form
\begin{equation}\label{Eq.Euler-monatomic-ansatz}
	\begin{aligned}
		u(t,x)&:=u_{\rm p}(t,x)\chi\left(\frac{|x|}{R_0}\right)+u_{\rm r}(t,x),\\
		\theta(t,x)&:=\theta_{\rm p}(t,x)\chi\left(\frac{|x|}{R_0}\right) +R_0^{1-r}\left(1-\chi\left(\frac{|x|}{R_0}\right)\right) +\theta_{\rm r}(t,x),
	\end{aligned}
\end{equation}
Here $R_0>10$, and $\chi\in C_c^\infty(\mathbb R;[0,1])$ is a smooth cutoff satisfying
\[\operatorname{supp}\chi\subset[-1,1], \qquad \chi\equiv1\quad\text{on }[-1/2,1/2].\]
The functions $u_{\rm r}$ and $\theta_{\rm r}$ are correction terms to be constructed backward from $t=T$. For convenience, we set
\[\chi_{R_0}(x):=\chi\left(\frac{|x|}{R_0}\right), \qquad \bar u_{\rm p}:=u_{\rm p}\chi_{R_0}, \qquad \bar\theta_{\rm p}:=\theta_{\rm p}\chi_{R_0} +R_0^{1-r}(1-\chi_{R_0}).\]
We define the truncation errors by
\begin{align*}
	\cE_u &:=\partial_t\bar u_{\rm p} +\bar u_{\rm p}\cdot\nabla\bar u_{\rm p} +5\bar\theta_{\rm p}\nabla\bar\theta_{\rm p},\\
	\cE_\theta &:=\partial_t\bar\theta_{\rm p} +\bar u_{\rm p}\cdot\nabla\bar\theta_{\rm p} +\frac13\bar\theta_{\rm p}\dive\bar u_{\rm p}.
\end{align*}
Since $(u_{\rm p},\theta_{\rm p})$ solves \eqref{Eq.Euler-u-theta}, both errors vanish outside the transition annulus:
\[\operatorname{supp}(\cE_u,\cE_\theta) \subset \left\{x\in\mathbb R^3: R_0/2\leq|x|\leq R_0\right\}.\]
Substituting \eqref{Eq.Euler-monatomic-ansatz} into \eqref{Eq.Euler-u-theta}, we obtain the following system for $(u_{\rm r},\theta_{\rm r})$:
\begin{equation}\label{Eq.u_r-theta_r-eq}
	\begin{cases}
		\begin{aligned}
			\partial_tu_{\rm r} &+(\bar u_{\rm p}+u_{\rm r})\cdot\nabla u_{\rm r} +u_{\rm r}\cdot\nabla\bar u_{\rm p}\\
			&\qquad\qquad +5(\bar\theta_{\rm p}+\theta_{\rm r}) \nabla\theta_{\rm r} +5\theta_{\rm r}\nabla\bar\theta_{\rm p} =-\cE_u,
		\end{aligned}
		\\[3mm] \begin{aligned} \partial_t\theta_{\rm r} &+(\bar u_{\rm p}+u_{\rm r})\cdot\nabla\theta_{\rm r} +u_{\rm r}\cdot\nabla\bar\theta_{\rm p}\\ &\qquad\qquad +\frac13(\bar\theta_{\rm p}+\theta_{\rm r}) \dive u_{\rm r} +\frac13\theta_{\rm r}\dive\bar u_{\rm p} =-\cE_\theta. \end{aligned} \end{cases}
\end{equation} 
\begin{proposition}[A truncated imploding Euler solution]\label{Thm.Euler-solution} 
	Let $r\in\{r_n\}_{n\geq1}$, and let $(u_{\rm s},\theta_{\rm s})$ be a profile from Theorem~\ref{Thm.profiles}. There exist $T_0\in(0,1]$ and $R_*>10$ such that, for every $0<T\leq T_0$ and $R_0\geq R_*$, system \eqref{Eq.u_r-theta_r-eq} has a smooth radially symmetric solution $(u_{\rm r},\theta_{\rm r})$ on $[0,T]$ with $(u_{\rm r},\theta_{\rm r})(T)=0$. For $0\leq t\leq T$ and $k\in\mathbb N_0$, 
	\begin{align} 
		\supp u_{\rm r}(t)\cup\supp\theta_{\rm r}(t)\subset\{R_0/4<|x|<2R_0\},\label{Eq.remainder-support}\\ \|\nabla_x^ku_{\rm r}(t)\|_{L^\infty_x}+\|\nabla_x^k\theta_{\rm r}(t)\|_{L^\infty_x}\leq C_k(T-t)R_0^{1-2r-k}.\label{Eq.remainder-pointwise-bound} 
	\end{align} 
	Moreover, the functions $(u,\theta)$ defined by \eqref{Eq.Euler-monatomic-ansatz} solve \eqref{Eq.Euler-u-theta} on $[0,T)\times\R^3$ and satisfy 
	\begin{align} 
		\theta(t,x)\geq c_0R_0^{1-r},\qquad\supp u(t)\cup\supp(\theta(t)-R_0^{1-r})\subset B_{2R_0},\label{eq:background}\\ (u,\theta)(t,x)=(T-t)^{-\alpha_{\rm ss}}(U,\vartheta)\left(\frac{x}{(T-t)^{1/r}}\right),\qquad |x|\leq R_0/4,\label{eq:profile} 
	\end{align} 
	where $(U,\vartheta)$ is defined in \eqref{eq:normalized-profile}. For all $j,k\in\mathbb N_0$, \begin{equation}\label{eq:background-hierarchy} 
		|\partial_t^j\nabla_x^ku|+|\partial_t^j\nabla_x^k\theta|\leq C_{j,k}\theta^{1+(jr+k)/(r-1)}. 
	\end{equation} 
	Here the constants $T_0,R_*,c_0,C_k,C_{j,k}$ are independent of $T$ and $R_0$. In particular, $(\rho,u):=(\theta^3,u)$ solves \eqref{Eq.Euler-monatomic}. 
\end{proposition} 
The proof of Proposition~\ref{Thm.Euler-solution} is given in Section~\ref{subsec:euler-correction-completion}; it uses the backward symmetric-hyperbolic construction in Section~\ref{subsec:euler-terminal} and the finite speed of propagation. In the kinetic construction, we fix $T$, $R_0$ and the Euler solution in Proposition~\ref{Thm.Euler-solution}. We define 
\[\theta_\infty:=R_0^{1-r},\qquad\theta_*:=c_0R_0^{1-r},\qquad R_{\rm e}:=R_0/4,\qquad d(t,x):=\big((T-t)^{2/r}+|x|^2\big)^{1/2}.\]
By the profile bounds, \eqref{Eq.remainder-pointwise-bound} and \eqref{eq:background-hierarchy}, we have
\begin{equation}\label{eq:background-bounds}
	|\nabla_xu|+|\nabla_x\theta|\leq C\theta^p,\qquad\theta^p\leq\frac{C}{T-t},\qquad d^j|\nabla_x^j\log\theta|\leq C_j\quad(j\geq1).
\end{equation}
At $x=0$, by \eqref{eq:profile}, we have $\rho(t,0)=\vartheta(0)^3(T-t)^{-3\alpha_{\rm ss}}\to\infty$ and
\[\theta(t,0)^{3+\gamma}\asymp(T-t)^{-1-\eta}.\]


\subsection{A coordinate transform and the microscopic decomposition}\label{Subsec:coordinate-transform}
The Euler solution fixed by Proposition \ref{Thm.Euler-solution} determines the expected mean velocity $u$ and the thermal speed $\theta$. We write
\begin{equation}\label{Eq.f-F-relation}
	f(t,x,v)=F\left(t,x,\frac{v-u(t,x)}{\theta(t,x)}\right),\qquad\widehat v:=\frac{v-u(t,x)}{\theta(t,x)}.
\end{equation}
At fixed $(t,x)$, we have $\dd v=\theta^3\dd\widehat v$. Differentiating \eqref{Eq.f-F-relation}, we obtain
\[\partial_tF+(u+\theta\widehat v)\cdot\nabla_xF+\mathsf H_{\rm tr}\cdot\nabla_{\widehat v}F=\frac{\theta^{3+\gamma}}{\eps}Q(F,F),\]
where 
\[\mathsf H_{\rm tr}:=-\frac{\partial_tu+u\cdot\nabla_xu+\theta\widehat v\cdot\nabla_xu}{\theta}-\frac{\partial_t\theta+u\cdot\nabla_x\theta+\theta\widehat v\cdot\nabla_x\theta}{\theta}\widehat v.\]
Using \eqref{Eq.Euler-u-theta}, we obtain
\begin{equation}\label{Eq.H-expression}
	\mathsf H_{\rm tr}=\left(\frac13\Div u-\widehat v\cdot\nabla_x\theta\right)\widehat v+5\nabla_x\theta-\widehat v\cdot\nabla_xu.
\end{equation}

For $F\geq0$ with positive density and temperature, we define its normalized density, mean velocity and temperature by
\begin{equation}\label{eq:kinetic-normalization}
	n:=\int F\dd\widehat v,\qquad V:=\frac1n\int\widehat vF\dd\widehat v,\qquad\Theta:=\frac1{3n}\int|\widehat v-V|^2F\dd\widehat v.
\end{equation}
The change of variables $v=u+\theta\widehat v$ gives the physical moments
\begin{equation}\label{eq:physical-parameters}
	\begin{aligned}
		\rho_f&:=\int f\dd v=\theta^3n,\qquad
		u_f:=\frac1{\rho_f}\int vf\dd v=u+\theta V,\\
		T_f&:=\frac1{3\rho_f}\int|v-u_f|^2f\dd v=\theta^2\Theta,\qquad P_f:=(\rho_f,u_f,T_f).
	\end{aligned}
\end{equation}
Let $M:=\mathcal M(\rho_f,u_f,T_f)$ be the local Maxwellian with the same density, mean velocity and temperature as $f$. In the normalized velocity, we have
\begin{equation}\label{eq:normalized-Maxwellian}
	\widetilde M(t,x,\widehat v):=M(t,x,u+\theta\widehat v)=\frac{n}{(2\pi\Theta)^{3/2}}\mathrm e^{-\frac{|\widehat v-V|^2}{2\Theta}}=\mathcal M(n, V, \Theta).
\end{equation}
We define the remainders by
\begin{equation}\label{eq:microscopic-remainders}
	H:=f-M,\qquad h(t,x,\widehat v):=H(t,x,u+\theta\widehat v)=F-\widetilde M.
\end{equation}
Since $f$ and $M$ have the same moments, we have
\begin{equation}\label{eq:zero-moments}
	\begin{aligned}
		\int H\dd v&=0,&\int vH\dd v&=0,&\int|v|^2H\dd v&=0,\\
		\int h\dd\widehat v&=0,&\int\widehat vh\dd\widehat v&=0,&\int|\widehat v|^2h\dd\widehat v&=0.
	\end{aligned}
\end{equation}
The three columns in each line are the mass, momentum and energy moment conditions, respectively.

We distinguish differentiation in physical velocity coordinates from evaluation in normalized velocity coordinates. The expression $(\partial_x^\alpha H)(t,x,u+\theta\widehat v)$ means that $H(t,x,v)$ is first differentiated in $x$ with $v$ fixed, and that $v=u+\theta\widehat v$ is substituted afterwards. For example,
\begin{equation}\label{eq:physical-derivative-convention}
	(\partial_{x_i}H)(t,x,u+\theta\widehat v)=\partial_{x_i}h-\frac{\partial_{x_i}u+\widehat v\,\partial_{x_i}\theta}{\theta}\cdot\nabla_{\widehat v}h,\qquad i=1,2,3,
\end{equation}
where the derivatives of $h$ are taken at fixed $\widehat v$. In velocity integrals, we abbreviate the substitution in physical functions. In particular,
\begin{equation}\label{eq:physical-integral-convention}
	\int\frac{|\partial_x^\alpha H|^2}{M}\dd\widehat v:=\int\frac{|(\partial_x^\alpha H)(t,x,u+\theta\widehat v)|^2}{\widetilde M(t,x,\widehat v)}\dd\widehat v=\theta^{-3}\int\frac{|\partial_x^\alpha H(t,x,v)|^2}{M(t,x,v)}\dd v.
\end{equation}
The same convention applies to other functions originally defined in $(t,x,v)$. Thus, a denominator $M$ under $\dd\widehat v$ denotes the function $\widetilde M$ in \eqref{eq:normalized-Maxwellian}. Differentiating the physical moment identity in \eqref{eq:zero-moments}, we get
\begin{equation}\label{eq:differentiated-zero-moments}
	\int(1,v,|v|^2)\partial_x^\alpha H\dd v=\partial_x^\alpha\int(1,v,|v|^2)H\dd v=0,\qquad\alpha\in\mathbb N_0^3.
\end{equation}

The prescribed Euler Maxwellian is $\bar M:=\mathcal M(\theta^3,u,\theta^2)$, and hence $\bar M(t,x,u+\theta\widehat v)=\mu(\widehat v)$. We define $D:=\partial_t+v\cdot\nabla_x$, with the physical velocity held fixed, and $S_0:=D\bar M$. By direct differentiation, we have
\begin{equation}\label{eq:source}
	S_0(t,x,u+\theta\widehat v)=\left[\left(\widehat v\otimes\widehat v-\frac13|\widehat v|^2\Id\right):\nabla_xu+(|\widehat v|^2-5)\widehat v\cdot\nabla_x\theta\right]\mu(\widehat v).
\end{equation}
For $i,j\in\{1,2,3\}$, by Gaussian integration, we have
\[\int\widehat v_i\widehat v_j\mu\dd\widehat v=\delta_{ij},\qquad\int|\widehat v|^2\widehat v_i\widehat v_j\mu\dd\widehat v=5\delta_{ij},\qquad\int|\widehat v|^4\mu\dd\widehat v=15.\]
Here $\delta_{ij}$ is the Kronecker symbol. Moreover, all Gaussian moments containing an odd power of one coordinate vanish. Applying these identities to \eqref{eq:source} and then substituting $v=u+\theta\widehat v$, we obtain
\begin{equation}\label{eq:source-moments}
	\begin{aligned}
		\int(1,\widehat v,|\widehat v|^2)S_0(t,x,u+\theta\widehat v)\dd\widehat v=0,\quad\int(1,v,|v|^2)S_0(t,x,v)\dd v=0.
	\end{aligned}
\end{equation}
As a consequence,
\begin{equation}\label{eq:source-moments-spatial}
	\int(1,v,|v|^2)\partial_x^\alpha S_0(t,x,v)\dd v=0,\qquad\alpha\in\mathbb N_0^3.
\end{equation}

The moments $(n,V,\Theta)$ are convenient for reconstructing $M$. For the Euler part of the argument, we use the pressure $\rho_fT_f$ and the logarithmic entropy variable $\log(T_f/\rho_f^{2/3})$. The Euler Maxwellian $\bar M$ has density $\theta^3$ and temperature $\theta^2$, and hence its pressure is $\theta^5$ and its logarithmic entropy variable is zero. Therefore, we define
\begin{equation}\label{eq:roadmap-macro}
	\theta_f:=(\rho_fT_f)^{1/5},\qquad B:=\frac{\theta_f}{\theta}-1,\qquad S:=\log\frac{T_f}{\rho_f^{2/3}},\qquad\mathbf z:=(V,B,S).
\end{equation}
The variable $V$ is the normalized mean velocity defined in \eqref{eq:kinetic-normalization}. By \eqref{eq:physical-parameters}, we have
\[B=(n\Theta)^{1/5}-1,\qquad S=\log\Theta-\frac23\log n.\]
Solving for $n$ and $\Theta$, we get
\begin{equation}\label{eq:normalized-reconstruction}
	n=(1+B)^3\mathrm e^{-3S/5},\qquad\Theta=(1+B)^2\mathrm e^{3S/5}.
\end{equation}
The transform \eqref{eq:normalized-reconstruction} reconstructs the same moments $n$ and $\Theta$ as in \eqref{eq:kinetic-normalization}. In particular, 
\begin{equation}\label{Eq.background}
	M=\bar M=\cM(\theta^3, u,\theta^2)\ \Longleftrightarrow\ (n,V,\Theta)=(1,0,1)\ \Longleftrightarrow\ \mathbf z=0.
\end{equation}
The change from $(n,V,\Theta)$ to $(V,B,S)$ is smooth and invertible for $n,\Theta>0$, or equivalently $B>-1$ and $S\in\R$.

The non-Maxwellian part of the momentum and energy fluxes is described by the stress tensor $\Pi$ and the heat flux $q$:
\begin{equation}\label{eq:roadmap-stress-heat}
	\begin{aligned}
		\Pi&:=\int(v-u_f)\otimes(v-u_f)f\dd v-\rho_fT_f\Id
		=\int(v-u_f)\otimes(v-u_f)H\dd v,\\
		q&:=\frac12\int|v-u_f|^2(v-u_f)f\dd v
		=\frac12\int|v-u_f|^2(v-u_f)H\dd v.
	\end{aligned}
\end{equation}
We use the notation $(\Div\Pi)_i:=\sum_j\partial_{x_j}\Pi_{ij}$ and $\Pi:\nabla_xu_f:=\sum_{i,j}\Pi_{ij}\partial_{x_j}u_{f,i}$. Integrating the Boltzmann equation \eqref{Eq.Boltzmann} against $1$, $v$ and $|v|^2$, we obtain, with $D_f:=\partial_t+u_f\cdot\nabla_x$,
\begin{equation}\label{eq:moment-system}
	\begin{cases}
		D_f\rho_f+\rho_f\Div u_f=0,\\
		D_fu_f+\rho_f^{-1}\nabla_x(\rho_fT_f)=-\rho_f^{-1}\Div\Pi,\\
		\frac32\rho_fD_fT_f+\rho_fT_f\Div u_f=-\Pi:\nabla_xu_f-\Div q.
	\end{cases}
\end{equation}
The use of the pressure variable removes the derivatives of $S$ from the acoustic coupling. The exact reformulation of \eqref{eq:moment-system} is recorded in the following lemma, whose proof by direct computation is given in Section~\ref{subsec:moment-microscopic}.

\begin{lemma}\label{prop:pressure}
	Let $(\rho_f,u_f,T_f,\Pi,q)$ be smooth functions satisfying \eqref{eq:moment-system}, with $\rho_f,T_f>0$. We define $\theta_f$ and $S$ by \eqref{eq:roadmap-macro}, and set $\mathcal C_f:=\Pi:\nabla_xu_f+\Div q$. Then
	\begin{equation}\label{eq:reconstruction}
		\rho_f=\theta_f^3\mathrm e^{-3S/5},\qquad T_f=\theta_f^2\mathrm e^{3S/5},
	\end{equation}
	and
	\begin{equation}\label{eq:pressure-system}
		\begin{cases}
			D_fu_f+5\theta_f\mathrm e^{3S/5}\nabla_x\theta_f=-\rho_f^{-1}\Div\Pi,\\
			D_f\theta_f+\frac13\theta_f\Div u_f=-\frac2{15}\theta_f^{-4}\mathcal C_f,\\
			D_fS=-\frac23\theta_f^{-5}\mathcal C_f.
		\end{cases}
	\end{equation}
\end{lemma}

\subsection{Velocity weights and collision estimates}
In this subsection, we first define a family of weights in an independent velocity variable $v\in\R^3$ and establish estimates for the collision operator linearized around $\mu$. In these estimates, the parameters $m$ and $\ell$ are fixed. Then we apply the weights to the local Maxwellian $M$ and choose $m$ as a function of $(t,x)$.

We fix a smooth nonincreasing function $\chi_{\rm vel}:[1,\infty)\to[0,1]$, which is equal to $1$ near $1$ and equal to $0$ on $[2,\infty)$. Let $\ell\geq1$ be a fixed tail-order parameter, and let $m\geq\ell$. We set
\begin{equation}\label{eq:J}
	I_{\rm vel}(q):=\int_1^q\frac{\chi_{\rm vel}(s)}s\dd s\quad(q\geq1), \qquad I_{{\rm vel},\infty}:=I_{\rm vel}(2).
\end{equation}
In particular, $0\leq I_{\rm vel}(q)\leq\log2$. We define
\begin{equation}\label{eq:Phi-fixed-tail}
	\Phi_{m,\ell}(a):=
	\begin{cases}
		a,&0\leq a\leq m,\\
		m+\ell\log(a/m)+(m-\ell)I_{\rm vel}(a/m),&a\geq m,
	\end{cases}
\end{equation}
and
\begin{equation}\label{eq:fixed-tail-weight}
	w_{m,\ell}(v):=(2\pi)^{-3/4}\mathrm e^{-\Phi_{m,\ell}(|v|^2)/4},\qquad v\in\R^3.
\end{equation}

\begin{lemma}\label{thm:fixed-tail}
	For fixed $\ell$ and $m\geq\ell$, the function $\Phi_{m,\ell}$ is increasing, concave, and $C^{1,1}$ on bounded intervals. It satisfies
	\begin{equation}\label{eq:weight-order}
		\sqrt\mu\leq w_{m,\ell},\qquad w_{m,\ell}=\sqrt\mu\quad\text{on }|v|^2\leq m.
	\end{equation}
	On $|v|^2\geq2m$, its decay in $|v|$ has the fixed order $\ell/2$. If $L\geq0$ is fixed and $\ell>2L+6$, then
	\begin{equation}\label{eq:new-moments}
		\sup_{m\geq\ell} \left\{\int\langle v\rangle^{2L}w_{m,\ell}^2\dd v +\int\langle v\rangle^Lw_{m,\ell}\dd v\right\}<\infty,
	\end{equation}
	and for some $c>0$,
	\begin{equation}\label{eq:new-tail-small}
		\norm{\langle v\rangle^L(w_{m,\ell}-\sqrt\mu)}_{L^2_v} +\norm{\langle v\rangle^L(\mu/w_{m,\ell}-w_{m,\ell})}_{L^2_v} \leq C_{L,\ell}\mathrm e^{-cm}.
	\end{equation}
	For an elastic collision $(v,v_*,v',v_*')$,
	\begin{equation}\label{eq:elastic-ratio}
		\frac{w_{m,\ell}(v')w_{m,\ell}(v_*')}{w_{m,\ell}(v)}\leq1, \qquad \frac{w_{m,\ell}(v)w_{m,\ell}(v_*)}{w_{m,\ell}(v')}\leq1.
	\end{equation}
	In addition,
	\begin{equation}\label{eq:bounded-m-derivative}
		0\leq\partial_m\Phi_{m,\ell}(a)\leq1+\log2, \qquad |\partial_m\log w_{m,\ell}(v)|\leq\frac{1+\log2}{4}.
	\end{equation}
	Here all constants may depend on the fixed $\chi_{\rm vel}$ and $\ell$.
\end{lemma}

\begin{figure}[htbp]
	\centering
	\includegraphics[width=0.92\textwidth]{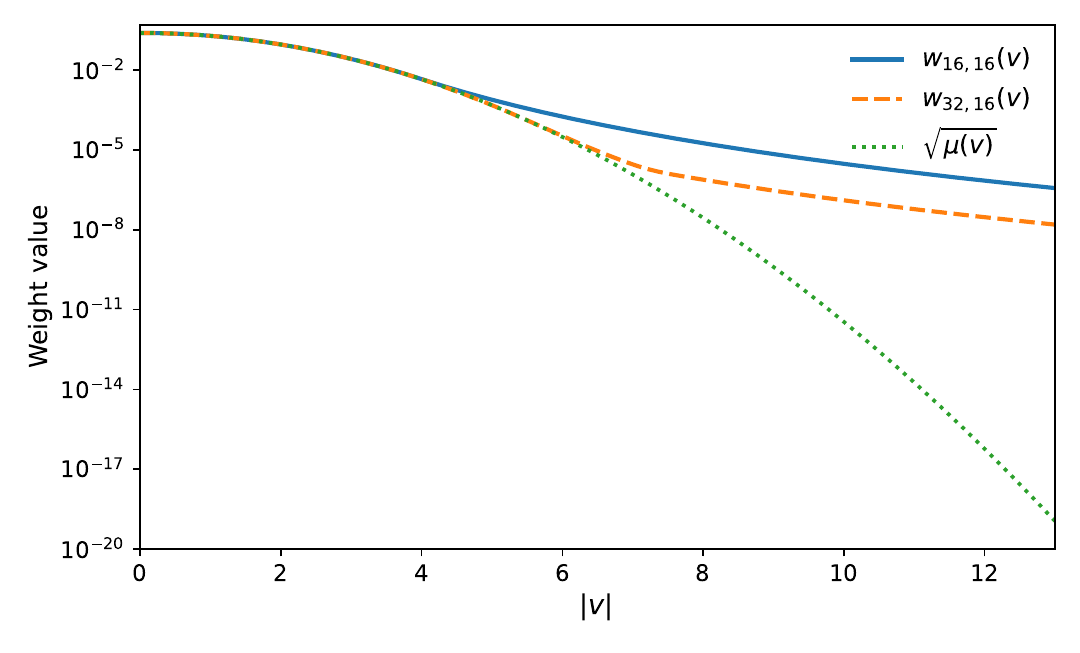}
	\caption{Comparison of $w_{m,\ell}(v)$ with $\sqrt{\mu(v)}$ for $\ell=16$ and $m=16,32$, using the same smooth cutoff $\chi_{\rm vel}$. The horizontal axis represents $|v|$, and the vertical axis is logarithmic. The weight $w_{m,\ell}(v)$ agrees with $\sqrt\mu$ for $|v|\leq\sqrt m$ and has polynomial decay of order $\ell/2$ for $|v|\geq\sqrt{2m}$.}
	\label{fig:weight-comparison}
\end{figure}

Figure~\ref{fig:weight-comparison} compares the radial profiles of $w_{m,\ell}$ with $\sqrt\mu$. As $m$ increases, the ball on which the two weights coincide is enlarged, while the polynomial decay order $\ell/2$ stays fixed.

Next we state our new weighted collision estimates. We write $w:=w_{m,\ell}$ and $w_\infty:=\sqrt\mu$, and define
\begin{align*}
	&\qquad\qquad\qquad\nu(v):=\int_{\R^3}|v-v_*|^\gamma\mu(v_*)\dd v_*,\\
	L_wg&:=-w^{-1}\{Q(\mu,wg)+Q(wg,\mu)\},\qquad Q_w(g_1,g_2):=w^{-1}Q(wg_1,wg_2).
\end{align*}
Recall \eqref{eq:postcollisional-velocities}. It is standard that
\begin{equation}\label{eq:normalized-frequency-derivatives}
	c\langle v\rangle^\gamma\leq\nu(v)\leq C\langle v\rangle^\gamma,\qquad
	|\nabla_v^j\nu(v)|\leq C_j\langle v\rangle^{\gamma-j}\quad(j=1,2).
\end{equation}

\begin{proposition}[Uniform collision bounds]\label{thm:new-collision}
	There exist $\ell_0$ and $c_0>0$ such that, after fixing any sufficiently large $\ell\geq\ell_0$ and choosing $m_0\geq\ell$ sufficiently large, for every $m\geq m_0$,
	\begin{equation}\label{eq:collision-missing}
		\int(1,v,|v|^2)w g\dd v=0 \quad\Longrightarrow\quad \operatorname{Re}\langle L_wg,g\rangle\geq c_0\|\nu^{1/2}g\|_{L^2_v}^2.
	\end{equation}
	Moreover,
	\begin{align}
		&\|\nu^{-1/2}Q_w(g_1,g_2)\|_{L^2_v}\leq C_\ell\big(\|\nu^{1/2}g_1\|_{L^2_v}\|g_2\|_{L^2_v}+\|g_1\|_{L^2_v}\|\nu^{1/2}g_2\|_{L^2_v}\big),\label{eq:new-bilinear}\\
		&\qquad\qquad\qquad\quad\|\nu^{-1/2}L_wg\|_{L^2_v}\leq C_\ell\|\nu^{1/2}g\|_{L^2_v}.\label{eq:new-linear-upper}
	\end{align}
	Here the constants are independent of $m$. The estimates first hold for smooth functions with convergent integrals, and then extend to the corresponding weighted spaces by density. 
\end{proposition}

We prove Proposition~\ref{thm:new-collision} in Subsection~\ref{subsec:uniform-collision}. For $w_\infty=\sqrt\mu$, we shall establish
\[|\langle(L_w-L_{w_\infty})g,g\rangle|\leq C_\gamma(\ell^{-1}+\mathrm e^{-cm})\|\nu^{1/2}g\|_{L^2_v}^2.\]
By first taking $\ell$ sufficiently large and then taking $m_0\geq\ell$ sufficiently large, we obtain \eqref{eq:collision-missing} with a positive constant $c_0$ independent of $m\geq m_0$. In the proof of \eqref{eq:new-bilinear}, we use the $L^1_v\times L^2_v\to L^2_v$ estimate for the gain operator of Maxwell molecules in \cite[Corollary~6]{ACG2010}. Since the constants in \eqref{eq:collision-missing}--\eqref{eq:new-linear-upper} are independent of $m$, we will choose $m=m_\eps(t,x)$ later and keep $\ell$ fixed.

We next introduce the velocity weight associated with the local Maxwellian $M$. In the energy estimates, we assume that the normalized macroscopic variables $(n,V,\Theta)$ defined in \eqref{eq:kinetic-normalization} belong to
\[\mathcal U:=\{(n,V,\Theta):3/4\leq n,\Theta\leq5/4,\ |V|\leq1/4\}.\]
For the physical velocity $v$, we define the centered normalized velocity and the corresponding weight by
\begin{equation}\label{eq:fitted-velocity-weight}
	\widetilde v:=\frac{v-u_f(t,x)}{\sqrt{T_f(t,x)}},\qquad W_{m,\ell}(t,x,v):=n(t,x)^{1/2}\Theta(t,x)^{-3/4}w_{m,\ell}(\widetilde v).
\end{equation}
Thus, $W_{m,\ell}^2\geq M$. By translation and dilation, Proposition~\ref{thm:new-collision} transfers to the local Maxwellian $M$, with constants uniform in $(n,V,\Theta)\in\mathcal U$. In the background velocity $\widehat v:=(v-u)/\theta$, the normalized collision frequency is
\begin{equation}\label{eq:normalized-collision-frequency}
	\widehat\nu_M(t,x,\widehat v):=\int|\widehat v-\widehat v_*|^\gamma\widetilde M(t,x,\widehat v_*)\dd\widehat v_*=n\Theta^{\gamma/2}\nu\left(\frac{\widehat v-V}{\sqrt\Theta}\right).
\end{equation}
We next choose $m$ as a function of $\theta(t,x)$ and $\eps$.

\begin{proposition}\label{thm:adaptive}
	Let $m_0\geq\ell$, $c_*>0$ and $\kappa_{\rm w}\geq4/c_*$. We define
	\begin{equation}\label{eq:adaptive-m}
		\Lambda_\eps(t,x):=\frac{\theta(t,x)^{\delta_{\rm gap}}}{\eps},\qquad m_\eps(t,x):=m_0+\kappa_{\rm w}\log(1+\Lambda_\eps(t,x)).
	\end{equation}
	Then on the whole interval $0\leq t<T$, we have
	\[\Lambda_\eps \mathrm e^{-c_*m_\eps/4}\leq \mathrm e^{-c_*m_0/4}.\]
	Moreover, for every fixed exponent $a\geq0$,
	\begin{equation}\label{eq:adaptive-absorption}
		\lim_{\eps\downarrow0} \sup_{0\leq t<T,\,x\in\R^3} \eps\theta^{-\delta_{\rm gap}}m_\eps^a=0.
	\end{equation}
\end{proposition}
In the remainder of the proof, we fix $\ell$ and take
\[W(t,x,v):=W_{m_\eps(t,x),\ell}(t,x,v),\]
where $m_\eps$ is given by \eqref{eq:adaptive-m}. By Proposition~\ref{thm:adaptive}, for any fixed $a\geq0$, we have
\[\theta^pm_\eps^a=o(1)\eps^{-1}\theta^{3+\gamma},\qquad \eps^{-1}\theta^{3+\gamma}\mathrm e^{-c_*m_\eps/4}\leq\mathrm e^{-c_*m_0/4}\theta^p,\]
where $o(1)\to0$ as $\eps\downarrow0$ uniformly in $(t,x)\in[0,T)\times\R^3$. See Subsection~\ref{subsec:velocity-weight-construction} for the proof of Proposition~\ref{thm:adaptive}.

\subsection{Finite-codimensional stability of the Euler linearization}\label{subsec:euler-linearization}
We will use the weighted energy estimates and the spectral decomposition for the linearized Euler operator in \cite{BCGJVY}. In addition, we also need to localize the unstable coordinates and construct a spatially localized correction to the Euler energy. See Subsections~\ref{subsec:unstable-localization} and \ref{subsec:localized-euler-correction}.

Let $(u,\theta)$ be the Euler solution constructed in Proposition~\ref{Thm.Euler-solution}. Recall that $u_f=u+\theta V$, $\theta_f=\theta(1+B)$ and $D_f=\partial_t+(u+\theta V)\cdot\nabla_x$. By \eqref{Eq.Euler-u-theta} and \eqref{eq:pressure-system}, we have
\begin{equation}\label{eq:normalized-macro-system}
	\begin{cases}
		D_fV+5\theta(1+B)\mathrm e^{3S/5}\nabla_xB=-(V\cdot\nabla_x)u+\tfrac13(\Div_xu)V-(V\cdot\nabla_x\theta)V\\
		\qquad\qquad\qquad\qquad\qquad\qquad\qquad\qquad-5\big((1+B)^2\mathrm e^{3S/5}-1\big)\nabla_x\theta-\frac{\Div_x\Pi}{\rho_f\theta},\\
		D_fB+\tfrac13\theta(1+B)\Div_xV=-\tfrac43(1+B)V\cdot\nabla_x\theta-\tfrac2{15}\theta^{-5}(1+B)^{-4}\mathcal C_f,\\
		D_fS=-\tfrac23\theta^{-5}(1+B)^{-5}\mathcal C_f,
	\end{cases}
\end{equation}
where $\mathcal C_f=\Pi:\nabla_xu_f+\Div_xq$. 

We next write \eqref{eq:normalized-macro-system} in self-similar coordinates. Let $y=x/(T-t)^{1/r}$. For a function $G=G(t,x)$, we denote
\[G^\sharp(t,y):=G(t,(T-t)^{1/r}y).\]
Throughout this subsection, the derivative $\partial_tG^\sharp$ is taken at fixed $y$. Recall from \eqref{eq:profile} that
\[(u^\sharp,\theta^\sharp)(t,y)=(T-t)^{-\alpha_{\rm ss}}(U(y),\vartheta(y)),\qquad |y|<R_{\rm e}(T-t)^{-1/r}.\]
We denote
\[b(y):=\frac yr+U(y),\qquad c_{\rm E}(y):=c_{\rm s}\vartheta(y),\qquad c_{\rm s}=\sqrt{\frac53}.\]
For $0\leq t<T$ and $|y|<R_{\rm e}(T-t)^{-1/r}$, we have
\[(T-t)(D_fG)^\sharp=(T-t)\partial_tG^\sharp+(b+\vartheta V^\sharp)\cdot\nabla_yG^\sharp.\]
We define
\begin{equation}\label{eq:euler-microscopic-source}
	\begin{aligned}
		F_V&:=-(T-t)\left(\frac{\Div_x\Pi}{\rho_f\theta}\right)^\sharp,\quad
		F_B:=-\tfrac2{15}(T-t)\big(\theta^{-5}(1+B)^{-4}\mathcal C_f\big)^\sharp,\\
		F_S&:=-\tfrac23(T-t)\big(\theta^{-5}(1+B)^{-5}\mathcal C_f\big)^\sharp,\qquad F_h:=(F_V,F_B,F_S).
	\end{aligned}
\end{equation}
By \eqref{eq:normalized-macro-system}, for $0\leq t<T$ and $|y|<R_{\rm e}(T-t)^{-1/r}$, we have
\begin{equation}\label{eq:profile-macro-system}
	\begin{cases}
		(T-t)\partial_tV^\sharp=-b\cdot\nabla_yV^\sharp-(V^\sharp\cdot\nabla_y)U+\tfrac13(\Div_yU)V^\sharp\\
		\qquad\qquad\qquad\qquad\qquad\qquad-5\vartheta\nabla_yB^\sharp-(10B^\sharp+3S^\sharp)\nabla_y\vartheta+N_V+F_V,\\
		(T-t)\partial_tB^\sharp=-b\cdot\nabla_yB^\sharp-\tfrac13\vartheta\Div_yV^\sharp-\tfrac43V^\sharp\cdot\nabla_y\vartheta+N_B+F_B,\\
		(T-t)\partial_tS^\sharp=-b\cdot\nabla_yS^\sharp+N_S+F_S,
	\end{cases}
\end{equation}
where the nonlinear terms $N:=(N_V,N_B,N_S)$ are given by
\begin{equation}\label{eq:nonlinear-N-explicit}
	\begin{aligned}
		N_V:={}&-\vartheta(V^\sharp\cdot\nabla_y)V^\sharp-(V^\sharp\cdot\nabla_y\vartheta)V^\sharp-5\big((1+B^\sharp)^2\mathrm e^{3S^\sharp/5}-1-2B^\sharp-3S^\sharp/5\big)\nabla_y\vartheta\\
		&-5\vartheta\big((1+B^\sharp)\mathrm e^{3S^\sharp/5}-1\big)\nabla_yB^\sharp,\\
		N_B:={}&-\vartheta V^\sharp\cdot\nabla_yB^\sharp-\tfrac43B^\sharp V^\sharp\cdot\nabla_y\vartheta-\tfrac13\vartheta B^\sharp\Div_yV^\sharp,\qquad
		N_S:=-\vartheta V^\sharp\cdot\nabla_yS^\sharp.
	\end{aligned}
\end{equation}

To symmetrize the acoustic terms in \eqref{eq:profile-macro-system} and to use the variables in \cite{BCGJVY}, we introduce the change of variables
\begin{equation}\label{eq:nonlinear-transfer}
	\begin{aligned}
		\mathsf T(y)&:=\vartheta(y)^3\operatorname{diag}\big(c_{\rm s}^{-1}I_3,3,-3/5\big),\qquad
		\mathsf W(t,y):=\mathsf T(y)\mathbf z^\sharp=(\mathsf W_V,\mathsf W_B,\mathsf W_S),\\
		&\qquad\mathsf W_V:=\vartheta^3V^\sharp/c_{\rm s},\qquad \mathsf W_B:=3\vartheta^3B^\sharp,\qquad \mathsf W_S:=-3\vartheta^3S^\sharp/5.
	\end{aligned}
\end{equation}
Since $\vartheta>0$ on $\R^3$, the matrix $\mathsf T(y)$ is invertible for all $y\in\R^3$. By the equation for the profile $\vartheta$, we have
\begin{equation}\label{eq:profile-transport}
	b\cdot\nabla_y\log\vartheta=-\alpha_{\rm ss}-\tfrac13\Div_yU.
\end{equation}
By direct computation,
\begin{align*}
	-\frac{\vartheta^3}{c_{\rm s}}b\cdot\nabla_yV^\sharp&=-b\cdot\nabla_y\mathsf W_V+3(b\cdot\nabla_y\log\vartheta)\mathsf W_V,\\
	-\frac{\vartheta^3}{c_{\rm s}}\big(5\vartheta\nabla_yB^\sharp+(10B^\sharp+3S^\sharp)\nabla_y\vartheta\big)&=-c_{\rm E}\nabla_y\mathsf W_B+(\mathsf W_B+3\mathsf W_S)\nabla_yc_{\rm E},\\
	3\vartheta^3\big(-\tfrac13\vartheta\Div_yV^\sharp-\tfrac43V^\sharp\cdot\nabla_y\vartheta\big)&=-c_{\rm E}\Div_y\mathsf W_V-\nabla_yc_{\rm E}\cdot\mathsf W_V.
\end{align*}
Multiplying \eqref{eq:profile-macro-system} by $\mathsf T(y)$ and using \eqref{eq:profile-transport}, we obtain
\begin{equation}\label{eq:actual-W-equation}
	\partial_t\mathsf W=\frac1{T-t}\big\{\mathcal L_{\rm E}\mathsf W+\mathsf T(N+F_h)\big\},\qquad |y|<R_{\rm e}(T-t)^{-1/r}.
\end{equation}
The linearized operator is given by
\begin{equation}\label{eq:explicit-euler-L}
	\begin{aligned}
		(\mathcal L_{\rm E}\mathsf W)_V:={}&-b\cdot\nabla_y\mathsf W_V-c_{\rm E}\nabla_y\mathsf W_B+(-3\alpha_{\rm ss}-\tfrac23\Div_yU)\mathsf W_V\\
		&\qquad-(\mathsf W_V\cdot\nabla_y)U+(\mathsf W_B+3\mathsf W_S)\nabla_yc_{\rm E},\\
		(\mathcal L_{\rm E}\mathsf W)_B:={}&-b\cdot\nabla_y\mathsf W_B-c_{\rm E}\Div_y\mathsf W_V+(-3\alpha_{\rm ss}-\Div_yU)\mathsf W_B-\nabla_yc_{\rm E}\cdot\mathsf W_V,\\
		(\mathcal L_{\rm E}\mathsf W)_S:={}&-b\cdot\nabla_y\mathsf W_S+(-3\alpha_{\rm ss}-\Div_yU)\mathsf W_S.
	\end{aligned}
\end{equation}
By the construction of $(u,\theta)$ in Proposition~\ref{Thm.Euler-solution}, equations \eqref{eq:profile-macro-system} and \eqref{eq:actual-W-equation} hold for $|y|<R_{\rm e}(T-t)^{-1/r}$. Linearizing \eqref{eq:actual-W-equation} at $\mathbf z=0$ with $F_h=0$, we obtain
\[(T-t)\partial_t\mathsf W=\mathcal L_{\rm E}\mathsf W.\]
The operator $\mathcal L_{\rm E}$ in \eqref{eq:explicit-euler-L} is defined on the whole space $\R^3_y$ and coincides with the linearized operator in \cite[(3.10)]{BCGJVY}, with $(\mathsf W_V,\mathsf W_B,\mathsf W_S)$ being the velocity, pressure and entropy components, respectively, and $c_{\rm E}=\theta_{\rm s}/(3r)$. To apply the spectral decomposition of $\mathcal L_{\rm E}$ to $\mathsf W$ in \eqref{eq:actual-W-equation}, we will introduce a spatial cutoff supported in $\{|y|<R_{\rm e}T^{-1/r}\}$. See Subsection~\ref{subsec:unstable-localization}.

\medskip
By the decay estimates \eqref{eq:profile-decay} for $U$ and $\vartheta$, we consider the following operator at infinity:
\[\mathcal L_\infty\mathsf W:=-\frac yr\cdot\nabla_y\mathsf W-3\alpha_{\rm ss}\mathsf W.\]
Let $D_*:=\operatorname{diag}(I_3,1,3/2)$ and $|Z|_{D_*}^2:=Z^{\mathsf T}D_*Z$ for $Z\in\R^5$. For $\varpi\in\R$ and $\mathsf W\in C_c^\infty(\R^3\setminus\{0\};\R^5)$, integrating by parts, we have
\begin{align*}
	\int |y|^\varpi\mathsf W^{\mathsf T}D_*\mathcal L_\infty\mathsf W\dd y&=-\frac1{2r}\int |y|^\varpi y\cdot\nabla_y|\mathsf W|_{D_*}^2\dd y-3\alpha_{\rm ss}\int |y|^\varpi|\mathsf W|_{D_*}^2\dd y\\
	&=\frac{\varpi+3-6r\alpha_{\rm ss}}{2r}\int |y|^\varpi|\mathsf W|_{D_*}^2\dd y.
\end{align*}
Hence we define the critical exponent $\bar\varpi$ and the decay rate $\lambda_\varpi$ by
\begin{equation}\label{eq:euler-critical-exponent}
	\bar\varpi:=6r\alpha_{\rm ss}-3=-3+6(r-1),\qquad \lambda_\varpi:=\frac{\bar\varpi-\varpi}{4r}>0\quad(\varpi<\bar\varpi).
\end{equation}
For the estimates of $\mathcal L_{\rm E}$, we replace $|y|^\varpi$ by $\langle y\rangle^\varpi$. To estimate the derivatives of order $2k$, we introduce the weight $\phi_1^{4k}$, where $\phi_1$ is the smooth positive radial function constructed in \cite[Lemma 4.1]{BCGJVY}, satisfying $\phi_1(y)\asymp\langle y\rangle$.

For an integer $k$ and $\varpi\in[-6,\bar\varpi]$, we define
\begin{equation}\label{eq:actual-X-norm}
	\|\mathsf W\|_{X^{2k}_\varpi}^2:=\int\left(\phi_1^{4k}|\Delta_y^k\mathsf W|_{D_*}^2+\varkappa_{k,\varpi}|\mathsf W|_{D_*}^2\right)\langle y\rangle^\varpi\dd y,
\end{equation}
where the constant $\varkappa_{k,\varpi}>0$ will be chosen in Proposition~\ref{prop:euler-coercivity}. Let $X^{2k}_\varpi$ be the completion of the space of smooth compactly supported radial fields under the norm \eqref{eq:actual-X-norm}, where the radial symmetry of $\mathsf W$ means that
\[\mathsf W_V(Oy)=O\mathsf W_V(y),\qquad \mathsf W_B(Oy)=\mathsf W_B(y),\qquad \mathsf W_S(Oy)=\mathsf W_S(y),\qquad O\in SO(3).\]
We denote the corresponding inner product by $\langle\cdot,\cdot\rangle_{X^{2k}_\varpi}$. For the complexification of $X^{2k}_\varpi$, we use $|Z|_{D_*}^2:=\overline Z^{\mathsf T}D_*Z$ for $Z\in\C^5$. We also use
\begin{equation}\label{eq:euler-zero-order-inner-product}
	\|\mathsf W\|_{X^0_\varpi}^2:=\int|\mathsf W|_{D_*}^2\langle y\rangle^\varpi\dd y,\qquad \langle\mathsf W,\mathsf Y\rangle_0:=\int\mathsf W^{\mathsf T}D_*\mathsf Y\dd y.
\end{equation}
The norms $\|\cdot\|_{X^{2k}_\varpi}$ and $\|\cdot\|_{X^0_\varpi}$ coincide with \cite[(4.6)]{BCGJVY}, with the three numerical coefficients written as the matrix $D_*$. For fixed $k$ and $\varpi$, we have
\begin{equation}\label{eq:euler-weighted-derivatives}
	\|\mathsf W\|_{X^{2k}_\varpi}^2\asymp\sum_{j=0}^{2k}\int\langle y\rangle^{\varpi+2j}|\nabla_y^j\mathsf W|^2\dd y.
\end{equation}
We verify \eqref{eq:euler-weighted-derivatives} in Section~\ref{subsec:unstable-localization}. In particular, the weighted norm controls the usual $H^{2k}$ norm on every fixed ball.

The construction of $\phi_1$ relies on the profile repulsivity \eqref{eq:normalized-repulsivity}. Condition \eqref{eq:normalized-repulsivity} gives a favorable sign for the highest-order acoustic terms in the bounded spatial region, while the outward transport and the profile decay control the exterior region. The resulting weight inequality is \cite[(4.2a)]{BCGJVY}. The calculation in \cite[Section 4.1]{BCGJVY} leaves only a localized zeroth-order error:
\begin{proposition}[Weighted dissipativity up to a localized error]\label{prop:euler-coercivity}
	There exists an integer $k_0\geq6$ depending only on the fixed profiles and $\phi_1$ with the following property. For every $\varpi\in[-6,\bar\varpi)$, there exists $R_\varpi>0$ such that for every $k\geq k_0$, we can choose $\varkappa_{k,\varpi}>0$ and $C_{k,\varpi}>0$ such that
	\begin{equation}\label{eq:euler-subcritical-coercivity}
		\langle\mathcal L_{\rm E}\mathsf W,\mathsf W\rangle_{X^{2k}_\varpi}\leq-\lambda_\varpi\|\mathsf W\|_{X^{2k}_\varpi}^2+C_{k,\varpi}\int_{|y|\leq R_\varpi}|\mathsf W|_{D_*}^2\dd y.
	\end{equation}
	For $\varpi=\bar\varpi$ and $k\geq k_0$, we can choose $\varkappa_{k,\bar\varpi}>0$ and $C_k>0$ such that
	\begin{equation}\label{eq:critical-Euler-input}
		\langle\mathcal L_{\rm E}\mathsf W,\mathsf W\rangle_{X^{2k}_{\bar\varpi}}\leq C_k\int\langle y\rangle^{\bar\varpi-r}|\mathsf W|^2\dd y.
	\end{equation}
	The inequalities \eqref{eq:euler-subcritical-coercivity} and \eqref{eq:critical-Euler-input} hold for smooth compactly supported radial fields and extend to the corresponding generator domains. All constants depend only on the fixed profiles, $\phi_1$, $\varpi$ and $k$, with $R_\varpi$ independent of $k$.
	\end{proposition}
	
	The inequalities \eqref{eq:euler-subcritical-coercivity} and \eqref{eq:critical-Euler-input} are given in \cite[Theorem 4.2, (4.5a)--(4.5c)]{BCGJVY}. For the construction of the weight $\phi_1$, see \cite[Lemma 4.1]{BCGJVY} and \cite[Lemma 3.1 and Section 3.1]{CCSV2024}.
	
	\medskip
	
	By \eqref{eq:euler-weighted-derivatives} and the Rellich compactness theorem, the restriction map
	\[\mathsf W\mapsto\mathsf W|_{B_{R_\varpi}}:X^{2k}_\varpi\to L^2(B_{R_\varpi};\R^5)\]
	is compact. Using \eqref{eq:euler-subcritical-coercivity}, we next decompose $\mathcal L_{\rm E}$ into a dissipative operator and a compact operator.
\begin{proposition}\label{prop:euler-compact}
	For $k\geq k_0$ and $\varpi\in[-6,\bar\varpi)$, there exists a real bounded linear map
	\[\mathcal K_{{\rm E},k,\varpi}:X^0_\varpi\longrightarrow X^{2k+6}_\varpi\]
	whose restriction to $X^{2k}_\varpi$ is compact as an operator from $X^{2k}_\varpi$ to itself. We define
	\[\mathcal D_{{\rm E},k,\varpi}:=\mathcal L_{\rm E}-\mathcal K_{{\rm E},k,\varpi},\qquad \mathcal D_{2k,\varpi}:=\{\mathsf W\in X^{2k}_\varpi:\mathcal L_{\rm E}\mathsf W\in X^{2k}_\varpi\},\]
	where $\mathcal L_{\rm E}\mathsf W$ is evaluated in the sense of distributions. On the domain $\mathcal D_{2k,\varpi}$, we have
	\begin{equation}\label{eq:euler-dissipative-part}
		\langle\mathcal D_{{\rm E},k,\varpi}\mathsf W,\mathsf W\rangle_{X^{2k}_\varpi}\leq-\lambda_\varpi\|\mathsf W\|_{X^{2k}_\varpi}^2.
	\end{equation}
	Both $\mathcal D_{{\rm E},k,\varpi}$ and $\mathcal L_{\rm E}$ generate strongly continuous semigroups on $X^{2k}_\varpi$ and on its complexification. Moreover, for $q\geq0$,
	\[\|\mathrm e^{q\mathcal D_{{\rm E},k,\varpi}}\|_{X^{2k}_\varpi\to X^{2k}_\varpi}\leq\mathrm e^{-\lambda_\varpi q},\qquad \|\mathrm e^{q\mathcal L_{\rm E}}\|_{X^{2k}_\varpi\to X^{2k}_\varpi}\leq\mathrm e^{C_{k,\varpi}q}.\]
\end{proposition}
The compact operator and \eqref{eq:euler-dissipative-part} are given in \cite[Proposition 4.6, (4.18)]{BCGJVY}, which refers to the proof of \cite[Proposition 3.4, Section 3.3]{CCSV2024}; there the compact operator is constructed by the Riesz representation, and its regularity and compactness are proved. For the generator properties, see \cite[Proposition 4.7, (4.24)--(4.26)]{BCGJVY}; see also \cite[Section 3.4]{CCSV2024} for the construction of the linear evolution.

We fix an even integer $K:=2k\geq\max\{14,2k_0+2\}$ and set
\begin{equation}\label{eq:euler-stability-space}
	\varpi:=\bar\varpi-\frac r2,\qquad \mathscr X:=X^K_\varpi,\qquad \lambda_\varpi=\frac18.
\end{equation}
We consider the spectrum of $\mathcal L_{\rm E}$ on the complexification of $\mathscr X$ with the complexified domain $\mathcal D_{K,\varpi}$. By the compact operator decomposition and \eqref{eq:euler-dissipative-part}, for every $\delta>0$, the set
\[\sigma(\mathcal L_{\rm E})\cap\left\{\zeta\in\C:\operatorname{Re}\zeta\geq-\frac18+\delta\right\}\]
consists of finitely many eigenvalues of finite algebraic multiplicity. Hence we can choose $d_*$ satisfying
\[0<d_*<\min\{\eta,1/4\},\qquad \sigma(\mathcal L_{\rm E})\cap\{\zeta\in\C:\operatorname{Re}\zeta=-d_*/2\}=\varnothing.\]
\begin{proposition}[Finite-codimensional stability]\label{prop:euler-spectrum}
	Let
	\[\sigma_{\rm u}:=\{\zeta\in\sigma(\mathcal L_{\rm E}):\operatorname{Re}\zeta>-d_*/2\}.\]
	The set $\sigma_{\rm u}$ is finite, and every $\zeta\in\sigma_{\rm u}$ is an eigenvalue of finite algebraic multiplicity. Let $\Gamma_{\rm u}$ be a finite union of positively oriented contours in the resolvent set enclosing precisely $\sigma_{\rm u}$, and we define
	\begin{equation}\label{eq:euler-riesz-projection}
		\Pu:=\frac1{2\pi\mathrm i}\int_{\Gamma_{\rm u}}(\zeta-\mathcal L_{\rm E})^{-1}\dd\zeta,\qquad E_{\rm u}:=\Pu\mathscr X,\qquad E_{\rm s}:=(I-\Pu)\mathscr X.
	\end{equation}
	If $\sigma_{\rm u}$ is empty, we set $\Pu:=0$. The projection $\Pu$ is bounded, commutes with the Euler evolution, and preserves real-valued fields. The space $E_{\rm u}$ is finite dimensional and consists of smooth fields. The space $E_{\rm s}$ is closed, invariant, and has codimension $\dim E_{\rm u}$. Thus, $\mathscr X=E_{\rm u}\oplus E_{\rm s}$. Moreover, there exist constants $\sigma_*>d_*/2$, $\kappa<d_*/2$ and $C>0$ such that
	\begin{equation}\label{eq:euler-spectral-evolution}
		\begin{aligned}
			\|\mathrm e^{q\mathcal L_{\rm E}}F\|_{\mathscr X}&\leq C\mathrm e^{-\sigma_*q}\|F\|_{\mathscr X},&&F\in E_{\rm s},\quad q\geq0,\\
			\|\mathrm e^{-q\mathcal L_{\rm E}}F\|_{\mathscr X}&\leq C\mathrm e^{\kappa q}\|F\|_{\mathscr X},&&F\in E_{\rm u},\quad q\geq0.
		\end{aligned}
	\end{equation}
	Moreover, there holds $E_{\rm u}\subset \bigcap_{j\geq k_0}X^{2j}_\varpi$.
\end{proposition}
The spectral decomposition and the exponential estimates in Proposition~\ref{prop:euler-spectrum} follow from \cite[Section 4.3.3, (4.27)--(4.31)]{BCGJVY}; the smoothness follows from \cite[Lemma 4.8, Section 4.3.4]{BCGJVY}. The underlying operator-theoretic argument and the regularity argument are given in \cite[Sections 3.5--3.6]{CCSV2024}. The strict separation from $\operatorname{Re}\zeta=-d_*/2$ gives the strict inequalities for $\sigma_*$ and $\kappa$ in \eqref{eq:euler-spectral-evolution}. We point out that the term unstable refers here to modes that fail the prescribed decay rate $d_*/2$, which may include modes with a small negative real part. Note that the projection need not be orthogonal.

\medskip

For $0\leq t_0\leq t<T$, we define
\[\mathcal U_E(t,t_0):=\exp\left(\log\frac{T-t_0}{T-t}\,\mathcal L_{\rm E}\right).\]
Then $(T-t)\partial_t\mathcal U_E(t,t_0)F=\mathcal L_{\rm E}\mathcal U_E(t,t_0)F$. In particular, \eqref{eq:euler-spectral-evolution} becomes
\begin{align}
	\|\mathcal U_E(t,t_0)F\|_{\mathscr X}\leq C\left(\frac{T-t}{T-t_0}\right)^{\sigma_*}\|F\|_{\mathscr X},\qquad F\in E_{\rm s},\nonumber\\
	\|\mathcal U_E(t_0,t)F\|_{\mathscr X}\leq C\left(\frac{T-t_0}{T-t}\right)^\kappa\|F\|_{\mathscr X},\qquad F\in E_{\rm u}.\label{eq:exceptional-backward-t}
\end{align}
For sufficiently large $R_{\rm sp}$, the outgoing-support estimate \cite[Proposition 4.9, Section 4.3.5]{BCGJVY} also gives
\begin{equation}\label{eq:exterior-decay-t}
	\|\mathcal U_E(t,t_0)F\|_{\mathscr X}\leq C\left(\frac{T-t}{T-t_0}\right)^{1/8}\|F\|_{\mathscr X},\qquad \supp F\subset\{|y|>R_{\rm sp}\}.
\end{equation}
Therefore, exterior data decay faster than the possible backward growth on $E_{\rm u}$, since $\kappa<d_*/2<1/8$. Comparing \eqref{eq:exceptional-backward-t} with \eqref{eq:exterior-decay-t}, we get the vanishing of the projection stated below.
\begin{proposition}[Compactly supported functionals for the unstable projection]\label{prop:projection-local}
	Every $F\in\mathscr X$ supported in $\{|y|>R_{\rm sp}\}$ satisfies $\Pu F=0$. For any radial cutoff $\chi_{\rm sp}\in C_c^\infty(\R^3)$ that equals one on a neighborhood of $\overline B_{R_{\rm sp}}$, we have
	\begin{equation}\label{eq:euler-projection-local}
		\Pu F=\Pu(\chi_{\rm sp}F),\qquad F\in\mathscr X.
	\end{equation}
	Let $\Phi_1,\ldots,\Phi_J$ be a real basis of $E_{\rm u}$, where $J:=\dim E_{\rm u}$ and each $\Phi_j$ is a five-component field on $\R^3$. Then there exist continuous real linear functionals $\ell_1,\ldots,\ell_J$ on $\mathscr X$ such that
	\[\Pu F=\sum_{j=1}^J\ell_j(F)\Phi_j,\qquad \ell_j(F)=\ell_j(\chi_{\rm sp}F).\]
	
	For every smooth compactly supported radial field $F$,
	\[\ell_j(F)=0\quad\text{if }\supp F\cap\overline B_{R_{\rm sp}}=\varnothing,\qquad j=1,\ldots,J.\]
	For every fixed $R>R_{\rm sp}$, the coordinate functionals also satisfy
	\begin{equation}\label{eq:unstable-functional-local-bound}
		|\ell_j(F)|\leq C_{j,R}\|F\|_{H^K(B_R)},\qquad F\in\mathscr X.
	\end{equation}
	In particular, if $F_1=F_2$ on a neighborhood of $\overline B_{R_{\rm sp}}$, then $\ell_j(F_1)=\ell_j(F_2)$ for $j=1,\ldots,J$. 
\end{proposition}

The proof of Proposition~\ref{prop:projection-local} in Section~\ref{subsec:unstable-localization} combines the backward estimate \eqref{eq:exceptional-backward-t} on $E_{\rm u}$ with the exterior decay estimate \eqref{eq:exterior-decay-t}.

For the source estimate, we also use the realization of $\mathcal L_{\rm E}$ on $X^{K-2}_\varpi$. We denote its unstable Riesz projection by $\Pu^{(K-2)}$, and write $\Pu^{(K)}:=\Pu$ for the projection on $X^K_\varpi$. In Lemma~\ref{lem:new-spectral-consistency}, we prove
\begin{equation}\label{eq:two-order-projection-agreement}
	\operatorname{Ran}\Pu^{(K-2)}=\operatorname{Ran}\Pu^{(K)}=E_{\rm u},\qquad \Pu^{(K-2)}F=\Pu^{(K)}F\quad(F\in X^K_\varpi).
\end{equation}
Relative to the same basis $\Phi_1,\ldots,\Phi_J$, we define the continuous coordinate functionals on $X^{K-2}_\varpi$ by
\[\Pu^{(K-2)}G=\sum_{j=1}^J\ell_j^{(K-2)}(G)\Phi_j,\qquad G\in X^{K-2}_\varpi.\]
Then $\ell_j^{(K-2)}(F)=\ell_j(F)$ for $F\in X^K_\varpi$. 

\medskip

We fix the real basis $\Phi_1,\ldots,\Phi_J$ from Proposition~\ref{prop:projection-local}. The matrix $A_{\rm u}\in\R^{J\times J}$ is defined by
\[\mathcal L_{\rm E}\Phi_j=\sum_{i=1}^J(A_{\rm u})_{ij}\Phi_i,\qquad j=1,\ldots,J.\]
For any fixed field $\mathsf W\in\mathscr X$, we define its coordinates by
\[\mathbf a:=(\ell_1(\mathsf W),\ldots,\ell_J(\mathsf W))^{\mathsf T},\qquad \Pu\mathsf W=\sum_{j=1}^J a_j\Phi_j.\]
We choose the decay constants such that
\begin{equation}\label{eq:new-spectral-choices}
	2\kappa<d_*,\qquad d_*/2<h_*<\min\{\sigma_*,1/8\},\qquad d_*<\lambda<\lambda_E<2h_*<1/4.
\end{equation}
If $E_{\rm u}=\{0\}$, we set $J:=0$, $\mathbf a:=0$, $|\mathbf a|_P:=0$ and $\kappa:=-1$. For $J>0$, let $B_{\rm u}:=A_{\rm u}+d_*I_J/2$ and define
\begin{equation}\label{eq:unstable-coordinate-metric}
	P:=\int_0^\infty\mathrm e^{-qB_{\rm u}^{\mathsf T}}\mathrm e^{-qB_{\rm u}}\dd q,\qquad \langle a_1,a_2\rangle_P:=a_1^{\mathsf T}Pa_2,\qquad |a_1|_P^2:=\langle a_1,a_1\rangle_P^2{}^{\!0}\langle a_1,a_1\rangle_P.
\end{equation}
By $\sigma(B_{\rm u})\subset\{\zeta\in\C:\operatorname{Re}\zeta>0\}$, the integral in \eqref{eq:unstable-coordinate-metric} converges and $P$ is symmetric positive definite. Moreover, we have
\[B_{\rm u}^{\mathsf T}P+PB_{\rm u}=-\int_0^\infty\frac{\dd}{\dd q}\left(\mathrm e^{-qB_{\rm u}^{\mathsf T}}\mathrm e^{-qB_{\rm u}}\right)\dd q=I_J.\]
Consequently,
\begin{equation}\label{eq:unstable-coordinate-lower-bound}
	2\langle A_{\rm u}\mathbf a,\mathbf a\rangle_P=|\mathbf a|^2-d_*|\mathbf a|_P^2\geq(\|P\|^{-1}-d_*)|\mathbf a|_P^2,
\end{equation}
where $\|P\|$ denotes the operator norm of $P$ on $\R^J$. For $J>0$, replacing $\kappa$ by $\max\{\kappa,(d_*-\|P\|^{-1})/2\}$, we have
\[2\langle A_{\rm u}\mathbf a,\mathbf a\rangle_P\geq-2\kappa|\mathbf a|_P^2,\qquad 2\kappa<d_*.\]
For both $J=0$ and $J>0$, we fix $\delta$ satisfying
\[\max\{0,2\kappa\}<\delta<d_*.\]
The constants $d_*,h_*,\sigma_*,\lambda,\lambda_E$ are unchanged. We fix $\kappa$ and $\delta$ before choosing $T$ and $\eps$.

\medskip

We next construct a quadratic correction to the Euler energy to replace the local term $C_{k,\varpi}\int_{|y|\leq R_\varpi}|\mathsf W|_{D_*}^2\dd y$ in \eqref{eq:euler-subcritical-coercivity} by $C|\mathbf a|_P^2$.
\begin{proposition}[A localized correction to the Euler energy]\label{thm:local-observation}
	There exist a bounded symmetric nonnegative bilinear form $\mathcal R_*$ on $\mathscr X$ and constants $h_*,c_*,C_*>0$ independent of $t$ and $\eps$ such that
	\[\mathcal E_*(\mathsf W):=\tfrac12\big(\|\mathsf W\|_{\mathscr X}^2+\mathcal R_*(\mathsf W,\mathsf W)\big)\]
	is equivalent to $\|\mathsf W\|_{\mathscr X}^2$. For every real $\mathsf W\in\mathcal D_{K,\varpi}$, with the notation $\mathbf a:=(\ell_1(\mathsf W),\ldots,\ell_J(\mathsf W))^{\mathsf T}$, we have
	\begin{equation}\label{eq:observation-stable-sign}
		\begin{aligned}
			\langle\mathcal L_{\rm E}\mathsf W,\mathsf W\rangle_{\mathscr X}+\mathcal R_*(\mathcal L_{\rm E}\mathsf W,\mathsf W) \leq{}&-2h_*\mathcal E_*(\mathsf W)-c_*\|\mathsf W\|_{\mathscr X}^2+C_*|\mathbf a|^2.
		\end{aligned}
	\end{equation}
	Moreover, there exist a fixed ball $B_{R_c}$ and a real operator $\mathscr R_*$ such that, for $\mathsf W,\mathsf Y\in\mathscr X$,
	\begin{equation}\label{eq:observation-local-map}
		\begin{aligned}
			\mathcal R_*(\mathsf W,\mathsf Y)&=\langle\mathscr R_*\mathsf W,\mathsf Y\rangle_0,\qquad \supp(\mathscr R_*\mathsf W)\subset B_{R_c},\\
			\|\mathscr R_*\mathsf W\|_{H^j}&\leq C_j\|\mathsf W\|_{H^j(B_{R_c})},\qquad j=0,1.
		\end{aligned}
	\end{equation}
	After increasing the fixed radius $R_c$, we can require $\supp(\mathscr R_*\mathsf W)\subset B_{R_c/2}$ for every $\mathsf W\in\mathscr X$. The locality estimate in \eqref{eq:observation-local-map} also gives
	\[\mathscr R_*\mathsf W_1=\mathscr R_*\mathsf W_2\quad\text{if }\mathsf W_1=\mathsf W_2\text{ on }B_{R_c}.\]
	The quadratic correction requires only local $L^2$ control:
	\begin{equation}\label{eq:observation-zero-order-bound}
		0\leq\mathcal R_*(\mathsf W,\mathsf W)\leq C\|\mathsf W\|_{L^2(B_{R_c})}^2.
	\end{equation}
\end{proposition}
Proposition~\ref{thm:local-observation} is proved in Section~\ref{subsec:localized-euler-correction} by integrating a localized Euler energy over a finite time interval and using the stable decay. The weighted energy estimate uses \eqref{eq:observation-stable-sign}. The evolution of the localized unstable coordinates is derived in Section~\ref{subsec:unstable-localization}.

For clarity, Proposition~\ref{prop:euler-coercivity} is the weighted dissipativity estimate of \cite[Theorem 4.2, (4.5a)--(4.5c)]{BCGJVY}. Proposition~\ref{prop:euler-compact} uses the compact perturbation and linear evolution statements of \cite[Propositions 4.6 and 4.7]{BCGJVY}. Proposition~\ref{prop:euler-spectrum} uses the finite-codimensional spectral decomposition of \cite[Section 4.3.3, (4.27)--(4.31)]{BCGJVY}, together with the regularity result of \cite[Lemma 4.8]{BCGJVY}. We do not repeat the proofs of these three propositions. Propositions~\ref{prop:projection-local} and \ref{thm:local-observation} are proved in Sections~\ref{subsec:unstable-localization} and \ref{subsec:localized-euler-correction}, respectively. The identities in \eqref{eq:two-order-projection-agreement} are proved in Lemma~\ref{lem:new-spectral-consistency}.

To apply the Euler results to the Boltzmann solution, we fix $R_2>R_1>\max\{R_{\rm sp},R_c\}$ and a smooth radial cutoff $0\leq\chi_{\rm sp}\leq1$ with
\begin{equation}\label{eq:fixed-spectral-cutoff}
	\chi_{\rm sp}=1\quad\text{on }B_{R_1},\qquad \supp\chi_{\rm sp}\subset B_{R_2},\qquad R_2T^{1/r}<R_{\rm e}/2.
\end{equation}
The radii and $\chi_{\rm sp}$ are fixed before the sufficiently small time $T$ is chosen. For every $0\leq t<T$, the physical points $x=(T-t)^{1/r}y$ with $y\in\supp\chi_{\rm sp}$ lie in $B_{R_{\rm e}/2}$. Hence, \eqref{eq:actual-W-equation} holds on $\supp\chi_{\rm sp}$. Since $\mathcal L_{\rm E}$ is a first-order differential operator, we have
\[\supp\big([\mathcal L_{\rm E},\chi_{\rm sp}]\mathsf W\big)\subset\{R_1\leq|y|\leq R_2\},\qquad [\mathcal L_{\rm E},\chi_{\rm sp}]\mathsf W:=\mathcal L_{\rm E}(\chi_{\rm sp}\mathsf W)-\chi_{\rm sp}\mathcal L_{\rm E}\mathsf W.\]
Therefore, by Proposition~\ref{prop:projection-local}, we have $\Pu([\mathcal L_{\rm E},\chi_{\rm sp}]\mathsf W)=0$. In Section~\ref{subsec:unstable-localization}, we compute the commutator and derive the equation for the coordinates of $\Pu(\chi_{\rm sp}\mathsf W(t,\cdot))$.

For the actual macroscopic perturbation $\mathsf W$ from \eqref{eq:nonlinear-transfer}, we define the unstable coordinates by
\begin{equation}\label{eq:roadmap-unstable-coordinates}
	\mathbf a(t):=\big(\ell_1(\chi_{\rm sp}\mathsf W(t,\cdot)),\ldots,\ell_J(\chi_{\rm sp}\mathsf W(t,\cdot))\big)^{\mathsf T}.
\end{equation}
When $J=0$, we set $\mathbf a(t):=0$. The evolution equation for $\mathbf a(t)$ is derived in Section~\ref{subsec:unstable-localization}.

The source $F_h:=(F_V,F_B,F_S)$ in \eqref{eq:euler-microscopic-source} contains the first spatial derivatives of the stress $\Pi$ and the heat flux $q$. In Proposition~\ref{thm:actual-spectral-forcing}, we estimate $\chi_{\rm sp}\mathsf T(N+F_h)$ in $X^{K-2}_\varpi$, using at most $K-1$ spatial derivatives of the microscopic moments. The energy $\mathcal E_*$ is defined on $X^K_\varpi$, whereas the source term $\chi_{\rm sp}\mathsf T(N+F_h)$ is estimated in $X^{K-2}_\varpi$. By \eqref{eq:two-order-projection-agreement}, we have
\[\ell_j^{(K-2)}(F)=\ell_j(F),\qquad F\in X^K_\varpi,\quad j=1,\ldots,J.\]
Thus, in the equation for $\mathbf a$, we use $\ell_j^{(K-2)}(\chi_{\rm sp}\mathsf T(N+F_h))$ as the source coordinates in the basis $\Phi_1,\ldots,\Phi_J$.

\subsection{The weighted energy estimates}\label{subsec:weighted-energies}
We seek an energy controlling the macroscopic perturbation in the norm \eqref{eq:actual-X-norm} and the microscopic perturbation in the velocity weights \eqref{eq:fitted-velocity-weight}. To avoid the loss of one spatial derivative, we combine the highest-order terms in the macroscopic and microscopic equations before performing the energy estimates.
\medskip

Let $M:=\mathcal M(\rho_f,u_f,T_f)$ be the local Maxwellian determined by the physical density $\rho_f$, the mean velocity $u_f$ and the temperature $T_f$ in \eqref{eq:physical-parameters}. The remainder $H:=f-M$ in \eqref{eq:microscopic-remainders} satisfies the five zero-moment identities \eqref{eq:zero-moments}. The stress tensor $\Pi$ and the heat flux $q$ are given by the velocity integrals in \eqref{eq:roadmap-stress-heat}.
Recall the moment equations \eqref{eq:moment-system}.The precise system for $\mathbf z:=(V,B,S)$ is \eqref{eq:normalized-macro-system}. At spatial derivative order $K$, the terms $\Div\Pi$ and $\Div q$ in \eqref{eq:moment-system} contain derivatives of $H$ of order $K+1$. 

We fix $(t,x)$ with $\rho_f(t,x)>0$ and $T_f(t,x)>0$, and regard $M=M(t,x,v)$ as a function of the physical velocity $v$. We define a projection operator $\cP_M$ acting on functions $g=g(v)$ of the physical velocity. Assume that $\int_{\R^3}(1+|v|^2)|g(v)|\dd v<\infty$. We define its five moments by
\begin{equation*}
	\begin{aligned}
		m_0[g]:=\int_{\R^3}g\dd v,\quad
		m_1[g]:=\int_{\R^3}(v-u_f)g\dd v\in\R^3,\quad
		m_2[g]:=\int_{\R^3}(|v-u_f|^2-3T_f)g\dd v.
	\end{aligned}
\end{equation*}
We define $\cP_Mg$ by the explicit formula
\begin{equation}\label{eq:roadmap-projection-definition}
	\mathcal P_Mg:=M\left(\frac{m_0[g]}{\rho_f}+\frac{(v-u_f)\cdot m_1[g]}{\rho_fT_f}+\frac{|v-u_f|^2-3T_f}{6\rho_fT_f^2}m_2[g]\right).
\end{equation}
Differentiating the Maxwellian with respect to its physical parameters, we get
\begin{equation}\label{eq:maxwellian-parameter-differential}
	\partial_{\rho_f}M=\frac{M}{\rho_f},\qquad\partial_{(u_f)_j}M=\frac{(v_j-(u_f)_j)M}{T_f},\qquad\partial_{T_f}M=\frac{|v-u_f|^2-3T_f}{2T_f^2}M.
\end{equation}
In particular, $\mathcal P_M|_{L^2(M^{-1}\,\mathrm dv)}$ is the orthogonal projection in $L^2(\R_v^3, M^{-1}\,\mathrm dv)$ onto the five-dimensional space spanned by the five partial derivatives in \eqref{eq:maxwellian-parameter-differential}. The moments in \eqref{eq:roadmap-projection-definition} satisfy
\begin{equation}\label{eq:projected-moment-preservation}
	\int (1,v,|v|^2)(I-\mathcal P_M)g\dd v=0.
\end{equation}

We write $D:=\partial_t+v\cdot\nabla_x$, with the physical velocity held fixed as in \eqref{eq:physical-derivative-convention}. By the Boltzmann equation \eqref{Eq.Boltzmann} and $f=M+H$, we have $DH=\eps^{-1}Q(f,f)-DM$. The microscopic identity is derived from the moment projection \eqref{eq:roadmap-projection-definition} and the stress and heat flux in \eqref{eq:roadmap-stress-heat}.
\begin{lemma}[The microscopic equation]\label{lem:roadmap-projection-microscopic}
	Let $f$ be a smooth solution to \eqref{Eq.Boltzmann} with $\rho_f,T_f>0$. Assume that $f$ is continuous in $(t,x)$ with values in $L^1(\R^3_v,\langle v\rangle^{2+\gamma}\dd v)$, and that $\partial_tf,\nabla_xf$ are continuous with values in $L^1(\R^3_v,\langle v\rangle^3\dd v)$. We define $M:=\mathcal M(\rho_f,u_f,T_f)$ by \eqref{eq:physical-parameters} and $H:=f-M$ by \eqref{eq:microscopic-remainders}. Then
	\begin{equation}\label{eq:roadmap-microscopic-equation}
		DH=\eps^{-1}Q(f,f)-(I-\mathcal P_M)(v\cdot\nabla_xM)+\mathcal P_M(v\cdot\nabla_xH).
	\end{equation}
	More explicitly, with $\widetilde v:=(v-u_f)/\sqrt{T_f}$ and $S_M:=(I-\mathcal P_M)(v\cdot\nabla_xM)$, $T_H:=\mathcal P_M(v\cdot\nabla_xH)$,
	\begin{align}\label{eq:roadmap-projected-sources}
		S_M={}&M\left[\left(\widetilde v\otimes\widetilde v-\frac13|\widetilde v|^2\Id\right):\nabla_xu_f+(|\widetilde v|^2-5)\widetilde v\cdot\nabla_x\sqrt{T_f}\right],\notag\\
		T_H={}&M\left[\frac{(v-u_f)\cdot\Div\Pi}{\rho_fT_f}+\frac{(|v-u_f|^2-3T_f)(\Pi:\nabla_xu_f+\Div q)}{3\rho_fT_f^2}\right].
	\end{align}
\end{lemma}
Lemma~\ref{lem:roadmap-projection-microscopic} is proved in Section~\ref{subsec:moment-microscopic} by using the moment identities and the Gaussian integration.

The physical parameters $(\rho_f,u_f,T_f)$ in \eqref{eq:physical-parameters} are precisely the density, the three components of the mean velocity and the temperature; equivalently, they can be expressed through the five entries of $\mathbf z:=(V,B,S)$ by \eqref{eq:roadmap-macro} and \eqref{eq:normalized-reconstruction}. Formula \eqref{eq:roadmap-projected-sources} shows that $S_M$ contains the first spatial derivatives of $u_f$ and $T_f$, but no derivative of $\rho_f$. After taking $K$ spatial derivatives of \eqref{eq:roadmap-microscopic-equation}, $\nabla_xu_f$ and $\nabla_xT_f$ produce derivatives of order $K+1$. The term $T_H$ contains $\Div\Pi$ and $\Div q$, whose $K$-th derivatives would similarly involve $K+1$ spatial derivatives of $H$. 

This shows the loss of a spatial derivative. Therefore, the differentiated moment equations \eqref{eq:moment-system} and the differentiated microscopic equation \eqref{eq:roadmap-microscopic-equation} must be estimated together, so that the cancellation of their highest-order terms is preserved.

\medskip\noindent\emph{The macroscopic Euler energy.}
We keep the even derivative order $K:=2k\geq\max\{14,2k_0+2\}$ chosen for the Euler space \eqref{eq:euler-stability-space}, and the critical exponent $\bar\varpi:=-3+6(r-1)$ in \eqref{eq:euler-critical-exponent}. The positive radial weight $\phi_1$ and the coefficient $\varkappa_{k,\varpi}$ occur in the Euler norm \eqref{eq:actual-X-norm}. We set $a:=r/2$ and $b\in\{0,a\}$.\footnote{The scalar $a$ is unrelated to the unstable-coordinate vector $\mathbf a(t)$ defined in \eqref{eq:roadmap-unstable-coordinates}.} For each $b\in\{0,a\}$, we define $\kappa_b^{\rm E}:=\varkappa_{k,\bar\varpi-b}>0$, where $\varkappa_{k,\bar\varpi-b}$ is the coefficient in the Euler norm \eqref{eq:actual-X-norm}.

Let $J_x\in C^\infty(\R^3; (0, +\infty))$ be a radial function such that
\begin{equation}\label{eq:spatial-weight-J}
	J_x(x)=1\ \text{for}\ |x|\leq R_{\rm e}/2,\quad J_x(x)=\langle x\rangle^{M_x}\ \text{for}\ |x|\geq R_{\rm e},\quad J_x\asymp\langle x\rangle^{M_x},
\end{equation}
where $M_x>3$ is fixed. Such a function can be chosen with $d^j|\nabla_x^jJ_x|\leq C_jJ_x$, where $d:=((T-t)^{2/r}+|x|^2)^{1/2}$ is the length in \eqref{eq:background-bounds}. By \eqref{eq:spatial-weight-J} and \eqref{eq:fixed-spectral-cutoff}, we have $J_x((T-t)^{1/r}y)=1$ for every $y\in\supp\chi_{\rm sp}$ and $0\leq t<T$. For $y:=x/(T-t)^{1/r}$ and $\mathbf z^\sharp(t,y):=\mathbf z(t,(T-t)^{1/r}y)$, we define
\begin{equation}\label{eq:weighted-macro-variable}
	\mathsf V(t,y):=J_x((T-t)^{1/r}y)\mathsf T(y)\mathbf z^\sharp(t,y),
\end{equation}
where $\mathsf T(y)$ is the explicit matrix in \eqref{eq:nonlinear-transfer}.\footnote{The five-component field $\mathsf V$ is distinct from the three-component velocity perturbation $V$.}

We use the symmetric matrix
\begin{equation}\label{eq:roadmap-energy-matrix}
	D(\mathbf z):=\operatorname{diag}\left((1+B)\mathrm e^{-6S/5}I_3,\ (1+B)\mathrm e^{-3S/5},\ \tfrac32(1+B)^3\mathrm e^{-3S/5}\right),
\end{equation}
and denote $|Z|_{D(\mathbf z)}^2:=Z^{\mathsf T}D(\mathbf z)Z$.
This matrix satisfies $D(0)=D_*=\operatorname{diag}(I_3,1,3/2)$, and it is uniformly positive when $(n,V,\Theta)$ belongs to $\mathcal U:=\{(n,V,\Theta):3/4\leq n,\Theta\leq5/4,\ |V|\leq1/4\}$, with $(n,V,\Theta)$ defined in \eqref{eq:kinetic-normalization}. For $b\in\{0,a\}$, we define the macroscopic energy
\begin{equation}\label{eq:roadmap-natural-Euler-energy}
	A_b^{\rm E}:=\frac{c_{\rm s}^2}{2}\int\langle y\rangle^{\bar\varpi-b}\left(\phi_1^{2K}|\Delta_y^k\mathsf V|_{D(\mathbf z^\sharp)}^2+\kappa_b^{\rm E}|\mathsf V|_{D(\mathbf z^\sharp)}^2\right)\dd y.
\end{equation}
For $(n,V,\Theta)\in\mathcal U$, by $D(\mathbf z)\asymp D_*$ and \eqref{eq:actual-X-norm}, we have
\[A_b^{\rm E}\asymp\|\mathsf V\|_{X^K_{\bar\varpi-b}}^2,\qquad b=0,a.\]

\medskip

To couple \eqref{eq:roadmap-natural-Euler-energy} with the microscopic energy, we consider the derivatives of
\[M=n\Theta^{-3/2}\mu\left(\frac{\widehat v-V}{\sqrt\Theta}\right),\qquad n=(1+B)^3\mathrm e^{-3S/5},\qquad \Theta=(1+B)^2\mathrm e^{3S/5}\]
with respect to $\mathbf z=(V,B,S)$. For $\zeta=(\zeta_V,\zeta_B,\zeta_S)\in\R^3\times\R\times\R$, we define
\begin{equation}\label{eq:macroscopic-Maxwellian-differential}
	D_{\mathbf z}M[\zeta]:=\sum_{j=1}^3\frac{\partial M}{\partial V_j}(\zeta_V)_j+\frac{\partial M}{\partial B}\zeta_B+\frac{\partial M}{\partial S}\zeta_S,
\end{equation}
where the derivatives are taken at fixed $(t,x,v)$.
Using \eqref{eq:physical-parameters} and \eqref{eq:normalized-reconstruction} in \eqref{eq:macroscopic-Maxwellian-differential}, we obtain
\begin{equation}\label{eq:explicit-Maxwellian-derivative}
	\frac{D_{\mathbf z}M[\zeta]}{M}=\frac{(\widehat v-V)\cdot\zeta_V}{\Theta}+\frac{|\widehat v-V|^2}{\Theta(1+B)}\zeta_B+\left(\frac{3|\widehat v-V|^2}{10\Theta}-\frac32\right)\zeta_S.
\end{equation}
In particular, $D_{\mathbf z}M[\zeta]/M$ belongs to $\operatorname{span}\{1,v_1,v_2,v_3,|v|^2\}$. Therefore,  \eqref{eq:differentiated-zero-moments} imply
\[\int\frac{\partial_x^\alpha H\,D_{\mathbf z}M[\zeta]}{M}\dd\widehat v=0\qquad (\alpha\in\mathbb N_0^3).\]
Here $\zeta$ may depend on $(t,x)$ but not on $v$. The substitution and the integration follow the conventions in \eqref{eq:physical-derivative-convention} and \eqref{eq:physical-integral-convention}. 

\begin{lemma}[Gaussian identities for the macroscopic variables]\label{prop:fitted-metric}
	For $B>-1$, $S\in\R$, and every $\zeta\in\R^5$, the derivative in \eqref{eq:explicit-Maxwellian-derivative} satisfies
	\begin{equation}\label{eq:fitted-metric-ABS}
		\int\frac{|D_{\mathbf z}M[\zeta]|^2}{M}\dd\widehat v=n\left(\frac{|\zeta_V|^2}{\Theta}+\frac{15\zeta_B^2}{(1+B)^2}+\frac9{10}\zeta_S^2\right).
	\end{equation}
	For $\mathsf T(y)$ in \eqref{eq:nonlinear-transfer} and $\Psi\in\R^5$, we also have
	\begin{equation}\label{eq:fitted-metric-W}
		\int\frac{|D_{\mathbf z}M[\mathsf T(y)^{-1}\Psi]|^2}{M}\dd\widehat v=c_{\rm s}^2\vartheta(y)^{-6}|\Psi|_{D(\mathbf z)}^2.
	\end{equation}
\end{lemma}
\begin{proof}
	The identity \eqref{eq:fitted-metric-ABS} follows by applying the standard Gaussian moments to \eqref{eq:explicit-Maxwellian-derivative}. Substituting \eqref{eq:nonlinear-transfer}, \eqref{eq:normalized-reconstruction} and \eqref{eq:roadmap-energy-matrix}, we then obtain \eqref{eq:fitted-metric-W}.
\end{proof}

The identity \eqref{eq:fitted-metric-W} identifies the quadratic form \eqref{eq:roadmap-natural-Euler-energy} with a Gaussian integral of Maxwellian derivatives. The chain rule for $M(t,x,\mathbf z(t,x);v)$ explains why the highest-order derivatives of $\mathbf z$ can be canceled.

\begin{lemma}[The highest-order derivatives of the Maxwellian]\label{lem:roadmap-Maxwellian-composition}
	Let $K:=2k\geq2$ and let $u,\theta$ and $\mathbf z=(V,B,S)$ be smooth functions with $\theta>0$ and $B>-1$. Let $M$ be given by \eqref{eq:normalized-Maxwellian} and $\bar M=\mathcal M(\theta^3,u,\theta^2)$. Then
	\begin{equation}\label{eq:roadmap-Maxwellian-chain-rule}
		\Delta_x^k(M-\bar M)=D_{\mathbf z}M[\Delta_x^k\mathbf z]+\mathcal R_K,
	\end{equation}
	where $\mathcal R_K$ contains at most $K-1$ spatial derivatives of $\mathbf z$ and at most $K$ spatial derivatives of $(u,\theta)$. All spatial derivatives are taken at fixed $v$.
\end{lemma}
\begin{proof}
	By \eqref{eq:normalized-Maxwellian}, we regard $M$ as a function of $(u,\theta,\mathbf z,v)$. For $j=1,2,3$, we have
	\[\partial_{x_j}M=D_{\mathbf z}M[\partial_{x_j}\mathbf z]+D_uM[\partial_{x_j}u]+\partial_\theta M\partial_{x_j}\theta.\]
	By induction on $|\alpha|\geq1$, we obtain
	\[\partial_x^\alpha M=D_{\mathbf z}M[\partial_x^\alpha\mathbf z]+D_uM[\partial_x^\alpha u]+\partial_\theta M\partial_x^\alpha\theta+\widetilde{\mathcal R}_\alpha,\]
	where all derivatives of $(u,\theta,\mathbf z)$ in $\widetilde{\mathcal R}_\alpha$ are of order at most $|\alpha|-1$. Since $\Delta_x^k$ is a linear combination of $\partial_x^\alpha$ with $|\alpha|=K$ and $\bar M$ depends only on $(u,\theta,v)$, we obtain \eqref{eq:roadmap-Maxwellian-chain-rule}.
\end{proof}

By \eqref{eq:roadmap-Maxwellian-chain-rule} and $H=f-M$, we have
\[\Delta_x^kH+D_{\mathbf z}M[\Delta_x^k\mathbf z]=\Delta_x^k(f-\bar M)-\mathcal R_K.\]
In the macroscopic energy \eqref{eq:roadmap-natural-Euler-energy}, however, we estimate $\Delta_y^k(J_x\mathsf T\mathbf z^\sharp)$. Since
\[\Delta_y^k(J_x\mathsf T\mathbf z^\sharp)=J_x\mathsf T\Delta_y^k\mathbf z^\sharp+[\Delta_y^k,J_x\mathsf T]\mathbf z^\sharp,\]
and the commutator contains at most $K-1$ derivatives of $\mathbf z^\sharp$, we can incorporate the spatial weights into the Maxwellian differential without affecting the cancellation of the highest-order derivatives. 

We use the untruncated Euler profiles
\[u_{\rm p}(t,x):=(T-t)^{-\alpha_{\rm ss}}U(y),\qquad\theta_{\rm p}(t,x):=(T-t)^{-\alpha_{\rm ss}}\vartheta(y),\qquad y:=\frac{x}{(T-t)^{1/r}}.\]
The spatial weights use $\theta_{\rm p}$ on all of $\R^3_x$, and $M$ in \eqref{eq:normalized-Maxwellian} uses the truncated background $(u,\theta)$. By \eqref{eq:profile}, the backgrounds agree on $|x|<R_{\rm e}$. Recall $d:=((T-t)^{2/r}+|x|^2)^{1/2}$ and \eqref{eq:background-bounds}. By the profile estimates \eqref{eq:profile-decay}, we have
\begin{equation}\label{eq:untruncated-profile-bounds}
	\theta_{\rm p}\asymp d^{-(r-1)},\quad d^j|\nabla_x^j\log\theta_{\rm p}|\leq C_j\ (j\geq1),\quad d^j|\nabla_x^j\partial_t\log\theta_{\rm p}|\leq C_j\theta/d\ (j\geq0).
\end{equation}
For $|x|\geq R_{\rm e}/2$, thanks to the lower bound $d\geq R_{\rm e}/2$, the coefficients in \eqref{eq:untruncated-profile-bounds} are bounded uniformly in $t<T$. We define
\begin{equation}\label{eq:new-global-weights}
	\Sigma_x:=J_x\theta_{\rm p}^3,\quad d_\phi:=(T-t)^{1/r}\phi_1(y),\quad \Omega_b^*:=(T-t)^{b/r}J_x^2\theta_{\rm p}^6d^{\bar\varpi-b}\ \text{for}\ b\in\{0,a\}.
\end{equation}
The weight in \cite[Lemma 4.1]{BCGJVY} has the form $\phi_1=\phi_{\rm b}^{c_2}(1+c_3\langle y\rangle)$, where $\phi_{\rm b}>0$ is radial and nonincreasing. Hence, $\phi_1\asymp\langle y\rangle$ and $y\cdot\nabla\log\phi_1\leq1$. As a consequence,
\begin{equation}\label{eq:good-length-time}
	d_\phi\asymp d,\qquad \partial_t\log d_\phi=-\frac{1-y\cdot\nabla_y\log\phi_1}{r(T-t)}\leq0,\qquad|\partial_t\log d_\phi|\leq C(T-t)^{-1+1/r}d^{-1}.
\end{equation}
The explicit expressions for the weights also give
\[\Omega_b^*\asymp J_x^2\langle y\rangle^{-b}d^{-3},\qquad \Omega_a^*/\Omega_0^*=\langle y\rangle^{-a}.\]
The weight $\Omega_0^*$ is used for global control, and $\Omega_a^*$ is used for the decaying Euler estimate.

Let $\mathsf T_0:=\operatorname{diag}(c_{\rm s}^{-1}I_3,3,-3/5)$, so that $\mathsf T(y)=\vartheta(y)^3\mathsf T_0$. By \eqref{eq:nonlinear-transfer}, \eqref{eq:weighted-macro-variable} and \eqref{eq:new-global-weights}, at $x=(T-t)^{1/r}y$, we have
\[\Sigma_x=(T-t)^{-3\alpha_{\rm ss}}J_x\vartheta^3,\qquad \Sigma_x\mathbf z=(T-t)^{-3\alpha_{\rm ss}}\mathsf T_0^{-1}\mathsf V.\]
By $d_\phi=(T-t)^{1/r}\phi_1$ and $\Delta_x^k=(T-t)^{-K/r}\Delta_y^k$, we obtain
\begin{equation}\label{eq:roadmap-weighted-derivative-identity}
	d_\phi^K\Sigma_x^{-1}\Delta_x^k(\Sigma_x\mathbf z)=\phi_1^KJ_x^{-1}\vartheta^{-3}\mathsf T_0^{-1}\Delta_y^k\mathsf V=\phi_1^KJ_x^{-1}\mathsf T(y)^{-1}\Delta_y^k\mathsf V.
\end{equation}

By the Leibniz rule, we have
\begin{equation}\label{eq:roadmap-conjugated-Laplacian}
	\Sigma_x^{-1}\Delta_x^k(\Sigma_x\mathbf z)=\Delta_x^k\mathbf z+\Sigma_x^{-1}[\Delta_x^k,\Sigma_x]\mathbf z,
\end{equation}
where $[\Delta_x^k,\Sigma_x]\mathbf z$ contains at most $K-1$ spatial derivatives of $\mathbf z$. We define
\begin{equation}\label{eq:roadmap-composition-remainders}
	\begin{aligned}
		\mathcal B_K:=\Delta_x^k(M-\bar M)-D_{\mathbf z}M[\Sigma_x^{-1}\Delta_x^k(\Sigma_x\mathbf z)],\quad\mathcal B_0:=M-\bar M-D_{\mathbf z}M[\mathbf z].
	\end{aligned}
\end{equation}
By \eqref{eq:roadmap-Maxwellian-chain-rule} and \eqref{eq:roadmap-conjugated-Laplacian}, we have
\[\mathcal B_K=\mathcal R_K-D_{\mathbf z}M[\Sigma_x^{-1}[\Delta_x^k,\Sigma_x]\mathbf z].\]
Hence $\mathcal B_K$ contains at most $K-1$ spatial derivatives of $\mathbf z$, and $\mathcal B_0$ contains no spatial derivatives of $\mathbf z$.

\medskip\noindent\emph{The microscopic energy.}
Since $\int(1,v,|v|^2)H\dd v=0$, by Proposition~\ref{thm:new-collision}, the collision term in \eqref{eq:roadmap-microscopic-equation} provides the dissipation
\[\frac1\eps\iint\Omega_b^*\theta^{3+\gamma}\widehat\nu_M|H/W|^2\dd\widehat v\dd x,\qquad b=0,a.\]
We take $W=W_{m_\eps,\ell}$ as in \eqref{eq:fitted-velocity-weight} and \eqref{eq:adaptive-m}. Recall that $W^2=M$ for $|\widetilde v|^2\leq m_\eps$.

In addition to $H$ and $d_\phi^K\Delta_x^kH$, we need to estimate the lower-order spatial derivatives of $H$. We define
\begin{equation}\label{eq:spatial-lower-hierarchy}
	H_\alpha:=d^{|\alpha|}\partial_x^\alpha H,\qquad |\alpha|=j<K,\qquad \widehat s_j:=2(K-j).
\end{equation}
To control the loss of velocity weights in the differentiated equations, for $|\alpha|=j<K$, we estimate
\[\iint\Omega_b^*\langle\widetilde v\rangle^{2\widehat s_j}|H_\alpha/W|^2\dd\widehat v\dd x.\]
The coefficients of the lower-order energies will be chosen in Proposition~\ref{prop:new-component-collision}. For the recovery of $\partial_x^\alpha H$ with $|\alpha|=K$ from $d_\phi^K\Delta_x^kH$ and \eqref{eq:spatial-lower-hierarchy}, see Lemma~\ref{lem:weighted-elliptic-recovery}.

In the energy estimate for $\Delta_x^k\mathbf z$, the source terms $\Delta_x^k\Div_x\Pi$ and $\Delta_x^k\Div_xq$ in \eqref{eq:moment-system} contain $K+1$ spatial derivatives of $H$. Similarly, differentiating \eqref{eq:roadmap-projected-sources} gives terms containing $K+1$ spatial derivatives of $\mathbf z$ in the equation for $\Delta_x^kH$. To avoid the loss of spatial derivatives, by \eqref{eq:roadmap-Maxwellian-chain-rule}, we consider
\[\Delta_x^kH+D_{\mathbf z}M[\Delta_x^k\mathbf z]=\Delta_x^k(f-\bar M)-\mathcal R_K,\]
where $\mathcal R_K$ contains at most $K-1$ spatial derivatives of $\mathbf z$. Since $[\partial_t+v\cdot\nabla_x,\Delta_x^k]=0,$ we can combine the order-$(K+1)$ terms before performing the energy estimates. We next incorporate the spatial weights in \eqref{eq:roadmap-natural-Euler-energy} into the Maxwellian differential by using \eqref{eq:fitted-metric-W}.

We define
\begin{equation}\label{eq:matched-unknowns}
	\begin{aligned}
		H_K:=d_\phi^K\Delta_x^kH,\quad &P_K:=D_{\mathbf z}M\!\left[d_\phi^K\Sigma_x^{-1}\Delta_x^k(\Sigma_x\mathbf z)\right],\quad H_0:=H,\quad P_0:=D_{\mathbf z}M[\mathbf z],\\
		&G_0:=H_0+P_0,\qquad G_K:=H_K+P_K.
	\end{aligned}
\end{equation}
Thus, $P_K$ and $P_0$ are explicit Maxwellian derivatives of the weighted highest-order perturbation $d_\phi^K\Sigma_x^{-1}\Delta_x^k(\Sigma_x\mathbf z)$ and the zero-order perturbation $\mathbf z$ in \eqref{eq:roadmap-composition-remainders}. They combine the Euler and microscopic derivatives without entering the definition \eqref{eq:roadmap-natural-Euler-energy} of $A_b^{\rm E}$. The Gaussian identity \eqref{eq:fitted-metric-W}, the derivative identity \eqref{eq:roadmap-weighted-derivative-identity}, and the change of variables $x=(T-t)^{1/r}y$ show that the Gaussian squares of $P_0,P_K$ reproduce $A_b^{\rm E}$ exactly.

\begin{proposition}[Matching and the highest-order equation]\label{prop:global-matching}
	Let $f$ be a smooth solution to \eqref{Eq.Boltzmann} with positive density and temperature. Assume that $f$ is continuous with values in $L^1(\R^3_v,\langle v\rangle^{2+\gamma}\dd v)$, and that $\partial_tf$ and the spatial derivatives of $f$ up to order $K$ are continuous with values in $L^1(\R^3_v,\langle v\rangle^3\dd v)$. Let $M$ be the corresponding local Maxwellian, let $H:=f-M$, and let $H_i,P_i,G_i$ be defined by \eqref{eq:matched-unknowns}. Then for $i\in\{0,K\}$, the five moments of $H_i$ against $1,v_1,v_2,v_3,|v|^2$ vanish. Moreover, we have
	\begin{align}\label{eq:new-Gaussian-match}
		A_b^{\rm E}=\frac12\iint\Omega_b^*\frac{|P_K|^2+\kappa_b^{\rm E}|P_0|^2}{M}\dd\widehat v\dd x,\quad \int\frac{H_iP_i}{M}\dd\widehat v=0
	\end{align}
	for $b\in\{0,a\}$ and $i\in\{0, K\}$.
	For the remainders in \eqref{eq:roadmap-composition-remainders}, the transport operator $D$ and the prescribed source $S_0:=D\bar M$ in \eqref{eq:source} satisfy
	\begin{equation}
		\begin{aligned}\label{eq:new-good-evolution}
			&DG_K=\eps^{-1}d_\phi^K\Delta_x^kQ(f,f)-d_\phi^K\Delta_x^kS_0-d_\phi^KD\mathcal B_K+K(D\log d_\phi)G_K,\\
			&\qquad\qquad\qquad DG_0=\eps^{-1}Q(f,f)-S_0-D\mathcal B_0.
		\end{aligned}
	\end{equation}
	Moreover, every spatial derivative of $\mathbf z$ occurring in $\mathcal B_K$ has order at most $K-1$. 
\end{proposition}

By the definitions of $G_K$, $G_0$, $\mathcal B_K$ and $\mathcal B_0$, we have
\begin{equation}\label{eq:coupled-perturbation-identities}
	G_K=d_\phi^K\left(\Delta_x^k(f-\bar M)-\mathcal B_K\right),\qquad G_0=f-\bar M-\mathcal B_0.
\end{equation}
Since $[D,\Delta_x^k]=0$ for $D:=\partial_t+v\cdot\nabla_x$, by \eqref{Eq.Boltzmann} and \eqref{eq:coupled-perturbation-identities}, we have
\begin{align*}
	DG_K&=K\frac{Dd_\phi}{d_\phi}G_K+d_\phi^K\left(\eps^{-1}\Delta_x^kQ(f,f)-\Delta_x^kD\bar M-D\mathcal B_K\right),\\
	DG_0&=\eps^{-1}Q(f,f)-D\bar M-D\mathcal B_0.
\end{align*}
After substituting the macroscopic equations into $D\mathcal B_K$, we only need to estimate spatial derivatives of the perturbations up to order $K$. See Subsection~\ref{subsec:exact-matching} for the proof of Proposition~\ref{prop:global-matching}.

By \eqref{eq:new-Gaussian-match}, we have
\[A_b^{\rm E}=\frac12\iint\Omega_b^*\frac{|P_K|^2+\kappa_b^{\rm E}|P_0|^2}{M}\dd\widehat v\dd x,\qquad \int\frac{H_iP_i}{M}\dd\widehat v=0\quad(i=0,K).\]
Recall that $G_i=H_i+P_i$ and $W^2=M$ for $|\widetilde v|^2\leq m_\eps$. Thus, for $i=0,K$, we have
\begin{align*}
	\int\frac{|G_i|^2}{W^2}\dd\widehat v-\int\frac{|H_i|^2}{W^2}\dd\widehat v-\int\frac{|P_i|^2}{M}\dd\widehat v=\int_{|\widetilde v|^2>m_\eps}(2H_iP_i+|P_i|^2)(W^{-2}-M^{-1})\dd\widehat v.
\end{align*}
Since $P_i/M$ is a quadratic polynomial in $v$ of degree at most two, the Gaussian decay of $P_i$ yields exponentially small coefficients in the estimates of the integral on the right-hand side. In Lemma~\ref{lem:single-spatial-allocation}, we estimate the integral using $W^{-1}H_i\in L^2_{\widehat v}$ without assuming $M^{-1/2}H_i\in L^2_{\widehat v}$. See Subsection~\ref{subsec:weighted-norm-recovery} for the proof.

\medskip\noindent\emph{The spatial energies and their estimates.}
Let $b_j^{\rm sp},\kappa_j>0$ ($0\leq j<K$) be chosen as in Proposition~\ref{prop:new-component-collision} to absorb the lower-order collision terms. We fix $b_j^{\rm sp}$, $\kappa_j$ and $\ell$ before choosing $\eps$. We define
\begin{align}\label{eq:new-kinetic-energies}
	I_b&:=\frac12\iint\Omega_b^*\left(\frac{|G_K|^2}{W^2}+\kappa_b^{\rm E}\frac{|G_0|^2}{W^2}\right)\dd\widehat v\dd x,\notag\\
	I_b^{\rm low}&:=\frac12\sum_{j<K}\sum_{|\alpha|=j}\iint\Omega_b^*\left(b_j^{\rm sp}+\kappa_j\langle\widetilde v\rangle^{2\widehat s_j}\right)\frac{|H_\alpha|^2}{W^2}\dd\widehat v\dd x,\notag\\
	I_b^{\rm mic}&:=\frac12\iint\Omega_b^*\left(\frac{|H_K|^2}{W^2}+\kappa_b^{\rm E}\frac{|H_0|^2}{W^2}\right)\dd\widehat v\dd x,\\
	\mathcal D_b&:=\iint\Omega_b^*\theta^{3+\gamma}\widehat\nu_M\left(\frac{|H_K|^2}{W^2}+\kappa_b^{\rm E}\frac{|H_0|^2}{W^2}\right)\dd\widehat v\dd x\notag\\
	&\quad+\sum_{j<K}\sum_{|\alpha|=j}\iint\Omega_b^*\theta^{3+\gamma}\widehat\nu_M(b_j^{\rm sp}+\kappa_j\langle\widetilde v\rangle^{2\widehat s_j})\frac{|H_\alpha|^2}{W^2}\dd\widehat v\dd x.\notag
\end{align}
Here $\widetilde v=(\widehat v-V)/\sqrt\Theta$, and $\widehat\nu_M$ is given by \eqref{eq:normalized-collision-frequency}. In the definition of $I_b$, we replace $(H_0,H_K)$ in $I_b^{\rm mic}$ by $(G_0,G_K)$, so that we can use \eqref{eq:new-good-evolution} in the estimates of the transport terms. The dissipation $\mathcal D_b$ involves only the microscopic quantities.

By Lemma~\ref{lem:single-spatial-allocation}, for sufficiently large $m_0$ and sufficiently small macroscopic perturbations with $(n,V,\Theta)\in\mathcal U$, we have
\[I_b+I_b^{\rm low}\asymp A_b^{\rm E}+I_b^{\rm mic}+I_b^{\rm low},\qquad b=0,a.\]
For the estimates of all spatial derivatives $\partial_x^\alpha H$ with $|\alpha|\leq K$ in terms of $I_b^{\rm mic}+I_b^{\rm low}$, see also Lemma~\ref{lem:single-spatial-allocation}.

Recall from \eqref{eq:new-global-weights} that
\[\Omega_a^*=\langle y\rangle^{-a}\Omega_0^*,\qquad y=\frac{x}{(T-t)^{1/r}}.\]
We use $\Omega_0^*$ to keep the perturbations small in the whole space and $\Omega_a^*$ to obtain the decay of the perturbations in self-similar variables. To use the Euler decay estimate in Proposition~\ref{thm:local-observation}, we add the quadratic form $\frac{c_{\rm s}^2}{2}\mathcal R_*(\mathsf V,\mathsf V)$ to $I_a+I_a^{\rm low}$. We define
\begin{equation}\label{eq:spatial-master-energies}
	\mathcal B:=I_0+I_0^{\rm low},\qquad \mathcal E:=I_a+I_a^{\rm low}+\frac{c_{\rm s}^2}{2}\mathcal R_*(\mathsf V,\mathsf V).
\end{equation}

\begin{proposition}[The spatial energy inequalities]\label{thm:new-closed-pair}
	Let $\gamma$, $r$, $\eta$ and $K=2k\geq\max\{14,2k_0+2\}$ be given as in \eqref{Eq.gamma-range}, \eqref{eq:exponents} and \eqref{eq:euler-stability-space}. Let the Euler background, the spectral parameters and the spatial weights be given by Proposition~\ref{Thm.Euler-solution}, \eqref{eq:new-spectral-choices} and \eqref{eq:new-global-weights}, respectively. We choose $\ell$ and $m_0$ as in Lemma~\ref{lem:polynomial-collision} with $s_*=2K+8$, and take the velocity weights in \eqref{eq:fitted-velocity-weight} and \eqref{eq:adaptive-m} and the coefficients in \eqref{eq:new-kinetic-energies} as in Proposition~\ref{prop:new-component-collision}.
	
	There exist constants $e_0>0$, $\eps_0\in(0,1]$, $T_{\rm kin}\in(0,\min\{T_0,1\}]$ and $c,c_0,C>0$, with $R_2T_{\rm kin}^{1/r}<R_{\rm e}/2$, such that the following holds. Let $0<\eps\leq\eps_0$, $0\leq S<T\leq T_{\rm kin}$ and let $f\geq0$ be a smooth solution to \eqref{Eq.Boltzmann} on $[0,S]\times\R^3_x\times\R^3_v$ satisfying $f(t,Ox,Ov)=f(t,x,v)$ for all $O\in SO(3)$. Assume that
	\begin{equation}\label{eq:spatial-bootstrap-hypotheses}
		(n,V,\Theta)(t,x)\in\mathcal U\quad(0\leq t\leq S,\ x\in\R^3),\qquad \sup_{0\leq t\leq S}\mathcal B(t)\leq e_0.
	\end{equation}
	Assume further that the energy and dissipation integrals in \eqref{eq:new-kinetic-energies} and \eqref{eq:spatial-master-energies} are finite and continuous on $[0,S]$ and that $\mathcal E,\mathcal B\in C^1([0,S])$. Then, for $t\in[0,S]$, we have
	\begin{align}
		\mathcal E'+\frac\lambda{T-t}\mathcal E+\frac c\eps\mathcal D_a&\leq\frac C{T-t}|\mathbf a|_P^2+C\eps(T-t)^{-1+\eta},\label{eq:new-closed-inner}\\
		\mathcal B'+\frac{c_0}\eps\mathcal D_0&\leq\frac C{T-t}\mathcal E+C\mathcal B+C\eps(T-t)^{-1+\eta},\label{eq:new-closed-global}
	\end{align}
	where $\mathbf a$ is defined in \eqref{eq:roadmap-unstable-coordinates}. The constants $e_0,\eps_0,T_{\rm kin},c,c_0,C$ depend only on the fixed profile, differentiation orders, weights, hierarchy coefficients and spectral gaps.
\end{proposition}

See Subsection~\ref{subsec:closed-energy-proofs} for the proof of Proposition~\ref{thm:new-closed-pair}. By Gr\"onwall's inequality applied to \eqref{eq:new-closed-global}, we have
\[\sup_{0\leq q\leq t}\mathcal B(q)+\frac{c_0}\eps\int_0^t\mathcal D_0(q)\dd q\lesssim\mathcal B(0)+\int_0^t\frac{\mathcal E(q)}{T-q}\dd q+\eps T^\eta,\qquad t\in[0,S].\]
Thus, to control $\mathcal B$ on $[0,T)$, we only need the time integrability of $\mathcal E/(T-t)$. To obtain the decay of $\mathcal E$ from \eqref{eq:new-closed-inner}, we will choose the initial data so that $|\mathbf a(t)|_P^2\lesssim(T-t)^\delta$ for some $\delta>0$.

\medskip\noindent\emph{Velocity derivatives and continuation.}
The (conditional) continuation result is given in Proposition~\ref{thm:B-restart}, which requires the positivity of the density and the temperature, bounded first spatial derivatives of the physical parameters, and the smallness of $\sup_{x,v}\langle\widetilde v\rangle^8|H/W|$; the higher derivatives are then propagated on a common local time interval. The macroscopic part of $\mathcal B$ controls the parameter bounds, but its microscopic part controls only spatial derivatives in a velocity $L^2$ norm. In particular, it does not control the velocity supremum. Therefore, we add physical velocity derivatives, so that the Sobolev embedding in both $x$ and $v$ gives the required pointwise bound, as proved in Lemma~\ref{lem:B-pointwise}.

For spatial and velocity multi-indices $\alpha,\beta\in\mathbb N_0^3$, we define $i:=|\alpha|$, $j:=|\beta|$, $N:=i+j$, and
\begin{equation}\label{eq:mixed-velocity-derivatives}
	\begin{gathered}
		\mathfrak A_K:=\{(\alpha,\beta):1\leq|\beta|,\ |\alpha|+|\beta|\leq K\},\quad
		Z_{\alpha,\beta}:=d^i\theta^j\partial_x^\alpha\partial_v^\beta H,\quad p_N:=2(K-N).
	\end{gathered}
\end{equation}
We choose positive coefficients $b_N^{\rm vel},k_N^{\rm vel}$ for $1\leq N<K$, and set $b_K^{\rm vel}:=1$, $k_K^{\rm vel}:=0$. The mixed-derivative energy and the corresponding collision dissipation are defined by
\begin{equation}\label{eq:mixed-velocity-energy}
	\begin{aligned}
		\mathcal V_K&:=\frac12\sum_{(\alpha,\beta)\in\mathfrak A_K}\iint\Omega_0^*\left(b_N^{\rm vel}+k_N^{\rm vel}\langle\widetilde v\rangle^{2p_N}\right)\frac{|Z_{\alpha,\beta}|^2}{W^2}\dd\widehat v\dd x,\\
		\mathcal D_K^{\rm vel}&:=\sum_{(\alpha,\beta)\in\mathfrak A_K}\iint\Omega_0^*\theta^{3+\gamma}\widehat\nu_M\left(b_N^{\rm vel}+k_N^{\rm vel}\langle\widetilde v\rangle^{2p_N}\right)\frac{|Z_{\alpha,\beta}|^2}{W^2}\dd\widehat v\dd x.
	\end{aligned}
\end{equation}
In \eqref{eq:mixed-velocity-energy}, we write $N=|\alpha|+|\beta|$ and make the substitution $v=u+\theta\widehat v$ after taking the derivatives $\partial_x^\alpha\partial_v^\beta H$ in the physical coordinates. As in the choice of $\widehat s_j$ in \eqref{eq:new-kinetic-energies}, we take $p_N=2(K-N)$, so that $p_N-p_{N+1}=2$ for $1\leq N<K$.

\begin{proposition}[Mixed spatial and velocity derivatives]\label{prop:mixed-velocity-estimate}
	Let the fixed data be as in Proposition~\ref{thm:new-closed-pair}. There exist positive coefficients $b_N^{\rm vel},k_N^{\rm vel}$ ($1\leq N<K$), with $b_K^{\rm vel}=1$ and $k_K^{\rm vel}=0$, and constants
	\[e_{\rm vel}\in(0,\min\{1,e_0\}],\qquad \eps_{\rm vel}\in(0,\eps_0],\qquad c_{\rm vel},C_{\rm vel}>0,\]
	such that the following holds. Let $0<\eps\leq\eps_{\rm vel}$, $0\leq S<T\leq T_{\rm kin}$ and let $f\geq0$ be a smooth simultaneous-rotation-invariant solution to \eqref{Eq.Boltzmann} on $[0,S]\times\R^3_x\times\R^3_v$ satisfying
	\begin{equation}\label{eq:mixed-bootstrap-hypotheses}
		(n,V,\Theta)(t,x)\in\mathcal U\quad(0\leq t\leq S,\ x\in\R^3),\qquad \sup_{0\leq t\leq S}(\mathcal B(t)+\mathcal V_K(t))\leq e_{\rm vel}.
	\end{equation}
	Assume further that the energy and dissipation integrals in \eqref{eq:new-kinetic-energies}, \eqref{eq:spatial-master-energies} and \eqref{eq:mixed-velocity-energy} are finite and continuous on $[0,S]$, and that $\mathcal E,\mathcal B,\mathcal V_K\in C^1([0,S])$. Then, for $t\in[0,S]$, we have
	\begin{equation}\label{eq:mixed-velocity-differential}
		\mathcal V_K'+\frac{c_{\rm vel}}\eps\mathcal D_K^{\rm vel}\leq\frac{C_{\rm vel}}\eps\mathcal D_0+C_{\rm vel}\eps(T-t)^{-1+\eta}(1+\mathcal B).
	\end{equation}
	The coefficients and constants depend only on the fixed data in Proposition~\ref{thm:new-closed-pair}. Moreover, if $f(0,x,\cdot)$ is a Maxwellian for every $x\in\R^3$, then $H(0)=0$ and $\mathcal V_K(0)=0$.
\end{proposition}

See Subsection~\ref{subsec:closed-energy-proofs} for the proof of Proposition~\ref{prop:mixed-velocity-estimate}. To absorb the term $C_{\rm vel}\mathcal D_0/\eps$ in \eqref{eq:mixed-velocity-differential}, we choose
\[0<\kappa_{\rm reg}\leq\min\{1,c_0/(2C_{\rm vel})\}\]
and define
\begin{equation}\label{eq:augmented-global}
	\mathcal B_{\rm reg}:=\mathcal B+\kappa_{\rm reg}\mathcal V_K,\qquad \mathcal D_{\rm reg}:=\mathcal D_0+\kappa_{\rm reg}\mathcal D_K^{\rm vel}.
\end{equation}
Adding $\kappa_{\rm reg}$ times \eqref{eq:mixed-velocity-differential} to \eqref{eq:new-closed-global} and using $c_0-\kappa_{\rm reg}C_{\rm vel}\geq c_0/2$ and $\mathcal B\leq e_{\rm vel}\leq1$, we obtain
\begin{equation}\label{eq:augmented-global-inequality}
	\mathcal B_{\rm reg}'+\frac{c_{\rm reg}}\eps\mathcal D_{\rm reg}\leq\frac C{T-t}\mathcal E+C\mathcal B_{\rm reg}+C\eps(T-t)^{-1+\eta},
\end{equation}
where $c_{\rm reg}:=\min\{c_0/2,c_{\rm vel}\}$. We fix all hierarchy coefficients and $\kappa_{\rm reg}$ before choosing $\eps$.

By \eqref{eq:augmented-global} and $\kappa_{\rm reg}\leq1$, we have
\begin{equation}\label{eq:augmented-bootstrap-comparison}
	\mathcal B+\mathcal V_K\leq\kappa_{\rm reg}^{-1}\mathcal B_{\rm reg}.
\end{equation}
Hence, to verify the energy bound in \eqref{eq:mixed-bootstrap-hypotheses}, we only need to prove
\[\sup_{0\leq t\leq S}\mathcal B_{\rm reg}(t)\leq\kappa_{\rm reg}e_{\rm vel}.\]
For the bounds on the macroscopic coefficients and the weighted pointwise bounds on $H$ required in Proposition~\ref{thm:B-local}, see Lemmas~\ref{lem:single-spatial-allocation} and \ref{lem:B-pointwise}. The continuation of the solution on compact time intervals in $[0,T)$ is proved in Proposition~\ref{thm:B-restart}.

\medskip\noindent\emph{The unstable coordinates.}
The decaying energy inequality still contains the finite-dimensional term $|\mathbf a(t)|_P^2$. For the coordinates defined in \eqref{eq:roadmap-unstable-coordinates}, by Proposition~\ref{thm:actual-spectral-forcing}, we have
\begin{equation}\label{eq:new-actual-mode-bound}
	\mathbf a'=\frac{A_{\rm u}\mathbf a}{T-t}+F_{\rm u},\qquad |F_{\rm u}|_P^2\leq C(T-t)^{-1+\eta}\mathcal D_a+C(T-t)^{-2}\mathcal E^2.
\end{equation}
The source estimate uses the time-integrated dissipation from \eqref{eq:new-closed-inner}. The inequalities \eqref{eq:new-closed-inner}, \eqref{eq:augmented-global-inequality} and \eqref{eq:new-actual-mode-bound} provide the energy and finite-dimensional estimates used in Theorem~\ref{thm:coupled-shooting} on the initial-data selection and continuation.

\subsection{Selection of the initial data and continuation}\label{subsec:initial-selection}
By \eqref{eq:new-closed-inner} and \eqref{eq:augmented-global-inequality}, to close the energy estimates, we seek initial data such that
\[|\mathbf a(t)|_P^2\lesssim\eps T^{\eta-d_*}(T-t)^{d_*},\qquad 0\leq t<T.\]
Due to the term $(T-t)^{-1}A_{\rm u}\mathbf a$ in \eqref{eq:new-actual-mode-bound}, the smallness of $\mathbf a(0)$ alone does not imply the required decay. We consider the family of local Maxwellian initial data constructed in Proposition~\ref{prop:positive-realization}, satisfying
\[\mathbf a(0)=\mathbf a_0,\qquad H(0,x,v)=0.\]
We will choose $\mathbf a_0$ by using \eqref{eq:new-closed-inner} and the source estimate \eqref{eq:new-actual-mode-bound}.

\medskip\noindent\emph{A positive initial family.} Let $J:=\dim E_{\rm u}$ and let $\Phi_j$ and $\ell_j$ ($1\leq j\leq J$) be given by Proposition~\ref{prop:projection-local}. For the radii $R_{\rm sp}$ and $R_1$ in \eqref{eq:fixed-spectral-cutoff}, we fix a radial function $\chi_{\rm init}\in C_c^\infty(B_{R_1})$ satisfying
\[0\leq\chi_{\rm init}\leq1,\qquad \chi_{\rm init}\equiv1\quad\text{in a neighborhood of }\overline B_{R_{\rm sp}}.\]

\begin{proposition}[A positive initial family]\label{prop:positive-realization}
	For $0<T\leq T_{\rm kin}$ and $\mathbf a_0\in\R^J$, we define
	\begin{equation}\label{eq:positive-initial-field}
		\mathbf z_0^\sharp(y):=\mathsf T(y)^{-1}\sum_{j=1}^J(\mathbf a_0)_j\chi_{\rm init}(y)\Phi_j(y),\qquad \mathbf z_0(x):=\mathbf z_0^\sharp(x/T^{1/r}).
	\end{equation}
	We denote $\mathbf z_0=(V_0,B_0,S_0^{\rm ent})$, and define
	\begin{equation}\label{eq:positive-initial-Maxwellian}
		\begin{gathered}
			n_0:=(1+B_0)^3\mathrm e^{-3S_0^{\rm ent}/5},\qquad \Theta_0:=(1+B_0)^2\mathrm e^{3S_0^{\rm ent}/5},\\
			f_0(x,v;\mathbf a_0):=\mathcal M\bigl(\theta(0,x)^3n_0(x),u(0,x)+\theta(0,x)V_0(x),\theta(0,x)^2\Theta_0(x)\bigr)(v).
		\end{gathered}
	\end{equation}
	There exists a constant $R_{\rm in}>0$ such that for all $|\mathbf a_0|_P\leq R_{\rm in}$, we have
	\[7/8\leq n_0,\Theta_0\leq9/8,\qquad |V_0|\leq1/8\quad\text{on }\R^3.\]
	In particular, $f_0$ is smooth and strictly positive. Moreover,
	\[H(0)=0,\qquad \mathbf a(0)=\mathbf a_0,\qquad \HH(f_0\mid M_\infty)<\infty,\]
	and
	\[f_0(Ox,Ov;\mathbf a_0)=f_0(x,v;\mathbf a_0)\quad(\forall\,O\in SO(3)),\quad f_0(x,v;\mathbf a_0)=\bar M(0,x,v)\quad(|x|\geq R_1T^{1/r}).\]
	For each fixed $T$, the map $(\mathbf a_0,x,v)\mapsto f_0(x,v;\mathbf a_0)$ is smooth. For every integer $m\geq0$, we have $\|\mathbf z_0^\sharp\|_{H^m_y}\leq C_m|\mathbf a_0|_P$, where $C_m$ is independent of $T$ and $\eps$. For $J=0$, we set $\mathbf z_0=0$ and $f_0=\bar M(0)$.
\end{proposition}

Note that by Proposition~\ref{prop:projection-local} and the choice of $\chi_{\rm init}$,
\[\ell_i(\chi_{\rm init}\Phi_j)=\ell_i(\Phi_j)=\delta_{ij},\qquad 1\leq i,j\leq J.\]
Since $\chi_{\rm sp}\chi_{\rm init}=\chi_{\rm init}$, by \eqref{eq:positive-initial-field} and \eqref{eq:roadmap-unstable-coordinates}, we have
\[(\mathbf a(0))_i=\ell_i(\chi_{\rm sp}\mathsf T\mathbf z_0^\sharp)=\sum_{j=1}^J(\mathbf a_0)_j\ell_i(\chi_{\rm init}\Phi_j)=(\mathbf a_0)_i.\]
See Subsection~\ref{subsec:selection-local} for the proof of Proposition~\ref{prop:positive-realization}.

\medskip\noindent\emph{The energy inequalities used for the selection.}
Recall that \eqref{eq:new-closed-inner} controls the decaying energy $\mathcal E$ and the time integral of $\mathcal D_a$, whereas \eqref{eq:augmented-global-inequality} controls the continuation energy $\mathcal B_{\rm reg}$ and the time integral of $\mathcal D_{\rm reg}$. The coordinate equation and the source estimate in \eqref{eq:new-actual-mode-bound} use the same dissipation $\mathcal D_a$. Thus, the initial-data selection uses the time integral of $\mathcal D_a$, rather than a pointwise-in-time upper bound for this quantity. Here the matrix $P$ is defined in \eqref{eq:unstable-coordinate-metric}; when $J>0$, by \eqref{eq:unstable-coordinate-lower-bound} and $2\kappa\geq d_*-\|P\|^{-1}$, we have the lower bound $2\langle A_{\rm u}\mathbf a,\mathbf a\rangle_P\geq-2\kappa|\mathbf a|_P^2$. If $J=0$, then every term involving $\mathbf a$ or $F_{\rm u}$ is zero.

We first specify one smallness condition under which all three estimates apply to the Maxwellian family \eqref{eq:positive-initial-Maxwellian}. We fix $C_{\rm rec}\geq1$ large enough for the recovery bounds
\begin{align}
	&\sup_x\left(|n-1|+|V|+|\Theta-1|+\sum_{j=1}^2d^j|\nabla_x^j(n,V,\Theta)|\right)\leq C_{\rm rec}\mathcal B_{\rm reg}^{1/2},\label{eq:roadmap-continuation-macroscopic}\\
	&\qquad\qquad\qquad\qquad\sup_{x,v}\langle\widetilde v\rangle^8|H/W|\leq C_{\rm rec}\mathcal B_{\rm reg}^{1/2},\label{eq:roadmap-continuation-pointwise}
\end{align}
proved in Lemmas~\ref{lem:single-spatial-allocation} and \ref{lem:B-pointwise}, respectively, for $(n,V,\Theta)\in\mathcal U$ and $\mathcal B\leq e_{\rm vel}$. We increase $C_{\rm rec}$, if necessary, to include the local macroscopic estimate from Lemma~\ref{lem:single-spatial-allocation}:
\begin{equation}\label{eq:selection-local-macro-bound}
	\|\mathbf z^\sharp(t)\|_{H^{K-1}(B_{R_2})}\leq C_{\rm rec}\mathcal E(t)^{1/2}.
\end{equation}
Let $\delta_*>0$ be the threshold in Proposition~\ref{thm:B-local} for the weighted pointwise norm with velocity exponent $8$. Let $b_0>0$ be the threshold in Proposition~\ref{thm:actual-spectral-forcing} for the norm on the left-hand side of \eqref{eq:selection-local-macro-bound}. We define
\begin{equation}\label{eq:roadmap-continuation-threshold}
	e_{\rm reg}:=\min\left\{\frac{\kappa_{\rm reg}e_{\rm vel}}2,\ \frac1{64C_{\rm rec}^2},\ \frac{\delta_*^2}{64C_{\rm rec}^2}\right\}>0,\quad e_{\rm sh}:=\min\left\{\frac{e_{\rm reg}}2,\ \frac{b_0^2}{4C_{\rm rec}^2}\right\}>0.
\end{equation}
The threshold $e_{\rm reg}$ gives the parameter bounds and the velocity supremum required for the continuation. 

\begin{lemma}\label{lem:selection-energy-validity}
	Let the fixed data be as in Propositions~\ref{thm:new-closed-pair} and \ref{prop:mixed-velocity-estimate}, and let $e_{\rm reg},e_{\rm sh}$ be given by \eqref{eq:roadmap-continuation-threshold}. For $0<\eps\leq\eps_{\rm vel}$, $0\leq S<T\leq T_{\rm kin}$ and $|\mathbf a_0|_P\leq R_{\rm in}$, let $f$ be a smooth, strictly positive, simultaneous-rotation-invariant solution to \eqref{Eq.Boltzmann} on $[0,S]\times\R^3_x\times\R^3_v$ with the initial datum \eqref{eq:positive-initial-Maxwellian}. Assume that $(n,V,\Theta)\in C([0,S];L^\infty_x)$, $\mathcal E,\mathcal B,\mathcal V_K\in C^1([0,S])$, and that all energy and dissipation integrals in \eqref{eq:new-kinetic-energies}, \eqref{eq:spatial-master-energies} and \eqref{eq:mixed-velocity-energy} are finite and continuous on $[0,S]$.
	
	If
	\begin{equation}\label{eq:shooting-bootstrap}
		\mathcal E(t)+\mathcal B_{\rm reg}(t)+|\mathbf a(t)|_P^2<e_{\rm sh},\qquad t\in[0,S],
	\end{equation}
	then, for every $t\in[0,S]$, we have
	\[\mathcal B+\mathcal V_K\leq e_{\rm vel}/4,\qquad \sup_x(|n-1|+|V|+|\Theta-1|)\leq1/8,\qquad \sup_{x,v}\langle\widetilde v\rangle^8|H/W|\leq\delta_*/8.\]
	In particular, $(n,V,\Theta)(t,x)\in\operatorname{int}\mathcal U$ for all $(t,x)\in[0,S]\times\R^3$, and the inequalities \eqref{eq:new-closed-inner}, \eqref{eq:augmented-global-inequality} and \eqref{eq:new-actual-mode-bound} hold on $[0,S]$ with constants independent of $T$, $\eps$, $S$ and $\mathbf a_0$.
\end{lemma}

See Subsection~\ref{subsec:selection-local} for the proof of Lemma~\ref{lem:selection-energy-validity}. For the local existence, regularity and continuous dependence of solutions with initial data \eqref{eq:positive-initial-Maxwellian} and the continuation of such solutions, see Lemmas~\ref{lem:selection-local-family} and \ref{lem:selection-continuation}.

\medskip\noindent\emph{The initial energies.}
By \eqref{eq:positive-initial-Maxwellian}, we have $H(0)=0$ for every $|\mathbf a_0|_P\leq R_{\rm in}$. In particular, by \eqref{eq:mixed-velocity-energy} and \eqref{eq:augmented-global},
\[\mathcal V_K(0)=0,\qquad \mathcal B_{\rm reg}(0)=\mathcal B(0).\]
We next estimate the initial energies in terms of $\mathbf a_0$.

\begin{lemma}[Initial energy bounds]\label{lem:selection-initial-energy}
	Let $0<T\leq T_{\rm kin}$, $0<\eps\leq\eps_{\rm vel}$ and let $f_0(x,v;\mathbf a_0)$ be given by \eqref{eq:positive-initial-Maxwellian} for $|\mathbf a_0|_P\leq R_{\rm in}$. The initial energies $\mathcal E(0;\mathbf a_0)$ and $\mathcal B_{\rm reg}(0;\mathbf a_0)$ defined by \eqref{eq:spatial-master-energies} and \eqref{eq:augmented-global} are continuously differentiable with respect to $\mathbf a_0$ on a neighborhood of $\{|\mathbf a_0|_P\leq R_{\rm in}\}$. Moreover, there exists a constant $C_I>0$, independent of $T$, $\eps$ and $\mathbf a_0$, such that
	\begin{equation}\label{eq:shooting-initial}
		\begin{aligned}
			&\qquad\qquad\qquad\qquad\quad\mathcal E(0;\mathbf 0)=\mathcal B_{\rm reg}(0;\mathbf 0)=0,\\
			&\mathcal E(0;\mathbf a_0)+\mathcal B_{\rm reg}(0;\mathbf a_0)\leq C_I|\mathbf a_0|_P^2,\qquad
			|\nabla_{\mathbf a_0}\mathcal E(0;\mathbf a_0)|_P\leq C_I|\mathbf a_0|_P,
		\end{aligned}
	\end{equation}
	where $\nabla_{\mathbf a_0}\mathcal E(0;\mathbf a_0)$ is defined by
	\[D_{\mathbf a_0}\mathcal E(0;\mathbf a_0)[h]=\langle\nabla_{\mathbf a_0}\mathcal E(0;\mathbf a_0),h\rangle_P,\qquad h\in\R^J.\]
	For $J=0$, we only consider $\mathbf a_0=0$ and omit the gradient estimate in \eqref{eq:shooting-initial}.
\end{lemma}

See Subsection~\ref{subsec:selection-local} for the proof of Lemma~\ref{lem:selection-initial-energy}.

\medskip\noindent\emph{Selection of the initial parameter.}
The spectral constants chosen in Section~\ref{subsec:euler-linearization} satisfy
\begin{equation}\label{eq:shooting-exponents}
	\max\{0,2\kappa\}<\delta<d_*<\min\{\eta,\lambda\}.
\end{equation}
The exponent $d_*$ specifies the decay of $\mathcal E+|\mathbf a|_P^2$. The inequalities $d_*<\lambda$ and $d_*<\eta$ place this decay below the damping and error exponents in \eqref{eq:new-closed-inner}. The exponent $\delta$ is used in the weighted dissipation integral; the strict inequality $\delta>2\kappa$ makes the backward coordinate integral absolutely convergent.

\begin{theorem}[Correction of the unstable initial data]\label{thm:coupled-shooting}
	Let $0<T\leq T_{\rm kin}$ and let the spectral exponents be chosen as in \eqref{eq:shooting-exponents}. There exist constants $\eps_{\rm sh}\in(0,\eps_{\rm vel}]$ and $C_{\rm sh}>0$ such that for every $0<\eps\leq\eps_{\rm sh}$, we can find $\mathbf a_0^*\in\R^J$ satisfying
	\[|\mathbf a_0^*|_P\leq\min\{R_{\rm in},C_{\rm sh}\eps^{1/2}T^{\eta/2}\}\]
	such that the Boltzmann equation \eqref{Eq.Boltzmann} with initial datum $f_0(x,v;\mathbf a_0^*)$ given by \eqref{eq:positive-initial-Maxwellian} admits a strictly positive simultaneous-rotation-invariant solution $f\in C^\infty([0,T)\times\R^3_x\times\R^3_v)$. For all $(t,x)\in[0,T)\times\R^3$, we have
	\[(n,V,\Theta)(t,x)\in\operatorname{int}\mathcal U,\qquad \mathcal E(t)+\mathcal B_{\rm reg}(t)+|\mathbf a(t)|_P^2<e_{\rm sh}.\]
	Moreover,
	\begin{align}
		&\qquad\qquad\quad\mathcal E(t)+|\mathbf a(t)|_P^2\leq C_{\rm sh}\eps T^{\eta-d_*}(T-t)^{d_*},\qquad t\in[0,T),\label{eq:shooting-pointwise-result}\\
		&\qquad\quad\sup_{0\leq t<T}\mathcal B_{\rm reg}(t)+\frac{c_{\rm reg}}\eps\int_0^T\mathcal D_{\rm reg}(t)\dd t\leq C_{\rm sh}\eps T^\eta,\label{eq:shooting-pointwise-result-2}\\
		&\sup_{0\leq t<T}(T-t)^{-\delta}\mathcal E(t)+\frac c\eps\int_0^T(T-t)^{-\delta}\mathcal D_a(t)\dd t\leq C_{\rm sh}\eps T^{\eta-\delta}.\label{eq:shooting-weighted-result}
	\end{align}
	For $J>0$, we have the integral representation
	\begin{equation}\label{eq:selected-terminal-integral}
		\mathbf a(t)=-\int_t^T\mathcal U_{\rm u}(t,q)F_{\rm u}(q)\dd q,\qquad \mathcal U_{\rm u}(t,q):=\exp\left(A_{\rm u}\log\frac{T-q}{T-t}\right),\quad t\leq q<T,
	\end{equation}
	with the integral converging absolutely. For $J=0$, we take $\mathbf a_0^*=0$. The constants $\eps_{\rm sh}$ and $C_{\rm sh}$ depend only on the fixed profile, differentiation orders, weights, spectral data and $T$.
\end{theorem}

In Subsection~\ref{subsec:selection-brouwer}, we prove Theorem~\ref{thm:coupled-shooting} by constructing an exit map and applying Brouwer's fixed point theorem.

\subsection{Proof of the main theorem}\label{subsec:main-proof}
In addition to the conclusions of Theorem~\ref{thm:main}, we will prove the relative entropy inequality
\begin{equation}\label{eq:new-relative-entropy-conclusion}
	\mathcal H(f(t)\mid M_\infty)\leq\mathcal H(f_0\mid M_\infty)<\infty,\qquad 0\leq t<T.
\end{equation}
\begin{proof}[Proof of Theorem~\ref{thm:main}]
	We fix the parameters in Propositions~\ref{thm:new-closed-pair} and \ref{prop:mixed-velocity-estimate} and choose $T\in(0,T_{\rm kin}]$ independently of $\eps$. Take $\eps_0\in(0,\eps_{\rm sh}]$. For each $\eps\in(0,\eps_0]$, let $\mathbf a_0^*(\eps)$ and $f$ be given by Theorem~\ref{thm:coupled-shooting}, with $f_0$ defined by \eqref{eq:positive-initial-Maxwellian} for $\mathbf a_0=\mathbf a_0^*(\eps)$. In particular,
	\[f>0,\qquad (n,V,\Theta)\in\operatorname{int}\mathcal U,\qquad 0\leq t<T,\]
	and the estimates \eqref{eq:shooting-pointwise-result}--\eqref{eq:shooting-weighted-result} hold.
	
	Recall that
	\[M=n\Theta^{-3/2}\mu(\widetilde v),\qquad W=n^{1/2}\Theta^{-3/4}w_{m_\eps,\ell}(\widetilde v),\qquad \widetilde v=\frac{v-u_f}{\sqrt{T_f}}.\]
	By $\Phi_{m,\ell}\geq0$ and \eqref{eq:fixed-tail-weight}, we have
	\[0<w_{m,\ell}(v)\leq(2\pi)^{-3/4},\qquad m\geq m_0,\quad v\in\R^3.\]
	Since $(n,V,\Theta)\in\mathcal U$, we get
	\[\sup_{0\leq t<T}\big(\|M(t)\|_{L^\infty_{x,v}}+\|W(t)\|_{L^\infty_{x,v}}\big)\leq C.\]
	By Lemma~\ref{lem:B-pointwise} and \eqref{eq:shooting-pointwise-result-2}, we have
	\[\|H(t)\|_{L^\infty_{x,v}}\leq C\|H(t)/W(t)\|_{L^\infty_{x,v}}\leq C\mathcal B_{\rm reg}(t)^{1/2}\leq C\eps^{1/2}T^{\eta/2},\qquad t\in[0,T).\]
	Thus, by $f=M+H$, we obtain \eqref{eq:main-uniform-bound}.
	
	By \eqref{eq:positive-initial-Maxwellian}, for every pair of multi-indices $\alpha,\beta$, there exist $c>0$ and $C_{\alpha,\beta}>0$ such that
	\[|\partial_x^\alpha\partial_v^\beta f_0(x,v)|\leq C_{\alpha,\beta}\mathrm e^{-c|v|^2},\qquad (x,v)\in\R^3\times\R^3.\]
	By Propositions~\ref{prop:B-time-regularity} and \ref{thm:B-restart}, the local existence time can be chosen independently of the differentiation order. Successively applying the two propositions on compact time intervals in $[0,T)$, we obtain
	\[f\in C^\infty([0,T)\times\R^3_x\times\R^3_v).\]
	
	Fix $R>0$ and take $t$ sufficiently close to $T$ such that $R(T-t)^{1/r}<R_{\rm e}$ and $J_x=1$ for $|x|\leq R(T-t)^{1/r}$. For $|y|\leq R$, we have
	\[\langle R\rangle^{-a}\leq\Omega_a^*/\Omega_0^*=\langle y\rangle^{-a}\leq1.\]
	By \eqref{eq:roadmap-natural-Euler-energy}, Lemma~\ref{lem:single-spatial-allocation} and the Sobolev embedding in $y$, we have
	\[\sup_{|x|\leq R(T-t)^{1/r}}\left(|n-1|+|V|+|\Theta-1|+\|h/W\|_{L^2_{\widehat v}}\right)\leq C_R\mathcal E(t)^{1/2}.\]
	Recall that
	\[\widetilde M=n\Theta^{-3/2}\mu\left(\frac{\widehat v-V}{\sqrt\Theta}\right),\qquad F=\widetilde M+h.\]
	By the mean-value formula and $W\leq C$, for $(n,V,\Theta)\in\mathcal U$, we have
	\[\|F-\mu\|_{L^2_{\widehat v}}\leq C\left(|n-1|+|V|+|\Theta-1|+\|h/W\|_{L^2_{\widehat v}}\right).\]
	Consequently, by \eqref{eq:shooting-pointwise-result}, we obtain
	\[\sup_{|x|\leq R(T-t)^{1/r}}\left(|n-1|+\|F-\mu\|_{L^2_{\widehat v}}\right)\leq C_R\eps^{1/2}T^{(\eta-d_*)/2}(T-t)^{d_*/2},\]
	which gives \eqref{eq:main-rescaled-limit}. Moreover, by $\rho_f=\theta^3n$ and \eqref{eq:profile}, we have
	\[\rho_f(t,0)=\vartheta(0)^3(T-t)^{-3\alpha_{\rm ss}}\left(1+O\left(\eps^{1/2}T^{(\eta-d_*)/2}(T-t)^{d_*/2}\right)\right).\]
	Hence \eqref{eq:main-density} holds.
	
	We next prove \eqref{eq:main-holder}. By the radial symmetry and smoothness of the profiles and \eqref{Eq.Euler-monatomic-ss}, we have $U(0)=0$ and $\nabla U(0)=-\alpha_{\rm ss}I$. Take a fixed $y_*\ne0$ sufficiently close to the origin such that $U(y_*)\ne0$, and let
	\[\mathfrak M_y(w):=\mu\left(\frac{\vartheta(0)w-U(y)}{\vartheta(y)}\right),\qquad y\in\{0,y_*\}.\]
	Since $U(0)=0\ne U(y_*)$, there exists $R>0$ such that $c_R:=\|\mathfrak M_{y_*}-\mathfrak M_0\|_{L^2(B_R)}>0$. For $y\in\{0,y_*\}$, by \eqref{eq:main-rescaled-limit} and the change of variables $\widehat v=(\vartheta(0)w-U(y))/\vartheta(y)$, we have
	\[f\big(t,(T-t)^{1/r}y,\vartheta(0)(T-t)^{-\alpha_{\rm ss}}w\big)\longrightarrow\mathfrak M_y(w)\quad\text{in }L^2(B_R)\]
	as $t\uparrow T$. Thus, for $t$ sufficiently close to $T$, we get
	\begin{align*}
		&\sup_{v\in\R^3}|f(t,(T-t)^{1/r}y_*,v)-f(t,0,v)|\\
		&\quad\geq\frac1{|B_R|^{1/2}}\left\|f\big(t,(T-t)^{1/r}y_*,\vartheta(0)(T-t)^{-\alpha_{\rm ss}}w\big)-f\big(t,0,\vartheta(0)(T-t)^{-\alpha_{\rm ss}}w\big)\right\|_{L^2_w(B_R)}\\
		&\quad\geq\frac{c_R}{2|B_R|^{1/2}}.
	\end{align*}
	For every $0<\delta_x\leq1$, by \eqref{eq:holder-definition}, we obtain
	\[[f(t)]_{C_x^{\delta_x}L_v^\infty}\geq\frac{c_R}{2|B_R|^{1/2}|y_*|^{\delta_x}}(T-t)^{-\delta_x/r}.\]
	
		To estimate the relative moments and entropy, we fix $S\in(0,T)$. By the estimates for the Euler background and the definition of $\Omega_0^*$, on $[0,S]\times\R^3$, we have
	\[|u|+\theta+\theta^{-1}\leq C_S,\qquad \Omega_0^*\asymp_S\langle x\rangle^{2M_x-3}.\]
	Since $M_x>3$, we have
	\[\sup_{0\leq t\leq S}\int_{\R^3}(\Omega_0^*)^{-1}\dd x\leq C_S.\]
	Moreover, by $(n,V,\Theta)\in\mathcal U$ and the polynomial decay of $w_{m_\eps,\ell}$, for $q=1,2,3$, we have
	\[\sup_{(t,x)\in[0,S]\times\R^3}\int\langle\widehat v\rangle^{2q}W^2\dd\widehat v\leq C_S.\]
	By the change of variables $v=u+\theta\widehat v$ and the Cauchy--Schwarz inequality, for $t\in[0,S]$, we obtain
	\begin{align*}
		\iint\langle v\rangle^q|H|\dd v\dd x
		&\leq C_S\iint\langle\widehat v\rangle^q|H|\dd\widehat v\dd x\\
		&\leq C_S\left(\iint\Omega_0^*|H/W|^2\dd\widehat v\dd x\right)^{1/2}
		\left(\iint(\Omega_0^*)^{-1}\langle\widehat v\rangle^{2q}W^2\dd\widehat v\dd x\right)^{1/2}\\
		&\leq C_S\mathcal B_{\rm reg}^{1/2}.
	\end{align*}
	Moreover, using $\theta^3\langle u+\theta\widehat v\rangle^qW^2\leq C_S\Omega_0^*$, we get
	\[\iint\langle v\rangle^q|H|^2\dd v\dd x\leq C_S\iint\Omega_0^*|H/W|^2\dd\widehat v\dd x\leq C_S\mathcal B_{\rm reg}.\]
	By the mean-value formula for $M$ in the macroscopic variables, we have
	\[\|\langle v\rangle^q(M-\bar M)\|_{L^1_v}+\|\langle v\rangle^{q/2}(M-\bar M)\|_{L^2_v}\leq C_S|\mathbf z|,\qquad q=1,2,3.\]
	Since
	\[\int\Omega_0^*|\mathbf z|^2\dd x\lesssim\mathcal B_{\rm reg},\qquad \int(\Omega_0^*)^{-1}\dd x\leq C_S,\]
	we obtain
	\[\|\langle v\rangle^q(M-\bar M)\|_{L^1_{x,v}}+\|\langle v\rangle^{q/2}(M-\bar M)\|_{L^2_{x,v}}\leq C_S\mathcal B_{\rm reg}^{1/2}.\]
	Writing $f-M_\infty=H+(M-\bar M)+(\bar M-M_\infty)$ and using the compact spatial support and Gaussian decay in velocity of $\bar M-M_\infty$, we obtain
	\[\sup_{0\leq t\leq S}\iint\left(\langle v\rangle^3|f-M_\infty|+\langle v\rangle|f-M_\infty|^2\right)\dd v\dd x<\infty.\]
	In particular, the relative mass and kinetic-energy moments in Theorem~\ref{thm:main} are absolutely finite.
	
	For $m>0$ and $g\geq-m$, we have the elementary inequality
	\[0\leq(m+g)\log\frac{m+g}{m}-g\leq |g|+g^2+|g|\log(1+m^{-1}).\]
	Taking $m=M_\infty$ and $g=f-M_\infty$, and using $\log(1+M_\infty^{-1})\lesssim\langle v\rangle^2$, we obtain
	\begin{equation}\label{eq:main-proof-entropy-flux}
		\sup_{0\leq t\leq S}\iint(1+|v|)\left(f\log\frac f{M_\infty}-f+M_\infty\right)\dd v\dd x<\infty.
	\end{equation}
	
		For $(t,x,v)\in[0,S]\times\R^3\times\R^3$, by $0<\gamma\leq2$ and the bounds on the mass and kinetic-energy densities, we have
	\[\nu_f(t,x,v)\lesssim\langle v\rangle^\gamma\int\langle v_*\rangle^\gamma f(t,x,v_*)\dd v_*\leq C_S\langle v\rangle^\gamma\leq C_S\langle v\rangle^2.\]
	By \eqref{eq:positive-initial-Maxwellian}, there exist constants $c,C>0$ such that $f_0(x,v)\geq c\mathrm e^{-C|v|^2}$ for all $(x,v)\in\R^3\times\R^3$. Integrating \eqref{Eq.Boltzmann} along characteristics and using $Q^+(f,f)\geq0$, for $t\in[0,S]$, we obtain
	\begin{align*}
		f(t,x,v)&\geq f_0(x-tv,v)\exp\left(-\frac1\eps\int_0^t\nu_f(q,x-(t-q)v,v)\dd q\right)\\
		&\geq c\exp\left(-C|v|^2-\frac{C_St}\eps\langle v\rangle^2\right)\geq c_{S,\eps}\mathrm e^{-C_{S,\eps}|v|^2}.
	\end{align*}
	By \eqref{eq:main-uniform-bound}, we get
	\[|\log f|+|\log(f/M_\infty)|\leq C_{S,\eps}\langle v\rangle^2\quad\text{on }[0,S]\times\R^3\times\R^3.\]
	Since $\log M_\infty\in\operatorname{span}\{1,v_1,v_2,v_3,|v|^2\}$, by the conservation properties of $Q$ and the $H$-theorem, we have
	\[\int Q(f,f)\log(f/M_\infty)\dd v=\int Q(f,f)\log f\dd v\leq0.\]
	
	Let $\chi_R^{\rm ent}\in C_c^\infty(B_{2R})$ satisfy
	\[0\leq\chi_R^{\rm ent}\leq1,\qquad \chi_R^{\rm ent}\equiv1\ \text{on }B_R,\qquad |\nabla_x\chi_R^{\rm ent}|\leq C/R.\]
	Multiplying \eqref{Eq.Boltzmann} by $\chi_R^{\rm ent}\log(f/M_\infty)$ and integrating on $[0,t]\times\R^3\times\R^3$, we obtain
	\begin{align*}
		&\iint\chi_R^{\rm ent}\left(f(t)\log\frac{f(t)}{M_\infty}-f(t)+M_\infty\right)\dd v\dd x\leq\iint\chi_R^{\rm ent}\left(f_0\log\frac{f_0}{M_\infty}-f_0+M_\infty\right)\dd v\dd x\\
		&\qquad\qquad\qquad\qquad\qquad\qquad\qquad\quad+\frac CR\int_0^t\iint|v|\left(f\log\frac f{M_\infty}-f+M_\infty\right)\dd v\dd x\dd q.
	\end{align*}
	Here the integrations in $v$ can be justified by first introducing velocity cutoffs and then removing the cutoffs using the polynomial decay of $f$ and $|\log(f/M_\infty)|\leq C_{S,\eps}\langle v\rangle^2$. By \eqref{eq:main-proof-entropy-flux}, we have
	\[\frac CR\int_0^t\iint|v|\left(f\log\frac f{M_\infty}-f+M_\infty\right)\dd v\dd x\dd q\leq\frac{C_{S,\eps}}R.\]
	Letting $R\to\infty$ and using the dominated convergence theorem, we obtain
	\[\mathcal H(f(t)\mid M_\infty)\leq\mathcal H(f_0\mid M_\infty),\qquad t\in[0,S].\]
	Since $S\in(0,T)$ is arbitrary, we obtain \eqref{eq:new-relative-entropy-conclusion}.
\end{proof}

For the Euler background, the relative entropy has a direct fluid interpretation. We define $\rho:=\theta^3$, $\rho_\infty:=\theta_\infty^3$ and $\mathfrak e(\rho):=\frac32\rho^{5/3}$. By direct Gaussian integration, we have
\[\theta_\infty^2\HH(\bar M\mid M_\infty)=\int e_{\rm rel}\dd x,\qquad e_{\rm rel}:=\frac12\rho|u|^2+\mathfrak e(\rho)-\mathfrak e(\rho_\infty)-\mathfrak e'(\rho_\infty)(\rho-\rho_\infty).\]
By the Euler mass and energy equations and a direct computation, we have
\[\partial_te_{\rm rel}+\Div\bigl((e_{\rm rel}+\rho^{5/3}-\rho_\infty^{5/3})u\bigr)=0.\]
The density $e_{\rm rel}$ and its flux have compact spatial support. Hence, $\HH(\bar M(t)\mid M_\infty)$ is finite and independent of $t<T$.

\section{Backward construction of the Euler correction}\label{Sec.backward-Euler-correction}

Fix a profile $(u_{\rm s},\theta_{\rm s})$ and $r\in\{r_n\}_{n\geq1}$ as in Theorem~\ref{Thm.profiles}. We construct the correction in Proposition~\ref{Thm.Euler-solution} by solving \eqref{Eq.backward-regularized-system} with terminal value zero. The constants in this section may depend on the profile and cutoffs and are independent of $T,R_0$.

\subsection{Profile estimates and the terminal problem}\label{subsec:euler-terminal}

\begin{lemma}[Profile derivatives and truncation errors]\label{lem:euler-profile-errors}
	Let $0<T\leq1$, $R_0>10$ and $d=((T-t)^{2/r}+|x|^2)^{1/2}$. For $t\in[0,T)$ and $j,k\in\mathbb N_0$,
	\begin{equation}\label{eq:background-hierarchy-base}
		c d^{1-r}\leq\theta_{\rm p}\leq C d^{1-r},\qquad |\nabla_x^k u_{\rm p}|+|\nabla_x^k\theta_{\rm p}|\leq C_k d^{1-r-k},
	\end{equation}
	and
	\begin{equation}\label{eq:background-hierarchy-profile}
		|\partial_t^j\nabla_x^k u_{\rm p}|+|\partial_t^j\nabla_x^k\theta_{\rm p}|\leq C_{j,k}d^{1-r-jr-k}.
	\end{equation}
	On $R_0/8\leq|x|\leq2R_0$, the profile $(u_{\rm p},\theta_{\rm p})$ extends smoothly to $t=T$ and satisfies
	\begin{equation}\label{Eq.profile-mixed-bound}
		|\partial_t^j\nabla_x^k u_{\rm p}|+|\partial_t^j\nabla_x^k\theta_{\rm p}|\leq C_{j,k}R_0^{1-(j+1)r-k}.
	\end{equation}
	Moreover, the errors extend smoothly to $[0,T]\times\R^3$ and satisfy
	\begin{align}
		&\supp\cE_u(t)\cup\supp\cE_\theta(t)\subset\{R_0/2\leq|x|\leq R_0\},\label{Eq.error-annulus-estimate}\\
		&\quad |\partial_t^j\nabla_x^k\cE_u|+|\partial_t^j\nabla_x^k\cE_\theta|\leq C_{j,k}R_0^{1-(j+2)r-k}.\label{Eq.error-mixed-estimate}
	\end{align}
\end{lemma}
\begin{proof}
	By \eqref{Eq.u_p-theta_p-def} and \eqref{eq:profile-decay}, we get \eqref{eq:background-hierarchy-base}, using
	\[(T-t)^{-\alpha_{\rm ss}-k/r}\left\langle\frac{x}{(T-t)^{1/r}}\right\rangle^{1-r-k}=d^{1-r-k}.\]
	To obtain the time derivatives, we use the equations for the untruncated fields
	\begin{equation}\label{Eq.profile-Euler-u-theta}
		\begin{cases}
			\partial_tu_{\rm p}+u_{\rm p}\cdot\nabla u_{\rm p}+5\theta_{\rm p}\nabla\theta_{\rm p}=0,\\
			\partial_t\theta_{\rm p}+u_{\rm p}\cdot\nabla\theta_{\rm p}+\dfrac13\theta_{\rm p}\Div u_{\rm p}=0.
		\end{cases}
	\end{equation}
	The case $j=0$ is \eqref{eq:background-hierarchy-base}. Assume the bounds in \eqref{eq:background-hierarchy-profile} up to time order $j$ for all $k$. Apply $\partial_t^j\nabla_x^k$ to \eqref{Eq.profile-Euler-u-theta} and use the product rule to get
	\begin{align*}
		&|\partial_t^{j+1}\nabla_x^ku_{\rm p}|+|\partial_t^{j+1}\nabla_x^k\theta_{\rm p}|\\
		&\quad\leq C_{j,k}\sum_{\substack{j_1+j_2=j\\k_1+k_2=k}}d^{1-r-j_1r-k_1}d^{1-r-j_2r-k_2-1}\leq C_{j,k}d^{1-r-(j+1)r-k}.
	\end{align*}
	This proves \eqref{eq:background-hierarchy-profile}. On $R_0/8\leq|x|\leq2R_0$ we have $d\asymp R_0$, hence \eqref{Eq.profile-mixed-bound}. For $0\leq t_1<t_2<T$, \eqref{Eq.profile-mixed-bound} at order $j+1$ gives
	\begin{equation}\label{eq:profile-terminal-lipschitz}
		\begin{aligned}
			&\|\partial_t^j\nabla_x^k(u_{\rm p},\theta_{\rm p})(t_2)-\partial_t^j\nabla_x^k(u_{\rm p},\theta_{\rm p})(t_1)\|_{L^\infty(R_0/8\leq|x|\leq2R_0)}\\
			&\qquad\leq C_{j,k}(t_2-t_1)R_0^{1-(j+2)r-k}.
		\end{aligned}
	\end{equation}
	Thus all the mixed derivatives have uniform limits at $t=T$ on this annulus.

	Put $c_{R_0}=R_0^{1-r}$, $\bar u_{\rm p}=\chi_{R_0}u_{\rm p}$ and $\bar\theta_{\rm p}=c_{R_0}+\chi_{R_0}(\theta_{\rm p}-c_{R_0})$. Subtracting $\chi_{R_0}$ times \eqref{Eq.profile-Euler-u-theta} from the equations for $(\bar u_{\rm p},\bar\theta_{\rm p})$ gives
	\begin{align*}
		\cE_u&=\chi_{R_0}(\chi_{R_0}-1)u_{\rm p}\cdot\nabla u_{\rm p}+5\chi_{R_0}(\bar\theta_{\rm p}-\theta_{\rm p})\nabla\theta_{\rm p}\\
		&\quad+\chi_{R_0}(u_{\rm p}\cdot\nabla\chi_{R_0})u_{\rm p}+5\bar\theta_{\rm p}(\theta_{\rm p}-c_{R_0})\nabla\chi_{R_0},\\
		\cE_\theta&=\chi_{R_0}(\chi_{R_0}-1)u_{\rm p}\cdot\nabla\theta_{\rm p}+\tfrac13\chi_{R_0}(\bar\theta_{\rm p}-\theta_{\rm p})\Div u_{\rm p}\\
		&\quad+\bigl(\chi_{R_0}(\theta_{\rm p}-c_{R_0})+\tfrac13\bar\theta_{\rm p}\bigr)u_{\rm p}\cdot\nabla\chi_{R_0}.
	\end{align*}
	Since $\chi_{R_0}(\bar\theta_{\rm p}-\theta_{\rm p})=\chi_{R_0}(\chi_{R_0}-1)(\theta_{\rm p}-c_{R_0})$, every term contains $\chi_{R_0}(\chi_{R_0}-1)$ or $\nabla_x\chi_{R_0}$. This gives \eqref{Eq.error-annulus-estimate}. On its support, \eqref{Eq.profile-mixed-bound} and $|\nabla_x^k\chi_{R_0}|\leq C_kR_0^{-k}$ give
	\begin{align*}
		|\partial_t^j\nabla_x^k\cE_u|+|\partial_t^j\nabla_x^k\cE_\theta|&\leq C_{j,k}\sum_{j_1+j_2=j}\sum_{k_0+k_1+k_2=k+1}R_0^{-k_0}R_0^{1-(j_1+1)r-k_1}R_0^{1-(j_2+1)r-k_2}\\
		&\leq C_{j,k}R_0^{1-(j+2)r-k}.
	\end{align*}
	The limits from \eqref{eq:profile-terminal-lipschitz} give the smooth extension of both errors to $[0,T]\times\R^3$.
\end{proof}

We choose a radial function $\zeta\in C^\infty(\R^3;[0,1])$ such that $\zeta=0$ on $|X|\leq1/8$ and $\zeta=1$ on $|X|\geq1/4$. We define
\begin{align*}
	U_{R_0}(t,X)&:=\zeta(X)\bar u_{\rm p}(t,R_0X),\\
	\Theta_{R_0}(t,X)&:=\sqrt{15}R_0^{1-r}+\sqrt{15}\zeta(X)\bigl(\bar\theta_{\rm p}(t,R_0X)-R_0^{1-r}\bigr),\\
	F_{R_0}(t,X)&:=\bigl(\cE_u(t,R_0X),\sqrt{15}\cE_\theta(t,R_0X)\bigr).
\end{align*}
Here the values at $t=T$ are given by Lemma~\ref{lem:euler-profile-errors}. These coefficients are smooth, agree with $(\bar u_{\rm p},\sqrt{15}\bar\theta_{\rm p})(t,R_0X)$ for $|X|\geq1/4$, and satisfy, for $j,m\in\mathbb N_0$,
\begin{equation}\label{Eq.modified-coefficient-bound}
	\begin{aligned}
		&\|\partial_t^jU_{R_0}(t)\|_{H^{m+1}_X}+\|\partial_t^j(\Theta_{R_0}(t)-\sqrt{15}R_0^{1-r})\|_{H^{m+1}_X}\leq C_{j,m}R_0^{1-r-jr},\\
		&\qquad\qquad\qquad\|\partial_t^jF_{R_0}(t)\|_{H^m_X}\leq C_{j,m}R_0^{1-(j+2)r}.
	\end{aligned}
\end{equation}
Here and below, $H^m_X:=H^m(\R^3_X)$, and a norm of a vector denotes the Sobolev norm of all its components. Moreover, by \eqref{Eq.error-annulus-estimate}, we have $\supp F_{R_0}(t)\subset\{1/2\leq|X|\leq1\}$. We consider the terminal problem
\begin{equation}\label{Eq.backward-regularized-system}
	\begin{cases}
		\partial_t W+\frac1{R_0}(U_{R_0}+W)\cdot\nabla_XW+\frac1{3R_0}(\Theta_{R_0}+P)\nabla_XP\\
		\qquad=-\frac1{R_0}W\cdot\nabla_XU_{R_0}-\frac1{3R_0}P\nabla_X\Theta_{R_0}-\cE_u(t,R_0X),\\
		\partial_t P+\frac1{R_0}(U_{R_0}+W)\cdot\nabla_XP+\frac1{3R_0}(\Theta_{R_0}+P)\Div_XW\\
		\qquad=-\frac1{R_0}W\cdot\nabla_X\Theta_{R_0}-\frac1{3R_0}P\Div_XU_{R_0}-\sqrt{15}\cE_\theta(t,R_0X),\\
		(W,P)(T,X)=(0,0).
	\end{cases}
\end{equation}

\begin{lemma}\label{lem:euler-terminal-problem}
	There exist $T_0\in(0,1]$ and $R_*>10$ such that, for $0<T\leq T_0$ and $R_0\geq R_*$, \eqref{Eq.backward-regularized-system} has a smooth solution $Z=(W,P)$ on $[0,T]\times\R^3$, unique in $C([0,T];H^4_X)\cap C^1([0,T];H^3_X)$. On the same interval for every $m\in\mathbb N_0$, we have
	\begin{equation}\label{Eq.backward-Hs-bound}
		\|Z(t)\|_{H^m_X}\leq C_m(T-t)R_0^{1-2r},\qquad 0\leq t\leq T.
	\end{equation}
\end{lemma}
\begin{proof}
	With $e_i$ denoting the $i$th coordinate vector, the matrices of $\partial_{X_i}Z$ in \eqref{Eq.backward-regularized-system} satisfy
	\begin{equation}\label{eq:terminal-principal-matrices}
		A_i=R_0^{-1}\begin{pmatrix}
			(U_{R_0,i}+W_i)I_3&(\Theta_{R_0}+P)e_i/3\\
			(\Theta_{R_0}+P)e_i^{\mathsf T}/3&U_{R_0,i}+W_i
		\end{pmatrix}=A_i^{\mathsf T}.
	\end{equation}
	Therefore we can apply \cite{Kato1975,Majda1984} at order $m_0=4$ with $Z(T)=0$. By \eqref{Eq.modified-coefficient-bound} and $H^4_X\hookrightarrow W^{1,\infty}_X$, we have
	\[-\frac{\mathrm d}{\mathrm dt}\|Z\|_{H^{m_0}_X}\leq C_{m_0}\bigl(R_0^{-r}+R_0^{-1}\|Z\|_{H^{m_0}_X}\bigr)\|Z\|_{H^{m_0}_X}+C_{m_0}R_0^{1-2r}.\]
	Under the bootstrap bound $\|Z(t)\|_{H^{m_0}_X}\leq2C_{m_0}(T-t)R_0^{1-2r}$, we have
	\begin{align*}
		R_0^{-r}+R_0^{-1}\|Z(t)\|_{H^{m_0}_X}&\leq CR_0^{-r},\\
		\|Z(t)\|_{H^{m_0}_X}&\leq C_{m_0}R_0^{1-2r}\int_t^T\mathrm e^{CR_0^{-r}(q-t)}\dd q\\
		&\leq C_{m_0}(T-t)R_0^{1-2r}\mathrm e^{CT_0R_*^{-r}}\leq\tfrac32C_{m_0}(T-t)R_0^{1-2r},
	\end{align*}
	after choosing $T_0\leq1$ and $R_*>10$. Therefore, by continuation, we obtain the solution on $[0,T]$ and \eqref{Eq.backward-Hs-bound} for $m\leq m_0$.
	
	For $m>m_0$, by the tame estimate and \eqref{Eq.modified-coefficient-bound}, we have
	\begin{align*}
		-\frac{\mathrm d}{\mathrm dt}\|Z\|_{H^m_X}&\leq C_m\bigl(R_0^{-r}+R_0^{-1}\|Z\|_{H^{m_0}_X}\bigr)\|Z\|_{H^m_X}+C_mR_0^{-r}\|Z\|_{H^{m_0}_X}+C_mR_0^{1-2r}\\
		&\leq C_mR_0^{-r}\|Z\|_{H^m_X}+C_mR_0^{1-2r},\\
		\|Z(t)\|_{H^m_X}&\leq C_mR_0^{1-2r}\int_t^T\mathrm e^{C_mR_0^{-r}(q-t)}\dd q\leq C_m(T-t)R_0^{1-2r}.
	\end{align*}
	The $H^4_X$ bound fixes $T_0,R_*$ for all $m$. Equation \eqref{Eq.backward-regularized-system} gives time smoothness up to the endpoints.
\end{proof}

\subsection{Localization and completion of the Euler construction}\label{subsec:euler-correction-completion}

\begin{lemma}\label{lem:euler-correction-support}
	Take $Z=(W,P)$ from Lemma~\ref{lem:euler-terminal-problem}. For $R_*$ sufficiently large, $R_0\geq R_*$ and $0\leq t\leq T$,
	\[\supp Z(t)\subset\{3/8\leq|X|\leq9/8\},\qquad (u_{\rm r},\theta_{\rm r})(t,x):=(W,P/\sqrt{15})(t,x/R_0).\]
	Then $(u_{\rm r},\theta_{\rm r})$ is a radial solution of \eqref{Eq.u_r-theta_r-eq} with zero terminal data and satisfies \eqref{Eq.remainder-support}, \eqref{Eq.remainder-pointwise-bound} and
	\begin{equation}\label{eq:euler-correction-mixed}
		\|\partial_t^j\nabla_x^ku_{\rm r}(t)\|_{L^\infty_x}+\|\partial_t^j\nabla_x^k\theta_{\rm r}(t)\|_{L^\infty_x}\leq C_{j,k}R_0^{1-(j+1)r-k}\qquad(j\geq1,\ k\geq0).
	\end{equation}
\end{lemma}
\begin{proof}
	By \eqref{eq:terminal-principal-matrices}, \eqref{Eq.modified-coefficient-bound} and \eqref{Eq.backward-Hs-bound},
	\begin{align*}
		\sup_{|\xi|=1}\Big\|\sum_i\xi_iA_i(t,X)\Big\|_{\mathrm{op}}
		&=R_0^{-1}\bigl(|U_{R_0}+W|+\tfrac13|\Theta_{R_0}+P|\bigr)
		\leq CR_0^{-r}.
	\end{align*}
	By finite propagation from $Z(T)=0$, with $\supp F_{R_0}(q)\subset\{1/2\leq|X|\leq1\}$ and $\int_t^T CR_0^{-r}\dd q\leq1/8$,
	\begin{equation}\label{eq:backward-correction-annulus}
		\supp Z(t)\subset\{1/2-C(T-t)R_0^{-r}\leq |X|\leq1+C(T-t)R_0^{-r}\}\subset\{3/8\leq|X|\leq9/8\}.
	\end{equation}
	By \eqref{eq:backward-correction-annulus}, $\zeta=1$ on $\supp Z(t)$. Use $x=R_0X$ and $P=\sqrt{15}\theta_{\rm r}$ to get \eqref{Eq.u_r-theta_r-eq} and \eqref{Eq.remainder-support}. Rotation invariance and uniqueness give radiality.

	By Sobolev embedding,
	\[\|\partial_t^j\nabla_x^k(u_{\rm r},\theta_{\rm r})\|_{L_x^\infty}\leq C_kR_0^{-k}\|\partial_t^jZ\|_{H_X^{k+2}}.\]
	For $j=0$ use \eqref{Eq.backward-Hs-bound} to get \eqref{Eq.remainder-pointwise-bound}. We prove
	\begin{equation}\label{eq:terminal-correction-time}
		\|\partial_t^jZ(t)\|_{H_X^m}\leq C_{j,m}R_0^{1-(j+1)r}\qquad(j,m\in\mathbb N_0).
	\end{equation}
	Start with \eqref{Eq.backward-Hs-bound} and $T\leq1$. Apply $\partial_t^j$ to \eqref{Eq.backward-regularized-system} and use \eqref{Eq.modified-coefficient-bound} and the $H_X^{m+2}$ product estimate to get
	\begin{align*}
		\|\partial_t^{j+1}Z\|_{H_X^m}
		&\leq C_{j,m}R_0^{-1}\sum_{\ell=0}^j R_0^{1-(\ell+1)r}R_0^{1-(j-\ell+1)r}+\|\partial_t^jF_{R_0}\|_{H_X^m}\\
		&\leq C_{j,m}R_0^{1-(j+2)r},
	\end{align*}
	where the constant part of $\Theta_{R_0}$ is $\sqrt{15}R_0^{1-r}$. Rescale \eqref{eq:terminal-correction-time} to get \eqref{eq:euler-correction-mixed}.
\end{proof}

\begin{lemma}\label{lem:euler-background-bounds}
	For the correction in Lemma~\ref{lem:euler-correction-support}, define $(u,\theta)$ by \eqref{Eq.Euler-monatomic-ansatz}. Then there are $R_*,c_0>0$ independent of $T,R_0$ such that for $R_0\geq R_*$,
	\[\inf_{0\leq t<T,\,x\in\R^3}\theta(t,x)\geq c_0R_0^{1-r},\]
	and \eqref{eq:background-hierarchy} holds for all $j,k\in\mathbb N_0$.
\end{lemma}
\begin{proof}
	On $\supp\chi_{R_0}$ we have $\theta_{\rm p}\geq cR_0^{1-r}$ by \eqref{eq:background-hierarchy-base}. By \eqref{Eq.Euler-monatomic-ansatz} and \eqref{Eq.remainder-pointwise-bound},
	\[\theta\geq R_0^{1-r}\big(c\chi_{R_0}+1-\chi_{R_0}-C(T-t)R_0^{-r}\big)\geq\tfrac12\min\{1,c\}R_0^{1-r}\]
	after increasing $R_*$ so that $CT_0R_*^{-r}\leq\tfrac12\min\{1,c\}$. For $R_0/4\leq |x|\leq 2R_0$, \eqref{Eq.profile-mixed-bound}, \eqref{Eq.remainder-pointwise-bound}, \eqref{eq:euler-correction-mixed} and $|\nabla_x^k\chi_{R_0}|\leq C_kR_0^{-k}$ give
	\[|\partial_t^j\nabla_x^ku|+|\partial_t^j\nabla_x^k\theta|
	\leq C_{j,k}R_0^{1-r-jr-k}\leq C_{j,k}\theta^{1+(jr+k)/(r-1)}.\]
	For $|x|\leq R_0/4$, use $(u,\theta)=(u_{\rm p},\theta_{\rm p})$ and \eqref{eq:background-hierarchy-profile}. For $|x|\geq2R_0$, $(u,\theta)=(0,R_0^{1-r})$ by \eqref{eq:backward-correction-annulus}. These bounds give \eqref{eq:background-hierarchy}.
\end{proof}

The support of $1-\chi_{R_0}$ in the $y$ variables starts at $|y|=R_0/(2(T-t)^{1/r})$. We use \eqref{eq:profile-decay} there and \eqref{Eq.remainder-pointwise-bound} on the correction.

\begin{lemma}\label{lem:euler-background-asymptotics}
	Let $(u,\theta)$ be as in \eqref{Eq.Euler-monatomic-ansatz}, using the correction of Lemma~\ref{lem:euler-correction-support}, and $\rho=\theta^3$. We have \eqref{eq:profile} and
	\begin{align*}
		&\|(T-t)^{\alpha_{\rm ss}}u(t,(T-t)^{1/r}\cdot)-U\|_{L^\infty_y}+\|(T-t)^{\alpha_{\rm ss}}\theta(t,(T-t)^{1/r}\cdot)-\vartheta\|_{L^\infty_y}\\
		&\hspace{35mm}\leq C(T-t)^{\alpha_{\rm ss}}R_0^{1-r}+C(T-t)^{1+\alpha_{\rm ss}}R_0^{1-2r}\longrightarrow0.
	\end{align*}
	As a consequence,
	\begin{align*}
		\|(T-t)^{3\alpha_{\rm ss}}\rho(t,(T-t)^{1/r}\cdot)-\vartheta^3\|_{L^\infty_y}&\longrightarrow0,\quad
		(T-t)^{\alpha_{\rm ss}}\|u(t)\|_{L^\infty_x}\longrightarrow\|U\|_{L^\infty_y}>0,\\
		(T-t)^{3\alpha_{\rm ss}}\|\rho(t)\|_{L^\infty_x}&\longrightarrow\|\vartheta\|_{L^\infty_y}^3>0.
	\end{align*}
\end{lemma}
\begin{proof}
	Since the correction vanishes and $\chi_{R_0}=1$ on $|x|\leq R_0/4$, we obtain \eqref{eq:profile}. We define $\chi_t(y):=\chi((T-t)^{1/r}|y|/R_0)$. By \eqref{Eq.Euler-monatomic-ansatz}, we have
	\begin{align*}
		(T-t)^{\alpha_{\rm ss}}u(t,(T-t)^{1/r}y)-U(y)&=(\chi_t(y)-1)U(y)+(T-t)^{\alpha_{\rm ss}}u_{\rm r}(t,(T-t)^{1/r}y),\\
		(T-t)^{\alpha_{\rm ss}}\theta(t,(T-t)^{1/r}y)-\vartheta(y)&=(\chi_t(y)-1)\vartheta(y)+(T-t)^{\alpha_{\rm ss}}R_0^{1-r}(1-\chi_t(y))\\
		&\quad+(T-t)^{\alpha_{\rm ss}}\theta_{\rm r}(t,(T-t)^{1/r}y).
	\end{align*}
	The profile decay and \eqref{Eq.remainder-pointwise-bound} imply
	\begin{align*}
		\|(1-\chi_t)(U,\vartheta)\|_{L^\infty_y}&\leq C\sup_{|y|\geq R_0/(2(T-t)^{1/r})}\langle y\rangle^{1-r}\leq C(T-t)^{\alpha_{\rm ss}}R_0^{1-r},\\
		(T-t)^{\alpha_{\rm ss}}\|(u_{\rm r},\theta_{\rm r})(t)\|_{L^\infty_x}&\leq C(T-t)^{1+\alpha_{\rm ss}}R_0^{1-2r}.
	\end{align*}
	By $\rho=\theta^3$ and uniform convergence of the rescaled temperature,
	\[\|(T-t)^{3\alpha_{\rm ss}}\rho(t,(T-t)^{1/r}\cdot)-\vartheta^3\|_{L^\infty_y}
	\leq C\|(T-t)^{\alpha_{\rm ss}}\theta(t,(T-t)^{1/r}\cdot)-\vartheta\|_{L^\infty_y}\longrightarrow0.\]
	Apply the reverse triangle inequality to both norm limits. If $U=0$, \eqref{Eq.Euler-monatomic-ss} gives $\nabla(\theta_{\rm s}^2)=0$, contrary to positivity and decay of $\theta_{\rm s}$.
\end{proof}

\begin{proof}[Proof of Proposition~\ref{Thm.Euler-solution}]
	Take $T_0,R_*$ from Lemmas~\ref{lem:euler-terminal-problem}-\ref{lem:euler-background-bounds}. Add the smooth radial correction with zero terminal data from Lemma~\ref{lem:euler-correction-support}. The residuals cancel by \eqref{Eq.u_r-theta_r-eq}, giving \eqref{Eq.Euler-u-theta}. Set $\rho=\theta^3$ to obtain
	\begin{align*}
		\partial_t\rho+\Div(\rho u)&=3\theta^2\left(\partial_t\theta+u\cdot\nabla\theta+\tfrac13\theta\Div u\right)=0,\\
		\rho(\partial_tu+u\cdot\nabla u)+\nabla(\rho^{5/3})&=\theta^3\left(\partial_tu+u\cdot\nabla u+5\theta\nabla\theta\right)=0.
	\end{align*}
	Lemmas~\ref{lem:euler-background-bounds} and \ref{lem:euler-background-asymptotics} give \eqref{eq:background}, \eqref{eq:profile} and \eqref{eq:background-hierarchy}. The thresholds depend only on the profile and cutoffs.
\end{proof}

\section{Velocity weights and collision estimates}\label{sec:velocity}
We keep the tail exponent $\ell$ fixed throughout the proof of Lemma~\ref{thm:fixed-tail} and Propositions~\ref{thm:new-collision}, \ref{thm:adaptive}. In particular, the constants in \eqref{eq:collision-missing} and \eqref{eq:new-bilinear} are uniform for $m\geq m_0$.

\subsection{Construction of the velocity weights}\label{subsec:velocity-weight-construction}

\begin{proof}[Proof of Lemma~\ref{thm:fixed-tail}]
	Differentiating the integral in \eqref{eq:Phi-fixed-tail}, for $a>m$ we get
	\begin{equation}\label{eq:Phi-a}
		\partial_a\Phi_{m,\ell}(a)=\frac{\ell+(m-\ell)\chi_{\rm vel}(a/m)}a.
	\end{equation}
	At $a=m$, $\partial_a\Phi(m+)=\partial_a\Phi(m-)=1$. Moreover,
	\[\partial_a^2\Phi_{m,\ell}(a)=\frac{(m-\ell)(a/m)\chi_{\rm vel}'(a/m)-\ell-(m-\ell)\chi_{\rm vel}(a/m)}{a^2}\leq0\qquad(a>m).\]
	Thus, $\Phi_{m,\ell}\in C^{1,1}_{\rm loc}([0,\infty))$ is increasing and concave, and
	\[0\leq\Phi_{m,\ell}(a)=\int_0^a\partial_b\Phi_{m,\ell}(b)\dd b\leq a,\]
	which proves \eqref{eq:weight-order}.
	
	For $a\geq2m$, by \eqref{eq:Phi-fixed-tail}, we have
	\begin{align*}
		\Phi_{m,\ell}(a)=(1+I_{{\rm vel},\infty})m+\ell\log(a/m)-\ell I_{{\rm vel},\infty},\quad w_{m,\ell}(v)=C_{m,\ell}|v|^{-\ell/2}\quad(|v|^2\geq2m).
	\end{align*}
	For $a\geq m$, we have $\Phi_{m,\ell}(a)\geq m+\ell\log(a/m)$, and hence
	\[w_{m,\ell}(v)\leq C\mathrm e^{-m/4}(|v|^2/m)^{-\ell/4}.\]
	By the substitution $|v|=\sqrt m\,R$, we obtain
	\begin{align*}
		\int_{|v|\geq\sqrt m}\langle v\rangle^{2L}w_{m,\ell}^2\dd v&\leq C_Lm^{L+3/2}\mathrm e^{-m/2}\int_1^\infty R^{2L+2-\ell}\dd R=\frac{C_Lm^{L+3/2}\mathrm e^{-m/2}}{\ell-2L-3},\\
		\int_{|v|\geq\sqrt m}\langle v\rangle^Lw_{m,\ell}\dd v&\leq C_Lm^{(L+3)/2}\mathrm e^{-m/4}\int_1^\infty R^{L+2-\ell/2}\dd R=\frac{C_Lm^{(L+3)/2}\mathrm e^{-m/4}}{\ell/2-L-3}.
	\end{align*}
	For $\ell>2L+6$, both bounds are uniform in $m\geq\ell$. Since $w_{m,\ell}=\sqrt\mu$ on $|v|^2\leq m$, this proves \eqref{eq:new-moments}. Moreover,
	\begin{align*}
		0\leq w_{m,\ell}-\sqrt\mu&\leq w_{m,\ell}{\bf1}_{\{|v|^2\geq m\}},&0\leq w_{m,\ell}-\mu/w_{m,\ell}&\leq w_{m,\ell}{\bf1}_{\{|v|^2\geq m\}}.
	\end{align*}
	\begin{align*}
		&\|\langle v\rangle^L(w_{m,\ell}-\sqrt\mu)\|_{L^2_v}^2+\|\langle v\rangle^L(w_{m,\ell}-\mu/w_{m,\ell})\|_{L^2_v}^2\leq C_{L,\ell}m^{L+3/2}\mathrm e^{-m/2}\leq C_{L,\ell}\mathrm e^{-m/4}.
	\end{align*}
	This proves \eqref{eq:new-tail-small}.
	
	The concavity and $\Phi(0)=0$ imply that, for $a+b>0$,
	\[\Phi(a)+\Phi(b)\geq\frac a{a+b}\Phi(a+b)+\frac b{a+b}\Phi(a+b)=\Phi(a+b).\]
	Therefore, by the conservation of energy in collisions, we have
	\[\Phi(|v'|^2)+\Phi(|v_*'|^2)\geq\Phi(|v|^2+|v_*|^2)\geq\Phi(|v|^2),\qquad\frac{w(v')w(v_*')}{w(v)}\leq(2\pi)^{-3/4}\leq1.\]
	Interchanging the pre- and post-collision velocities, we get the other inequality in \eqref{eq:elastic-ratio}.
	
	At fixed $a>m$, with $q:=a/m$,
	\begin{equation}\label{eq:Phi-m}
		\partial_m\Phi_{m,\ell}(a)=I_{\rm vel}(q)+(1-\ell/m)(1-\chi_{\rm vel}(q)).
	\end{equation}
	As a consequence,
	\begin{align*}
		0\leq\partial_m\Phi_{m,\ell}\leq1+\log2,\quad\partial_m\Phi(m+)=\partial_m\Phi(m-)=0,\quad
		|\partial_m\log w|=\tfrac14|\partial_m\Phi|\leq\tfrac14(1+\log2).
	\end{align*}
	Here we have used $I_{\rm vel}(1)=0$ and $\chi_{\rm vel}=1$ near $1$ for the boundary limits. This proves \eqref{eq:bounded-m-derivative}.
\end{proof}

\begin{proof}[Proof of Proposition~\ref{thm:adaptive}]
	Since $c_*\kappa_{\rm w}/4\geq1$, we have
	\[\Lambda_\eps\mathrm e^{-c_*m_\eps/4}=\mathrm e^{-c_*m_0/4}\frac{\Lambda_\eps}{(1+\Lambda_\eps)^{c_*\kappa_{\rm w}/4}}\leq\mathrm e^{-c_*m_0/4}.\]
	Moreover, $\Lambda_\eps\geq\theta_*^{\delta_{\rm gap}}/\eps$, and hence
	\begin{align*}
		\sup_{t,x}\eps\theta^{-\delta_{\rm gap}}m_\eps^a&=\sup_{t,x}\frac{(m_0+\kappa_{\rm w}\log(1+\Lambda_\eps))^a}{\Lambda_\eps}\\
		&\leq\sup_{z\geq\theta_*^{\delta_{\rm gap}}/\eps}\frac{(m_0+\kappa_{\rm w}\log(1+z))^a}{z}\longrightarrow0\qquad(\eps\downarrow0).
	\end{align*}
	This proves \eqref{eq:adaptive-absorption}.
\end{proof}

For $D:=\partial_t+v\cdot\nabla_x$, by the chain rule, we have
\begin{equation}\label{eq:Dm}
	Dm_\eps=\kappa_{\rm w}\delta_{\rm gap}\frac{\Lambda_\eps}{1+\Lambda_\eps}D\log\theta,\quad |Dm_\eps|\le C\theta^p\langle \widehat v\rangle,\quad |\nabla_xm_\eps|\le C d^{-1},\ |\nabla_x^2m_\eps|\le C d^{-2}.
\end{equation}
Indeed, $D\log\theta=-\frac13\Div u+\widehat v\cdot\nabla\theta$ and
\begin{align*}
	\partial_{x_i}m_\eps&=\kappa_{\rm w}\delta_{\rm gap}\frac{\Lambda_\eps}{1+\Lambda_\eps}\partial_{x_i}\log\theta,\\
	\partial_{x_i}\partial_{x_j}m_\eps&=\kappa_{\rm w}\delta_{\rm gap}\frac{\Lambda_\eps}{1+\Lambda_\eps}\partial_{x_i}\partial_{x_j}\log\theta+\kappa_{\rm w}\delta_{\rm gap}^2\frac{\Lambda_\eps}{(1+\Lambda_\eps)^2}\partial_{x_i}\log\theta\,\partial_{x_j}\log\theta.
\end{align*}
The bounds on $m_\eps$ in \eqref{eq:Dm} follow from \eqref{eq:background-bounds}. To differentiate $W_{m,\ell}$ from \eqref{eq:fitted-velocity-weight}, we use $\Phi_{m,\ell}$ at $\widetilde v=(v-u_f)/\sqrt{T_f}$.

\begin{lemma}\label{prop:weight-derivatives}
	For $m\geq\ell\geq1$, we have almost everywhere
	\begin{equation}\label{eq:Phi-derivatives}
		\begin{split}
			&\langle v\rangle|\nabla_v\Phi_{m,\ell}(|v|^2)| +\langle v\rangle^2|\nabla_v^2\Phi_{m,\ell}(|v|^2)|\leq C m,\\
			&|\partial_m^2\Phi_{m,\ell}(|v|^2)|\leq C/m, \qquad \langle v\rangle|\nabla_v\partial_m\Phi_{m,\ell}(|v|^2)|\leq C.
		\end{split}
	\end{equation}
	For $W_{m,\ell}$ in \eqref{eq:fitted-velocity-weight}, let $m=m(t,x)\geq\ell$ and $\widetilde v=(v-u_f)/\sqrt{T_f}$. Assume the bounds at fixed $v$ for $i,j\in\{1,2,3\}$
	\begin{equation}\label{eq:coefficient-derivative-assumptions}
		\begin{split}
			&|\partial_{x_i}\widetilde v|\leq C d^{-1}\langle \widetilde v\rangle, \quad |\partial_{x_i}\partial_{x_j}\widetilde v| \leq Cd^{-2}\langle \widetilde v\rangle,\quad|\nabla_xm|\leq C/d,\quad |\nabla_x^2m|\leq C/d^2,\\
			&\qquad|\nabla_x\log n|+|\nabla_x\log\Theta|\leq C/d,\quad |\nabla_x^2\log n|+|\nabla_x^2\log\Theta|\leq C/d^2.
		\end{split}
	\end{equation}
	Then
	\begin{equation}\label{eq:W-x}
		d|\nabla_x\log W_{m,\ell}|+d^2|\nabla_x^2\log W_{m,\ell}|\leq Cm.
	\end{equation}
	For $m=m_\eps$ and $D=\partial_t+v\cdot\nabla_x$, the radius term in $D\log W_{m,\ell}$ satisfies
	\[
	 R_{m,\ell}:=-\tfrac14(Dm)\partial_m\Phi_{m,\ell}(|\widetilde v|^2),\qquad
	 |R_{m,\ell}|\leq C\theta^p\langle\widehat v\rangle.
	\]
\end{lemma}
\begin{proof}
	Set $a=|v|^2$ and $q=a/m$. By \eqref{eq:Phi-a} and \eqref{eq:Phi-m}, for $\Phi=\Phi_{m,\ell}$ and $a>m$,
	\begin{align}\label{eq:Phi-differentiated}
		a|\partial_a\Phi|&\leq m,\qquad a^2|\partial_a^2\Phi|\leq(m-\ell)|q\chi_{\rm vel}'(q)|+\ell+(m-\ell)\chi_{\rm vel}(q)\leq Cm,\notag\\
		\partial_m^2\Phi&=-\frac{\chi_{\rm vel}(q)}m+\frac{\ell}{m^2}\big(1-\chi_{\rm vel}(q)\big)+(1-\ell/m)\frac q m\chi_{\rm vel}'(q),\\
		\nabla_v\partial_m\Phi&=\frac{2v}m\left(\frac{\chi_{\rm vel}(q)}q-(1-\ell/m)\chi_{\rm vel}'(q)\right).\notag
	\end{align}
	For $a<m$ we have $\nabla_v\Phi=2v$, $\nabla_v^2\Phi=2\Id$ and $\partial_m\Phi=0$. For $a>m$, \eqref{eq:Phi-differentiated} gives
	\begin{align*}
		|\nabla_v\Phi|&=2|v\partial_a\Phi|\leq 2m/\sqrt a,\\
		|\nabla_v^2\Phi|&=|2\partial_a\Phi\,\Id+4\partial_a^2\Phi\,v\otimes v|\leq Cm/a,\\
		m|\partial_m^2\Phi|+\langle v\rangle|\nabla_v\partial_m\Phi|&\leq C.
	\end{align*}
	In the last line $\nabla_v\partial_m\Phi=0$ for $a>2m$. The matching traces $(\partial_a\Phi,\partial_m\Phi)|_{a=m\pm}=(1,0)$ give \eqref{eq:Phi-derivatives} for the weak derivatives.

	At fixed $(t,v)$ put $F(x)=\Phi_{m(t,x),\ell}(|\widetilde v|^2)$. By \eqref{eq:fitted-velocity-weight},
	\begin{equation}\label{eq:fitted-log-weight}
		\log W_{m,\ell}+\tfrac34\log(2\pi)=\tfrac12\log n-\tfrac34\log\Theta-\tfrac14F.
	\end{equation}
	Differentiating $F$ gives
	\begin{equation}\label{eq:fitted-Phi-chain}
		\begin{aligned}
			\partial_{x_i}F&=\nabla_{\widetilde v}\Phi\cdot\partial_{x_i}\widetilde v+\partial_m\Phi\,\partial_{x_i}m,\\
			\partial_{x_i}\partial_{x_j}F&=\nabla_{\widetilde v}\Phi\cdot\partial_{x_i}\partial_{x_j}\widetilde v
				+\nabla_{\widetilde v}^2\Phi[\partial_{x_i}\widetilde v,\partial_{x_j}\widetilde v]\\
			&\quad+\nabla_{\widetilde v}\partial_m\Phi\cdot(\partial_{x_i}\widetilde v\,\partial_{x_j}m+\partial_{x_j}\widetilde v\,\partial_{x_i}m)\\
			&\quad+\partial_m^2\Phi\,\partial_{x_i}m\,\partial_{x_j}m+\partial_m\Phi\,\partial_{x_i}\partial_{x_j}m.
		\end{aligned}
	\end{equation}
	In the order of the terms in \eqref{eq:fitted-Phi-chain}, \eqref{eq:coefficient-derivative-assumptions}, \eqref{eq:Phi-derivatives} and \eqref{eq:bounded-m-derivative} give
	\begin{align*}
		|\partial_{x_i}F|&\leq C(m/d+1/d),\\
		|\partial_{x_i}\partial_{x_j}F|&\leq C(m/d^2+m/d^2+1/d^2+1/(md^2)+1/d^2).
	\end{align*}
	Substituting these bounds in \eqref{eq:fitted-log-weight} proves \eqref{eq:W-x}, since $m\geq1$. For $R_{m,\ell}$ use \eqref{eq:Dm} and \eqref{eq:bounded-m-derivative} to get
	\[|R_{m,\ell}|\leq\tfrac14|Dm|\,|\partial_m\Phi|\leq C\theta^p\langle\widehat v\rangle.\qedhere\]
\end{proof}

Equation \eqref{eq:Dm} supplies the bounds on $m$ in \eqref{eq:coefficient-derivative-assumptions}. Lemma~\ref{lem:single-spatial-allocation} supplies the bounds on $n,V,\Theta$ on the bootstrap interval.

For $m\geq\max\{1,\ell\}$ and $1<\gamma\leq2$, \eqref{eq:Phi-a} gives
\begin{equation}\label{eq:cubic-weight-bound}
 |v|^3\partial_a\Phi_{m,\ell}(|v|^2)
 \leq
 \begin{cases}
  1,&|v|\leq1,\\
  m^{(3-\gamma)/2}|v|^\gamma,&1<|v|\leq\sqrt m,\\
  m|v|\leq m^{(3-\gamma)/2}|v|^\gamma,&|v|>\sqrt m.
 \end{cases}
\end{equation}
In the last case of \eqref{eq:cubic-weight-bound}, $a\partial_a\Phi\leq m$ and $m^{(\gamma-1)/2}\leq |v|^{\gamma-1}$. It follows from \eqref{eq:cubic-weight-bound} and $\widehat\nu_M\asymp\langle\widetilde v\rangle^\gamma$ that
\begin{equation}\label{eq:fitted-cubic-ratio}
 \frac{\theta^p\langle\widetilde v\rangle^3\partial_a\Phi_{m_\eps,\ell}(|\widetilde v|^2)}
      {\eps^{-1}\theta^{3+\gamma}\widehat\nu_M}
 \leq C\eps\theta^{-\delta_{\rm gap}}\bigl(1+m_\eps^{(3-\gamma)/2}\bigr).
\end{equation}
The supremum over $(t,x,\widetilde v)$ in \eqref{eq:fitted-cubic-ratio} tends to zero by \eqref{eq:adaptive-absorption}.

\subsection{Uniform collision estimates for the fixed-order family}\label{subsec:uniform-collision}

\begin{proof}[Proof of Proposition~\ref{thm:new-collision}]
	Set $\Phi=\Phi_{m,\ell}$ and $q(a)=\mathrm e^{-\Phi(a)/4}$. By \eqref{eq:Phi-fixed-tail},
	\begin{align*}
		\partial_\ell\Phi_{m,\ell}(a)&=0 &&(a\leq m),\\
		\partial_\ell\Phi_{m,\ell}(a)&=\log(a/m)-I_{\rm vel}(a/m)=\int_1^{a/m}\frac{1-\chi_{\rm vel}(s)}s\dd s\geq0 &&(a>m).
	\end{align*}
	By \eqref{eq:bounded-m-derivative}, $\partial_m\Phi\geq0$ also. Thus for $m\geq\ell\geq\ell_1>2L+6$,
	\begin{equation}\label{eq:q-moment-uniform}
		\begin{aligned}
			&\Phi_{m,\ell}(a)\geq\Phi_{m,\ell_1}(a)\geq\Phi_{\ell_1,\ell_1}(a),\\
			&\int\langle v\rangle^Lq(|v|^2)\dd v\leq\int\langle v\rangle^L\mathrm e^{-\Phi_{\ell_1,\ell_1}(|v|^2)/4}\dd v<\infty.
		\end{aligned}
	\end{equation}
	Fix $\ell_1$ for the largest $L$ used below. Integrating $\ell/(\ell+a)\leq\Phi'(a)\leq1$ gives
	\begin{equation}\label{eq:new-ratio-bounds}
		\frac{q(A)}{q(B)}\leq\mathrm e^{(B-A)/4},\qquad\frac{q(B)}{q(A)}\leq\left(\frac{\ell+A}{\ell+B}\right)^{\ell/4},\qquad0\leq A\leq B.
	\end{equation}
	For $L_\infty=\nu-K$ and $\zeta_w:=w/w_\infty$, $L_w=\nu-\zeta_w^{-1}K\zeta_w$. We shall prove
	\begin{equation}\label{eq:new-conjugation}
		|\langle(L_w-L_\infty)g,g\rangle|\leq C_\gamma\bigl(\ell^{-1}+\mathrm e^{-cm}\bigr)\|\nu^{1/2}g\|_{L^2_v}^2,\qquad m\geq\ell\geq\ell_1(\gamma).
	\end{equation}
	For $v\ne v'$, we define
	\begin{align*}
		\rho&:=|v-v'|,\quad e:=(v-v')/\rho,\quad a:=v\cdot e,\quad a':=v'\cdot e,\\
		y&:=v-ae=v'-a'e,\qquad E:=|v|^2,\qquad E':=|v'|^2.
	\end{align*}
	With $a-a'=\rho$, $E=|y|^2+a^2$ and $E'=|y|^2+(a')^2$, the starred loss is bounded by $C\rho^\gamma\mathrm e^{-(E+E')/4}$. The sum $K^+$ of the gains satisfies
	\begin{equation}\label{eq:new-grad-kernel}
		|K^+(v,v')|\leq C\rho^{-1}(\rho+\langle y\rangle)^{\gamma-1}\mathrm e^{-(a^2+(a')^2)/4}.
	\end{equation}
	With $\rho=(v-v_*)\cdot e>0$, the change in \eqref{eq:postcollisional-velocities} is
	\begin{equation}\label{eq:gain-collision-coordinates}
		\begin{aligned}
			\omega&=(v-v_*-2\rho e)/|v-v_*|,\qquad \dd\omega=4\frac{\rho}{|v-v_*|}\dd e,\\
			v'&=v-\rho e,\qquad v_*'=v_*+\rho e,\qquad
			v_*=y_*+a'e,\quad y_*\in e^\perp,\quad a'=a-\rho .
		\end{aligned}
	\end{equation}
	By $\dd v_*=\dd y_*\,\dd a'$ and $\dd v'=\rho^2\dd\rho\dd e$,
	\[B_\gamma(v-v_*)\dd\omega\dd v_*=\frac{|\rho e+y-y_*|^{\gamma-1}}{\pi\rho}\dd y_*\,\dd v'.\]
	By $\mu\mu_*=\mu'\mu_*'$, the gain with $g(v')$ has factor
	\[\frac{\mu(v_*')\sqrt{\mu(v')}}{\sqrt{\mu(v)}}=\sqrt{\mu(v_*)\mu(v_*')}=(2\pi)^{-3/2}\mathrm e^{-(a^2+(a')^2)/4-|y_*|^2/2}.\]
	The change $\omega\mapsto-\omega$ exchanges the gains. Thus
	\begin{align*}
		|K^+(v,v')|
		&\leq \frac{C\mathrm e^{-(a^2+(a')^2)/4}}{\rho}
		\int_{e^\perp}|\rho e+y-y_*|^{\gamma-1}\mathrm e^{-|y_*|^2/2}\dd y_*\\
		&\leq \frac{C_\gamma\mathrm e^{-(a^2+(a')^2)/4}}{\rho}(\rho+\langle y\rangle)^{\gamma-1}.
	\end{align*}
	The last bound uses $0<\gamma-1\leq1$.

	To estimate the kernel difference in \eqref{eq:new-conjugation}, we have
	\[\mathrm e^{-(a^2+(a')^2)/4}\frac{\zeta_w(v')}{\zeta_w(v)}=\mathrm e^{-a^2/2}\frac{q(E')}{q(E)},\qquad \frac{\zeta_w(v')}{\zeta_w(v)}=1\quad(\max\{E,E'\}\leq m).\]
	The transposed kernel exchanges $(a,E)$ and $(a',E')$. Using $\dd v'=\rho^2\dd\rho\dd e$ in \eqref{eq:new-grad-kernel} and enlarging $a'<a$ to $\R$, we bound the sum of the absolute row and column by $C\mathfrak R(v)$, where
	\begin{align}\label{eq:new-row}
	 \mathfrak R(v):=\int_{\mathbb S^2}\int_\R
	  &\bigl(|a-a'|^\gamma+\langle y\rangle^{\gamma-1}|a-a'|\bigr)
	      \mathbf1_{\{\max(E,E')\geq m\}}\notag\\
	  &\times\left(\mathrm e^{-a^2/2}\frac{q(E')}{q(E)}
	     +\mathrm e^{-(a')^2/2}\frac{q(E)}{q(E')}
	     +2\mathrm e^{-(a^2+(a')^2)/4}\right)\dd a'\dd e .
	\end{align}
	Put $A=\ell+E$ and $k=\ell/4$. By \eqref{eq:new-ratio-bounds},
	\[
	 \mathrm e^{-a^2/2}\frac{q(E')}{q(E)}
	 \leq\begin{cases}
	  \mathrm e^{-(a^2+(a')^2)/4},&|a'|\leq|a|,\\
	  \mathrm e^{-a^2/2}\bigl(1+(|a'|-|a|)^2/A\bigr)^{-k},&|a'|>|a|.
	 \end{cases}
	\]
	For $\tau\in\{0,1,\gamma\}$ and $k\geq\tau+2$, by the substitution $u:=s^2/A$, we obtain
	\begin{equation}\label{eq:rw34-beta}
		\begin{aligned}
			\int_0^\infty s^\tau(1+s^2/A)^{-k}\dd s&=\tfrac12A^{(\tau+1)/2}\int_0^\infty u^{(\tau-1)/2}(1+u)^{-k}\dd u\\
			&\leq C_\tau A^{(\tau+1)/2}\left(\int_0^1u^{(\tau-1)/2}\mathrm e^{-ku/2}\dd u+\frac{2^{1-k/2}}{k-\tau-1}\right)\\
			&\leq C_\tau A^{(\tau+1)/2}k^{-(\tau+1)/2}.
		\end{aligned}
	\end{equation}
	Thus,
	\begin{align*}
		&\int_\R\bigl(|a-a'|^\gamma+\langle y\rangle^{\gamma-1}|a-a'|\bigr)\mathrm e^{-a^2/2}\frac{q(E')}{q(E)}\dd a'\\
		&\quad\leq C\mathrm e^{-a^2/4}\left[1+|a|^{\gamma+1}+\langle y\rangle^{\gamma-1}(1+a^2)+A^{(\gamma+1)/2}k^{-(\gamma+1)/2}+\langle y\rangle^{\gamma-1}A/k\right].
	\end{align*}
	For $|a'|\leq|a|$,
	\[\log\frac{\ell+E}{\ell+E'}\geq\frac{a^2-(a')^2}{A},\qquad\frac kA\leq\frac14,\qquad\mathrm e^{-(a')^2/2}\frac{q(E)}{q(E')}\leq\mathrm e^{-ka^2/A-(a')^2/4}.\]
	For $|a'|\geq|a|$, by \eqref{eq:new-ratio-bounds}, we have $\mathrm e^{-(a')^2/2}q(E)/q(E')\leq\mathrm e^{-(a^2+(a')^2)/4}$. Hence,
	\begin{align*}
		\int_\R\bigl(|a-a'|^\gamma+\langle y\rangle^{\gamma-1}|a-a'|\bigr)\mathrm e^{-(a')^2/2}\frac{q(E)}{q(E')}\dd a'\leq C\mathrm e^{-ka^2/A}\left[1+|a|^\gamma+\langle y\rangle^{\gamma-1}(1+|a|)\right].
	\end{align*}
	For $a=Rt$ and $R=|v|>0$, use the variables $Rt$ and $\sqrt{k}t$ to obtain
	\begin{equation}\label{eq:gain-angular-integrals}
		\begin{aligned}
			\int_{\mathbb S^2}\mathrm e^{-c(v\cdot e)^2}\dd e&=\frac{2\pi}{R}\int_{-R}^R\mathrm e^{-cu^2}\dd u\leq C/R,\\
			\int_{-1}^1\mathrm e^{-kt^2/3}|t|^\tau\dd t&\leq k^{-(\tau+1)/2}\int_{\R}\mathrm e^{-u^2/3}|u|^\tau\dd u\leq C_\tau k^{-(\tau+1)/2}\qquad(\tau=0,1,\gamma).
		\end{aligned}
	\end{equation}
	For $E\geq m/2$, use $A\leq3R^2$, $R^2\geq\ell/2$ and \eqref{eq:gain-angular-integrals} in \eqref{eq:new-row}. Since $\gamma>1$ and $k=\ell/4$,
	\begin{equation}\label{eq:conjugation-gain-large}
		\frac{\mathfrak R(v)}{\langle v\rangle^\gamma}\leq C\bigl(k^{-(\gamma+1)/2}+k^{-1}+R^{-2}+k^{-1/2}R^{-1}\bigr)\leq C\ell^{-1}.
	\end{equation}

	For $E<m/2$, \eqref{eq:new-row} is restricted to $E'\geq m$. Then $\Phi(E)=E$ and $(a')^2=E'-|y|^2\geq m/2$. Integrating \eqref{eq:Phi-a} yields
	\begin{equation}\label{eq:gain-small-input-tail}
		q(E')\leq\mathrm e^{-m/4}(E'/m)^{-k},\qquad \mathrm e^{-a^2/2}\frac{q(E')}{q(E)}\leq\mathrm e^{-m/8}(E'/m)^{-k}.
	\end{equation}
	Use \eqref{eq:gain-small-input-tail} in the first term of \eqref{eq:new-row} and put $a'=\sqrt m\,u$. By $|a|,|y|\leq\sqrt m$ and $E'/m\geq\max\{1,u^2\}$, for fixed $k_{\min}>\gamma+3$ and $k\geq k_{\min}$ we get
	\begin{align*}
	 &\int_{\mathbb S^2}\int_{E'\geq m}
	   \bigl(|a-a'|^\gamma+\langle y\rangle^{\gamma-1}|a-a'|\bigr)
	       \mathrm e^{-a^2/2}\frac{q(E')}{q(E)}\dd a'\dd e\\
	 &\qquad\leq C(1+m^{(\gamma+1)/2})\mathrm e^{-m/8}
	   \int_\R(1+|u|)^\gamma\max\{1,u^2\}^{-k_{\min}}\dd u
	   \leq C\mathrm e^{-cm}.
	\end{align*}
	For the sum of the other two factors, \eqref{eq:new-ratio-bounds} gives the integrable majorant
	\[
	 \mathrm e^{-(a')^2/2}\frac{q(E)}{q(E')}+2\mathrm e^{-(a^2+(a')^2)/4}
	   \leq3\mathrm e^{-m/16}\mathrm e^{-(a^2+(a')^2)/8}
	   \qquad(E<m/2,\ E'\geq m).
	\]
	Since $|y|^2<m/2$, the two remaining terms in \eqref{eq:new-row} integrate to $O(\mathrm e^{-cm})$. Use \eqref{eq:gain-small-input-tail} and \eqref{eq:conjugation-gain-large} to obtain
	\begin{equation}\label{eq:conjugation-gain-rows}
		\mathfrak R(v)\leq
		\begin{cases}
			C\mathrm e^{-cm},&E<m/2,\\
			C\ell^{-1}\langle v\rangle^\gamma,&E\geq m/2.
		\end{cases}
	\end{equation}

	For the starred loss kernel put $b(E)=\mathrm e^{-E/2}/q(E)$. By \eqref{eq:weight-order}, $b(E)\leq\mathrm e^{-E/4}$. From \eqref{eq:gain-small-input-tail} and $v'=\sqrt m\,u$, for $\ell\geq\ell_1>2L+6$ we have
	\begin{align*}
		\sup_{E\geq m}\{b(E)+q(E)\}&\leq2\mathrm e^{-m/4},\\
		\int_{E'\geq m}\langle v'\rangle^L\bigl(b(E')+q(E')\bigr)\dd v'
		&\leq C_L m^{(L+3)/2}\mathrm e^{-m/4}\leq C_L\mathrm e^{-cm}.
	\end{align*}
	Take $L=\gamma$. Use $\mathbf1_{\{\max(E,E')\geq m\}}\leq\mathbf1_{\{E\geq m\}}+\mathbf1_{\{E'\geq m\}}$ and $|v-v'|^\gamma\leq C\langle v\rangle^\gamma\langle v'\rangle^\gamma$. By the preceding tail bounds and \eqref{eq:q-moment-uniform},
	\begin{align}\label{eq:conjugation-loss-row}
		&\int |v-v'|^\gamma {\bf1}_{\{\max(E,E')\geq m\}}
		\bigl(b(E)q(E')+b(E')q(E)+2\mathrm e^{-(E+E')/4}\bigr)\dd v'\notag\\
		&\qquad\leq C\mathrm e^{-cm}\langle v\rangle^\gamma .
	\end{align}
	By $2|g(v)g(v')|\leq|g(v)|^2+|g(v')|^2$, the row-plus-column bounds \eqref{eq:conjugation-gain-rows} and \eqref{eq:conjugation-loss-row} give \eqref{eq:new-conjugation}.

	Let $\mathbb P$ be the orthogonal projection in $L^2_v$ onto $\sqrt\mu\,\operatorname{span}\{1,v_i,|v|^2\}$. We prove
	\begin{equation}\label{eq:gaussian-gap-new}
		\langle L_\infty g,g\rangle\geq\lambda_0\|\nu^{1/2}(I-\mathbb P)g\|_{L^2_v}^2\qquad\text{for some }\lambda_0>0.
	\end{equation}
	With $\psi=g/\sqrt\mu$, use $\mu\mu_*=\mu'\mu_*'$ to write
	\begin{equation}\label{eq:gaussian-collision-form}
		\langle L_\infty g,g\rangle=\frac1{16\pi}\iiint |v-v_*|^\gamma\mu\mu_*|\psi'+\psi_*'-\psi-\psi_*|^2\dd\omega\dd v_*\dd v.
	\end{equation}
	When \eqref{eq:gaussian-collision-form} vanishes, mollify $\psi'+\psi_*'=\psi+\psi_*$ by common translations. By \eqref{eq:postcollisional-velocities},
	\[\psi_\delta(c+se)+\psi_\delta(c-se)=\psi_\delta(c+s\omega)+\psi_\delta(c-s\omega),\qquad s\geq0,\quad |e|=|\omega|=1.\]
	Differentiate twice at $s=0$. Then $D^2\psi_\delta=\lambda_\delta\Id$ and $\partial_j\lambda_\delta=\partial_i\partial_{ij}\psi_\delta=0$ for $i\ne j$. The limit in distributions belongs to $\operatorname{span}\{1,v_i,|v|^2\}$. By \eqref{eq:gaussian-collision-form},
	\begin{equation}\label{eq:gaussian-collision-nullspace}
		\ker L_\infty=\operatorname{ran}\mathbb P .
	\end{equation}

	Let $\mathcal K=\nu^{-1/2}K\nu^{-1/2}$ with gain part $\mathcal K^+$. For $R=|v|\gg1$ and $|a|+\rho\leq R/2$, use $|v'|\asymp R$, $\rho+\langle y\rangle\lesssim R$ and $a^2+(a-\rho)^2\geq(a^2+\rho^2)/4$ in \eqref{eq:new-grad-kernel} to get
	\begin{align*}
		\int_{|a|+\rho\leq R/2}|\mathcal K^+(v,v')|\dd v'
		&\leq CR^{-1}\int_{\mathbb S^2}\mathrm e^{-ca^2}\dd e\int_0^\infty \rho\mathrm e^{-c\rho^2}\dd\rho
		\leq CR^{-2}
	\end{align*}
	by \eqref{eq:gain-angular-integrals}. The integral on $|a|+\rho>R/2$ and the starred loss are $O(\mathrm e^{-cR^2})$. Hence
	\begin{equation}\label{eq:gaussian-normalized-row}
		\int|\mathcal K(v,v')|\dd v'\leq C\langle v\rangle^{-2}.
	\end{equation}
	Let $\mathcal K_{R,N}$ be the restriction to $|v|,|v'|\leq R$, $\rho\geq N^{-1}$. By symmetry, Schur's test with \eqref{eq:gaussian-normalized-row} and the local $\rho^{-1}$ singularity in \eqref{eq:new-grad-kernel} gives
	\[\|\mathcal K-\mathcal K_{R,N}\|_{L^2\to L^2}\leq C/R+C_R\int_0^{N^{-1}}\rho\dd\rho,\qquad \mathcal K_{R,N}\in L^2_{v,v'}.\]
	If \eqref{eq:gaussian-gap-new} fails, take $\mathbb Pg_j=0$, $\|\nu^{1/2}g_j\|_2=1$ and $\langle L_\infty g_j,g_j\rangle\to0$. For a subsequence $h_j:=\nu^{1/2}g_j\rightharpoonup h$. By compactness of $\mathcal K$ and \eqref{eq:gaussian-collision-form},
	\[1=\lim_j\langle\mathcal Kh_j,h_j\rangle=\langle\mathcal Kh,h\rangle\leq\|h\|_2^2\leq 1,\]
	Thus $\|h\|_2=1$ and $\langle L_\infty\nu^{-1/2}h,\nu^{-1/2}h\rangle=0$. The weak limit satisfies $\mathbb P(\nu^{-1/2}h)=0$ by testing with $\nu^{-1/2}\sqrt\mu\,(1,v,|v|^2)$, contradicting \eqref{eq:gaussian-collision-nullspace}. See also \cite{Mouhot}.

	For $g$ in \eqref{eq:collision-missing} write
	\[g=g^\perp+\sqrt\mu\,\boldsymbol\psi_{\rm inv}\cdot b,\qquad g^\perp=(I-\mathbb P)g,\qquad \boldsymbol\psi_{\rm inv}=(1,v_1,v_2,v_3,|v|^2-3).\]
	By orthogonality of $g^\perp$ to $\sqrt\mu\,\boldsymbol\psi_{\rm inv}$,
	\[\left(\int\boldsymbol\psi_{\rm inv}\otimes\boldsymbol\psi_{\rm inv}\,w\sqrt\mu\dd v\right)b=-\int\boldsymbol\psi_{\rm inv}(w-\sqrt\mu)g^\perp\dd v.\]
	The matrix is $\operatorname{diag}(1,I_3,6)+O(\mathrm e^{-cm})$ by \eqref{eq:new-tail-small}. Invert the matrix for $m$ large and apply the Cauchy-Schwarz inequality to get
	\begin{equation}\label{eq:gaussian-projection-error}
		|b|+\|\nu^{1/2}\mathbb P g\|_{L^2_v}\leq C\mathrm e^{-cm}\|\nu^{1/2}g^\perp\|_{L^2_v}.
	\end{equation}
	Combining \eqref{eq:new-conjugation}, \eqref{eq:gaussian-gap-new} and \eqref{eq:gaussian-projection-error} gives
	\[\frac{\operatorname{Re}\langle L_wg,g\rangle}{\|\nu^{1/2}g\|_{L^2_v}^2}\geq\frac{\lambda_0}{(1+C\mathrm e^{-cm})^2}-C_\gamma(\ell^{-1}+\mathrm e^{-cm})\qquad(g\ne0).\]
	Choose $\ell$ and then $m_0\geq\ell$ to make the right-hand side at least $c_0>0$ for $m\geq m_0$, proving \eqref{eq:collision-missing}.
	
	For the quadratic operator $Q_w$, take $L>\gamma+5$ and $\ell\geq2L$. By \eqref{eq:Phi-a},
	\[\tfrac14\Phi'(a)\geq\tfrac14\min\{1,\ell/a\}\geq\frac{L}{2(1+a)}\qquad(a\geq2L).\]
	Integrating yields $\frac14(\Phi(B)-\Phi(A))\geq \frac L2\log((1+B)/(1+A))-C_L$ for $0\leq A\leq B$. To estimate the gain, we use
	\begin{equation}\label{eq:two-channel-new}
		R_w:=\frac{w(v')w(v_*')}{w(v)}\leq C_L\left(\langle v_*\rangle^{-L}+\frac{\langle v\rangle^L}{\langle v'\rangle^L\langle v_*'\rangle^L}\right).
	\end{equation}
	We define $E_*:=|v_*|^2$ and $E_*':=|v_*'|^2$, and order the post-collisional energies so that $E'\leq E_*'$. Since $E+E_*=E'+E_*'$, we have
	\begin{align*}
		E_*'\geq E&\quad\Longrightarrow\quad\frac{q(E')q(E_*')}{q(E)}\leq C_L(1+E')^{-L/2}\left(\frac{1+E}{1+E_*'}\right)^{L/2},\\
		E_*'<E&\quad\Longrightarrow\quad E_*\leq E'\leq E_*'\leq E,\\
		&\hspace{20mm}\Phi(E')+\Phi(E_*')\geq\Phi(E)+\Phi(E_*),\quad\frac{q(E')q(E_*')}{q(E)}\leq C_L(1+E_*)^{-L/2}.
	\end{align*}
	The concavity applies when $E_*' < E$, which proves \eqref{eq:two-channel-new}.
	
	For the collision rule \eqref{eq:postcollisional-velocities}, we define the Maxwell-molecule gain operator
	\[Q_0^+(f,g)(v):=\frac1{4\pi}\int_{\R^3}\int_{\mathbb S^2}f(v_*')g(v')\dd\omega\dd v_*.\]
	Alonso, Carneiro and Gamba \cite[Corollary~6]{ACG2010} proved
	\begin{equation}\label{eq:maxwell-gain-L1-L2}
		\|Q_0^+(f,g)\|_{L^2_v}\leq C\|f\|_{L^1_v}\|g\|_{L^2_v},\quad\|Q_0^+(f,g)\|_{L^2_v}\leq C\|f\|_{L^2_v}\|g\|_{L^1_v}.
	\end{equation}
	Let $Q^\pm$ denote the gain and loss terms in \eqref{Eq.collision-operator}. We define $Q_w^\pm(g_1,g_2):=w^{-1}Q^\pm(wg_1,wg_2)$ and
	\begin{align*}
		I_1(v)&:=\iint B_\gamma(v-v_*)\langle v_*\rangle^{-L}|g_1(v_*')|\,|g_2(v')|\dd\omega\dd v_*,\\
		I_2(v)&:=\langle v\rangle^LQ^+(\langle\cdot\rangle^{-L}|g_1|,\langle\cdot\rangle^{-L}|g_2|)(v).
	\end{align*}
	By \eqref{eq:two-channel-new}, we have $|Q_w^+(g_1,g_2)|\leq C_L(I_1+I_2)$. Since $L>\gamma+3$, we have
	\[\int|v-v_*|^\gamma\langle v_*\rangle^{-L}\dd v_*\leq C_\gamma\langle v\rangle^\gamma\int\langle v_*\rangle^{-L}\dd v_*+C_\gamma\int\langle v_*\rangle^{\gamma-L}\dd v_*\leq C_L\nu(v).\]
	Therefore, by the Cauchy--Schwarz inequality and the pre--post collision change of variables, we have
	\begin{align*}
		\|\nu^{-1/2}I_1\|_{L^2_v}^2
		&\leq C_L\iiint B_\gamma\langle v_*\rangle^{-L}|g_1(v_*')|^2|g_2(v')|^2\dd\omega\dd v_*\dd v\\
		&=C_L\iiint B_\gamma\langle v_*'\rangle^{-L}|g_1(v_*)|^2|g_2(v)|^2\dd\omega\dd v_*\dd v\\
		&\leq C_L\iint|v-v_*|^\gamma|g_1(v_*)|^2|g_2(v)|^2\dd v_*\dd v\\
		&\leq C_L\bigl(\|\nu^{1/2}g_1\|_{L^2_v}^2\|g_2\|_{L^2_v}^2+\|g_1\|_{L^2_v}^2\|\nu^{1/2}g_2\|_{L^2_v}^2\bigr).
	\end{align*}
	Here $B_\gamma=B_\gamma(v-v_*)$. To separate the two input norms we used $|v-v_*|^\gamma\leq C_\gamma(\langle v\rangle^\gamma+\langle v_*\rangle^\gamma)$.
	
	For $I_2$, by the conservation of kinetic energy and relative speed, we have
	\[\langle v\rangle\leq\langle v'\rangle+\langle v_*'\rangle,\qquad |v-v_*|=|v'-v_*'|\leq\langle v'\rangle+\langle v_*'\rangle.\]
	Since $L\geq\gamma/2$, we have
	\begin{align*}
		\frac{|v-v_*|^\gamma\langle v\rangle^{L-\gamma/2}}{\langle v'\rangle^L\langle v_*'\rangle^L}
		&\leq\frac{(\langle v'\rangle+\langle v_*'\rangle)^{L+\gamma/2}}{\langle v'\rangle^L\langle v_*'\rangle^L}\leq C_{L,\gamma}\bigl(\langle v_*'\rangle^{-L}\langle v'\rangle^{\gamma/2}+\langle v_*'\rangle^{\gamma/2}\langle v'\rangle^{-L}\bigr).
	\end{align*}
	As a consequence, \eqref{eq:maxwell-gain-L1-L2}, $\nu\asymp\langle v\rangle^\gamma$ and the Cauchy--Schwarz inequality imply
	\begin{align*}
		\|\nu^{-1/2}I_2\|_{L^2_v}
		&\leq C_{L,\gamma}\bigl(\|\langle v\rangle^{-L}g_1\|_{L^1_v}\|\langle v\rangle^{\gamma/2}g_2\|_{L^2_v}+\|\langle v\rangle^{\gamma/2}g_1\|_{L^2_v}\|\langle v\rangle^{-L}g_2\|_{L^1_v}\bigr)\\
		&\leq C_{L,\gamma}\bigl(\|g_1\|_{L^2_v}\|\nu^{1/2}g_2\|_{L^2_v}+\|\nu^{1/2}g_1\|_{L^2_v}\|g_2\|_{L^2_v}\bigr),
	\end{align*}
	since $\langle v\rangle^{-L}\in L^2_v$.
	
	For the loss, set $\nu[h](v)=\int|v-v_*|^\gamma h(v_*)\dd v_*$. Then $Q_w^-(g_1,g_2)=g_2\nu[wg_1]$ and the Cauchy-Schwarz inequality gives
	\[\frac{|\nu[wg_1](v)|}{\nu(v)}\leq C\|\langle v\rangle^\gamma w\|_{L^2_v}\|g_1\|_{L^2_v}\leq C_\ell\|g_1\|_{L^2_v}\]
	by \eqref{eq:new-moments} and \eqref{eq:normalized-frequency-derivatives}. Hence
	\begin{equation}\label{eq:weighted-collision-loss}
		\|\nu^{-1/2}Q_w^-(g_1,g_2)\|_{L^2_v}\leq C_\ell\|g_1\|_{L^2_v}\|\nu^{1/2}g_2\|_{L^2_v}.
	\end{equation}
	Adding the $I_1,I_2$ estimates gives \eqref{eq:new-bilinear}. Since $0\leq\mu/w\leq\sqrt\mu$ and $\nu\geq c>0$, applying \eqref{eq:new-bilinear} to $L_wg=-Q_w(\mu/w,g)-Q_w(g,\mu/w)$ proves \eqref{eq:new-linear-upper}.
\end{proof}

The transfer to $M,W$ uses $\widehat v=V+\sqrt\Theta\,\widetilde v$ in \eqref{eq:fitted-velocity-weight} and \eqref{eq:normalized-collision-frequency}.
\[\dd\widehat v=\Theta^{3/2}\dd\widetilde v,\qquad \widehat\nu_M=n\Theta^{\gamma/2}\nu(\widetilde v),\qquad W=n^{1/2}\Theta^{-3/4}w_{m,\ell}(\widetilde v).\]
The factors are bounded above and below on $\mathcal U$ by $3/4\leq n,\Theta\leq5/4$.

\subsection{Polynomially weighted velocity estimates}\label{subsec:polynomial-velocity}

\begin{lemma}\label{lem:polynomial-collision}
	Fix $s_*\geq0$ and choose $\ell,m_0$ sufficiently large. For all $s\in[0,s_*]$ and $m\geq m_0$,
	\begin{equation}\label{eq:polynomial-gap}
		\operatorname{Re}\langle L_wg,\langle v\rangle^{2s}g\rangle \geq\tfrac12\|\nu^{1/2}\langle v\rangle^sg\|_{L^2_v}^2 -C_s\|\nu^{1/2}g\|_{L^2_v}^2.
	\end{equation}
	For the same $s,m$, we have
	\begin{align}\label{eq:polynomial-bilinear}
		&\|\nu^{-1/2}\langle v\rangle^sQ_w(g_1,g_2)\|_{L^2_v}\notag\\
		&\quad\leq C_s\sum_{(i,j)=(1,2),(2,1)} \big(\|\nu^{1/2}\langle v\rangle^sg_i\|_{L^2_v}\|g_j\|_{L^2_v} +\|\langle v\rangle^sg_i\|_{L^2_v}\|\nu^{1/2}g_j\|_{L^2_v}\big).
	\end{align}
	The constants in \eqref{eq:polynomial-gap} and \eqref{eq:polynomial-bilinear} are independent of $m\geq m_0$. For $M,W$, they are also uniform on $\mathcal U$.
\end{lemma}
\begin{proof}
	Let $q_s(E)=\mathrm e^{-\Phi_{m,\ell}(E)/4}(1+E)^{-s/2}$. Then
	\[
	 \langle v\rangle^sL_wg=L_{w/\langle v\rangle^s}(\langle v\rangle^sg),
	 \qquad \langle v\rangle^s\nu g=\nu\langle v\rangle^sg.
	\]
	Let $k_s^\pm\geq0$ be the magnitudes of the gain and starred loss kernels for $w/\langle v\rangle^s$. By \eqref{eq:new-ratio-bounds},
	\begin{equation}\label{eq:polynomial-gain-ratio}
		\mathrm e^{-a^2/2}\frac{q_s(E')}{q_s(E)}\leq
		\begin{cases}
			\mathrm e^{-a^2/2}q_0(E')/q_0(E),& E'\geq E,\\
			C_{s_*}\mathrm e^{-(a^2+(a')^2)/8},& E'<E .
		\end{cases}
	\end{equation}
	In the second line use $q_s(E')/q_0(E')\cdot q_0(E)/q_s(E)\leq(1+a^2)^{s/2}$. With $b=\big||a'|-|a|\big|$, $A=\ell+E$ and $k=\ell/4$, the first branch of \eqref{eq:polynomial-gain-ratio} is bounded by $\mathrm e^{-a^2/2}(1+b^2/A)^{-k}$. By \eqref{eq:rw34-beta},
	\begin{equation}\label{eq:polynomial-tail-integral}
		\int_0^\infty (b^\gamma+\langle y\rangle^{\gamma-1}b)(1+b^2/A)^{-k}\dd b
		\leq C_\gamma\bigl((A/k)^{(\gamma+1)/2}+\langle y\rangle^{\gamma-1}A/k\bigr).
	\end{equation}
	For $R=|v|\geq\sqrt\ell$, we have $A\leq2R^2$. The Gaussian branch contributes at most $C_{s_*}R^{\gamma-2}$ by \eqref{eq:gain-angular-integrals}. For the first branch use $|a-a'|\leq b+2|a|$ in the row and exchange $(a,E),(a',E')$ in the column, where the factor is bounded by $\mathrm e^{-ka^2/A-(a')^2/4}$. The two bounds are
	\begin{align*}
		&\frac{C}{R}\big((A/k)^{(\gamma+1)/2}+R^{\gamma-1}A/k\big)
			\leq C\nu(v)\big(k^{-(\gamma+1)/2}+k^{-1}\big),\\
		&C\int_{-1}^1\mathrm e^{-kt^2/2}(1+R^\gamma|t|^\gamma+R^{\gamma-1}(1+R|t|))\dd t
			\leq C\nu(v)\big(k^{-(\gamma+1)/2}+k^{-1}+k^{-1/2}R^{-1}\big).
	\end{align*}
	Thus
	\begin{equation}\label{eq:polynomial-absolute-row}
		\int\big(k_s^+(v,v')+k_s^+(v',v)\big)\dd v'\leq C_{s_*}\nu(v)(\ell^{-1}+R^{-2}+\ell^{-1/2}R^{-1}),\qquad R\geq\sqrt\ell.
	\end{equation}
	For the starred loss, by \eqref{eq:weight-order} and the monotonicity in $m$,
	\[\frac{\mathrm e^{-E/2}}{q_s(E)}\leq\langle v\rangle^s\mathrm e^{-E/4},\qquad q_s(E)\leq\mathrm e^{-\Phi_{m_0,\ell}(E)/4}\qquad(m\geq m_0).\]
	Using \eqref{eq:q-moment-uniform} for the integrals in $v'$ gives
	\begin{equation}\label{eq:polynomial-loss-row}
		\frac1{\nu(v)}\int\big(k_s^-(v,v')+k_s^-(v',v)\big)\dd v'\leq C_{s_*}\big((1+E)^{s/2}\mathrm e^{-E/4}+\mathrm e^{-\Phi_{m_0,\ell}(E)/4}\big).
	\end{equation}
	Choose $C_{s_*}/\ell<1/4$. By \eqref{eq:polynomial-absolute-row} and \eqref{eq:polynomial-loss-row}, each row-plus-column integral is at most $\nu(v)/2$ on $|v|\geq R_{\rm v}$ for $R_{\rm v}$ large. On $B_{R_{\rm v}}$ use \eqref{eq:polynomial-gain-ratio}, \eqref{eq:rw34-beta} and \eqref{eq:polynomial-loss-row}. Thus for $k_s=k_s^++k_s^-$,
	\begin{equation}\label{eq:absolute-row-sum}
		\int\bigl(k_s(v,v')+k_s(v',v)\bigr)\dd v'\leq\nu(v)\bigl(1+C_{s_*}{\bf1}_{B_{R_{\rm v}}}(v)\bigr).
	\end{equation}
	Writing $G=\langle v\rangle^sg$, \eqref{eq:absolute-row-sum} gives
	\begin{align*}
		\big|\langle(\nu-L_w)g,\langle v\rangle^{2s}g\rangle\big|
		&\leq\tfrac12\iint (k_s(v,v')+k_s(v',v))|G(v)|^2\dd v'\dd v\\
		&\leq\tfrac12\|\nu^{1/2}G\|_2^2+C_{s_*}\langle R_{\rm v}\rangle^{2s}\|\nu^{1/2}g\|_2^2.
	\end{align*}
	Subtract from $\|\nu^{1/2}G\|_2^2$ to get \eqref{eq:polynomial-gap}.

	By energy conservation, $\langle v\rangle^s\leq C_s(\langle v'\rangle^s+\langle v_*'\rangle^s)$. Hence
	\begin{align*}
		\langle v\rangle^s|Q_w^+(g_1,g_2)|
		&\leq C_sQ_w^+(\langle v\rangle^s|g_1|,|g_2|)+C_sQ_w^+(|g_1|,\langle v\rangle^s|g_2|),\\
		\langle v\rangle^s Q_w^-(g_1,g_2)&=Q_w^-(g_1,\langle v\rangle^sg_2).
	\end{align*}
	Use \eqref{eq:new-bilinear} on the two gains and \eqref{eq:weighted-collision-loss} on the loss. This proves \eqref{eq:polynomial-bilinear}. Substituting $\widehat v=V+\sqrt\Theta\,\widetilde v$ transfers the estimates to $M,W$ with uniform constants on $\mathcal U$.
\end{proof}

\begin{lemma}\label{lem:log-algebra}
	Let $a\geq0$ and $b>0$. For $0<\eps\leq1$,
	\begin{equation}\label{eq:log-spatial}
		\sup_{t<T,\,x}m_\eps(t,x)^a\theta(t,x)^{-b}\leq C_{a,b,\theta_*}(1+|\log\eps|)^a.
	\end{equation}
	If $|A|\leq C\eps^{1/2}\theta^{-b}$, then $\sup_{t<T,\,x}m_\eps^a|A|\leq C'\eps^{1/2}(1+|\log\eps|)^a\to0$ as $\eps\downarrow0$.
\end{lemma}
\begin{proof}
	By \eqref{eq:adaptive-m} and $0<\eps\leq1$,
	\[m_\eps\leq m_0+\kappa_{\rm w}|\log\eps|+\kappa_{\rm w}\log(1+\theta^{\delta_{\rm gap}})\leq C(1+|\log\eps|)(1+\log(1+\theta^{\delta_{\rm gap}})).\]
	Using $\theta\geq\theta_*$, we obtain
	\begin{align*}
		\sup_{t<T,\,x}m_\eps^a\theta^{-b}
		&\leq C(1+|\log\eps|)^a\sup_{\theta\geq\theta_*}\theta^{-b}(1+\log(1+\theta^{\delta_{\rm gap}}))^a\\
		&\leq C_{a,b,\theta_*}(1+|\log\eps|)^a .
	\end{align*}
	Multiplying \eqref{eq:log-spatial} by $C\eps^{1/2}$ gives the last assertion.
\end{proof}

\section{The moment equations and the Euler linearization}\label{sec:moments}
We derive the moment and microscopic equations \eqref{eq:moment-system} and \eqref{eq:roadmap-microscopic-equation}. The localization of the Euler evolution then gives \eqref{eq:actual-a-equation} and \eqref{eq:observation-local-map}.

\subsection{The moment and microscopic equations}\label{subsec:moment-microscopic}

The $1$ and $v$ moments of \eqref{Eq.Boltzmann} give the first two equations of \eqref{eq:moment-system}. The $|v|^2/2$ moment and the momentum equation multiplied by $\rho_fu_f$ give
\begin{equation}\label{eq:moment-energy-subtraction}
	\begin{aligned}
		&\partial_t\left(\tfrac12\rho_f|u_f|^2+\tfrac32\rho_fT_f\right)
		+\Div\left[\left(\tfrac12\rho_f|u_f|^2+\tfrac52\rho_fT_f\right)u_f+\Pi u_f+q\right]=0,\\
		&\partial_t\left(\tfrac12\rho_f|u_f|^2\right)
		+\Div\left(\tfrac12\rho_f|u_f|^2u_f\right)
		+u_f\cdot\nabla(\rho_fT_f)+u_f\cdot\Div\Pi=0 .
	\end{aligned}
\end{equation}
	Subtract the second line of \eqref{eq:moment-energy-subtraction} from the first. By the mass equation and $\Pi=\Pi^{\mathsf T}$,
	\begin{align*}
		\partial_t(\rho_fT_f)+\Div(\rho_fT_fu_f)&=\rho_fD_fT_f,\\
		\Div(\rho_fT_fu_f)-u_f\cdot\nabla_x(\rho_fT_f)&=\rho_fT_f\Div u_f,\\
		\Div(\Pi u_f)-u_f\cdot\Div\Pi&=\Pi:\nabla_xu_f .
	\end{align*}
	These three identities give the temperature equation in \eqref{eq:moment-system}.

\begin{proof}[Proof of Lemma~\ref{prop:pressure}]
	The identities $\rho_fT_f=\theta_f^5$ and $T_f/\rho_f^{2/3}=\mathrm e^S$ give \eqref{eq:reconstruction}. Substitution and \eqref{eq:moment-system} yield
	\begin{align*}
		\rho_f^{-1}\nabla_x(\rho_f T_f)&=\theta_f^{-3}\mathrm e^{3S/5}\nabla_x(\theta_f^5)=5\theta_f\mathrm e^{3S/5}\nabla_x\theta_f,\\
		5D_f\log\theta_f&=\frac{D_f\rho_f}{\rho_f}+\frac{D_fT_f}{T_f}
			=-\tfrac53\Div u_f-\frac{2\mathcal C_f}{3\theta_f^5},\\
		D_fS&=\frac{D_fT_f}{T_f}-\frac{2D_f\rho_f}{3\rho_f}=-\frac{2\mathcal C_f}{3\theta_f^5}.
	\end{align*}
	Multiply the second line by $\theta_f/5$ to get \eqref{eq:pressure-system}.
\end{proof}

\begin{proof}[Proof of Lemma~\ref{lem:roadmap-projection-microscopic}]
	Differentiate the absolute moments in \eqref{eq:physical-parameters} at fixed $v$. By $\int(v-u_f)f\dd v=0$,
	\[\partial_t\rho_f=m_0[\partial_tf],\qquad \rho_f\partial_tu_f=m_1[\partial_tf],\qquad 3\rho_f\partial_tT_f=m_2[\partial_tf].\]
	Substitution in \eqref{eq:maxwellian-parameter-differential} yields
	\begin{equation}\label{eq:roadmap-Maxwellian-time-derivative}
		\begin{aligned}
			\partial_tM&=M\left(\frac{m_0[\partial_tf]}{\rho_f}+\frac{(v-u_f)\cdot m_1[\partial_tf]}{\rho_fT_f}
			+\frac{|v-u_f|^2-3T_f}{6\rho_fT_f^2}m_2[\partial_tf]\right)\\
			&=\mathcal P_M(\partial_tf).
		\end{aligned}
	\end{equation}
	The absolute $(2+\gamma)$-moment ensures $m_j[Q(f,f)]=0$ for $j=0,1,2$. Substituting \eqref{Eq.Boltzmann} into \eqref{eq:roadmap-Maxwellian-time-derivative}, we get
	\begin{align*}
		DM&=v\cdot\nabla_xM-\mathcal P_M(v\cdot\nabla_x(M+H))\\
		&=(I-\mathcal P_M)(v\cdot\nabla_xM)-\mathcal P_M(v\cdot\nabla_xH),\\
		DH&=\eps^{-1}Q(f,f)-DM.
	\end{align*}
	This is \eqref{eq:roadmap-microscopic-equation}. To compute $S_M$ in \eqref{eq:roadmap-projected-sources}, differentiate before substituting $v=u_f+\sqrt{T_f}\widetilde v$. Then
	\begin{align*}
		M^{-1}v\cdot\nabla_x M
		&=v\cdot\nabla_x\log\rho_f+\frac1{\sqrt{T_f}}\sum_{i,j=1}^3v_i\widetilde v_j\partial_{x_i}(u_f)_j
			+\frac{|\widetilde v|^2-3}{2T_f}v\cdot\nabla_xT_f .
	\end{align*}
	In \eqref{eq:roadmap-projection-definition} use $\int\widetilde v_i\widetilde v_j\mu\dd\widetilde v=\delta_{ij}$ and $\int\widetilde v_i\widetilde v_j|\widetilde v|^2\mu\dd\widetilde v=5\delta_{ij}$. We have
	\begin{equation}\label{eq:projected-maxwellian-polynomials}
		\begin{aligned}
			(I-\mathcal P_M)(M\widetilde v_i\widetilde v_j)&=M(\widetilde v_i\widetilde v_j-\tfrac13\delta_{ij}|\widetilde v|^2),\\
			(I-\mathcal P_M)\bigl(M\widetilde v_i(|\widetilde v|^2-3)\bigr)&=M\widetilde v_i(|\widetilde v|^2-5).
		\end{aligned}
	\end{equation}
	Since $I-\mathcal P_M$ vanishes on $M\operatorname{span}\{1,\widetilde v_i,|\widetilde v|^2\}$, the terms in \eqref{eq:projected-maxwellian-polynomials} are multiplied by $\partial_{x_i}(u_f)_j$ and $\partial_{x_i}\sqrt{T_f}$, respectively. This gives $S_M$.

	For $T_H$, the zero moments \eqref{eq:zero-moments} in the definitions of $\Pi,q$ give
	\begin{equation}\label{eq:roadmap-raw-moment-identities}
		\Pi=\int v\otimes vH\dd v,\qquad \tfrac12\int|v|^2vH\dd v=q+\Pi u_f.
	\end{equation}
	Taking the divergence and using $\Pi=\Pi^{\mathsf T}$, we have
	\begin{equation}\label{eq:roadmap-raw-heat-source}
		\mathcal C_f=\Pi:\nabla_xu_f+\Div q=\tfrac12\Div\int|v|^2vH\dd v-u_f\cdot\Div\Pi.
	\end{equation}
	Therefore
	\begin{equation}\label{eq:microscopic-transport-moments}
		\begin{aligned}
			m_0[v\cdot\nabla_xH]&=\Div\int vH\dd v=0,\qquad m_1[v\cdot\nabla_xH]=\Div\Pi,\\
			m_2[v\cdot\nabla_xH]&=\Div\int |v|^2vH\dd v-2u_f\cdot\Div\Pi=2\mathcal C_f .
		\end{aligned}
	\end{equation}
	Insert \eqref{eq:microscopic-transport-moments} in \eqref{eq:roadmap-projection-definition} to get the second line of \eqref{eq:roadmap-projected-sources}.
\end{proof}

	At fixed $v$, $[D,\partial_{x_i}]=0$. For the moments use \eqref{eq:source-moments-spatial}. In \eqref{eq:normalized-collision-frequency}, the change $v_*=u+\theta\widehat v_*$ gives $\dd v_*=\theta^3\dd\widehat v_*$ and $|v-v_*|^\gamma=\theta^\gamma|\widehat v-\widehat v_*|^\gamma$, hence
	\begin{equation}\label{eq:physical-collision-frequency}
		\nu_M^{\rm phys}(t,x,u+\theta\widehat v)=\theta^{3+\gamma}\widehat\nu_M(t,x,\widehat v).
	\end{equation}
	See Lemma~\ref{lem:matched-source} for the derivatives of $S_0$.

\subsection{Localization of the unstable coordinates}\label{subsec:unstable-localization}

\begin{proof}[Proof of \eqref{eq:euler-weighted-derivatives}]
	For $F\in C_c^\infty$, Plancherel and interpolation give
	\[\sum_{|\alpha|=2k}\|\partial_z^\alpha F\|_2^2\leq C\|\Delta_z^k F\|_2^2,\qquad \|\nabla_z^jF\|_2\leq\eps\|\nabla_z^{2k}F\|_2+C_{\eps,j}\|F\|_2\quad(0<j<2k).\]
	Use the interior estimate for $\Delta_z^k$ on $1<|z|<2$ in $1/2<|z|<4$. Rescaling by $y=Rz$ gives
	\begin{equation}\label{eq:rw56-annulus-elliptic}
		\sum_{j=0}^{2k}R^{2j}\|\nabla_y^j\mathsf W\|_{L^2(R<|y|<2R)}^2
		\leq C\bigl(R^{4k}\|\Delta_y^k\mathsf W\|_{L^2(R/2<|y|<4R)}^2+\|\mathsf W\|_{L^2(R/2<|y|<4R)}^2\bigr),\quad R\geq1.
	\end{equation}
	Sum \eqref{eq:rw56-annulus-elliptic} with weight $R^\varpi$ over $R=2^n$, $n\geq0$, and apply the same interior estimate on $B_2$ in $B_4$. By bounded overlap and $\phi_1\asymp\langle y\rangle$,
	\[\sum_{j=0}^{2k}\|\langle y\rangle^{\varpi/2+j}\nabla_y^j\mathsf W\|_2^2\leq C\Big(\|\langle y\rangle^{\varpi/2+2k}\Delta_y^k\mathsf W\|_2^2+\|\langle y\rangle^{\varpi/2}\mathsf W\|_2^2\Big)\leq C\|\mathsf W\|_{X_\varpi^{2k}}^2.\]
	Conversely expand $\Delta_y^k=\sum_{|\alpha|=k}(k!/\alpha!)\partial_y^{2\alpha}$ in \eqref{eq:actual-X-norm}. This yields \eqref{eq:euler-weighted-derivatives} by density.
\end{proof}

\begin{proof}[Proof of Proposition~\ref{prop:projection-local}]
	Let $\supp F\subset\{|y|>R_{\rm sp}\}$. Invert the evolution on $E_{\rm u}$ and use \eqref{eq:exceptional-backward-t} and \eqref{eq:exterior-decay-t} for $t_0\leq t<T$ to get
	\begin{align}\label{eq:projection-exterior-bound}
		\|\Pu F\|_{\mathscr X}
		&=\|\mathcal U_E(t_0,t)\Pu\mathcal U_E(t,t_0)F\|_{\mathscr X}\notag\\
		&\leq C\left(\frac{T-t_0}{T-t}\right)^\kappa\|\mathcal U_E(t,t_0)F\|_{\mathscr X}
		\leq C\left(\frac{T-t}{T-t_0}\right)^{1/8-\kappa}\|F\|_{\mathscr X}.
	\end{align}
	Let $t\uparrow T$, using $1/8>\kappa$. Applying $\Pu F=0$ to $(1-\chi_{\rm sp})F$ gives \eqref{eq:euler-projection-local}.

	Fix $R>R_{\rm sp}$ and a radial $\chi\in C_c^\infty(B_R)$ with $\chi=1$ near $\overline B_{R_{\rm sp}}$. For a general field $F=(F_V,F_B,F_S)$ we use the extension $\ell_j(F):=\ell_j(\overline F)$, where
	\[\overline F(y):=\int_{SO(3)}(O^{-1}F_V(Oy),F_B(Oy),F_S(Oy))\dd O\]
	and $\dd O$ is Haar probability measure. By \eqref{eq:projection-exterior-bound}, $\ell_j(\overline F)=\ell_j(\chi\overline F)$. Then \eqref{eq:euler-weighted-derivatives} and $OB_R=B_R$ give
	\[
		|\ell_j(F)|\leq C_j\|\chi\overline F\|_{\mathscr X}\leq C_{j,R}\|\overline F\|_{H^K(B_R)}\leq C_{j,R}'\|F\|_{H^K(B_R)}.
	\]
	Let $R\downarrow R_{\rm sp}$ in the support inclusion $\supp\ell_j\subset\overline B_R$. This proves \eqref{eq:unstable-functional-local-bound}.
\end{proof}

\begin{lemma}\label{lem:new-spectral-consistency}
	Let $\varpi=\bar\varpi-r/2$ and $K\geq2k_0+2$ be even. On $X^j_\varpi$, $j=K,K-2$, write $\mathcal L_{{\rm E},j}$ for \eqref{eq:explicit-euler-L} with domain $\mathcal D_{j,\varpi}$. In $\Re\zeta>-1/8$ the spectra and generalized eigenspaces coincide. In \eqref{eq:two-order-projection-agreement} the coefficients in the common basis satisfy $\ell_j^{(K-2)}(F)=\ell_j(F)$ for $F\in X^K_\varpi$ and $1\leq j\leq J$.
\end{lemma}
\begin{proof}
	Uniqueness on $X^{K-2}_\varpi$ gives $\mathrm e^{q\mathcal L_{{\rm E},K}}F=\mathrm e^{q\mathcal L_{{\rm E},K-2}}F$ for $F\in X^K_\varpi$. Taking Laplace transforms to the right of both growth bounds gives
	\begin{equation}\label{eq:rw56-resolvent-agreement}
		(\zeta-\mathcal L_{{\rm E},K})^{-1}F=(\zeta-\mathcal L_{{\rm E},K-2})^{-1}F\quad\text{in }X^{K-2}_\varpi,\qquad F\in X^K_\varpi .
	\end{equation}
	Propositions~\ref{prop:euler-compact} and \ref{prop:euler-spectrum} give discrete spectra of finite algebraic multiplicity in $\Re\zeta>-1/8$. By \cite[Lemma 4.8]{BCGJVY}, each chain field $F$ at order $K-2$ has $F,(\mathcal L_{\rm E}-\zeta)F\in X^K_\varpi$. Thus $F\in\mathcal D_{K,\varpi}$. Since $\mathcal D_{K,\varpi}\subset\mathcal D_{K-2,\varpi}$,
	\begin{equation}\label{eq:euler-generalized-eigenspaces}
		\ker(\mathcal L_{{\rm E},K}-\zeta)^m=\ker(\mathcal L_{{\rm E},K-2}-\zeta)^m\qquad(\Re\zeta>-1/8,\ m\geq1).
	\end{equation}
	By \eqref{eq:euler-generalized-eigenspaces}, the two pole sets agree in the half-plane. Extend \eqref{eq:rw56-resolvent-agreement} to their connected complement by the identity theorem. Integrate over $\Gamma_{\rm u}$ to obtain \eqref{eq:two-order-projection-agreement} in the common basis.
\end{proof}

On $|x|<R_{\rm e}$, \eqref{eq:profile} gives $(T-t)\theta\nabla_x=\vartheta\nabla_y$ and $(T-t)\nabla_xu=\nabla_yU$. In \eqref{eq:expanded-flux-vector} this yields
\[(G_{\rm mom},G_{\rm heat})=(T-t)(G_{{\rm mom},x}^\sharp,C_H^\sharp).\]
Write the pullbacks without $\sharp$ in the formulas below.
\begin{align}\label{eq:dimensionless-forces}
	G_{\rm mom}=\vartheta\Div_y\pi+5\pi\nabla_y\vartheta,\quad
	G_{\rm heat}=\pi:\nabla_yU+\vartheta\Div_yq_H+6q_H\cdot\nabla_y\vartheta-V\cdot G_{\rm mom}.
\end{align}
In \eqref{eq:dimensionless-forces}, $(\Div_y\pi)_i=\sum_j\partial_{y_j}\pi_{ij}$ and $\pi:\nabla_yU=\sum_{i,j}\pi_{ij}\partial_{y_j}U_i$. By \eqref{eq:euler-microscopic-source}, the components of $F_h=(F_V,F_B,F_S)$ are
\begin{align*}
	F_V=-(1+B)^{-3}\mathrm e^{3S/5}G_{\rm mom},\quad
	F_B=-\tfrac2{15}(1+B)^{-4}G_{\rm heat},\quad F_S=-\tfrac23(1+B)^{-5}G_{\rm heat},
\end{align*}

By \eqref{eq:fixed-spectral-cutoff}, $|(T-t)^{1/r}y|\leq T^{1/r}R_2<R_{\rm e}/2$ on $\supp\chi_{\rm sp}$ for $0\leq t<T$. Thus \eqref{eq:actual-W-equation} gives on $\R^3_y$
\begin{equation}\label{eq:cutoff-actual-euler}
 (T-t)\partial_t(\chi_{\rm sp}\mathsf W)
 =\mathcal L_{\rm E}(\chi_{\rm sp}\mathsf W)-[\mathcal L_{\rm E},\chi_{\rm sp}]\mathsf W
   +\chi_{\rm sp}\mathsf T(N+F_h).
\end{equation}
The commutator from \eqref{eq:explicit-euler-L} is
\begin{align*}
	[\mathcal L_{\rm E},\chi_{\rm sp}]\mathsf W=-(b\cdot\nabla_y\chi_{\rm sp})\mathsf W-c_{\rm E}\big(\mathsf W_B\nabla_y\chi_{\rm sp},\,\nabla_y\chi_{\rm sp}\cdot\mathsf W_V,\,0\big).
\end{align*}
As a consequence,
\[\supp\big([\mathcal L_{\rm E},\chi_{\rm sp}]\mathsf W\big)\subset\{R_1\leq|y|\leq R_2\}\subset\{|y|>R_{\rm sp}\}.\]
Proposition~\ref{prop:projection-local} gives $\Pu[\mathcal L_{\rm E},\chi_{\rm sp}]\mathsf W=0$. Applying $\Pu$ to \eqref{eq:cutoff-actual-euler} and commuting $\Pu$ with $\mathcal L_{\rm E}$ gives
\[(T-t)\partial_t\Pu(\chi_{\rm sp}\mathsf W)=\mathcal L_{\rm E}\Pu(\chi_{\rm sp}\mathsf W)+\Pu\{\chi_{\rm sp}\mathsf T(N+F_h)\}.\]
Recall that the coordinate vector $\mathbf a(t)$ is defined in \eqref{eq:roadmap-unstable-coordinates}. The projected equation gives
\begin{equation}\label{eq:actual-a-equation}
	\begin{aligned}
		\mathbf a'(t)&=\frac{A_{\rm u}\mathbf a(t)}{T-t}+F_{\rm u}(t),\\
		(F_{\rm u}(t))_j&:=\frac1{T-t}\ell_j^{(K-2)}\!\left(\chi_{\rm sp}\mathsf T\big(N(t,\cdot)+F_h(t,\cdot)\big)\right),\qquad 1\leq j\leq J.
	\end{aligned}
\end{equation}
The coordinate vector of the projected equation has this form, since $A_{\rm u}$ represents $\mathcal L_{\rm E}|_{E_{\rm u}}$ in the basis $\Phi_1,\ldots,\Phi_J$. For smooth solutions, Lemma~\ref{lem:new-spectral-consistency} identifies the source coordinates obtained from $\ell_j$ and $\ell_j^{(K-2)}$. The expression $(T-t)^{-1}\ell_j^{(K-2)}(\chi_{\rm sp}\mathsf T(N+F_h))$ remains meaningful when the localized source has only $K-2$ spatial derivatives.

In Proposition~\ref{thm:actual-spectral-forcing}, we estimate \eqref{eq:actual-a-equation} by using the norm recovery in Lemma~\ref{lem:single-spatial-allocation}. Since \eqref{eq:fixed-spectral-cutoff} gives $J_x((T-t)^{1/r}y)=1$ on $\supp\chi_{\rm sp}$, the coordinates of $\mathsf V$ in the energy framework are also $\mathbf a(t)$.

\subsection{A localized correction to the Euler energy}\label{subsec:localized-euler-correction}
	Keep $\mathscr X=X^K_\varpi$, $\varpi=\bar\varpi-r/2$ from \eqref{eq:euler-stability-space} and $\lambda_0=\lambda_\varpi=1/8$. Choose a smooth nonnegative radial cutoff $\chi_c$ of compact support with $\chi_c=1$ on $\overline B_{R_\varpi}$ from \eqref{eq:euler-subcritical-coercivity}. Then for a sufficiently large $C_c$,
	\begin{equation}\label{eq:local-observation-coercivity}
		2\langle\mathcal L_{\rm E}\mathsf W,\mathsf W\rangle_{\mathscr X}\leq -2\lambda_0\|\mathsf W\|_{\mathscr X}^2+C_c\int\chi_c|\mathsf W|_{D_*}^2\dd y.
	\end{equation}
	Here $D_*=\operatorname{diag}(I_3,1,3/2)$ and the real pairing $\langle\cdot,\cdot\rangle_0$ is as in \eqref{eq:euler-zero-order-inner-product}.

	Put $S_{\rm E}(q)=\mathrm e^{q\mathcal L_{\rm E}}$, $\mathcal L_{{\rm E},*}=\mathcal L_{\rm E}+h_*I$ and $S_*(q)=\mathrm e^{h_*q}S_{\rm E}(q)$. By \eqref{eq:euler-spectral-evolution},
	\begin{equation}\label{eq:shifted-stable-decay}
		\|S_*(q)(I-\Pu)\mathsf W\|_{\mathscr X}\leq C_{\rm E}\mathrm e^{-(\sigma_*-h_*)q}\|\mathsf W\|_{\mathscr X},\qquad q\geq0,\quad d_*/2<h_*<\min\{\sigma_*,\lambda_0\}.
	\end{equation}
	For $t_{\rm obs}>0$ to be fixed, set
	\begin{equation}\label{eq:observation-correction-form}
		\mathcal R_*(\mathsf W,\mathsf Y):=C_c\int_0^{t_{\rm obs}}\langle\chi_c S_*(q)\mathsf W,S_*(q)\mathsf Y\rangle_0\dd q .
	\end{equation}

\begin{proof}[Proof of Proposition~\ref{thm:local-observation}]
	Write $\mathcal L_{\rm E}=-\sum_j A_j\partial_{y_j}+B_{\rm E}$. By \eqref{eq:explicit-euler-L},
	\[A_j=b_jI_5+c_{\rm E}\begin{pmatrix}0_{3\times3}&e_j&0_{3\times1}\\e_j^{\mathsf T}&0&0\\0_{1\times3}&0&0\end{pmatrix},\qquad D_*A_j=A_j^{\mathsf T}D_* .\]
	By $b=y/r+U$ and \eqref{eq:profile-decay},
	\begin{equation}\label{eq:observation-coefficient-bounds}
		|A_j(y)|\leq C\langle y\rangle,\qquad
		\|\nabla_y^mA_j\|_\infty+\|\nabla_y^{m-1}B_{\rm E}\|_\infty\leq C_m\qquad(m\geq1).
	\end{equation}
	Integration by parts in $\langle\cdot,\cdot\rangle_0$ gives
	\[\mathcal L_{\rm E}^\dagger=\sum_j A_j\partial_{y_j}+B_{\rm E}^\dagger,\qquad B_{\rm E}^\dagger=D_*^{-1}\left(\sum_j\partial_{y_j}(D_*A_j)+B_{\rm E}^{\mathsf T}D_*\right).\]
	For $\partial_q\mathsf Y=\mathcal A\mathsf Y$, $\mathcal A=\mathcal L_{{\rm E},*}$ or $\mathcal L_{{\rm E},*}^\dagger$,
	\begin{align*}
		\frac{\dd}{\dd q}\int|\mathsf Y|_{D_*}^2\dd y
		&=\langle(\mathcal L_{{\rm E},*}+\mathcal L_{{\rm E},*}^\dagger)\mathsf Y,\mathsf Y\rangle_0\\
		&=\int \mathsf Y^{\mathsf T}\left(\sum_j\partial_{y_j}(D_*A_j)+D_*B_{\rm E}
			+B_{\rm E}^{\mathsf T}D_*+2h_*D_*\right)\mathsf Y\dd y .
	\end{align*}
	By \eqref{eq:observation-coefficient-bounds}, the matrix is bounded, and
	\[(\partial_q-\mathcal A)\partial_{y_i}\mathsf Y=[\partial_{y_i},\mathcal A]\mathsf Y,\qquad \|[\partial_{y_i},\mathcal A]\mathsf Y\|_2\leq C\|\mathsf Y\|_{H^1}.\]
	Consequently
	\begin{equation}\label{eq:rw56-euler-adjoint-energy}
		\frac{\dd}{\dd q}\sum_{|\alpha|\leq j}\int|\partial_y^\alpha\mathsf Y|_{D_*}^2\dd y
		\leq C\sum_{|\alpha|\leq j}\int|\partial_y^\alpha\mathsf Y|_{D_*}^2\dd y\qquad(j=0,1).
	\end{equation}
	Also
	\[\frac{\dd}{\dd s}\langle\mathrm e^{s\mathcal L_{\rm E}^\dagger}F,\mathrm e^{(q-s)\mathcal L_{\rm E}}G\rangle_0=0\qquad(F,G\in C_c^\infty).\]
	Evaluate at $s=0,q$ and use \eqref{eq:rw56-euler-adjoint-energy} with $j=0$ to get $\mathrm e^{q\mathcal L_{\rm E}^\dagger}=S_{\rm E}(q)^\dagger$ on $L^2$.

	On $B_{R(q)}$, \eqref{eq:observation-coefficient-bounds} gives
	\begin{equation}\label{eq:observation-boundary-flux}
		\int_{|y|=R(q)}\left(R'(q)|\mathsf Y|_{D_*}^2\mp\mathsf Y^{\mathsf T}D_*\sum_i A_i\frac{y_i}{|y|}\mathsf Y\right)\dd S_y
		\leq\bigl(R'(q)+C(1+R(q))\bigr)\int_{|y|=R(q)}|\mathsf Y|_{D_*}^2\dd S_y ,
	\end{equation}
	with the upper sign for $\mathcal L_{{\rm E},*}$. Take $R(q)=(1+R)\mathrm e^{C(q_{\max}-q)}-1$, so $R'(q)+C(1+R(q))=0$. Apply \eqref{eq:observation-boundary-flux} in \eqref{eq:rw56-euler-adjoint-energy} for $|\alpha|\leq j$ and $R_+=R(0)$. Since $R(q)\geq R$,
	\begin{equation}\label{eq:rw56-local-evolution}
		\|S_*(q)\mathsf Y\|_{H^j(B_R)}+\|S_*(q)^\dagger\mathsf Y\|_{H^j(B_R)}
		\leq C_{q_{\max},R,j}\|\mathsf Y\|_{H^j(B_{R_+})}\qquad(0\leq q\leq q_{\max},\ j=0,1).
	\end{equation}
	For the exterior energy, $R(q)=(1+R)\mathrm e^{Cq}-1$ gives $-R'(q)+C(1+R(q))=0$. Thus
	\begin{equation}\label{eq:observation-propagation}
		\supp\mathsf Y(0)\subset B_R\quad\Longrightarrow\quad
		\supp S_*(q)\mathsf Y(0),\ \supp S_*(q)^\dagger\mathsf Y(0)\subset\overline B_{(1+R)\mathrm e^{Cq}-1}.
	\end{equation}
	By density these estimates hold in $H^j$. The local estimate \eqref{eq:rw56-local-evolution} also applies on $B_R$ to $\mathsf W\in\mathscr X$.

	By \eqref{eq:observation-correction-form}, $\mathcal R_*(\mathsf W,\mathsf Y)=\langle\mathscr R_*\mathsf W,\mathsf Y\rangle_0$ with
	\begin{equation}\label{eq:observation-operator}
		\mathscr R_*=C_c\int_0^{t_{\rm obs}}S_*(q)^\dagger\chi_cS_*(q)\dd q.
	\end{equation}
	Choose $R$ and $R_c$ with $\supp\chi_c\subset B_R$ and $R_c>2((1+R)\mathrm e^{Ct_{\rm obs}}-1)$. Equations \eqref{eq:rw56-local-evolution} and \eqref{eq:observation-propagation} give, for $0\leq q\leq t_{\rm obs}$,
	\begin{align*}
		\|S_*(q)^\dagger\chi_cS_*(q)\mathsf W\|_{H^j}&\leq C\|\chi_cS_*(q)\mathsf W\|_{H^j}\leq C_j\|\mathsf W\|_{H^j(B_{R_c})},\\
		\supp(S_*(q)^\dagger\chi_cS_*(q)\mathsf W)&\subset B_{R_c/2}.
	\end{align*}
	Integrate in \eqref{eq:observation-operator} to obtain \eqref{eq:observation-local-map}. In particular, with $j=0$,
	\[0\leq\mathcal R_*(\mathsf W,\mathsf W)=\langle\mathscr R_*\mathsf W,\mathsf W\rangle_0\leq C\|\mathsf W\|_{L^2(B_{R_c})}^2,\qquad
	\tfrac12\|\mathsf W\|_{\mathscr X}^2\leq\mathcal E_*(\mathsf W)\leq C\|\mathsf W\|_{\mathscr X}^2.\]
	
	For $\mathsf W\in D(\mathcal L_{\rm E})$, we have $\partial_qS_*(q)\mathsf W=S_*(q)\mathcal L_{{\rm E},*}\mathsf W$ in $\mathscr X$. By \eqref{eq:observation-correction-form},
	\begin{equation}\label{eq:rw56-correction-endpoints}
		\begin{aligned}
			2\mathcal R_*(\mathcal L_{{\rm E},*}\mathsf W,\mathsf W)&=C_c\int_0^{t_{\rm obs}}\frac\dd{\dd q}\int\chi_c|S_*(q)\mathsf W|_{D_*}^2\dd y\dd q\\
			&=C_c\int\chi_c|S_*(t_{\rm obs})\mathsf W|_{D_*}^2\dd y-C_c\int\chi_c|\mathsf W|_{D_*}^2\dd y.
		\end{aligned}
	\end{equation}
	Apply \eqref{eq:shifted-stable-decay} on $(I-\Pu)\mathscr X$ and use finite dimensionality on $E_{\rm u}$ to obtain
	\[\begin{aligned}
		C_c\int\chi_c|S_*(t_{\rm obs})(I-\Pu)\mathsf W|_{D_*}^2\dd y&\leq CC_cC_{\rm E}^2\mathrm e^{-2(\sigma_*-h_*)t_{\rm obs}}\|\mathsf W\|_{\mathscr X}^2,\\
		C_c\int\chi_c|S_*(t_{\rm obs})\Pu\mathsf W|_{D_*}^2\dd y&\leq C_{t_{\rm obs}}|\mathbf a|^2.
	\end{aligned}\]
	We choose $t_{\rm obs}$ such that $2CC_cC_{\rm E}^2\mathrm e^{-2(\sigma_*-h_*)t_{\rm obs}}\leq\lambda_0-h_*$. Then
	\[C_c\int\chi_c|S_*(t_{\rm obs})\mathsf W|_{D_*}^2\dd y\leq(\lambda_0-h_*)\|\mathsf W\|_{\mathscr X}^2+C_{t_{\rm obs}}|\mathbf a|^2.\]
	Adding \eqref{eq:rw56-correction-endpoints} to \eqref{eq:local-observation-coercivity}, with the contribution of $h_*I$ included, we get
	\[2\big\{\langle\mathcal L_{{\rm E},*}\mathsf W,\mathsf W\rangle_{\mathscr X}+\mathcal R_*(\mathcal L_{{\rm E},*}\mathsf W,\mathsf W)\big\}\leq-(\lambda_0-h_*)\|\mathsf W\|_{\mathscr X}^2+C_{t_{\rm obs}}|\mathbf a|^2.\]
	Since $\|\mathsf W\|_{\mathscr X}^2+\mathcal R_*(\mathsf W,\mathsf W)=2\mathcal E_*(\mathsf W)$, subtracting the contribution of $h_*I$, we obtain
	\[\langle\mathcal L_{\rm E}\mathsf W,\mathsf W\rangle_{\mathscr X}+\mathcal R_*(\mathcal L_{\rm E}\mathsf W,\mathsf W)\leq-2h_*\mathcal E_*(\mathsf W)-\frac{\lambda_0-h_*}{2}\|\mathsf W\|_{\mathscr X}^2+C_*|\mathbf a|^2.\]
	This proves \eqref{eq:observation-stable-sign} with $c_*:=(\lambda_0-h_*)/2$.
\end{proof}

\section{Matched energies and norm recovery}\label{sec:energies}
\subsection{Proof of the exact matching identities}\label{subsec:exact-matching}

\begin{proof}[Proof of Proposition~\ref{prop:global-matching}]
	By \eqref{eq:differentiated-zero-moments} at fixed $v$ and \eqref{eq:explicit-Maxwellian-derivative},
	\[P_i/M\in\operatorname{span}\{1,v_1,v_2,v_3,|v|^2\},\qquad \int H_iP_i/M\dd\widehat v=\theta^{-3}\int H_iP_i/M\dd v=0,\quad i=0,K.\]
	At $x=(T-t)^{1/r}y$, use \eqref{eq:roadmap-weighted-derivative-identity}, \eqref{eq:weighted-macro-variable} and $\bar\varpi+3=6r\alpha_{\rm ss}$ to get
	\begin{equation}\label{eq:matched-change-of-variables}
		\begin{aligned}
			P_K&=D_{\mathbf z}M[\mathsf T^{-1}\phi_1^KJ_x^{-1}\Delta_y^k\mathsf V],\qquad P_0=D_{\mathbf z}M[\mathsf T^{-1}J_x^{-1}\mathsf V],\\
			\Omega_b^*\dd x&=(T-t)^{b/r-6\alpha_{\rm ss}+(\bar\varpi-b+3)/r}J_x^2\vartheta^6\langle y\rangle^{\bar\varpi-b}\dd y
			=J_x^2\vartheta^6\langle y\rangle^{\bar\varpi-b}\dd y.
		\end{aligned}
	\end{equation}
	Using \eqref{eq:matched-change-of-variables} in \eqref{eq:fitted-metric-W}, we have
	\begin{align*}
		\iint\Omega_b^*\frac{|P_K|^2}{M}\dd\widehat v\dd x&=c_{\rm s}^2\int\langle y\rangle^{\bar\varpi-b}\phi_1^{2K}|\Delta_y^k\mathsf V|_{D(\mathbf z^\sharp)}^2\dd y,\\
		\iint\Omega_b^*\frac{|P_0|^2}{M}\dd\widehat v\dd x&=c_{\rm s}^2\int\langle y\rangle^{\bar\varpi-b}|\mathsf V|_{D(\mathbf z^\sharp)}^2\dd y.
	\end{align*}
	This gives \eqref{eq:new-Gaussian-match} by \eqref{eq:roadmap-natural-Euler-energy}.

	Since $[D,\Delta_x^k]=0$ at fixed $v$, differentiating \eqref{eq:coupled-perturbation-identities} gives
	\begin{align*}
		DG_0+D\mathcal B_0&=D(f-\bar M),\\
		D(d_\phi^{-K}G_K)+D\mathcal B_K&=\Delta_x^kD(f-\bar M),\qquad D(f-\bar M)=\eps^{-1}Q(f,f)-S_0 .
	\end{align*}
	Multiplying the second line by $d_\phi^K$ yields \eqref{eq:new-good-evolution}. In $\mathcal B_K$, \eqref{eq:roadmap-Maxwellian-chain-rule} and \eqref{eq:roadmap-conjugated-Laplacian} leave
	\[\mathcal B_K=\mathcal R_K-D_{\mathbf z}M[\Sigma_x^{-1}[\Delta_x^k,\Sigma_x]\mathbf z].\]
	Both terms have order at most $K-1$ in $\mathbf z$. By \eqref{eq:roadmap-raw-moment-identities} and \eqref{eq:roadmap-raw-heat-source},
	\begin{equation}\label{eq:physical-flux-derivatives}
		\begin{aligned}
			\partial_x^\alpha\Div\Pi&=\sum_{j=1}^3\int vv_j\,\partial_x^{\alpha+e_j}H\dd v,\\
			\partial_x^\alpha\mathcal C_f&=\tfrac12\sum_{j=1}^3\int|v|^2v_j\partial_x^{\alpha+e_j}H\dd v
				-\sum_{\beta\leq\alpha}\binom\alpha\beta\partial_x^\beta u_f\cdot\partial_x^{\alpha-\beta}\Div\Pi .
		\end{aligned}
	\end{equation}
	By \eqref{eq:physical-flux-derivatives} in \eqref{eq:normalized-macro-system}, $D\partial_x^\alpha\mathbf z$ has order at most $|\alpha|+1$ in $(\mathbf z,H)$. Since $\mathcal B_0$ has order zero by \eqref{eq:roadmap-composition-remainders}, the orders of $(D\mathcal B_K,D\mathcal B_0)$ are at most $(K,1)$.
\end{proof}

By \eqref{eq:new-global-weights},
\[
 D\log\Omega_b^*+\frac b{r(T-t)}
 =2D\log\Sigma_x+(\bar\varpi-b)D\log d .
\]
The profile bounds \eqref{eq:untruncated-profile-bounds}, $|\partial_t\log d|\leq C\theta^p$ and $|v|\leq C\theta\langle\widetilde v\rangle$ imply
\begin{equation}\label{eq:transport-weight-log-bounds}
 |D\log\Sigma_x|+\left|D\log\Omega_b^*+\frac{b}{r(T-t)}\right|
  \leq C\theta^p\langle\widetilde v\rangle,\qquad b\in\{0,a\}.
\end{equation}
By \eqref{eq:good-length-time},
\[
 D(d_\phi^K)-K d_\phi^K v\cdot\nabla_x\log d_\phi
 =K d_\phi^K\partial_t\log d_\phi\leq0 .
\]

Define $\partial_t^E P_f$ by \eqref{eq:moment-system} with $\Pi=q=0$ for $P_f=(\rho_f,u_f,T_f)$. On explicit $(t,x,v)$ arguments, $\partial_t^E$ is $\partial_t$. Let $D_E=\partial_t^E+v\cdot\nabla_x$ and $\mathscr D_{\rm mic}=\partial_t-\partial_t^E$. Then
\begin{equation}\label{eq:microscopic-moment-time-difference}
 \mathscr D_{\rm mic}(\rho_f,u_f,T_f)
 =\left(0,-\rho_f^{-1}\Div\Pi,-\frac{2\mathcal C_f}{3\rho_f}\right).
\end{equation}
Use $\mathcal Z_K^2(t)=\sum_{j=0}^K\int\Omega_0^*|d^j\nabla_x^j\mathbf z|^2\dd x$. For $(n,V,\Theta)\in\mathcal U$, we have
\begin{equation}\label{eq:macro-equivalence-new}
 \begin{gathered}
  A_b^{\rm E}\asymp\|\mathsf V\|_{X^K_{\bar\varpi-b}}^2
    \asymp\sum_{j=0}^K\int\Omega_b^*|d^j\nabla_x^j\mathbf z|^2\dd x,\qquad b\in\{0,a\},\\
  \sup_x J_x\sum_{j\leq K-2}d^j|\nabla_x^j\mathbf z|\leq C\mathcal Z_K .
 \end{gathered}
\end{equation}
	By \eqref{eq:roadmap-energy-matrix}, $D(\mathbf z)\asymp D_*$. Equations \eqref{eq:spatial-weight-J}, \eqref{eq:profile-decay} and \eqref{eq:matched-change-of-variables} give
	\begin{align*}
		\mathsf V&=J_x\vartheta^3\mathsf T_0\mathbf z^\sharp,\qquad
		|\nabla_y^j(J_x\vartheta^3)^{\pm1}|\leq C_j\langle y\rangle^{-j}(J_x\vartheta^3)^{\pm1},\\
		d^j\nabla_x^j\mathbf z&=\langle y\rangle^j\nabla_y^j\mathbf z^\sharp,\qquad
		\Omega_b^*\dd x=J_x^2\vartheta^6\langle y\rangle^{\bar\varpi-b}\dd y .
	\end{align*}
	Use the product rule in both directions and apply \eqref{eq:euler-weighted-derivatives}. Then
	\[\|\mathsf V\|_{X^K_{\bar\varpi-b}}^2\asymp\sum_{j=0}^K\int J_x^2\vartheta^6\langle y\rangle^{\bar\varpi-b+2j}|\nabla_y^j\mathbf z^\sharp|^2\dd y=\sum_{j=0}^K\int\Omega_b^*|d^j\nabla_x^j\mathbf z|^2\dd x.\]
	For $d_0=d(t,x_0)$, we have $d\asymp d_0$ and $\Omega_0^*\gtrsim J_x(x_0)^2d_0^{-3}$ on $B_{2cd_0}(x_0)$. Apply the $H^2$ embedding with $x=x_0+d_0X$ to obtain, for $j\leq K-2$,
	\[J_x(x_0)^2d_0^{2j}|\nabla_x^j\mathbf z(t,x_0)|^2\leq C\sum_{l=0}^2\int_{B_{2cd_0}(x_0)}\Omega_0^*|d^{j+l}\nabla_x^{j+l}\mathbf z|^2\dd x\leq C\mathcal Z_K^2.\]
	Take the supremum in $x_0$.

\subsection{Weighted elliptic recovery for the microscopic variables}\label{subsec:weighted-norm-recovery}

\begin{lemma}\label{lem:weighted-elliptic-recovery}
	Fix $t$, even $N\geq0$ and $b=0,a$. Suppose $(n,V,\Theta)\in\mathcal U$ and $\mathbf r=(n-1,V,\Theta-1)$ satisfies
	\[d(\theta^{-1}|\nabla_xu|+|\nabla_x\log\theta|+|\nabla_x\mathbf r|)\leq L.\]
	For $\sigma=0,1$, use the physical frequency \eqref{eq:physical-collision-frequency} and set $\mathfrak w=\Omega_b^*\theta^{-3}W^{-2}(\nu_M^{\rm phys})^\sigma$ with $W=W_{m_\eps,\ell}$. Then, with $\partial_x$ at fixed $v$,
	\begin{align}\label{eq:weighted-elliptic-recovery}
		\sum_{|\alpha|=N}\iint\mathfrak w\,|d^N\partial_x^\alpha H|^2\dd v\dd x
		&\leq C_L\iint\mathfrak w\,|d^N\Delta_x^{N/2}H|^2\dd v\dd x\notag\\
		&\quad+C_L\sum_{j<N}\sum_{|\alpha|=j}\iint\mathfrak w\,\langle\widetilde v\rangle^{4(N-j)}|d^j\partial_x^\alpha H|^2\dd v\dd x .
	\end{align}
	The constant $C_L$ is independent of $m\geq m_0$.
\end{lemma}
\begin{proof}
	For even $N\geq2$ set $\rho_v=d\langle\widetilde v\rangle^{-2}$. Differentiating $\theta\sqrt{\Theta}\,\widetilde v=v-u-\theta V$ at fixed $v$ gives
	\[\partial_{x_i}\widetilde v=-\frac{\partial_{x_i}u+\partial_{x_i}(\theta V)}{\theta\sqrt\Theta}
	-\left(\partial_{x_i}\log\theta+\tfrac12\partial_{x_i}\log\Theta\right)\widetilde v.\]
	Hence $d|\nabla_x\widetilde v|\leq C_L\langle\widetilde v\rangle$. By \eqref{eq:new-global-weights}, \eqref{eq:physical-collision-frequency} and \eqref{eq:normalized-collision-frequency},
	\begin{align*}
		\nabla_x\log\mathfrak w
		&=\nabla_x\log\Omega_b^*-3\nabla_x\log\theta-2\nabla_x\log W
			+\sigma\nabla_x\log\nu_M^{\rm phys},\\
		\nabla_x\log\Omega_b^*
		&=2\nabla_x\log J_x+6\nabla_x\log\theta_{\rm p}+(\bar\varpi-b)\nabla_x\log d,\\
		\partial_{x_i}\log\nu_M^{\rm phys}
		&=(3+\gamma)\partial_{x_i}\log\theta+\partial_{x_i}\log n
			+\tfrac\gamma2\partial_{x_i}\log\Theta+(\nabla\log\nu)(\widetilde v)\cdot\partial_{x_i}\widetilde v .
	\end{align*}
	Use \eqref{eq:spatial-weight-J}, \eqref{eq:untruncated-profile-bounds} and \eqref{eq:normalized-frequency-derivatives} for $\Omega_b^*$ and $\nu_M^{\rm phys}$. For $W$, \eqref{eq:Dm} and \eqref{eq:fitted-log-weight} give $d|\nabla_x\log W|\leq C_L\langle\widetilde v\rangle^2$ by \eqref{eq:Phi-derivatives} and \eqref{eq:bounded-m-derivative}. Differentiating $\rho_v=d\langle\widetilde v\rangle^{-2}$ now gives
	\begin{equation}\label{eq:elliptic-first-weight-bounds}
		|\nabla_x\widetilde v|\leq C_Ld^{-1}\langle\widetilde v\rangle,\qquad |\nabla_x\log\mathfrak w|\leq C_L\rho_v^{-1},\qquad |\nabla_x\rho_v|\leq C_L .
	\end{equation}
	Choose $c_L>0$ with $10C_Lc_L<1/2$. On $B_{10c_L\rho_v(x_0)}(x_0)$, \eqref{eq:elliptic-first-weight-bounds} gives
	\begin{equation}\label{eq:elliptic-ball-comparison}
		\begin{aligned}
			\rho_v(x)&\asymp\rho_v(x_0),\qquad d(x)\asymp d(x_0),\qquad \langle\widetilde v(x,v)\rangle\asymp\langle\widetilde v(x_0,v)\rangle,\\
			\left|\log\frac{\mathfrak w(x,v)}{\mathfrak w(x_0,v)}\right|&\leq C_L\frac{|x-x_0|}{\rho_v(x_0)} .
		\end{aligned}
	\end{equation}
	Set $r(x)=c_L\rho_v(x)$ and $\lambda=1/4$. Take a maximal disjoint family $B_{\lambda r_i}(x_i)$ with $r_i=r(x_i)$. For each $x$, some $B_{\lambda r_i}(x_i)$ intersects $B_{\lambda r(x)}(x)$. By \eqref{eq:elliptic-first-weight-bounds},
	\begin{align*}
		|x-x_i|\leq\lambda(r(x)+r_i),\quad
		|r(x)-r_i|\leq C_Lc_L|x-x_i|
		&\quad\Longrightarrow\quad r(x)\leq2r_i,\quad |x-x_i|\leq3r_i/4,\\
		x\in2B_i\cap2B_j
		&\quad\Longrightarrow\quad |r_i-r_j|\leq2C_Lc_L(r_i+r_j).
	\end{align*}
	Thus $B_i=B_{r_i}(x_i)$ cover $\R^3$. For fixed $x\in2B_i$, each $j$ with $x\in2B_j$ has $r_i/2\leq r_j\leq2r_i$ and $B_{\lambda r_j}(x_j)\subset B_{6r_i}(x)$. By disjointness,
	\begin{equation}\label{eq:elliptic-recovery-overlap}
		\#\{j:x\in2B_j\}\leq \frac{|B_{6r_i}|}{|B_{\lambda r_i/2}|}=(12/\lambda)^3 .
	\end{equation}
	Use a cutoff $\chi_i=1$ on $B_i$, supported in $2B_i$, with $|\nabla_x^j\chi_i|\leq C_jr_i^{-j}$. Plancherel and the product rule give
	\begin{equation}\label{eq:elliptic-recovery-ball}
	 \begin{aligned}
	  \sum_{|\alpha|=N}\|\partial_x^\alpha H\|_{L^2(B_i)}^2
	  &\leq C\|\Delta_x^{N/2}(\chi_iH)\|_{L^2_x}^2\\
	  &\leq C\|\Delta_x^{N/2}H\|_{L^2(2B_i)}^2
	   +C\sum_{j<N}\sum_{|\alpha|=j}r_i^{-2(N-j)}\|\partial_x^\alpha H\|_{L^2(2B_i)}^2.
	 \end{aligned}
	\end{equation}
	Multiply \eqref{eq:elliptic-recovery-ball} by $\mathfrak w(x_i,v)d(x_i)^{2N}$ and sum. On $2B_i$, \eqref{eq:elliptic-ball-comparison} and $r_i=c_Ld(x_i)\langle\widetilde v(x_i,v)\rangle^{-2}$ give
	\[\mathfrak w(x_i,v)d(x_i)^{2N}r_i^{-2(N-j)}
	\asymp_L \mathfrak w(x,v)d(x)^{2j}\langle\widetilde v(x,v)\rangle^{4(N-j)}.\]
	With \eqref{eq:elliptic-recovery-overlap}, the summed bound is
	\begin{align*}
		\sum_{|\alpha|=N}\int\mathfrak w\,|d^N\partial_x^\alpha H|^2\dd x
		&\leq C_L\int\mathfrak w\,|d^N\Delta_x^{N/2}H|^2\dd x\\
		&\quad+C_L\sum_{j<N}\sum_{|\alpha|=j}\int\mathfrak w\,\langle\widetilde v\rangle^{4(N-j)}|d^j\partial_x^\alpha H|^2\dd x .
	\end{align*}
	Integrate in $v$ to get \eqref{eq:weighted-elliptic-recovery}. For $N=0$ the bound is immediate.
\end{proof}

\begin{lemma}\label{lem:single-spatial-allocation}
	For $\mathcal B$ sufficiently small and $m_0$ sufficiently large,
	\begin{equation}\label{eq:new-norm-comparison}
		I_b+I_b^{\rm low}\asymp A_b^{\rm E}+I_b^{\rm mic}+I_b^{\rm low},\qquad \mathcal E+\mathcal Z_K^2\leq C\mathcal B .
	\end{equation}
	The $L^2_{x,\widehat v}$ bounds for $d^{|\alpha|}\partial_x^\alpha H$, $|\alpha|\leq K$, are stated in \eqref{eq:physical-spatial-recovery}. Let $\zeta_b=(T-t)^{b/r}d^{-b}$, so that $\zeta_0=1$. For $|\alpha|=i\leq K-3$ and $0\leq q\leq2(K-i-2)$,
	\begin{align}\label{eq:single-weight-pointwise}
		\|\langle \widetilde v\rangle^qH_\alpha/W\|_{L^2_{\widehat v}} &\leq C J_x^{-1}\zeta_b^{-1/2}(I_b+I_b^{\rm low})^{1/2},\notag\\
		\theta^{(3+\gamma)/2} \|\widehat\nu_M^{1/2}\langle \widetilde v\rangle^qH_\alpha/W\|_{L^2_{\widehat v}} &\leq C J_x^{-1}\zeta_b^{-1/2}\mathcal D_b^{1/2}, \qquad |\alpha|=i.
	\end{align}
	For $J_H:=\|H_0/W\|_{L^2_{\widehat v}}+\sum_{|\alpha|=1}\|H_\alpha/W\|_{L^2_{\widehat v}}$, we have
	\begin{equation}\label{eq:single-weight-low-flux}
		\sup_x\theta^p J_H\leq C(T-t)^{(-1+\eta)/2}\mathcal D_b^{1/2}, \qquad b=0,a.
	\end{equation}
\end{lemma}
\begin{proof}
	By \eqref{eq:explicit-Maxwellian-derivative} and the compact parameter range, $P_i/M\in\operatorname{span}\{1,\widetilde v_1,\widetilde v_2,\widetilde v_3,|\widetilde v|^2\}$ and
	\[\left|P_i/M\right|\leq C\langle\widetilde v\rangle^2\|P_i/\sqrt M\|_{L^2_{\widehat v}},\qquad i=0,K.\]
	By \eqref{eq:new-Gaussian-match}, $\int H_iP_i/W^2\dd\widehat v=\int H_iP_i(W^{-2}-M^{-1})\dd\widehat v$. Both $W^{-2}-M^{-1}$ and $M/W-W$ vanish for $|\widetilde v|^2\leq m$. By \eqref{eq:new-tail-small} and $W^2\geq M$, we have
	\begin{align*}
		\left|\int H_iP_i(W^{-2}-M^{-1})\dd\widehat v\right|
		&\leq\|H_i/W\|_{L^2_{\widehat v}}\|(P_i/M)(M/W-W)\|_{L^2_{\widehat v}}\\
		&\leq C\mathrm e^{-cm}\|H_i/W\|_{L^2_{\widehat v}}\|P_i/\sqrt M\|_{L^2_{\widehat v}},\\
		0\leq\int|P_i|^2(M^{-1}-W^{-2})\dd\widehat v
		&\leq\int_{|\widetilde v|^2\geq m}|P_i|^2/M\dd\widehat v
			\leq C\mathrm e^{-cm}\|P_i/\sqrt M\|_{L^2_{\widehat v}}^2.
	\end{align*}
	Expand $G_i=H_i+P_i$ in \eqref{eq:new-kinetic-energies} and use \eqref{eq:new-Gaussian-match} to obtain
	\begin{equation}\label{eq:matched-energy-tail-comparison}
		|I_b-I_b^{\rm mic}-A_b^{\rm E}|\leq C\mathrm e^{-cm_0}(I_b^{\rm mic}+A_b^{\rm E}).
	\end{equation}
	Choose $m_0$ with $C\mathrm e^{-cm_0}\leq1/2$ in \eqref{eq:matched-energy-tail-comparison}. By \eqref{eq:macro-equivalence-new}, \eqref{eq:observation-zero-order-bound} and $\Omega_a^*\leq\Omega_0^*$,
	\begin{align*}
		\mathcal Z_K^2+\mathcal R_*(\mathsf V,\mathsf V)
		&\leq CA_0^{\rm E}\leq C\mathcal B,\\
		I_a+I_a^{\rm low}
		&\leq \max\{1,\kappa_a^{\rm E}/\kappa_0^{\rm E}\}\mathcal B .
	\end{align*}
	Thus $\mathcal E\leq C\mathcal B$ by \eqref{eq:spatial-master-energies}, proving \eqref{eq:new-norm-comparison}.

	For $\mathbf r=(n-1,V,\Theta-1)$, \eqref{eq:normalized-reconstruction} and the pointwise bound in \eqref{eq:macro-equivalence-new} give
	\[\sup_x\sum_{j=0}^2d^j\bigl(|\nabla_x^j\mathbf z|+|\nabla_x^j\mathbf r|\bigr)\leq C\mathcal B^{1/2}.\]
	For small $\mathcal B$ this gives $(n,V,\Theta)\in\operatorname{int}\mathcal U$ and \eqref{eq:coefficient-derivative-assumptions}. Use
	\[\theta^{-3}\dd v=\dd\widehat v,\qquad d^K\Delta_x^kH=(d/d_\phi)^KH_K,\qquad 4(K-j)=2\widehat s_j\quad(j<K)\]
	in \eqref{eq:weighted-elliptic-recovery} with $N=K$, $d_\phi\asymp d$ and $\sigma=0,1$ to get
	\begin{equation}\label{eq:physical-spatial-recovery}
		\begin{aligned}
			\sum_{|\alpha|\leq K}\iint\Omega_b^*\left|\frac{d^{|\alpha|}\partial_x^\alpha H}{W}\right|^2\dd\widehat v\dd x&\leq C(I_b^{\rm mic}+I_b^{\rm low}),\\
			\sum_{|\alpha|\leq K}\iint\Omega_b^*\theta^{3+\gamma}\widehat\nu_M\left|\frac{d^{|\alpha|}\partial_x^\alpha H}{W}\right|^2\dd\widehat v\dd x&\leq C\mathcal D_b .
		\end{aligned}
	\end{equation}

	On $B_{2cd_0}(x_0)$, where $d_0=d(t,x_0)$, we have $(d,J_x,\zeta_b)\asymp(d_0,J_x(x_0),\zeta_b(x_0))$. For $A=1,(\theta^{3+\gamma}\widehat\nu_M)^{1/2}$, \eqref{eq:normalized-collision-frequency}, \eqref{eq:normalized-frequency-derivatives} and \eqref{eq:coefficient-derivative-assumptions} give $d^j|\nabla_x^j A|\leq CA$, $j\leq2$. By \eqref{eq:fitted-Phi-chain} and \eqref{eq:Phi-derivatives},
	\begin{equation}\label{eq:rw56-pointwise-multiplier}
		\begin{aligned}
			|\nabla_x\log W|&\leq Cd^{-1}\langle\widetilde v\rangle^2,\qquad
			|\nabla_x^2\log W|\leq Cd^{-2}\langle\widetilde v\rangle^2,\\
			d^j\left|\nabla_x^j\left(\frac{A\theta^{-3/2}\langle\widetilde v\rangle^q}{W}\right)\right|
			&\leq C\frac{A\theta^{-3/2}\langle\widetilde v\rangle^{q+2j}}{W}\qquad(j=0,1,2).
		\end{aligned}
	\end{equation}
	The second line uses $\nabla_x^2W^{-1}=W^{-1}(\nabla_x\log W\otimes\nabla_x\log W-\nabla_x^2\log W)$.

	Take $i=|\alpha|\leq K-3$ and $0\leq q\leq2(K-i-2)$. When $l$ out of $j$ derivatives fall on $\partial_x^\alpha H$, we have
	\begin{equation}\label{eq:pointwise-order-allocation}
		i+l\leq K-1,\qquad q+2(j-l)\leq\widehat s_{i+l}=2(K-i-l)\qquad(0\leq l\leq j\leq2).
	\end{equation}
	Use $H^2_X(L^2_v)\hookrightarrow L^\infty_X(L^2_v)$ for $X=(x-x_0)/d_0$. Then $\Omega_b^*\asymp J_x^2\zeta_b d^{-3}$, $\theta^{-3}\dd v=\dd\widehat v$ and \eqref{eq:rw56-pointwise-multiplier} give
	\begin{align*}
		&J_x(x_0)^2\zeta_b(x_0)\left\|A\langle\widetilde v\rangle^q H_\alpha/W\right\|_{L^2_{\widehat v}}^2(x_0)\\
		&\quad\leq C\sum_{j=0}^2\sum_{l=0}^j\iint_{B_{2cd_0}(x_0)\times\R^3}
			\Omega_b^*A^2\left|\frac{\langle\widetilde v\rangle^{q+2(j-l)}d^{i+l}\nabla_x^{i+l}H}{W}\right|^2\dd\widehat v\dd x .
	\end{align*}
	By \eqref{eq:pointwise-order-allocation}, the right-hand sides for the two choices of $A$ are bounded by $C(I_b+I_b^{\rm low})$ and $C\mathcal D_b$. This proves \eqref{eq:single-weight-pointwise}.

	Take $q=0$ and $|\alpha|\leq1$ in the second bound of \eqref{eq:single-weight-pointwise}. Since $\widehat\nu_M\geq c$, this gives $\theta^pJ_H\leq C J_x^{-1}\zeta_b^{-1/2}\theta^{p-(3+\gamma)/2}\mathcal D_b^{1/2}$. By $(r-1)(2p-3-\gamma)=r(1-\eta)$ and $\theta\asymp d^{1-r}$ in the profile region,
	\begin{equation}\label{eq:pointwise-flux-time-factor}
		J_x^{-1}\zeta_b^{-1/2}\theta^{p-(3+\gamma)/2}\leq C
		\begin{cases}
			(T-t)^{-b/(2r)}d^{(b-r(1-\eta))/2},&|x|\leq R_{\rm e},\\
			(T-t)^{-b/(2r)},&|x|\geq R_{\rm e}.
		\end{cases}
	\end{equation}
	For $|x|\geq R_{\rm e}$ use $\theta\asymp1$ and $J_x^{-1}d^{b/2}\leq C$ in \eqref{eq:pointwise-flux-time-factor}. Use $b\leq r/2<r(1-\eta)$ and $(T-t)^{1/r}\leq d$ with $T\leq1$ to get
	\[(T-t)^{-b/(2r)}d^{(b-r(1-\eta))/2}\leq(T-t)^{(\eta-1)/2},\qquad (T-t)^{-b/(2r)}\leq(T-t)^{(\eta-1)/2}.\]
	This gives \eqref{eq:single-weight-low-flux}.
\end{proof}

\section{The differentiated collision estimate}\label{sec:collision-hierarchy}

\begin{proposition}\label{prop:new-component-collision}
	Let $f=M+H$ be a smooth solution of \eqref{Eq.Boltzmann}, with $M=\mathcal M(\rho_f,u_f,T_f)$ from \eqref{eq:physical-parameters}. Assume $(n,V,\Theta)\in\mathcal U$ and $\mathcal B$ sufficiently small. There are a fixed tail order and coefficients $b_j^{\rm sp},\kappa_j>0$, $0\leq j<K$, for which the dissipation-side collision pairing $\mathcal C_b$ in $(I_b+I_b^{\rm low})'$ satisfies
	\begin{equation}\label{eq:new-component-collision}
		\mathcal C_b\geq\frac{c-C\sqrt{\mathcal B}}{\eps}\mathcal D_b
			-C\eps(T-t)^{-1+\eta}(1+\mathcal B)(I_b+I_b^{\rm low})\qquad(b=0,a).
	\end{equation}
	The coefficients and constants are uniform in $m,t,\eps$.
\end{proposition}
\begin{proof}
	Take $\partial_x$ at fixed physical $v$. For $|\alpha|=j<K$ set $H_\alpha=d^j\partial_x^\alpha H$ and $\widehat s_j=2(K-j)$, and use $H_0=H$, $H_K=d_\phi^K\Delta_x^kH$. By \eqref{eq:physical-collision-frequency}, the two velocity norms are
	\begin{align*}
		\mathfrak n_s(R)&:=\|\langle\widetilde v\rangle^sR/W\|_{L^2_{\widehat v}},\\
		\mathfrak d_s(R)&:=\|(\nu_M^{\rm phys})^{1/2}\langle\widetilde v\rangle^sR/W\|_{L^2_{\widehat v}}
			=\theta^{(3+\gamma)/2}\|\widehat\nu_M^{1/2}\langle\widetilde v\rangle^sR/W\|_{L^2_{\widehat v}} .
	\end{align*}
	The moment constraints are differentiated before changing $v=u+\theta\widehat v$. Proposition~\ref{thm:new-collision} and Lemma~\ref{lem:polynomial-collision} give
	\begin{align}\label{eq:rw78-collision-coercivity}
		&\qquad\qquad\qquad\int(1,v,|v|^2)R\dd v=0,\quad R\in\{H_0,H_K,H_\alpha\},\notag\\
		&\qquad\qquad-\operatorname{Re}\int\frac{\{Q(M,R)+Q(R,M)\}\overline R}{W^2}\dd\widehat v\geq c\mathfrak d_0(R)^2,\\
		&-\operatorname{Re}\int\frac{\{Q(M,H_\alpha)+Q(H_\alpha,M)\}\langle\widetilde v\rangle^{2\widehat s_j}\overline{H_\alpha}}{W^2}\dd\widehat v
		\geq\tfrac12\mathfrak d_{\widehat s_j}(H_\alpha)^2-C_{\widehat s_j}\mathfrak d_0(H_\alpha)^2\notag
	\end{align}
	for $|\alpha|=j<K$.
	
	By the Leibniz formula, we have, for $|\alpha|=j$,
	\begin{align*}
		d^j\partial_x^\alpha Q(M,H)&=Q(M,H_\alpha)+\sum_{0<\lambda\leq\alpha}\binom\alpha\lambda Q(d^{|\lambda|}\partial_x^\lambda M,H_{\alpha-\lambda}),\\
		d^j\partial_x^\alpha Q(H,M)&=Q(H_\alpha,M)+\sum_{0<\lambda\leq\alpha}\binom\alpha\lambda Q(H_{\alpha-\lambda},d^{|\lambda|}\partial_x^\lambda M).
	\end{align*}
	By the chain rule, $W\geq\sqrt M$ and the uniform Gaussian moments, we have
	\begin{align*}
		\|\widehat\nu_M^{1/2}\langle\widetilde v\rangle^s d^{|\lambda|}\partial_x^\lambda M/W\|_{L^2_{\widehat v}}
		\leq C_{\lambda,s}\left(1+\sum_{0<\lambda_0\leq\lambda}\sum_{h=1}^{|\lambda_0|}\sum_{\substack{\lambda_1+\cdots+\lambda_h=\lambda_0\\|\lambda_q|\geq1}}\prod_{q=1}^h|d^{|\lambda_q|}\partial_x^{\lambda_q}\mathbf z|\right).
	\end{align*}
	If all the unknown derivatives have order at most $K-2$, then \eqref{eq:macro-equivalence-new} bounds the right-hand side uniformly. Hence, by \eqref{eq:polynomial-bilinear}, we have, for $s=0,\widehat s_j$ and $i:=j-|\lambda|<j$,
	\begin{align}\label{eq:rw78-collision-commutator}
		&\left|\int\frac{\{Q(d^{|\lambda|}\partial_x^\lambda M,H_{\alpha-\lambda})+Q(H_{\alpha-\lambda},d^{|\lambda|}\partial_x^\lambda M)\}\langle\widetilde v\rangle^{2s}\overline{H_\alpha}}{W^2}\dd\widehat v\right|\notag\\
		&\qquad\leq C\mathfrak d_s(H_{\alpha-\lambda})\mathfrak d_s(H_\alpha)
		\leq\delta\mathfrak d_s(H_\alpha)^2+C_\delta\mathfrak d_s(H_{\alpha-\lambda})^2.
	\end{align}
	In \eqref{eq:rw78-collision-commutator}, $\widehat s_{j-|\lambda|}=\widehat s_j+2|\lambda|$, so the last term is bounded by the dissipation at order $j-|\lambda|$. At order $K=2k$ the identity
	\[d_\phi^K\Delta_x^kQ(M,H)=(d_\phi/d)^K\sum_{|\beta|=k}\frac{k!}{\beta!}\sum_{\lambda\leq2\beta}\binom{2\beta}{\lambda}Q(d^{|\lambda|}\partial_x^\lambda M,H_{2\beta-\lambda})\]
	gives $Q(M,H_K)$ at $\lambda=0$. For $A=d^{|\lambda|}\partial_x^\lambda M$ with unknown derivatives of order at most $K-2$, use \eqref{eq:polynomial-bilinear} and $d_\phi\asymp d$ to obtain
	\begin{align*}
		&\left|\int(d_\phi/d)^K\frac{\{Q(A,H_{2\beta-\lambda})+Q(H_{2\beta-\lambda},A)\}\overline{H_K}}{W^2}\dd\widehat v\right|\\
		&\qquad\leq\delta\mathfrak d_0(H_K)^2+C_\delta\mathfrak d_0(H_{2\beta-\lambda})^2,\qquad 0<\lambda\leq2\beta .
	\end{align*}
	Choose $\delta$ so the tested-order terms use at most one quarter of \eqref{eq:rw78-collision-coercivity}. Fix the top coefficient to one. If $L_j,U_j$ are the losses of $\mathfrak d_{\widehat s_j}(H_\alpha)^2,\mathfrak d_0(H_\alpha)^2$ from orders above $j$, take
	\begin{equation}\label{eq:spatial-collision-hierarchy}
		\kappa_j\geq C(1+L_j),\qquad b_j^{\rm sp}\geq C(1+U_j+C_{\widehat s_j}\kappa_j),\qquad j=K-1,\ldots,0 .
	\end{equation}
	The unweighted loss at order $j$ is $C_{\widehat s_j}\kappa_j$ by \eqref{eq:rw78-collision-coercivity}. For $C$ large, these terms leave $c\mathcal D_b/\eps$, $b=0,a$.

	For $A_1,A_2,A_3\geq0$, use $\Omega_b^*A_1A_2A_3=(J_x\zeta_b^{1/2}A_1)(J_x^{-1}(\Omega_0^*)^{1/2}A_2)((\Omega_b^*)^{1/2}A_3)$. The Cauchy-Schwarz inequality gives
	\begin{equation}\label{eq:collision-spatial-allocation}
		\int\Omega_b^*A_1A_2A_3\dd x\leq\|J_x\zeta_b^{1/2}A_1\|_{L^\infty_x}
			\|J_x^{-1}A_2\|_{L^2(\Omega_0^*\dd x)}\|A_3\|_{L^2(\Omega_b^*\dd x)} .
	\end{equation}
	For $j_{\mathbf z}\geq K-1$, let $A_{\rm hi}$ be the scaled derivative of $\mathbf z$ times its bounded coefficients. The other input $R=H_\lambda$ has $|\lambda|\leq K-j_{\mathbf z}\leq1$. In \eqref{eq:collision-spatial-allocation} take
	\[(A_1,A_2,A_3)=(\mathfrak d_s(R),|A_{\rm hi}|,G),\qquad G=\mathfrak d_0(H_K)\ \text{or}\ \mathfrak d_s(H_\alpha).\]
	Here $s=0$ for $H_K$, and $s=0,\widehat s_{K-1}$ for $|\alpha|=K-1$. By \eqref{eq:single-weight-pointwise}, \eqref{eq:macro-equivalence-new} and $J_x^{-1}\leq C$,
	\begin{align*}
		\int\Omega_b^*|A_{\rm hi}|\mathfrak d_s(R)G\dd x
		&\leq\|J_x\zeta_b^{1/2}\mathfrak d_s(R)\|_\infty\|J_x^{-1}A_{\rm hi}\|_{L^2(\Omega_0^*\dd x)}
		\|G\|_{L^2(\Omega_b^*\dd x)}\\
		&\leq C\mathcal D_b^{1/2}\sqrt{\mathcal B}\mathcal D_b^{1/2}.
	\end{align*}
	
	For $Q(H,H)$, we use
	\[d^j\partial_x^\alpha Q(H,H)=\sum_{\lambda\leq\alpha}\binom\alpha\lambda Q(H_\lambda,H_{\alpha-\lambda}).\]
	Assume that $|\lambda|\leq|\alpha-\lambda|$. Then $|\lambda|\leq\lfloor K/2\rfloor\leq K-3$. By \eqref{eq:new-bilinear}, \eqref{eq:single-weight-pointwise} and \eqref{eq:collision-spatial-allocation}, we have
	\begin{align}\label{eq:rw78-quadratic-collision}
		\int\Omega_b^*\mathfrak n_0(H_\lambda)\mathfrak d_0(H_{\alpha-\lambda})G\dd x
		&\leq\|\mathfrak n_0(H_\lambda)\|_{L^\infty_x}\|\mathfrak d_0(H_{\alpha-\lambda})\|_{L^2(\Omega_b^*\dd x)}\|G\|_{L^2(\Omega_b^*\dd x)}\notag\\
		&\leq C\sqrt{\mathcal B}\mathcal D_b,\notag\\
		\int\Omega_b^*\mathfrak d_0(H_\lambda)\mathfrak n_0(H_{\alpha-\lambda})G\dd x
		&\leq\|J_x\zeta_b^{1/2}\mathfrak d_0(H_\lambda)\|_{L^\infty_x}\|J_x^{-1}\mathfrak n_0(H_{\alpha-\lambda})\|_{L^2(\Omega_0^*\dd x)}\|G\|_{L^2(\Omega_b^*\dd x)}\notag\\
		&\leq C\sqrt{\mathcal B}\mathcal D_b.
	\end{align}
	In \eqref{eq:rw78-quadratic-collision} the lower spatial order is in $L_x^\infty$, with the frequency on either input as in \eqref{eq:polynomial-bilinear}. At $j<K$, an input of order $i$ in $L^2_x$ can carry $\widehat s_j\leq\widehat s_i$. For a lower-order input in $L_x^\infty$, if the other input has at least two derivatives, then
	\[
	 j-i\geq2\quad\Longrightarrow\quad \widehat s_j=2(K-j)\leq 2(K-i-2).
	\]
	This is the exponent allowed by \eqref{eq:single-weight-pointwise}. For pairs of orders $(0,0),(0,1),(1,1)$, put the input with zero extra velocity power in $L_x^\infty$. Both frequency placements in \eqref{eq:polynomial-bilinear}, combined with \eqref{eq:collision-spatial-allocation}, give the bound $C\sqrt{\mathcal B}\mathcal D_b$ for the product pairing. For order $K$ the $L_x^2$ input is controlled by \eqref{eq:physical-spatial-recovery}.

	Let $\mathcal Q_0=Q(f,f)$ and $\mathcal Q_K=d_\phi^K\Delta_x^kQ(f,f)$. Differentiate the collision invariants at fixed $v$ to get
	\[
	 \int(1,v,|v|^2)\mathcal Q_i\dd v=0,\qquad i\in\{0,K\}.
	\]
	As $P_i/M$ is a linear combination of $(1,v,|v|^2)$ and $\dd\widehat v=\theta^{-3}\dd v$, we have $\int \mathcal Q_iP_i/M\dd\widehat v=0$. Since $W^2=M$ for $|\widetilde v|^2\leq m_\eps$, the pairing is
	\[
	 \int\frac{\mathcal Q_iP_i}{W^2}\dd\widehat v
	 =\int_{|\widetilde v|^2> m_\eps}\mathcal Q_iP_i(W^{-2}-M^{-1})\dd\widehat v .
	\]
	Apply \eqref{eq:new-tail-small} to this pairing to get
	\begin{align*}
		\left|\int\mathcal Q_iP_i/W^2\dd\widehat v\right|
		&=\left|\int\frac{\mathcal Q_i}{W}\frac{P_i}{M}(M/W-W)\dd\widehat v\right|\\
		&\leq C\mathrm e^{-cm}\|\widehat\nu_M^{-1/2}\mathcal Q_i/W\|_{L^2_{\widehat v}}\left(\int|P_i|^2/M\dd\widehat v\right)^{1/2}.
	\end{align*}
	The estimates for the bounded, high-parameter and quadratic terms above, along with \eqref{eq:polynomial-bilinear}, yield
	\[\sum_{i=0,K}\iint\Omega_b^*\theta^{-3-\gamma}|\widehat\nu_M^{-1/2}\mathcal Q_i/W|^2\dd\widehat v\dd x\leq C\mathcal D_b+C\mathcal B\mathcal D_b.\]
	Here the fixed coefficient $\kappa_b^{\rm E}$ at order zero is included in $C$. Therefore, by Young's inequality, we have
	\begin{align}\label{eq:rw78-collision-macro-tail}
		\frac1\eps\sum_{i=0,K}\left|\iint\Omega_b^*\mathcal Q_iP_i/W^2\dd\widehat v\dd x\right|
		\leq\frac{C\delta}\eps(1+\mathcal B)\mathcal D_b+\frac{C_\delta}\eps\sum_{i=0,K}\iint\Omega_b^*\theta^{3+\gamma}\mathrm e^{-2cm}|P_i|^2/M\dd\widehat v\dd x.
	\end{align}
	Take $\delta$ small to absorb the term $C\delta(1+\mathcal B)\mathcal D_b/\eps$ in \eqref{eq:rw78-collision-macro-tail}. Choose $2c\kappa_{\rm w}\geq2$ in \eqref{eq:adaptive-m}. Since $p-\delta_{\rm gap}=2p-3-\gamma>0$ and $\alpha_{\rm ss}\delta_{\rm gap}=\eta$,
	\begin{equation}\label{eq:collision-tail-time}
		\eps^{-1}\theta^{3+\gamma}\mathrm e^{-2cm}
		\leq C\eps\theta^{3+\gamma-2\delta_{\rm gap}}
		=C\eps\theta^{p-\delta_{\rm gap}}\leq C\eps(T-t)^{-1+\eta}.
	\end{equation}
	By \eqref{eq:collision-tail-time} and \eqref{eq:new-Gaussian-match}, the unabsorbed term in \eqref{eq:rw78-collision-macro-tail} satisfies
	\[
	 \frac{C_\delta}{\eps}\sum_{i=0,K}\iint
	    \Omega_b^*\theta^{3+\gamma}\mathrm e^{-2cm}\frac{|P_i|^2}{M}\dd\widehat v\dd x
	 \leq C\eps(T-t)^{-1+\eta}A_b^{\rm E}
	 \leq C\eps(T-t)^{-1+\eta}(1+\mathcal B)(I_b+I_b^{\rm low})
	\]
	using \eqref{eq:new-norm-comparison}. Combining \eqref{eq:rw78-collision-coercivity}, \eqref{eq:spatial-collision-hierarchy} and \eqref{eq:rw78-quadratic-collision} with the Gaussian-tail estimate proves \eqref{eq:new-component-collision}.
\end{proof}

\section{The nonlinear transport estimates}\label{sec:expanded-critical}

\subsection{Microscopic increments and localized products}
For the microscopic terms in the macroscopic equations, take the moments about the background velocity $u$.
\begin{equation}\label{eq:expanded-moments}
 \pi_{ij}:=\theta^{-5}\int(v_i-u_i)(v_j-u_j)H\dd v,\qquad
 (q_H)_i:=\frac1{2\theta^6}\int|v-u|^2(v_i-u_i)H\dd v,\quad 1\leq i,j\leq3.
\end{equation}
Using \eqref{eq:zero-moments} in the expansion about $u$,
\[
 \pi^{\mathsf T}=\pi,\qquad \tr\pi=0,\qquad \Pi=\theta^5\pi,\qquad
 \tfrac12\int |v|^2vH\dd v=\theta^6q_H+\Pi u .
\]
Differentiate the last two identities and use $u_f=u+\theta V$ to obtain
\begin{align}\label{eq:expanded-flux-vector}
	G_{{\rm mom},x}&:=\theta\Div\pi+5\pi\nabla\theta,\notag\\
	C_H&:=\theta\Div q_H+6q_H\cdot\nabla\theta +\pi:\nabla u-V\cdot G_{{\rm mom},x},\\
	\Div\Pi&=\theta^4G_{{\rm mom},x},\qquad \mathcal C_f=\tfrac12\Div m_3-u_f\cdot\Div\Pi=\theta^5C_H. \notag
\end{align}
Here $(\Div_x\pi)_i=\sum_j\partial_{x_j}\pi_{ij}$ and $\pi:\nabla_xu=\sum_{i,j}\pi_{ij}\partial_{x_j}u_i$. In the last line of \eqref{eq:expanded-flux-vector}, use \eqref{eq:roadmap-raw-heat-source} and
\[
 \tfrac12\Div m_3-u_f\cdot\Div\Pi
 =\Div(\theta^6q_H)+\Pi:\nabla u-\theta V\cdot\Div\Pi .
\]
On $|(T-t)^{1/r}y|<R_{\rm e}$, \eqref{eq:profile} gives $(T-t)G_{{\rm mom},x}^\sharp=G_{\rm mom}$ and $(T-t)C_H^\sharp=G_{\rm heat}$ from \eqref{eq:dimensionless-forces}.

\begin{lemma}\label{lem:expanded-increment}
	For $\mathbf z=(V,B,S)$ and $\mathscr D_{\rm mic}$ as in \eqref{eq:microscopic-moment-time-difference},
	\begin{equation}\label{eq:expanded-kinetic-increment}
		\mathscr D_{\rm mic} V=-n^{-1}G_{{\rm mom},x},\qquad
		\mathscr D_{\rm mic} B=-\frac{2C_H}{15(1+B)^4},\qquad
		\mathscr D_{\rm mic} S=-\frac{2C_H}{3(1+B)^5}.
	\end{equation}
	Linearity in $(\pi,q_H,\nabla_x\pi,\nabla_xq_H)$ follows from \eqref{eq:expanded-flux-vector}.
\end{lemma}
\begin{proof}
	By \eqref{eq:microscopic-moment-time-difference} and \eqref{eq:expanded-flux-vector},
	\[\mathscr D_{\rm mic}\rho_f=0,\qquad \mathscr D_{\rm mic}u_f=-\frac{\theta^4}{\rho_f}G_{{\rm mom},x},\qquad
	\mathscr D_{\rm mic}(\rho_fT_f)=-\tfrac23\theta^5 C_H.\]
	Apply $\mathscr D_{\rm mic}(u,\theta)=0$ to $V=(u_f-u)/\theta$, $(1+B)^5=\rho_fT_f/\theta^5$ and $S=\log(\rho_fT_f)-\tfrac53\log\rho_f$ to get \eqref{eq:expanded-kinetic-increment}.
\end{proof}

	At fixed $t$, cover $|x|\leq(T-t)^{1/r}$ and the dyadic annuli outside by $B_{cd_c}(x_c)$, $d_c=d(t,x_c)$, with bounded overlap after enlargement. Put $x=x_c+d_cX$ and $(\theta_c,J_c,\zeta_{b,c})=(\theta,J_x,\zeta_b)(t,x_c)$. For $c$ small,
	\begin{equation}\label{eq:local-product-measure}
		(d,\theta,J_x,\zeta_b)\asymp(d_c,\theta_c,J_c,\zeta_{b,c}),\qquad
		\Omega_b^*\dd x\asymp J_c^2\zeta_{b,c}d_c^{-3}\dd x=J_c^2\zeta_{b,c}\dd X .
	\end{equation}

\begin{lemma}\label{lem:expanded-products}
	For $K\geq3$, $\ell\geq2$, $1\leq j_i\leq K$ and $\sum_i j_i\leq K+1$, the smaller ball of \eqref{eq:local-product-measure} satisfies
	\begin{equation}\label{eq:localized-tame-product}
		\left\|\prod_{i=1}^{\ell}\nabla^{j_i}f_i\right\|_{L^2}\leq C\sum_{i=1}^{\ell}\|f_i\|_{H^K}\prod_{h\ne i}\|f_h\|_{W^{1,\infty}}.
	\end{equation}
	The norms on the right are on the enlarged ball. Any additional order-zero factor contributes its $L^\infty$ norm.
\end{lemma}
\begin{proof}
	Cut off $f_i$ in the enlarged ball. For $\alpha_i=(j_i-1)/(K-1)$, $\sum_i\alpha_i\leq(K+1-\ell)/(K-1)\leq1$. H\"older and Gagliardo-Nirenberg give
	\begin{align*}
		\left\|\prod_i\nabla^{j_i}f_i\right\|_{L^2}
		&\leq C\prod_i\|\nabla^{j_i}f_i\|_{L^{2/\alpha_i}}
		\leq C\prod_i\|f_i\|_{W^{1,\infty}}^{1-\alpha_i}\|f_i\|_{H^K}^{\alpha_i}\\
		&\leq C\sum_i\|f_i\|_{H^K}\prod_{h\ne i}\|f_h\|_{W^{1,\infty}},
	\end{align*}
	where $2/\alpha_i=\infty$ for $j_i=1$. The last inequality follows from the weighted arithmetic-geometric mean and $H^K\hookrightarrow W^{1,\infty}$.
\end{proof}

	At fixed $\mathbf z,v$ write $\partial_x^{\alpha,\mathrm{exp}}$. For $|\alpha|=K$ the chain rule gives
\begin{align}\label{eq:expanded-composition}
	\partial_x^\alpha M(t,x,\mathbf z(t,x);v)
	&=\sum_{q=0}^{K}\sum_{\substack{\alpha_0+\cdots+\alpha_q=\alpha\\|\alpha_i|\geq1\ (i\geq1)}}
	\frac{\alpha!}{q!\alpha_0!\cdots\alpha_q!}
	\partial_x^{\alpha_0,\mathrm{exp}}D_{\mathbf z}^qM
	[\partial_x^{\alpha_1}\mathbf z,\ldots,\partial_x^{\alpha_q}\mathbf z].
\end{align}
In \eqref{eq:roadmap-composition-remainders}, the terms $(q,\alpha_0)=(1,0)$ sum to $D_{\mathbf z}M[\Delta_x^k\mathbf z]$. The rest of the conjugated input in \eqref{eq:roadmap-conjugated-Laplacian} comes from the following terms for $|\alpha|=K$.
\begin{equation}\label{eq:expanded-sigma-commutator}
 D_{\mathbf z}M\left[\sum_{0<\lambda\leq\alpha}\binom\alpha\lambda
           \Sigma_x^{-1}\partial_x^\lambda\Sigma_x\,\partial_x^{\alpha-\lambda}\mathbf z\right].
\end{equation}
After the $(q,\alpha_0)=(1,0)$ cancellation in \eqref{eq:expanded-composition}, every derivative of $\mathbf z$ there and in \eqref{eq:expanded-sigma-commutator} has order at most $K-1$. For $\mathcal B_0$, using \eqref{eq:roadmap-composition-remainders} and $\mathscr D_{\rm mic}\bar M=0$ gives
\begin{equation}\label{eq:expanded-B0-increment}
 \mathscr D_{\rm mic}\mathcal B_0=-D_{\mathbf z}^2M[\mathbf z,\mathscr D_{\rm mic}\mathbf z].
\end{equation}

By \eqref{eq:explicit-Maxwellian-derivative} and \eqref{eq:background-hierarchy}, $d^{|\lambda|}\partial_x^{\lambda,\mathrm{exp}}D_{\mathbf z}^qM/M$ is a polynomial in $\widetilde v$ of degree at most $2(|\lambda|+q)$ with bounded coefficients. Since $W\geq\sqrt M$,
\begin{equation}\label{eq:expanded-gaussian-coefficient}
 \left\|\langle\widetilde v\rangle^N\widehat\nu_M^{-1/2}
       \frac{d^{|\lambda|}\partial_x^{\lambda,\mathrm{exp}}D_{\mathbf z}^qM}{W}
 \right\|_{L^2_{\widehat v}}\leq C_{N,K},\qquad |\lambda|+q\leq K+2.
\end{equation}
For the terms with $q=0$ in $M-\bar M$, use
\[
 M(t,x,\mathbf z;v)-M(t,x,0;v)=\int_0^1D_{\mathbf z}M(t,x,s\mathbf z;v)[\mathbf z]\dd s .
\]
For $\mathbf z$ small and $0\leq s\leq1$, the bound
\[
 \frac{M(t,x,s\mathbf z;v)^2}{W^2}
 \leq\frac{M(t,x,s\mathbf z;v)^2}{M(t,x,\mathbf z;v)}
 \leq C\mathrm e^{-c|\widetilde v|^2}
\]
gives the norm in \eqref{eq:expanded-gaussian-coefficient} with $C_{N,K}|\mathbf z|$ for $d^{|\lambda|}\partial_x^{\lambda,\mathrm{exp}}(M-\bar M)$ and $|\lambda|+1\leq K+2$. On $|\widetilde v|^2\geq m$ the Gaussian norm has an additional factor $C\mathrm e^{-cm}$.

	On the enlarged $X$ ball put $Z_c=\|\mathbf z\|_{H^K_X}$, $N_c=\|(\pi,q_H)\|_{H^K_X}$ and $j_c=\|(\pi,q_H)\|_{W^{1,\infty}_X}$. By \eqref{eq:background-hierarchy}, \eqref{eq:untruncated-profile-bounds} and \eqref{eq:expanded-full-macro-euler},
	\begin{equation}\label{eq:scaled-euler-increment}
		\begin{aligned}
			d^l|\nabla_x^l\partial_tu|&\leq C_l\theta^2/d,\qquad d^l|\nabla_x^l\partial_t\log\theta|\leq C_l\theta/d,\\
			\|\partial_t^E\mathbf z\|_{H_X^{K-1}}&\leq C(\theta_c/d_c)Z_c,\qquad
			\|\mathbf z\|_{W_X^{1,\infty}}\leq C\sqrt{\mathcal B}.
		\end{aligned}
	\end{equation}
	Use $\nabla_x=d_c^{-1}\nabla_X$ in \eqref{eq:expanded-flux-vector} and apply \eqref{eq:localized-tame-product} to \eqref{eq:expanded-kinetic-increment}.
	\begin{equation}\label{eq:scaled-kinetic-increment}
		\frac{d_c}{\theta_c}\|\mathscr D_{\rm mic}\mathbf z\|_{H_X^{K-1}}\leq C(N_c+Z_cj_c),\qquad
		\frac{d_c}{\theta_c}\|\mathscr D_{\rm mic}\mathbf z\|_{L_X^\infty}\leq Cj_c .
	\end{equation}
	In the remaining terms of \eqref{eq:expanded-composition}, \eqref{eq:expanded-sigma-commutator}, use $\mathscr D_{\rm mic}(u,\theta,\Sigma_x)=0$. By \eqref{eq:scaled-euler-increment} and $|v|\leq C\theta_c\langle\widetilde v\rangle$, $D_E$ on the prescribed coefficients contributes $\theta_c/d_c$ in \eqref{eq:expanded-gaussian-coefficient}. The derivative orders on $(\mathbf z,\pi,q_H)$ satisfy
	\[\max_i j_i\leq K,\qquad \sum_i j_i\leq K+1,\qquad \sum_i j_i=K+1\ \Longrightarrow\ \#\{i:j_i\geq1\}\geq2.\]
	For instance, \eqref{eq:localized-tame-product} gives
	\[\|\nabla_X^{K-1}\mathbf z\,\nabla_X^2\pi\|_{L^2_X}\leq C\bigl(\|\mathbf z\|_{W_X^{1,\infty}}\|\pi\|_{H_X^K}+Z_c\|\pi\|_{W_X^{1,\infty}}\bigr)\leq C(N_c+Z_cj_c).\]
	With \eqref{eq:expanded-gaussian-coefficient} and the increments \eqref{eq:scaled-euler-increment}, \eqref{eq:scaled-kinetic-increment}, this yields
	\begin{align}\label{eq:expanded-B-local}
		\left\|\frac{\widehat\nu_M^{-1/2}d_\phi^KD_E\mathcal B_K}{W}\right\|_{L^2_{X,\widehat v}}&\leq C(\theta_c/d_c)Z_c,\notag\\
		\left\|\frac{\widehat\nu_M^{-1/2}d_\phi^K\mathscr D_{\rm mic}\mathcal B_K}{W}\right\|_{L^2_{X,\widehat v}}&\leq C(\theta_c/d_c)(N_c+Z_cj_c).
	\end{align}
	For $\mathcal B_0$, use
	\[\mathcal B_0=-\int_0^1sD_{\mathbf z}^2M(t,x,s\mathbf z;v)[\mathbf z,\mathbf z]\dd s,\qquad \mathscr D_{\rm mic}\mathcal B_0=-D_{\mathbf z}^2M[\mathbf z,\mathscr D_{\rm mic}\mathbf z]\]
	from Taylor's formula and \eqref{eq:expanded-B0-increment}. By \eqref{eq:scaled-euler-increment} and \eqref{eq:scaled-kinetic-increment},
	\begin{align*}
		\|\widehat\nu_M^{-1/2}D_E\mathcal B_0/W\|_{L^2_{X,\widehat v}}&\leq C(\theta_c/d_c)Z_c,\quad
		\|\widehat\nu_M^{-1/2}\mathscr D_{\rm mic}\mathcal B_0/W\|_{L^2_{X,\widehat v}}\leq C(\theta_c/d_c)Z_cj_c.
	\end{align*}

\subsection{Composition errors and macroscopic principal terms}

Set $H_\beta=d^{|\beta|}\partial_x^\beta H$. Apply $\partial_x^\alpha$ at fixed $v$ to \eqref{eq:expanded-moments} and then set $v=u+\theta\widehat v$. We get
\[
 d^{|\alpha|}\partial_x^\alpha(\pi,q_H)
  =\sum_{\beta\leq\alpha}\int P_{\alpha\beta}(t,x,\widehat v)H_\beta\dd\widehat v,
 \qquad \deg_{\widehat v}P_{\alpha\beta}\leq3 ,
\]
with coefficients bounded by \eqref{eq:background-hierarchy}. By the uniform bound on
$\|P_{\alpha\beta}W/\widehat\nu_M^{1/2}\|_{L^2_{\widehat v}}$ and the Cauchy-Schwarz inequality,
\begin{align}\label{eq:rw78-moment-derivatives}
 d^j|\nabla_x^j(\pi,q_H)|
 &\leq C_K\sum_{i\leq j}\sum_{|\alpha|=i}
                        \|\widehat\nu_M^{1/2}H_\alpha/W\|_{L^2_{\widehat v}},\qquad j\leq K .
\end{align}
At $|\alpha|=K$, use \eqref{eq:physical-spatial-recovery} for $H_\alpha$.

By the finite-overlap summation, the norm recovery and the pointwise embedding, we have
\begin{equation}\label{eq:expanded-block-bounds}
	\begin{aligned}
		&\qquad\qquad\quad\sum_c J_c^2\zeta_{b,c}Z_c^2\leq C A_b^{\rm E},\quad \sum_c J_c^2\zeta_{b,c}\theta_c^{3+\gamma}N_c^2\leq C\mathcal D_b,\\
		&j_c\leq C\sqrt{\mathcal B},\quad \theta_c^{(3+\gamma)/2}j_c \leq C J_c^{-1}\zeta_{b,c}^{-1/2}\mathcal D_b^{1/2},\quad \sup_c\theta_c^p j_c\leq C(T-t)^{(-1+\eta)/2}\mathcal D_b^{1/2}.
	\end{aligned}
\end{equation}
The summed bounds in \eqref{eq:expanded-block-bounds} follow from \eqref{eq:local-product-measure}, \eqref{eq:rw78-moment-derivatives} and bounded overlap. The $j_c$ bounds follow from \eqref{eq:single-weight-pointwise} and \eqref{eq:pointwise-flux-time-factor}.

For a velocity function $F(t,x,v)$, we define its source norm
\begin{equation}\label{eq:expanded-source-norm}
	\|F\|_{b,\mathrm{src}}^2 :=\int\Omega_b^*\theta^{-3-\gamma} \|\widehat\nu_M^{-1/2}F/W\|_{L^2_{\widehat v}}^2\dd x.
\end{equation}
For $R=H_0$ or $R=H_K$, by the definition of $\mathcal D_b$ and the Cauchy--Schwarz inequality, we have
\begin{equation}\label{eq:expanded-source-pair}
	\left|\iint\Omega_b^*RF/W^2\dd \widehat v\dd x\right| \leq C\mathcal D_b^{1/2}\|F\|_{b,\mathrm{src}}.
\end{equation}
From \eqref{eq:expanded-B-local} and \eqref{eq:expanded-block-bounds}, we obtain
\begin{align}\label{eq:expanded-B-global}
	\|d_\phi^K D_E\mathcal B_K\|_{b,\mathrm{src}}^2 \leq C(T-t)^{-1+\eta}A_b^{\rm E},\quad
	\|d_\phi^K\mathscr D_{\rm mic}\mathcal B_K\|_{b,\mathrm{src}}^2 \leq C(1+\mathcal B)\mathcal D_b.
\end{align}
For the first inequality, $\theta/d\leq C\theta^p$ and \eqref{eq:expanded-block-bounds} give
\[\|d_\phi^KD_E\mathcal B_K\|_{b,\mathrm{src}}^2\leq C\sum_cJ_c^2\zeta_{b,c}\theta_c^{2p-3-\gamma}Z_c^2\leq C(T-t)^{-1+\eta}A_b^{\rm E}.\]
For the second, $2p-6-2\gamma=-2\delta_{\rm gap}$ and the pointwise bound on $j_c$ give
\begin{align*}
	\sum_cJ_c^2\zeta_{b,c}\theta_c^{2p-3-\gamma}N_c^2&\leq\theta_*^{-2\delta_{\rm gap}}\sum_cJ_c^2\zeta_{b,c}\theta_c^{3+\gamma}N_c^2\leq C\mathcal D_b,\\
	\sum_cJ_c^2\zeta_{b,c}\theta_c^{2p-3-\gamma}Z_c^2j_c^2&\leq C\mathcal D_b\sum_c\theta_c^{-2\delta_{\rm gap}}Z_c^2\leq C\mathcal B\mathcal D_b.
\end{align*}
The last bound uses $\theta_c\geq\theta_*$ and $\sum_c Z_c^2\leq CA_0^E\leq C\mathcal B$. Applying the same summation to the order-zero local bounds gives
\[
 \|D_E\mathcal B_0\|_{b,\mathrm{src}}^2\leq C(T-t)^{-1+\eta}A_b^{\rm E},\qquad
 \|\mathscr D_{\rm mic}\mathcal B_0\|_{b,\mathrm{src}}^2\leq C\mathcal B\mathcal D_b,\qquad b=0,a .
\]

We also have a stronger estimate for a macroscopic test:
\begin{align}\label{eq:expanded-B-macro-pair}
	&\left|\iint\Omega_b^*\frac{P_K}{W^2} d_\phi^K\mathscr D_{\rm mic}\mathcal B_K\dd \widehat v\dd x\right| +\kappa_b^{\rm E}\left|\iint\Omega_b^*\frac{P_0}{W^2} \mathscr D_{\rm mic}\mathcal B_0\dd \widehat v\dd x\right| \notag\\
	&\hspace{12mm}\leq C(T-t)^{(-1+\eta)/2} \big((A_b^{\rm E})^{1/2}+A_b^{\rm E}\big)\mathcal D_b^{1/2}.
\end{align}
By \eqref{eq:expanded-gaussian-coefficient}, the local pairings are bounded by $CJ_c^2\zeta_{b,c}(\theta_c/d_c)Z_c(N_c+Z_cj_c)$. By the Cauchy--Schwarz inequality and \eqref{eq:expanded-block-bounds}, we have
\begin{align*}
	\sum_cJ_c^2\zeta_{b,c}\theta_c^pZ_cN_c&\leq\left(\sum_cJ_c^2\zeta_{b,c}\theta_c^{2p-3-\gamma}Z_c^2\right)^{1/2}\left(\sum_cJ_c^2\zeta_{b,c}\theta_c^{3+\gamma}N_c^2\right)^{1/2}\\
	&\leq C(T-t)^{(-1+\eta)/2}(A_b^{\rm E})^{1/2}\mathcal D_b^{1/2}.
\end{align*}
The contribution containing $Z_cj_c$ is bounded by
\[\sum_cJ_c^2\zeta_{b,c}\theta_c^pZ_c^2j_c\leq\sup_c(\theta_c^pj_c)\sum_cJ_c^2\zeta_{b,c}Z_c^2\leq C(T-t)^{(-1+\eta)/2}A_b^{\rm E}\mathcal D_b^{1/2}.\]
Summing the contributions containing $Z_cN_c$ and $Z_c^2j_c$, we get \eqref{eq:expanded-B-macro-pair}.

To compare the macroscopic principal terms with $\mathcal L_{\rm E}$, we use the quadratic form $D(\mathbf z)$ of \eqref{eq:roadmap-energy-matrix}. With the microscopic terms of \eqref{eq:normalized-macro-system} set to zero and $D_f=\partial_t+(u+\theta V)\cdot\nabla_x$, the equations at fixed $x$ are
\begin{align}\label{eq:expanded-full-macro-euler}
 \partial_t^EV+(u+\theta V)\cdot\nabla V
      +5\theta(1+B)\mathrm e^{3S/5}\nabla B
 &=-V\cdot\nabla u+\tfrac13(\Div u)V-(V\cdot\nabla\theta)V\notag\\
 &\quad-5\big((1+B)^2\mathrm e^{3S/5}-1\big)\nabla\theta,\\
 \partial_t^EB+(u+\theta V)\cdot\nabla B+\tfrac13\theta(1+B)\Div V
 &=-\tfrac43(1+B)V\cdot\nabla\theta,\notag\\
 \partial_t^ES+(u+\theta V)\cdot\nabla S&=0.\notag
\end{align}
Put $\mathbf w=\mathsf T_0\mathbf z=(V/c_{\rm s},3B,-3S/5)$, with $\mathsf T_0=\operatorname{diag}(c_{\rm s}^{-1}I_3,3,-3/5)$. In the equation for $\mathbf w$, the principal matrix in the $x_j$ direction is
\begin{equation}\label{eq:expanded-principal-matrices}
	\mathcal A_j(\mathbf z):=(u_j+\theta V_j)I_5+ c_{\rm s}\theta(1+B)
	\begin{pmatrix}
		0_{3\times3}&\mathrm e^{3S/5}e_j&0_{3\times1}\\
		e_j^{\mathsf T}&0&0\\
		0_{1\times3}&0&0
	\end{pmatrix}
	.
\end{equation}
Here $e_j$ is the $j$th coordinate vector and $V_j$ is the $j$th component of $V$. With $D(\mathbf z)$ from \eqref{eq:roadmap-energy-matrix}, we have
\[D(\mathbf z)\mathcal A_j(\mathbf z)=(u_j+\theta V_j)D(\mathbf z)+c_{\rm s}\theta(1+B)^2\mathrm e^{-3S/5}\begin{pmatrix}0_{3\times3}&e_j&0_{3\times1}\\e_j^{\mathsf T}&0&0\\0_{1\times3}&0&0\end{pmatrix}.\]
As a consequence, we have
\begin{equation}\label{eq:expanded-symmetric-difference}
	D(\mathbf z)\mathcal A_j(\mathbf z)-g_eD_*\mathcal A_j(0) \quad\hbox{is symmetric},\qquad g_e:=1-Ce.
\end{equation}
The symmetry in \eqref{eq:expanded-symmetric-difference} allows integration by parts in the difference of the two principal forms.

\begin{proposition}\label{prop:expanded-macro}
	\label{prop:new-macro-sign}
	Assume that $\mathcal Z_K\leq e$ with $e$ sufficiently small. Then there exists a constant $C$ such that for $b=0,a$, with $g_e:=1-Ce>0$,
	\begin{equation}\label{eq:new-macro-comparison}
		\partial_t^EA_b^{\rm E}\leq \frac{c_{\rm s}^2g_e}{T-t} \langle \mathcal L_{\rm E}\mathsf V,\mathsf V\rangle_{X^K_{\bar\varpi-b}} +\frac{Ce}{T-t}\|\mathsf V\|_{X^K_{\bar\varpi-b-r}}^2 +C A_b^{\rm E}.
	\end{equation}
	In $\|\mathsf V\|_{X^K_{\bar\varpi-b-r}}$, the zero-order coefficient $\varkappa_{k,\bar\varpi-b-r}$ may be replaced by $\kappa_b^{\rm E}=\varkappa_{k,\bar\varpi-b}$, since both are fixed positive constants. In particular, we have
	\begin{equation}\label{eq:new-critical-macro}
		\partial_t^EA_0^E\leq C(T-t)^{-1}A_a^E+C A_0^E.
	\end{equation}
	Let $\mathcal R_*$ be the local correction from Proposition~\ref{thm:local-observation}, constructed at order $K$, exponent $\bar\varpi-a$ and rate $h_*<1/8$, and let $A_*^E:=A_a^E+c_{\rm s}^2\mathcal R_*(\mathsf V,\mathsf V)/2$. For every fixed $\lambda_E<2h_*$, we can choose $e,T>0$ small enough such that
	\begin{equation}\label{eq:new-stable-macro}
		\partial_t^EA_*^E+\frac{\lambda_E}{T-t}A_*^E \leq\frac C{T-t}|\mathbf a(t)|_P^2.
	\end{equation}
	Here the vector $\mathbf a(t)$ is defined in \eqref{eq:roadmap-unstable-coordinates}. The condition \eqref{eq:fixed-spectral-cutoff} gives $J_x((T-t)^{1/r}y)=1$ on $B_{R_2}$, which contains the supports of the coordinate functionals $\ell_j$ and the input and output regions of $\mathscr R_*$.
\end{proposition}
\begin{proof}
	For this proof, we define
	\[\varpi_b:=\bar\varpi-b,\qquad V_K:=\Delta_y^k\mathsf V,\qquad V_0:=\mathsf V,\qquad \omega_K:=\langle y\rangle^{\varpi_b}\phi_1^{2K},\qquad \omega_0:=\kappa_b^{\rm E}\langle y\rangle^{\varpi_b}.\]
	Use $\mathbf z$ in the $(t,y)$ variables and differentiate at fixed $y$.
	\begin{align}\label{eq:rw78-macro-variation}
	 A_b^{\rm E}&=\frac{c_{\rm s}^2}{2}\sum_{i=0,K}\int\omega_iV_i^{\mathsf T}D(\mathbf z)V_i\dd y,\notag\\
	 \partial_t^EA_b^{\rm E}
	 &=c_{\rm s}^2\sum_{i=0,K}\int\omega_i
	   \left((\partial_t^EV_i)^{\mathsf T}D(\mathbf z)V_i
	         +\tfrac12 V_i^{\mathsf T}(\partial_t^ED(\mathbf z))V_i\right)\dd y .
	\end{align}
	Write $\mathcal L_{{\rm E},\infty}=-r^{-1}y\cdot\nabla_y-3\alpha_{\rm ss}$. Since $[\Delta_y^k,y\cdot\nabla_y]=K\Delta_y^k$, integration by parts in the dilation term of \eqref{eq:rw78-macro-variation} gives
	\begin{align*}
		q_{i,b}&=\frac{3+y\cdot\nabla_y\log\omega_i-2i}{2r}-3\alpha_{\rm ss},\\
		q_{0,b}&=-\frac b{2r}-\frac{\varpi_b}{2r}\langle y\rangle^{-2},\qquad
		q_{K,b}=q_{0,b}+\frac Kr(y\cdot\nabla_y\log\phi_1-1).
	\end{align*}
	Choose $g_e=1-Ce$ so that $0\leq D(\mathbf z)-g_eD_*\leq CeI_5$. By \eqref{eq:good-length-time}, $q_{K,b}\leq q_{0,b}\leq C\langle y\rangle^{-2}$, and thus
	\begin{equation}\label{eq:macro-dilation-square}
		q_{i,b}V_i^{\mathsf T}(D(\mathbf z)-g_eD_*)V_i
		\leq Ce\langle y\rangle^{-2}|V_i|^2\leq Ce\langle y\rangle^{-r}|V_i|^2.
	\end{equation}
	The two derivatives of the metric in \eqref{eq:rw78-macro-variation} and in the dilation pairing satisfy
	\begin{align}\label{eq:macro-dilation-metric}
		\left((T-t)\left.\partial_t^E\right|_y+\frac yr\cdot\nabla_y\right)D(\mathbf z)
		&=(T-t)D_{\mathbf z}D(\mathbf z)[\left.\partial_t^E\mathbf z\right|_x],\notag\\
		(T-t)|D_{\mathbf z}D(\mathbf z)[\left.\partial_t^E\mathbf z\right|_x]|
		&\leq\begin{cases}
			Ce(T-t)\theta/d\leq Ce\langle y\rangle^{-r},& |x|<R_{\rm e}/2,\\
			Ce(T-t),& |x|\geq R_{\rm e}/2,
		\end{cases}
	\end{align}
	where the second line follows from \eqref{eq:expanded-full-macro-euler}.

	For the remaining principal terms, let $\widetilde{\mathcal A}_j$ be \eqref{eq:expanded-principal-matrices} with $(u,\theta)=(U,\vartheta)$. Set
	\[C_j:=D(\mathbf z)\widetilde{\mathcal A}_j(\mathbf z)-g_eD_*\widetilde{\mathcal A}_j(0).\]
	By the profile bounds and $\mathcal Z_K\leq e$, on $\langle y\rangle\asymp R$,
	\[\|C_j\|_{L^\infty}+R\|\nabla C_j\|_{L^\infty}\leq CeR^{1-r},\qquad |\nabla\omega_i|\leq CR^{-1}\omega_i.\]
	By $C_j=C_j^{\mathsf T}$ in \eqref{eq:expanded-symmetric-difference}, for $i=0,K$,
	\begin{align}\label{eq:rw78-macro-principal}
		-\int\omega_iV_i^{\mathsf T}C_j\partial_j V_i\dd y
		=\tfrac12\int V_i^{\mathsf T}\partial_j(\omega_iC_j)V_i\dd y
		\leq Ce\int\omega_i\langle y\rangle^{-r}|V_i|^2\dd y .
	\end{align}
	
	For $|\alpha|=K$, subtract the principal commutators as follows.
	\begin{align*}
	 \mathfrak C_{\alpha,j}
	 &:=D(\mathbf z)[\partial^\alpha,\widetilde{\mathcal A}_j(\mathbf z)]\partial_j\mathsf V
	      -g_eD_*[\partial^\alpha,\widetilde{\mathcal A}_j(0)]\partial_j\mathsf V\\
	 &=\sum_{0<\lambda\leq\alpha}\binom\alpha\lambda
	   \big[(D(\mathbf z)-g_eD_*)\partial^\lambda\widetilde{\mathcal A}_j(0)
	      +D(\mathbf z)\partial^\lambda(\widetilde{\mathcal A}_j(\mathbf z)
	                                       -\widetilde{\mathcal A}_j(0))\big]
	       \partial_j\partial^{\alpha-\lambda}\mathsf V .
	\end{align*}
	On $\langle y\rangle\asymp R$, put $y=RY$, $x=(T-t)^{1/r}RY$ and $\lambda_R=J_x(x)\vartheta(RY)^3$. By \eqref{eq:profile-decay}, \eqref{eq:spatial-weight-J} and $\mathsf V=\lambda_R\mathsf T_0\mathbf z$,
	\begin{align*}
		\|\widetilde{\mathcal A}_j(0)\|_{W_Y^{K,\infty}}&\leq CR^{1-r},\qquad \|\lambda_R\|_\infty\|\lambda_R^{-1}\|_\infty\leq C,\\
		\|\mathsf V\|_{W_Y^{1,\infty}}&\leq Ce\|\lambda_R\|_\infty,\qquad
		\|\mathbf z\|_{H_Y^K}\leq C\|\lambda_R^{-1}\|_\infty\|\mathsf V\|_{H_Y^K}.
	\end{align*}
	Use $D(\mathbf z)-g_eD_*=O(e)$, $\|\mathbf z\|_{W_Y^{1,\infty}}\leq Ce$ and \eqref{eq:localized-tame-product}, with the input norms on an enlarged annulus. Since $\mathfrak C_{\alpha,j}$ has $K+1$ derivatives,
	\begin{align}\label{eq:annular-principal-commutator}
		\|R^K\mathfrak C_{\alpha,j}\|_{L^2_Y}
		&\leq CR^{-r}\bigl(e\|\mathsf V\|_{H_Y^K}+\|\mathbf z\|_{H_Y^K}\|\mathsf V\|_{W_Y^{1,\infty}}\bigr)\leq CeR^{-r}\|\mathsf V\|_{H_Y^K}.
	\end{align}
	Pair \eqref{eq:annular-principal-commutator} with $\omega_K V_K$ and sum over the annuli using \eqref{eq:euler-weighted-derivatives}.
	\begin{equation}\label{eq:rw78-macro-commutator}
	 \begin{aligned}
	  \left|\int\omega_K V_K^{\mathsf T}\sum_{j=1}^3\sum_{|\beta|=k}
	                    \frac{k!}{\beta!}\mathfrak C_{2\beta,j}\dd y\right|
	  &\leq Ce\int\langle y\rangle^{\varpi_b-r}\sum_{j\leq K}
	                           \langle y\rangle^{2j}|\nabla_y^j\mathsf V|^2\dd y\leq Ce\|\mathsf V\|_{X^K_{\varpi_b-r}}^2 .
	 \end{aligned}
	\end{equation}
	The terms of order zero have the right-hand side of \eqref{eq:rw78-macro-commutator} by \eqref{eq:localized-tame-product}. For the change from $(u_{\rm p},\theta_{\rm p})$ to $(u,\theta)$, the cutoff and \eqref{eq:untruncated-profile-bounds} give
	\begin{align*}
		\sup_{|x|\geq R_{\rm e}/2}d^{j-1}|\nabla_x^j(u-u_{\rm p},\theta-\theta_{\rm p})|&\leq C_j,\qquad 0\leq j\leq K+1,\\
		\left((T-t)\partial_t+\tfrac yr\cdot\nabla_y\right)J_x((T-t)^{1/r}y)&=0.
	\end{align*}
	For $|x|\geq R_{\rm e}/2$, use \eqref{eq:rw78-macro-principal} and \eqref{eq:localized-tame-product} to bound the principal and commutator pairings by
	\[C\sum_{j=0}^K\int\Omega_b^*|d^j\nabla_x^j\mathbf z|^2\dd x\leq CA_b^{\rm E}\]
	by \eqref{eq:macro-equivalence-new}. Using this bound and \eqref{eq:macro-dilation-square}, \eqref{eq:macro-dilation-metric}, \eqref{eq:rw78-macro-principal}, \eqref{eq:rw78-macro-commutator} in \eqref{eq:rw78-macro-variation}, we get \eqref{eq:new-macro-comparison}.

	For $b=0$, \eqref{eq:critical-Euler-input} and $\bar\varpi-r<\bar\varpi-a$ give
	\begin{align*}
		\langle\mathcal L_{\rm E}\mathsf V,\mathsf V\rangle_{X^K_{\bar\varpi}}\leq C\int\langle y\rangle^{\bar\varpi-r}|\mathsf V|^2\dd y\leq CA_a^E,\quad
		\|\mathsf V\|_{X^K_{\bar\varpi-r}}^2\leq CA_a^E.
	\end{align*}
	These bounds in \eqref{eq:new-macro-comparison} give \eqref{eq:new-critical-macro}.
	
	On $B_{R_c}$, \eqref{eq:fixed-spectral-cutoff} and \eqref{eq:cutoff-actual-euler} give $\mathsf V=\mathsf W$, $\chi_{\rm sp}=1$ and $(T-t)\partial_t^E\mathsf V=\mathcal L_{\rm E}\mathsf V+\mathsf TN$ at fixed $y$. By the symmetry and locality in \eqref{eq:observation-local-map},
	\begin{equation}\label{eq:observation-euler-time}
		\frac{T-t}{2}\partial_t^E\mathcal R_*(\mathsf V,\mathsf V)
		=\mathcal R_*(\mathcal L_{\rm E}\mathsf V,\mathsf V)+\mathcal R_*(\chi_{\rm sp}\mathsf TN,\mathsf V).
	\end{equation}
	On $B_{R_c}$, \eqref{eq:nonlinear-N-explicit} and \eqref{eq:observation-coefficient-bounds} imply
	\[\|\mathsf TN\|_{L^2(B_{R_c})}+|1-g_e|\|\mathcal L_{\rm E}\mathsf V\|_{L^2(B_{R_c})}\leq Ce\|\mathsf V\|_{H^1(B_{R_c})}.\]
	With \eqref{eq:observation-local-map}, this yields
	\[|\mathcal R_*(\chi_{\rm sp}\mathsf TN,\mathsf V)|+|(1-g_e)\mathcal R_*(\mathcal L_{\rm E}\mathsf V,\mathsf V)|\leq Ce\|\mathsf V\|_{\mathscr X}^2.\]
	Add $c_{\rm s}^2$ times \eqref{eq:observation-euler-time} to $(T-t)$ times \eqref{eq:new-macro-comparison} at $b=a$. By \eqref{eq:observation-stable-sign}, with $c_1>0$,
	\[\partial_t^EA_*^E+\frac{2h_*}{T-t}A_*^E+\frac{c_1-Ce}{T-t}\|\mathsf V\|_{\mathscr X}^2\leq\frac{C_1}{T-t}|\mathbf a|_P^2+CA_a^E.\]
	Use $Ce\leq c_1/2$, $A_a^E\leq A_*^E$ and $CT\leq2h_*-\lambda_E$ to get \eqref{eq:new-stable-macro}.
\end{proof}

\subsection{Maxwellian source and transport terms}
\begin{lemma}\label{lem:matched-source}
	Let $b\in\{0,a\}$, $N\geq0$ and $j\in\{0,\ldots,K\}$. Then
	\begin{equation}\label{eq:matched-source}
		\sum_{|\alpha|=j}\|\langle\widetilde v\rangle^N d^j\partial_x^\alpha S_0\|_{b,\mathrm{src}}^2
		\leq C_{N,K}(T-t)^{-1+\eta}.
	\end{equation}
	The same bound holds for $\|\langle\widetilde v\rangle^N d_\phi^K\Delta_x^{K/2}S_0\|_{b,\mathrm{src}}^2$.
\end{lemma}
\begin{proof}
	Write $S_0/\mu(\widehat v)=P$ by \eqref{eq:source}. The degree of $P$ is at most three and its coefficients are bounded by $C\theta/d$. At fixed $v$,
	\begin{align*}
		\mu^{-1}\partial_{x_i}(\mu P)
		&=\partial_{x_i}P|_{\widehat v}-\theta^{-1}(\partial_{x_i}u+\widehat v\partial_{x_i}\theta)\cdot(\nabla_{\widehat v}P-\widehat v P).
	\end{align*}
	Iterate the last identity. The polynomial $d^j\partial_x^\alpha S_0/\mu$ has degree at most $2j+3$, with coefficients bounded by $C_j\theta/d$ by \eqref{eq:background-hierarchy}. Use $\theta/d\leq C\theta^p$ and the support in $|x|\leq2R_0$ to get
	\begin{equation}\label{eq:source-polynomial-order}
		|d^j\partial_x^\alpha S_0|\leq C_j\theta^p\langle\widehat v\rangle^{2j+3}\mu(\widehat v)\mathbf 1_{|x|\leq2R_0},\qquad |\alpha|=j\leq K.
	\end{equation}
	Use $W^2\geq M$, $\widehat\nu_M\geq c$ and $\mu(\widehat v)^2/M\leq C\mathrm e^{-c|\widehat v|^2}$ in \eqref{eq:source-polynomial-order}.
	\begin{align}\label{eq:rw78-source-gaussian}
		\left\|\widehat\nu_M^{-1/2}\langle\widetilde v\rangle^N d^j\partial_x^\alpha S_0/W\right\|_{L^2_{\widehat v}}^2&\leq C_{N,K}\theta^{2p}\mathbf 1_{|x|\leq2R_0}\int\langle\widehat v\rangle^{2N+4j+6}\mathrm e^{-c|\widehat v|^2}\dd\widehat v\notag\\
		&\leq C_{N,K}\theta^{2p}\mathbf 1_{|x|\leq2R_0}.
	\end{align}
	By $\theta\asymp d^{1-r}$ and $\Omega_b^*\lesssim d^{-3}$ on $|x|<R_{\rm e}$, with $(r-1)(2p-3-\gamma)=r(1-\eta)>0$,
	\begin{align*}
		\int_{|x|\leq2R_0}\Omega_b^*\theta^{2p-3-\gamma}\dd x
		&\leq C+C\int_{|x|<R_{\rm e}}d^{-3-r(1-\eta)}\dd x\\
		&\leq C+C(T-t)^{-1+\eta}\int_{\R^3}\langle y\rangle^{-3-r(1-\eta)}\dd y
		\leq C(T-t)^{-1+\eta}.
	\end{align*}
	This proves \eqref{eq:matched-source}. For $\Delta_x^{K/2}S_0$, expand the Laplacian power and use $d_\phi\asymp d$.
\end{proof}

For the test of the $G_i$ equation take $d_i=1$ if $i=0$ and $d_i=d_\phi^K$ if $i=K$. In the sum over $i=0,K$ the $i=0$ term has the factor $\kappa_b^{\rm E}$. Set
\begin{align*}
	&\qquad\qquad\qquad\qquad\qquad q_i^{\rm time}:=-\frac b{r(T-t)}+2i\partial_t\log d_\phi, \\
	\Lambda_i^E&:=D\log\Omega_b^*+\frac b{r(T-t)}-3D\log\theta -2D_E\log W+2iv\cdot\nabla\log d_\phi, &\Lambda_i^H:=-2\mathscr D_{\rm mic}\log W .
\end{align*}
Substituting $G_i=H_i+P_i$ and $D=D_E+\mathscr D_{\rm mic}$ into \eqref{eq:new-good-evolution}, we obtain
\begin{align}\label{eq:expanded-main-inventory}
	&\frac{\mathrm d}{\mathrm dt}\frac12\iint\Omega_b^*\frac{G_i^2}{W^2}-\frac1\eps\iint\Omega_b^*\frac{G_i}{W^2}d_i\partial_iQ(f,f)\notag\\
	&\quad=-\iint\Omega_b^*\frac{G_i}{W^2}d_i\partial_i S_0 -\iint\Omega_b^*\frac{H_i}{W^2}d_iD_E\mathcal B_i -\iint\Omega_b^*\frac{P_i}{W^2}d_iD_E\mathcal B_i\notag\\
	&\qquad-\iint\Omega_b^*\frac{H_i+P_i}{W^2}d_i\mathscr D_{\rm mic}\mathcal B_i +\frac12\iint\Omega_b^*\frac{H_i^2+2H_iP_i+P_i^2}{W^2} (q_i^{\rm time}+\Lambda_i^E+\Lambda_i^H).
\end{align}
Here $\partial_0:=\mathrm{Id}$, $\partial_K:=\Delta_x^k$, and the integrals are taken with respect to $\dd\widehat v\dd x$. The physical integration density is $\Omega_b^*\theta^{-3}W^{-2}$, which accounts for $-3D\log\theta$ in $\Lambda_i^E$.

For $\mathscr A\in\{D_E,\mathscr D_{\rm mic}\}$ and $a_{\rm vel}:=|\widetilde v|^2$, by differentiation, we have
\[\mathscr A\log W=\tfrac12\mathscr A\log n-\tfrac34\mathscr A\log\Theta-\tfrac14\partial_a\Phi\,\mathscr A a_{\rm vel}-\tfrac14\partial_m\Phi\,\mathscr A m.\]
Use \eqref{eq:transport-weight-log-bounds}, \eqref{eq:background-bounds}, \eqref{eq:macro-equivalence-new} and \eqref{eq:expanded-kinetic-increment}. Since $\mathscr D_{\rm mic}m=0$, we have
\begin{align*}
	|\Lambda_i^E|&\leq C\theta^p\bigl(\langle\widetilde v\rangle+\langle\widetilde v\rangle^3\partial_a\Phi\bigr),\\
	|\Lambda_i^H|&\leq C\theta^pJ_H\bigl(1+\langle\widetilde v\rangle^2\partial_a\Phi\bigr).
\end{align*}
Here $J_H$ is defined in Lemma~\ref{lem:single-spatial-allocation}. The bounds
\[|\partial_m\Phi|\leq C,\qquad 0\leq\partial_a\Phi(a_{\rm vel})\leq\min\{1,Cm/a_{\rm vel}\},\qquad a_{\rm vel}>0,\]
give, with $h_\gamma:=(3-\gamma)/2>0$,
\begin{align}\label{eq:expanded-multiplier-bounds}
	\int |\Lambda_i^E+\Lambda_i^H|\,|H_i/W|^2\dd \widehat v &\leq C\theta^p(1+m^{h_\gamma})(1+J_H) \|\widehat\nu_M^{1/2}H_i/W\|_{L^2_{\widehat v}}^2,\notag\\
	\left|\int \Lambda_i^E H_iP_i/W^2\dd \widehat v\right| &\leq C\theta^p p_i^{\rm mac}\|\widehat\nu_M^{1/2}H_i/W\|_{L^2_{\widehat v}},\\
	\left|\int \Lambda_i^H H_iP_i/W^2\dd \widehat v\right| &\leq C\theta^p J_H p_i^{\rm mac}\|\widehat\nu_M^{1/2}H_i/W\|_{L^2_{\widehat v}},\notag\\
	\left|\int \Lambda_i^H P_i^2/W^2\dd \widehat v\right| &\leq C\theta^pJ_H (p_i^{\rm mac})^2.\notag
\end{align}
Set $p_0^{\rm mac}=|\mathbf z|$ and $p_K^{\rm mac}=|d_\phi^K\Sigma_x^{-1}\Delta_x^k(\Sigma_x\mathbf z)|$. Applying \eqref{eq:new-Gaussian-match} and \eqref{eq:fitted-metric-W}, we get
\begin{equation}\label{eq:macroscopic-gaussian-size}
 \begin{aligned}
  \int\Omega_b^*\big((p_K^{\rm mac})^2+\kappa_b^{\rm E}(p_0^{\rm mac})^2\big)\dd x&\leq CA_b^{\rm E},\\
  \int\langle\widetilde v\rangle^j |P_i|^2/M\dd\widehat v&\leq C_j(p_i^{\rm mac})^2,\qquad j\geq0,\ i=0,K.
 \end{aligned}
\end{equation}
	In the $H_i^2$ term use \eqref{eq:cubic-weight-bound} in the form
	\[(\langle\widetilde v\rangle^3+\langle\widetilde v\rangle^2)\partial_a\Phi\leq C(1+m^{h_\gamma})\widehat\nu_M.\]
	In \eqref{eq:expanded-multiplier-bounds} use \eqref{eq:fitted-cubic-ratio} for $H_i^2$ and \eqref{eq:macroscopic-gaussian-size} for $P_i/W$, since $W^2\geq M$.

\subsection{Gaussian weights and the transport inequality}

	At the reference instant take $H=0$, $W=\sqrt M$ in \eqref{eq:expanded-main-inventory}. By \eqref{eq:roadmap-microscopic-equation},
	\[G_i=P_i,\qquad \partial_tH=-S_M,\qquad \mathscr D_{\rm mic}\mathbf z=0,\qquad \partial_t M=\partial_t^EM.\]
	Use \eqref{eq:new-Gaussian-match} for $A_b^{\rm E}$. For the differentiated cross term, \eqref{eq:projected-moment-preservation} and \eqref{eq:source-moments-spatial} give
	\[\int\frac{P_i}{M}\partial_tH_i\dd\widehat v=-\int\frac{P_i}{M}d_i\partial_i S_M\dd\widehat v=-\int\frac{P_i}{M}d_i\partial_i S_0\dd\widehat v=0,\qquad i=0,K.\]
	Replacing $\Lambda_i^E$ by $\Lambda_i^{E,G}:=\Lambda_i^E+2D_E\log W-D_E\log M$ in \eqref{eq:expanded-main-inventory} gives
	\begin{equation}\label{eq:rw78-gaussian-reference}
		\partial_t^E A_b^{\rm E}=\sum_{i=0,K}\left(-\iint\Omega_b^*\frac{P_i}{M}d_iD_E\mathcal B_i
		+\frac12\iint\Omega_b^*\frac{P_i^2}{M}(q_i^{\rm time}+\Lambda_i^{E,G})\right).
	\end{equation}
	The sums use $\dd\widehat v\dd x$ and the coefficient $\kappa_b^{\rm E}$ at $i=0$.

	Fix $\zeta>0$ independently of $S<T$. By \eqref{eq:expanded-source-pair}, \eqref{eq:expanded-B-global}, \eqref{eq:expanded-B-macro-pair} and Young's inequality,
\begin{align}\label{eq:rw78-transport-composition}
	\sum_{i=0,K}\left|\iint\Omega_b^*\frac{H_i}{W^2}d_iD_E\mathcal B_i\right|
	&\leq\frac\zeta\eps\mathcal D_b+C_\zeta\eps(T-t)^{-1+\eta}A_b^{\rm E},\notag\\
	\sum_{i=0,K}\left|\iint\Omega_b^*\frac{H_i}{W^2}d_i\mathscr D_{\rm mic}\mathcal B_i\right|
	&\leq C(1+\mathcal B)^{1/2}\mathcal D_b,\notag\\
	\sum_{i=0,K}\left|\iint\Omega_b^*\frac{P_i}{W^2}d_i\mathscr D_{\rm mic}\mathcal B_i\right|
	&\leq\frac\zeta\eps\mathcal D_b+C_\zeta\eps(T-t)^{-1+\eta}\bigl(A_b^{\rm E}+(A_b^{\rm E})^2\bigr).
\end{align}

For the velocity multipliers, by \eqref{eq:expanded-multiplier-bounds} and \eqref{eq:single-weight-low-flux}, we have
\begin{align}\label{eq:rw78-transport-multipliers}
	\sum_{i=0,K}\iint\Omega_b^*|\Lambda_i^E+\Lambda_i^H||H_i/W|^2
	&\leq C\sup_{t,x}\{\eps\theta^{-\delta_{\rm gap}}(1+m^{h_\gamma})(1+C\sqrt{\mathcal B})\}\frac{\mathcal D_b}\eps,\notag\\
	\sum_{i=0,K}\left|\iint\Omega_b^*(\Lambda_i^E+\Lambda_i^H)H_iP_i/W^2\right|
	&\leq C(T-t)^{(-1+\eta)/2}(1+\mathcal B)^{1/2}(A_b^{\rm E})^{1/2}\mathcal D_b^{1/2}\notag\\
	&\leq\frac\zeta\eps\mathcal D_b+C_\zeta\eps(T-t)^{-1+\eta}(1+\mathcal B)A_b^{\rm E},\notag\\
	\sum_{i=0,K}\left|\iint\Omega_b^*\Lambda_i^HP_i^2/W^2\right|
	&\leq CA_b^{\rm E}\sup_x(\theta^pJ_H)\notag\\
	&\leq\frac\zeta\eps\mathcal D_b+C_\zeta\eps(T-t)^{-1+\eta}(A_b^{\rm E})^2.
\end{align}
In the first line of \eqref{eq:rw78-transport-multipliers}, $1+J_H\leq1+C\sqrt{\mathcal B}$ by \eqref{eq:single-weight-pointwise}. Then \eqref{eq:adaptive-absorption} makes the coefficient of $\mathcal D_b/\eps$ small uniformly on $[0,T)$.

	By Lemma~\ref{lem:matched-source}, $\sum_{i=0,K}\|d_i\partial_iS_0\|_{b,\mathrm{src}}^2\leq C(T-t)^{-1+\eta}$. Since $\int P_id_i\partial_iS_0/M\dd\widehat v=0$ by \eqref{eq:source-moments-spatial},
	\[\int\frac{P_id_i\partial_iS_0}{W^2}\dd\widehat v=\int_{|\widetilde v|^2>m}\frac{d_i\partial_iS_0}{W}\frac{P_i}{M}(M/W-W)\dd\widehat v.\]
	Apply \eqref{eq:expanded-source-pair} and Young's inequality to $H_i$. For $P_i$, use the Cauchy-Schwarz inequality with \eqref{eq:new-tail-small}, \eqref{eq:collision-tail-time} and \eqref{eq:macroscopic-gaussian-size}. Then
	\begin{align}\label{eq:rw78-transport-source}
		\sum_{i=0,K}\left|\iint\Omega_b^*H_id_i\partial_iS_0/W^2\right|
		&\leq\frac\zeta\eps\mathcal D_b+C_\zeta\eps(T-t)^{-1+\eta},\notag\\
		\sum_{i=0,K}\left|\iint\Omega_b^*P_id_i\partial_iS_0/W^2\right|
		&\leq C(T-t)^{(-1+\eta)/2}
		\left(\sum_i\int\Omega_b^*\theta^{3+\gamma}\mathrm e^{-2cm}(p_i^{\rm mac})^2\dd x\right)^{1/2}\notag\\
		&\leq C\eps(T-t)^{-1+\eta}(A_b^{\rm E})^{1/2}.
	\end{align}

	On $|\widetilde v|^2<m$, $W^2=M$ and $\Lambda_i^E=\Lambda_i^{E,G}$. By \eqref{eq:Dm}, $|m-m_c|\leq C$ on $B_{cd_c}(x_c)$ with $m_c=m_\eps(t,x_c)$. On the tail apply \eqref{eq:expanded-B-local} and \eqref{eq:macroscopic-gaussian-size}. This yields
\begin{align}\label{eq:rw78-transport-gaussian-tail}
	&\sum_{i=0,K}\left|\iint\Omega_b^*P_id_iD_E\mathcal B_i
	                    \left(\frac1{W^2}-\frac1M\right)\right|
	 +\sum_{i=0,K}\left|\iint\Omega_b^*P_i^2
	                    \left(\frac{\Lambda_i^E}{W^2}-\frac{\Lambda_i^{E,G}}M\right)\right|\notag\\
	&\qquad\leq C\sum_cJ_c^2\zeta_{b,c}\mathrm e^{-cm_c}\theta_c^p Z_c^2
	 \leq C\eps(T-t)^{-1+\eta}A_b^{\rm E}.
\end{align}
In the last inequality, we used $\mathrm e^{-cm}\theta^p\leq C\eps\theta^{p-\delta_{\rm gap}}\leq C\eps(T-t)^{-1+\eta}$ and $\sum_cJ_c^2\zeta_{b,c}Z_c^2\leq CA_b^{\rm E}$ from \eqref{eq:expanded-block-bounds}.

For $q_i^{\rm time}\leq0$, by the cancellation $\int H_iP_i/M\dd\widehat v=0$, we have
\[\left|\int H_iP_i/W^2\dd\widehat v\right|\leq C\mathrm e^{-cm}\|H_i/W\|_{L^2_{\widehat v}}p_i^{\rm mac}.\]
Combining Young's inequality with the Gaussian tail bound for $P_i$, we obtain, pointwise in $x$,
\begin{align}\label{eq:rw78-transport-time}
	&\frac{q_i^{\rm time}}2\int\frac{(H_i+P_i)^2}{W^2}\dd\widehat v-\frac{q_i^{\rm time}}2\int\frac{P_i^2}{M}\dd\widehat v\notag\\
	&\qquad\leq\frac{q_i^{\rm time}}4\|H_i/W\|_{L^2_{\widehat v}}^2+C|q_i^{\rm time}|\mathrm e^{-cm}(p_i^{\rm mac})^2.
\end{align}
We choose $L\geq\max\{2,a/(r\eta),1/(r\eta)\}$ and then $\kappa_{\rm w}$ such that $\mathrm e^{-cm}\leq C\eps^L\theta^{-\delta_{\rm gap}L}$. By the profile bounds, we have
\[d^{-a}\theta^{-\delta_{\rm gap}L}+d^{-1}\theta^{-\delta_{\rm gap}L}\leq C.\]
Using \eqref{eq:good-length-time}, $a/r=1/2>\eta$ and $1/r>\eta$, we obtain
\begin{align*}
	\sum_{i=0,K}\int\Omega_b^*|q_i^{\rm time}|\mathrm e^{-cm}(p_i^{\rm mac})^2\dd x
	&\leq C\eps^L\left((T-t)^{a/r-1}\mathcal B+(T-t)^{-1+1/r}(I_b+I_b^{\rm low})\right)\\
	&\leq C\eps(T-t)^{-1+\eta}\mathcal B.
\end{align*}
The negative square in \eqref{eq:rw78-transport-time} retains $bI_b^{\rm mic}/[2r(T-t)]$ of damping.

For the auxiliary terms with $j<K$ and $|\alpha|=j$, by \eqref{eq:roadmap-microscopic-equation}, we have
\[DH_\alpha=\eps^{-1}d^j\partial_x^\alpha Q(f,f)-d^j\partial_x^\alpha S_M+d^j\partial_x^\alpha T_H+j(D\log d)H_\alpha.\]
	In \eqref{eq:roadmap-projected-sources}, $S_M-S_0$ has first derivatives of $\mathbf z$, and the fluxes in $T_H$ satisfy \eqref{eq:physical-flux-derivatives}. Thus $d^j\partial_x^\alpha$ uses at most $j+1\leq K$ derivatives of the unknowns. By \eqref{eq:localized-tame-product} and \eqref{eq:expanded-gaussian-coefficient},
\begin{align*}
	\|\widehat\nu_M^{-1/2}\langle\widetilde v\rangle^{\widehat s_j}d^j\partial_x^\alpha(S_M-S_0)/W\|_{L^2_{X,\widehat v}}&\leq C(\theta_c/d_c)Z_c,\\
	\|\widehat\nu_M^{-1/2}\langle\widetilde v\rangle^{\widehat s_j}d^j\partial_x^\alpha T_H/W\|_{L^2_{X,\widehat v}}&\leq C(\theta_c/d_c)(N_c+Z_cj_c).
\end{align*}
By Lemma~\ref{lem:matched-source} and the localized bounds \eqref{eq:expanded-block-bounds},
\begin{align}\label{eq:rw78-transport-low-sources}
	\sum_{j<K,|\alpha|=j}\|\langle\widetilde v\rangle^{\widehat s_j}d^j\partial_x^\alpha S_M\|_{b,\mathrm{src}}^2
	&\leq C(T-t)^{-1+\eta}+C\sum_cJ_c^2\zeta_{b,c}\theta_c^{2p-3-\gamma}Z_c^2\notag\\
	&\leq C(T-t)^{-1+\eta}(1+\mathcal B),\notag\\
	\sum_{j<K,|\alpha|=j}\|\langle\widetilde v\rangle^{\widehat s_j}d^j\partial_x^\alpha T_H\|_{b,\mathrm{src}}^2
	&\leq C\sum_cJ_c^2\zeta_{b,c}\theta_c^{2p-3-\gamma}(N_c^2+Z_c^2j_c^2)\notag\\
	&\leq C(1+\mathcal B)\mathcal D_b.
\end{align}
The $N_c^2$ and $Z_c^2j_c^2$ terms in \eqref{eq:rw78-transport-low-sources} are summed as in \eqref{eq:expanded-B-global}. Applying \eqref{eq:expanded-source-pair} and Young's inequality to \eqref{eq:rw78-transport-low-sources} gives
\begin{equation}\label{eq:auxiliary-source-pairing}
 \begin{aligned}
 &\sum_{0\leq j<K}\sum_{|\alpha|=j}
 \left|\iint\Omega_b^*(b_j^{\rm sp}+\kappa_j\langle\widetilde v\rangle^{2\widehat s_j})
       \frac{H_\alpha}{W^2}d^j\partial_x^\alpha(T_H-S_M)\dd\widehat v\dd x\right|\\
 &\qquad\leq\frac{\zeta}{\eps}\mathcal D_b
   +C_\zeta\eps(T-t)^{-1+\eta}(1+\mathcal B)
   +C_\zeta\eps(1+\mathcal B)\mathcal D_b .
 \end{aligned}
\end{equation}
For $0\leq j<K$, the multiplier from the lower-order weight satisfies
\[
 \left|D\log\bigl(b_j^{\rm sp}+\kappa_j\langle\widetilde v\rangle^{2\widehat s_j}\bigr)\right|
 \leq C\theta^p(1+J_H)\langle\widetilde v\rangle
 \leq C\theta^p(1+J_H)\widehat\nu_M .
\]
	Use \eqref{eq:rw78-transport-multipliers}. In the lower-order test,
	\[-\frac{b}{r(T-t)}+2j\partial_t\log d=-\frac{b}{r(T-t)}-\frac{2j(T-t)^{2/r-1}}{rd^2}\leq-\frac{b}{r(T-t)}.\]
	This gives $bI_b^{\rm low}/[2r(T-t)]$ on the left.

\begin{proposition}\label{prop:expanded-transport}
	\label{prop:new-transport-reduction}
	Use the background of Proposition~\ref{Thm.Euler-solution} and the coefficients of Proposition~\ref{prop:new-component-collision}. For $\kappa_{\rm w},m_0$ sufficiently large and $\mathcal B,\eps$ sufficiently small, we have, for $b=0,a$,
	\begin{equation}\label{eq:new-transport-reduction}
		(I_b+I_b^{\rm low})'+\frac c\eps\mathcal D_b+\frac b{2r(T-t)}(I_b^{\rm mic}+I_b^{\rm low})
		\leq\partial_t^EA_b^{\rm E}+C\eps(T-t)^{\eta-1}(1+\mathcal B)^2 .
	\end{equation}
	The constants $C,c$ and the smallness threshold for $\eps$ are independent of $S<T$.
\end{proposition}
\begin{proof}
	Add the lower-order tests to \eqref{eq:expanded-main-inventory} and subtract \eqref{eq:rw78-gaussian-reference}. Use \eqref{eq:new-component-collision} for the collision bound and \eqref{eq:rw78-transport-time} for the damping. Using \eqref{eq:rw78-transport-composition}--\eqref{eq:rw78-transport-time}, \eqref{eq:auxiliary-source-pairing} and $A_b^{\rm E}\leq C\mathcal B$, we get
	\begin{align}\label{eq:transport-before-absorption}
		&(I_b+I_b^{\rm low})'+\frac b{2r(T-t)}(I_b^{\rm mic}+I_b^{\rm low})-\partial_t^EA_b^{\rm E}\notag\\
		&\quad\leq-\frac{\mathcal D_b}{\eps}\Bigl[c-C\sqrt{\mathcal B}-C\zeta-C_\zeta\eps(1+\mathcal B)^{1/2}
		-C_\zeta\eps^2(1+\mathcal B)\notag\\
		&\hspace{28mm}{}-C\sup_{t,x}\big\{\eps\theta^{-\delta_{\rm gap}}(1+m^{h_\gamma})(1+J_H)\big\}\Bigr]
		+C_\zeta\eps(T-t)^{-1+\eta}(1+\mathcal B)^2 .
	\end{align}
	Choose $\zeta$, then the bound on $\mathcal B$, then the bound on $\eps$ by
	\begin{equation}\label{eq:transport-absorption-choice}
	 \begin{aligned}
	  C\zeta&\leq c/8,\qquad C\sqrt{\mathcal B}\leq c/8,\\
	  C_\zeta\eps(1+\mathcal B)^{1/2}+C_\zeta\eps^2(1+\mathcal B)&\leq c/8,\\
	  C\sup_{0\leq t<T,\,x}\{\eps\theta^{-\delta_{\rm gap}}(1+m^{h_\gamma})(1+J_H)\}&\leq c/8 .
	 \end{aligned}
	\end{equation}
	The last inequality in \eqref{eq:transport-absorption-choice} follows from \eqref{eq:adaptive-absorption} and $J_H\leq C\sqrt{\mathcal B}$ by \eqref{eq:single-weight-pointwise}. Substituting \eqref{eq:transport-absorption-choice} in \eqref{eq:transport-before-absorption} leaves at least $c/2$ multiplying $\mathcal D_b/\eps$ on the left of \eqref{eq:new-transport-reduction}.

\end{proof}

\section{The macroscopic correction and the closed inequalities}\label{sec:closed}

\subsection{The projected source estimate}

For even $K\geq2k_0+2$, set $Y^\alpha=(T-t)^{|\alpha|/r}\partial_x^\alpha|_v$ and define
\begin{align*}
	Z_{K-1,R_2}(t)&:=\|\mathbf z^\sharp(t)\|_{H^{K-1}(B_{R_2})},\\
	D_{K-1,R_2}(t)&:=(T-t)^{-3/r}\sum_{|\alpha|\leq K-1}
	\int_{|x|<R_2(T-t)^{1/r}}\int\theta^{3+\gamma}\widehat\nu_M|Y^\alpha H/W|^2\dd\widehat v\dd x.
\end{align*}
Here $B_R=\{|y|<R\}$, $W=W_{m_\eps,\ell}$ and $(n,V,\Theta)\in\mathcal U$. The substitution $v=u+\theta\widehat v$ is made after applying $Y^\alpha$.

\begin{proposition}\label{thm:actual-spectral-forcing}
	Let $Z_{K-1,R_2}\leq b_0$ with $b_0>0$ small. Fix the tail order such that
	\[\sup_{(n,V,\Theta)\in\mathcal U,\,m\geq m_0}\int\langle\widetilde v\rangle^6W^2/\widehat\nu_M\dd\widehat v<\infty.\]
	Then $F_{\rm u}$ from \eqref{eq:actual-a-equation} satisfies
	\begin{equation}\label{eq:local-mode-forcing}
		|F_{\rm u}(t)|^2\leq C(T-t)^{-1+\eta}D_{K-1,R_2}(t)+C(T-t)^{-2}Z_{K-1,R_2}(t)^4.
	\end{equation}
	Here $C$ is uniform in $m,\eps,t$.
\end{proposition}
\begin{proof}
	Let $J_P=\theta^{-3}\int P(\widehat v)H\dd v$ and $Y_i=Y^{e_i}$. At fixed $v$,
	\[Y_i\widehat v=-\theta^{-1}(Y_i u+\widehat v Y_i\theta),\qquad Y_i(\theta^{-3})=-3\theta^{-4}Y_i\theta.\]
	Thus
	\begin{equation}\label{eq:moment-physical-derivative}
		Y_iJ_P=\int P(\widehat v)Y_iH\dd\widehat v-\frac{Y_i\theta}{\theta}J_{3P+\widehat v\cdot\nabla_{\widehat v}P}-\frac{Y_i u}{\theta}\cdot J_{\nabla_{\widehat v}P}.
	\end{equation}
	On $B_{R_2}$, \eqref{eq:profile} gives
	\[\partial_{y_i}J_P^\sharp=(Y_iJ_P)^\sharp,\qquad \left(\frac{Y_i\theta}{\theta},\frac{Y_i u}{\theta}\right)^\sharp=\frac{(\partial_{y_i}\vartheta,\partial_{y_i}U)}{\vartheta}.\]
	Iterate \eqref{eq:moment-physical-derivative} to obtain
	\[\partial_y^\alpha J_P^\sharp=\sum_{|\lambda|\leq|\alpha|}\int P_{\alpha\lambda}(y,\widehat v)Y^\lambda H\dd\widehat v,\qquad \deg_{\widehat v}P_{\alpha\lambda}\leq\deg P.\]
	The coefficients of $P_{\alpha\lambda}$ are bounded on $B_{R_2}$ for $|\alpha|\leq K-1$. For $P=\widehat v_i\widehat v_j,|\widehat v|^2\widehat v_i/2$, we use $\langle\widehat v\rangle\asymp\langle\widetilde v\rangle$ and obtain
	\begin{align*}
		|\partial_y^\alpha J_P^\sharp|^2
		&\leq C\left(\int\frac{\langle\widetilde v\rangle^6W^2}{\widehat\nu_M}\dd\widehat v\right)
		\sum_{|\lambda|\leq|\alpha|}\int\widehat\nu_M|Y^\lambda H/W|^2\dd\widehat v
		\leq C\sum_{|\lambda|\leq|\alpha|}\int\widehat\nu_M|Y^\lambda H/W|^2\dd\widehat v .
	\end{align*}
	Use \eqref{eq:expanded-moments}, $\dd y=(T-t)^{-3/r}\dd x$ and $\theta^{3+\gamma}\asymp(T-t)^{-1-\eta}$ to get
	\begin{equation}\label{eq:rw9-local-moments}
		\|\pi^\sharp\|_{H^{K-1}(B_{R_2})}^2+\|q_H^\sharp\|_{H^{K-1}(B_{R_2})}^2\leq C(T-t)^{1+\eta}D_{K-1,R_2}.
	\end{equation}
	
	In the following coefficient estimates, $B,S$ denote their pullbacks to $B_{R_2}$. Since $K-2>3/2$ and $Z_{K-1,R_2}\leq b_0$, by the Sobolev product estimate and Taylor's formula, we have
	\begin{align*}
		\|(1+B)^2\mathrm e^{3S/5}-1-2B-3S/5\|_{H^{K-2}(B_{R_2})}&\leq CZ_{K-1,R_2}^2,\\
		\|(1+B)\mathrm e^{3S/5}-1\|_{H^{K-2}(B_{R_2})}&\leq CZ_{K-1,R_2},\\
		\|\big((1+B)\mathrm e^{3S/5}-1\big)\nabla B\|_{H^{K-2}(B_{R_2})}&\leq CZ_{K-1,R_2}\|\nabla B\|_{H^{K-2}(B_{R_2})}\leq CZ_{K-1,R_2}^2.
	\end{align*}
	Hence \eqref{eq:nonlinear-N-explicit} and \eqref{eq:dimensionless-forces} give
	\begin{align*}
		\|N\|_{H^{K-2}(B_{R_2})}
		&\leq C\|\mathbf z^\sharp\|_{H^{K-2}(B_{R_2})}\left(\|\mathbf z^\sharp\|_{H^{K-2}(B_{R_2})}+\|\nabla\mathbf z^\sharp\|_{H^{K-2}(B_{R_2})}\right)\leq CZ_{K-1,R_2}^2,\\
		\|F_h\|_{H^{K-2}(B_{R_2})}
		&\leq C\left(\|\pi^\sharp\|_{H^{K-1}(B_{R_2})}+\|q_H^\sharp\|_{H^{K-1}(B_{R_2})}\right).
	\end{align*}
	Multiplication by $\chi_{\rm sp}\mathsf T$ maps $H^{K-2}(B_{R_2})$ to $X^{K-2}_\varpi$ boundedly. By \eqref{eq:actual-a-equation} and \eqref{eq:rw9-local-moments},
	\begin{align*}
		|F_{\rm u}|^2&\leq C(T-t)^{-2}\|\chi_{\rm sp}\mathsf T(N+F_h)\|_{X^{K-2}_\varpi}^2\\
		&\leq C(T-t)^{-2}Z_{K-1,R_2}^4+C(T-t)^{-1+\eta}D_{K-1,R_2}.
	\end{align*}
	On $B_x=B_{R_2(T-t)^{1/r}}$, $J_x=1$, $d\asymp(T-t)^{1/r}$ and $\Omega_a^*\asymp(T-t)^{-3/r}$. Thus
	\begin{equation}\label{eq:projected-source-local-comparison}
		\begin{aligned}
			Z_{K-1,R_2}^2
			&\leq C\sum_{|\alpha|\leq K-1}\int_{B_x}\Omega_a^*|d^{|\alpha|}\partial_x^\alpha\mathbf z|^2\dd x
				\leq CA_a^{\rm E}\leq C\mathcal E,\\
			D_{K-1,R_2}
			&\leq C\sum_{|\alpha|\leq K-1}\iint_{B_x\times\R^3}\Omega_a^*\theta^{3+\gamma}\widehat\nu_M
				|d^{|\alpha|}\partial_x^\alpha H/W|^2\dd\widehat v\dd x\leq C\mathcal D_a .
		\end{aligned}
	\end{equation}
	Here we used \eqref{eq:macro-equivalence-new}, \eqref{eq:new-norm-comparison} and \eqref{eq:physical-spatial-recovery}. Apply \eqref{eq:local-mode-forcing} and the norm equivalence on $E_{\rm u}$ to get \eqref{eq:new-actual-mode-bound}.
\end{proof}

\subsection{Pairing the correction with the actual microscopic flux}
	Let $\mathsf W_c=\chi_{\rm sp}\mathsf W$ and $\psi=\mathscr R_*\mathsf W_c$ with $\mathscr R_*$ in \eqref{eq:observation-operator}. Equations \eqref{eq:observation-local-map} and \eqref{eq:fixed-spectral-cutoff} give
\begin{equation}\label{eq:observation-test-support}
	\supp\psi\subset B_{R_c}\subset\{|y|<R_1\},\qquad
	\chi_{\rm sp}\psi=\psi,\qquad \langle\psi,[\mathcal L_{\rm E},\chi_{\rm sp}]\mathsf W\rangle_0=0 .
\end{equation}
For the pairing of \eqref{eq:actual-W-equation} with $\psi$, define
	\begin{equation}\label{eq:observation-local-dissipation}
		D_{0,R_c}:=(T-t)^{-3/r}\int_{|x|<R_c(T-t)^{1/r}}\int\theta^{3+\gamma}\widehat\nu_M|H/W|^2\dd\widehat v\dd x,\qquad W=W_{m_\eps,\ell}.
	\end{equation}

\begin{proposition}\label{prop:correction-flux}
	Assume that the parameters of $\mathbf z^\sharp=(V,B,S)$ lie in the compact range and $\|\mathbf z^\sharp\|_{W^{1,\infty}(B_{R_c})}\leq b_c$, with $b_c$ fixed. For every $\zeta>0$,
	\begin{align}\label{eq:correction-flux}
		\left|\frac1{T-t}\mathcal R_* (\mathsf W_c,\chi_{\rm sp}\mathsf T F_h)\right| &\leq C(T-t)^{(-1+\eta)/2} \|\mathsf W_c\|_{H^1(B_{R_c})}D_{0,R_c}^{1/2}\notag\\
		&\leq\frac\zeta\eps D_{0,R_c} +C_\zeta\eps(T-t)^{-1+\eta} \|\mathsf W_c\|_{H^1(B_{R_c})}^2.
	\end{align}
	For the source $N$ we have
	\begin{equation}\label{eq:correction-nonlinear}
		\left|\frac1{T-t}\mathcal R_*(\mathsf W_c,\chi_{\rm sp}\mathsf T N)\right| \leq\frac C{T-t}\|\mathbf z^\sharp\|_{W^{1,\infty}(B_{R_c})} \|\mathsf W_c\|_{H^1(B_{R_c})}^2.
	\end{equation}
	The radius $R_c$ is fixed before $T,\eps$. The microscopic norm in \eqref{eq:correction-flux} is $D_{0,R_c}$ and only involves $H$ itself.
\end{proposition}
\begin{proof}
	Let $m_h=(\pi,q_H)^\sharp$ and $\partial_j=\partial_{y_j}$. Evaluate the coefficients $V,B,S$ at $x=(T-t)^{1/r}y$. Expanding \eqref{eq:dimensionless-forces} in \eqref{eq:euler-microscopic-source} gives
	\begin{equation}\label{eq:rw9-flux-divergence}
		\mathsf T F_h=\sum_j\partial_j\{C_j(y,\mathbf z^\sharp)m_h\}+C_0(y,\mathbf z^\sharp,\nabla\mathbf z^\sharp)m_h\qquad(|y|<R_c),
	\end{equation}
	where the matrices are specified by
	\begin{align*}
		(C_jm_h)_V&:=-\frac{\vartheta^4n^{-1}}{c_{\rm s}}\pi_{\cdot j},\quad (C_jm_h)_B:=-\frac{2\vartheta^4}{5(1+B)^4}(q_{H,j}-V\cdot\pi_{\cdot j}),\\
		(C_jm_h)_S&:=\frac{2\vartheta^4}{5(1+B)^5}(q_{H,j}-V\cdot\pi_{\cdot j}),\quad
		(C_0m_h)_V:=-\frac{5\vartheta^3n^{-1}}{c_{\rm s}}\pi\nabla\vartheta-\sum_j((\partial_jC_j)m_h)_V,\\
		(C_0m_h)_B&:=-\frac{2\vartheta^3}{5(1+B)^4}(\pi:\nabla U+6q_H\cdot\nabla\vartheta-5V\cdot\pi\nabla\vartheta)-\sum_j((\partial_jC_j)m_h)_B,\\
		(C_0m_h)_S&:=\frac{2\vartheta^3}{5(1+B)^5}(\pi:\nabla U+6q_H\cdot\nabla\vartheta-5V\cdot\pi\nabla\vartheta)-\sum_j((\partial_jC_j)m_h)_S.
	\end{align*}
	Here $\pi_{\cdot j}$ is the $j$th column of $\pi$ and $n=(1+B)^3\mathrm e^{-3S/5}$. The compact parameter range and $\|\mathbf z^\sharp\|_{W^{1,\infty}(B_{R_c})}\leq b_c$ imply
	\[\sum_j\|C_j\|_{W^{1,\infty}(B_{R_c})}+\|C_0\|_{L^\infty(B_{R_c})}\leq C.\]
	Integrate \eqref{eq:rw9-flux-divergence} against $D_*\psi$, using $\chi_{\rm sp}\psi=\psi$ from \eqref{eq:observation-test-support}. Then
	\[\mathcal R_*(\mathsf W_c,\chi_{\rm sp}\mathsf T F_h)=-\sum_j\int(\partial_j\psi)^{\mathsf T}D_*C_jm_h\dd y+\int\psi^{\mathsf T}D_*C_0m_h\dd y.\]
	By the Cauchy-Schwarz inequality in \eqref{eq:expanded-moments},
	\[\|m_h\|_{L^2(B_{R_c})}^2\leq C\int_{B_{R_c}}\int\widehat\nu_M|H/W|^2\dd\widehat v\dd y\leq C(T-t)^{1+\eta}D_{0,R_c}.\]
	By \eqref{eq:observation-local-map}, $\|\psi\|_{H^j}\leq C\|\mathsf W_c\|_{H^j(B_{R_c})}$, $j=0,1$. The pairing divided by $T-t$ is bounded by
	\[\frac C{T-t}\|\psi\|_{H^1}\|m_h\|_{L^2(B_{R_c})}\leq C(T-t)^{(-1+\eta)/2}\|\mathsf W_c\|_{H^1(B_{R_c})}D_{0,R_c}^{1/2}.\]
	Young's inequality yields \eqref{eq:correction-flux}. By \eqref{eq:nonlinear-N-explicit},
	\[\|\mathsf T N\|_{L^2(B_{R_c})}\leq C\|\mathbf z^\sharp\|_{W^{1,\infty}(B_{R_c})}\|\mathsf W_c\|_{H^1(B_{R_c})}.\]
	Pairing with $\psi/(T-t)$ gives \eqref{eq:correction-nonlinear}.
\end{proof}

\subsection{The closed energy inequalities}\label{subsec:closed-energy-proofs}

\begin{proof}[Proof of Proposition~\ref{thm:new-closed-pair}]
	For $b=0$, apply \eqref{eq:new-critical-macro} in \eqref{eq:new-transport-reduction}. By \eqref{eq:new-norm-comparison},
	\begin{align*}
		\mathcal B'+\frac{c_0}\eps\mathcal D_0
		&\leq\frac C{T-t}A_a^E+CA_0^E+C\eps(T-t)^{-1+\eta}(1+\mathcal B)^2
		\leq\frac C{T-t}\mathcal E+C\mathcal B+C\eps(T-t)^{-1+\eta}.
	\end{align*}
	This is \eqref{eq:new-closed-global}.

	For $b=a$, $\mathsf V=\mathsf W_c=\mathsf W$ on $B_{R_c}$ by \eqref{eq:fixed-spectral-cutoff}. Differentiate using \eqref{eq:cutoff-actual-euler} and \eqref{eq:observation-test-support}. Then
	\[\frac{T-t}{2}\frac{\dd}{\dd t}\mathcal R_*(\mathsf V,\mathsf V)=\langle\psi,\mathcal L_{\rm E}\mathsf W_c+\chi_{\rm sp}\mathsf T(N+F_h)\rangle_0.\]
	Subtract \eqref{eq:observation-euler-time} to get
	\begin{equation}\label{eq:observation-microscopic-time}
		\left(\frac\dd{\dd t}-\partial_t^E\right)\frac{c_{\rm s}^2}{2}\mathcal R_*(\mathsf V,\mathsf V)
		=\frac{c_{\rm s}^2}{T-t}\mathcal R_*(\mathsf W_c,\chi_{\rm sp}\mathsf T F_h).
	\end{equation}
	By \eqref{eq:observation-local-dissipation}, $D_{0,R_c}\leq D_{K-1,R_2}\leq C\mathcal D_a$. Using \eqref{eq:correction-flux} and \eqref{eq:macro-equivalence-new}, we get
	\begin{align*}
		\left|\frac{c_{\rm s}^2}{T-t}\mathcal R_*(\mathsf W_c,\chi_{\rm sp}\mathsf T F_h)\right|&\leq\frac\zeta\eps D_{0,R_c}+C_\zeta\eps(T-t)^{-1+\eta}\|\mathsf W_c\|_{H^1(B_{R_c})}^2\\
		&\leq\frac{C\zeta}\eps\mathcal D_a+C_\zeta\eps(T-t)^{-1+\eta}A_a^E.
	\end{align*}
	Take $C\zeta\leq c/2$ with $c$ from \eqref{eq:new-transport-reduction}. By \eqref{eq:observation-microscopic-time}, \eqref{eq:new-stable-macro} and $A_a^E\leq C\mathcal B\leq Ce_0$, with $M_a=I_a^{\rm mic}+I_a^{\rm low}$ we obtain
	\begin{equation}\label{eq:rw9-inner-before-damping}
		\mathcal E'+\frac c\eps\mathcal D_a+\frac{a}{2r(T-t)}M_a+\frac{\lambda_E}{T-t}A_*^E\leq\frac C{T-t}|\mathbf a|_P^2+C\eps(T-t)^{-1+\eta}.
	\end{equation}
	By \eqref{eq:matched-energy-tail-comparison},
	\[|\mathcal E-(A_*^E+M_a)|=|I_a-I_a^{\rm mic}-A_a^{\rm E}|\leq C\mathrm e^{-cm_0}(A_*^E+M_a).\]
	Since $a/(2r)=1/4>\lambda_E>\lambda$, take $m_0$ with $C\lambda\mathrm e^{-cm_0}\leq\min\{\lambda_E-\lambda,1/4-\lambda\}$. Then
	\begin{equation}\label{eq:closed-inner-damping-comparison}
		\lambda\mathcal E\leq\lambda(1+C\mathrm e^{-cm_0})(A_*^E+M_a)\leq\lambda_E A_*^E+\tfrac14M_a.
	\end{equation}
	Apply \eqref{eq:closed-inner-damping-comparison} in \eqref{eq:rw9-inner-before-damping} to get \eqref{eq:new-closed-inner}. The thresholds are independent of $S<T$, and \eqref{eq:new-actual-mode-bound} follows from Proposition~\ref{thm:actual-spectral-forcing}.
\end{proof}

\begin{proof}[Proof of Proposition~\ref{prop:mixed-velocity-estimate}]
	Set $(i,j,N)=(|\alpha|,|\beta|,|\alpha|+|\beta|)$ with $j\geq1$, $Z=Z_{\alpha,\beta}$ from \eqref{eq:mixed-velocity-derivatives}, and $Z_{\alpha,0}=d^{|\alpha|}\partial_x^\alpha H$. Let $\mathcal D_{\rm all}=\mathcal D_0+\mathcal D_K^{\rm vel}$.

	For $T_hg(v)=g(v+h)$, \eqref{eq:postcollisional-velocities} and translation of $v_*$ in \eqref{Eq.collision-operator} give $T_hQ(g,f)=Q(T_hg,T_hf)$. Differentiating in $h$ and $x$ gives
	\begin{equation}\label{eq:rw9-collision-derivatives}
		\partial_x^\alpha\partial_v^\beta Q(g,f)
		=\sum_{\alpha_1\leq\alpha}\sum_{\beta_1\leq\beta}\binom\alpha{\alpha_1}\binom\beta{\beta_1}
		Q(\partial_x^{\alpha_1}\partial_v^{\beta_1}g,\partial_x^{\alpha-\alpha_1}\partial_v^{\beta-\beta_1}f).
	\end{equation}
	For $\phi\in\operatorname{span}\{1,v_i,|v|^2\}$, $\partial_v^\beta\phi$ is in the same span. Integration by parts and \eqref{eq:zero-moments} give
	\[\int\phi Z_{\alpha,\beta}\dd v=(-1)^j d^i\theta^j\partial_x^\alpha\int(\partial_v^\beta\phi)H\dd v=0.\]
	This is \eqref{eq:collision-missing} for $Z/W$. The commutator $[\partial_{v_l},D]=\partial_{x_l}$ gives
	\[d^i\theta^j[\partial_x^\alpha\partial_v^\beta,D]H=\sum_{l:\beta_l>0}\beta_l(\theta/d)Z_{\alpha+e_l,\beta-e_l}.\]
	Apply $d^i\theta^j\partial_x^\alpha\partial_v^\beta$ to \eqref{eq:roadmap-microscopic-equation}.
	\begin{equation}\label{eq:mixed-exact-equation}
	 \begin{aligned}
	  D Z_{\alpha,\beta}
	  &={}\eps^{-1}d^i\theta^j\partial_x^\alpha\partial_v^\beta Q(f,f)
	      +d^i\theta^j\partial_x^\alpha\partial_v^\beta(T_H-S_M)\\
	  &\quad+(iD\log d+jD\log\theta)Z_{\alpha,\beta}
	      -\sum_{l:\beta_l>0}\beta_l\frac\theta d Z_{\alpha+e_l,\beta-e_l}.
	 \end{aligned}
	\end{equation}

	Take $N\leq K-3$, $0\leq q\leq2(K-N-2)$ and $A=1$ or $A=\theta^{(3+\gamma)/2}\widehat\nu_M^{1/2}$. By \eqref{eq:background-hierarchy}, $d^h|\nabla_x^h(d^i\theta^j)|\leq Cd^i\theta^j$ for $h\leq2$. With $l$ derivatives on $\partial_x^\alpha\partial_v^\beta H$, the loss in \eqref{eq:rw56-pointwise-multiplier} is covered by
	\[q+2(2-l)\leq p_{N+l}=2(K-N-l),\qquad 0\leq l\leq2,\quad N+l\leq K-1.\]
	Apply $H_X^2(L_v^2)\hookrightarrow L_X^\infty(L_v^2)$ as in Lemma~\ref{lem:single-spatial-allocation}. With $\dd\widehat v=\theta^{-3}\dd v$,
	\begin{align*}
		J_x(x_0)^2\|A\langle\widetilde v\rangle^q Z/W\|_{L^2_{\widehat v}}^2(x_0)\leq C\sum_{|\lambda|\leq2}\iint\Omega_0^*A^2\left|\langle\widetilde v\rangle^{p_{N+|\lambda|}} Z_{\alpha+\lambda,\beta}/W\right|^2\dd\widehat v\dd x.
	\end{align*}
	Using \eqref{eq:mixed-velocity-energy}, we have
	\begin{equation}\label{eq:rw9-mixed-pointwise}
	 \begin{aligned}
	  \|\langle\widetilde v\rangle^q Z/W\|_{L^2_{\widehat v}}
	     \leq CJ_x^{-1}(\mathcal B+\mathcal V_K)^{1/2},\quad
	  \theta^{(3+\gamma)/2}\|\widehat\nu_M^{1/2}\langle\widetilde v\rangle^q Z/W\|_{L^2_{\widehat v}}
	     \leq CJ_x^{-1}\mathcal D_{\rm all}^{1/2}.
	 \end{aligned}
	\end{equation}
	For $\beta=0$ the needed weights are $\widehat s_N=p_N$ when $N<K$, and \eqref{eq:physical-spatial-recovery} when $N=K$.

	Using $f=M+H$ and $Q(M,M)=0$, the terms in \eqref{eq:rw9-collision-derivatives} linear in $H$ with all derivatives on $H$ give $Q(M,Z)+Q(Z,M)$ after multiplying by $d^i\theta^j$. Each other linear term has inputs
	\[d^{|\alpha_1|}\theta^{|\beta_1|}\partial_x^{\alpha_1}\partial_v^{\beta_1}M,\qquad Z_{\alpha-\alpha_1,\beta-\beta_1},\qquad p_{N-N_1}=p_N+2N_1,\quad N_1=|\alpha_1|+|\beta_1|>0.\]
	When the $\mathbf z$ derivatives have order at most $K-2$, apply \eqref{eq:macro-equivalence-new} and \eqref{eq:expanded-gaussian-coefficient} with $\theta\partial_v=\Theta^{-1/2}\partial_{\widetilde v}$. By \eqref{eq:polynomial-bilinear}, the pairing for $R=Z_{\alpha-\alpha_1,\beta-\beta_1}$ and $s=0,p_N$ is bounded by
	\[\delta\theta^{3+\gamma}\|\widehat\nu_M^{1/2}\langle\widetilde v\rangle^s Z/W\|_2^2+C_\delta\theta^{3+\gamma}\|\widehat\nu_M^{1/2}\langle\widetilde v\rangle^s R/W\|_2^2.\]
	Use $s\leq p_{N-N_1}$ and $\mathcal D_0$ when $\beta-\beta_1=0$. As in \eqref{eq:spatial-collision-hierarchy}, start with $(b_K^{\rm vel},k_K^{\rm vel})=(1,0)$ and choose
	\begin{equation}\label{eq:mixed-weight-hierarchy}
		k_N^{\rm vel}\geq C(1+L_N),\qquad b_N^{\rm vel}\geq C(1+U_N+C_{p_N}k_N^{\rm vel}),\qquad N=K-1,\ldots,1.
	\end{equation}
	Here $L_N,U_N$ are the weighted and unweighted losses from larger orders. The contribution $C_{p_N}k_N^{\rm vel}$ is at order $N$ by \eqref{eq:polynomial-gap}.

	If the Maxwellian contains a derivative of $\mathbf z$ of order at least $K-1$, then
	\[K-1\leq|\alpha_1|\leq|\alpha|\leq K-j\leq K-1,\qquad j=1,\quad N=K,\quad p_N=0,\quad N-N_1\leq1.\]
	Write $A_{\rm hi}$ for the scaled $\mathbf z$ derivative times its bounded coefficients. By \eqref{eq:macro-equivalence-new},
	\[\|J_x^{-1}A_{\rm hi}\|_{L^2(\Omega_0^*\dd x)}\leq C\sqrt{\mathcal B}.\]
	Use \eqref{eq:polynomial-bilinear} with $q=0$ in \eqref{eq:rw9-mixed-pointwise} or \eqref{eq:single-weight-pointwise} for the low input. Thus
	\begin{align}\label{eq:mixed-high-Maxwellian}
		&\mathcal D_{\rm all}^{1/2}\int\Omega_0^*J_x^{-1}|A_{\rm hi}|
			\theta^{(3+\gamma)/2}\|\widehat\nu_M^{1/2}Z/W\|_{L^2_{\widehat v}}\dd x\leq C\mathcal D_{\rm all}\|J_x^{-1}A_{\rm hi}\|_{L^2(\Omega_0^*\dd x)}\leq C\sqrt{\mathcal B}\mathcal D_{\rm all}.
	\end{align}

	For $Q(H,H)$ take $N_1\leq N_2$, $N_1+N_2=N$. If $N_2\geq2$,
	\[N_1\leq\lfloor K/2\rfloor\leq K-3,\qquad p_N=2(K-N_1-N_2)\leq2(K-N_1-2).\]
	Thus \eqref{eq:rw9-mixed-pointwise} applies to the order-$N_1$ input with $q=0,p_N$. If $N_2\leq1$, place the input carrying $p_N$ in $L_x^2$ and use $q=0$ on the other input. For either placement, $p_N\leq p_{N_i}$ on the $L_x^2$ input, with \eqref{eq:physical-spatial-recovery} at purely spatial order $K$. The two terms in \eqref{eq:polynomial-bilinear} place $\widehat\nu_M^{1/2}$ on either factor, so \eqref{eq:rw9-mixed-pointwise} gives
	\begin{equation}\label{eq:mixed-quadratic-collision}
		\begin{aligned}
			\left|\sum_{\mathfrak A_K}\iint\Omega_0^*(b_N^{\rm vel}+k_N^{\rm vel}\langle\widetilde v\rangle^{2p_N})
			\frac Z{W^2}d^i\theta^j\partial_x^\alpha\partial_v^\beta Q(H,H)\dd\widehat v\dd x\right|\leq C\sqrt{\mathcal B+\mathcal V_K}\mathcal D_{\rm all}.
		\end{aligned}
	\end{equation}
	Add \eqref{eq:mixed-high-Maxwellian} and \eqref{eq:mixed-quadratic-collision} to the coercivity estimate from \eqref{eq:mixed-weight-hierarchy}.
	\begin{align}\label{eq:mixed-collision-lower-bound}
		-\eps^{-1}\sum_{\mathfrak A_K}\iint &\Omega_0^*(b_N^{\rm vel}+k_N^{\rm vel}\langle\widetilde v\rangle^{2p_N})\frac{Z_{\alpha,\beta}}{W^2}\,d^i\theta^j\partial_x^\alpha\partial_v^\beta Q(f,f)\dd\widehat v\dd x\notag\\
		&\geq\frac{c_1}\eps\mathcal D_K^{\rm vel}-\frac{C_1}\eps\mathcal D_0-\frac{C_1\sqrt{\mathcal B+\mathcal V_K}}\eps\mathcal D_{\rm all}.
	\end{align}
	
	For the source terms, we have $\theta\partial_v=\Theta^{-1/2}\partial_{\widetilde v}$ and $i+1\leq K$. Therefore, by \eqref{eq:expanded-composition}, \eqref{eq:expanded-kinetic-increment}, \eqref{eq:expanded-gaussian-coefficient} and Lemma~\ref{lem:matched-source}, we have
	\begin{align*}
		\sum_{\mathfrak A_K}\|\langle\widetilde v\rangle^{p_N}d^i\theta^j\partial_x^\alpha\partial_v^\beta S_M\|_{0,\rm src}^2&\leq C(T-t)^{-1+\eta}+C\sum_cJ_c^2\theta_c^{2p-3-\gamma}Z_c^2\\
		&\leq C(T-t)^{-1+\eta}(1+\mathcal B),\\
		\sum_{\mathfrak A_K}\|\langle\widetilde v\rangle^{p_N}d^i\theta^j\partial_x^\alpha\partial_v^\beta T_H\|_{0,\rm src}^2&\leq C\sum_cJ_c^2\theta_c^{2p-3-\gamma}(N_c^2+Z_c^2j_c^2)\leq C(1+\mathcal B)\mathcal D_0.
	\end{align*}
	In the bound involving $N_c$ in \eqref{eq:expanded-block-bounds}, $T_H$ uses at most $|\alpha|+1\leq K$ spatial derivatives of the moments. The $v$ derivatives hit only the Gaussian polynomial in \eqref{eq:roadmap-projected-sources}. By the Cauchy-Schwarz inequality and Young's inequality, for any $\zeta_0>0$,
	\begin{equation}\label{eq:mixed-source-pairing}
		\begin{aligned}
			&\left|\sum_{\mathfrak A_K}\iint\Omega_0^*(b_N^{\rm vel}+k_N^{\rm vel}\langle\widetilde v\rangle^{2p_N})\frac Z{W^2}d^i\theta^j\partial_x^\alpha\partial_v^\beta(T_H-S_M)\dd\widehat v\dd x\right|\\
			&\qquad\leq\frac{\zeta_0}\eps\mathcal D_K^{\rm vel}+C_{\zeta_0}\eps(T-t)^{-1+\eta}(1+\mathcal B)+C_{\zeta_0}\eps(1+\mathcal B)\mathcal D_0.
		\end{aligned}
	\end{equation}
	
	The measure in \eqref{eq:mixed-exact-equation} is $\Omega_0^*\theta^{-3}W^{-2}\dd v\dd x$. The coefficient of the square from the measure, $d^i\theta^j$ and the polynomial factor is
	\begin{align*}
	 \Lambda_{\alpha,\beta}
	 &=D\log(\Omega_0^*\theta^{-3}W^{-2})+2iD\log d+2jD\log\theta\\
	 &\quad+D\log(b_N^{\rm vel}+k_N^{\rm vel}\langle\widetilde v\rangle^{2p_N}).
	\end{align*}
	By \eqref{Eq.Euler-u-theta} and the definition of $d$,
	\[D\log\theta=-\tfrac13\Div u+\widehat v\cdot\nabla\theta,\qquad 2i\partial_t\log d=-\frac{2i(T-t)^{2/r-1}}{rd^2}\leq0.\]
	Also, at fixed $v$,
	\[\left|D\log(b_N^{\rm vel}+k_N^{\rm vel}\langle\widetilde v\rangle^{2p_N})\right|\leq C\frac{|D\widetilde v|}{\langle\widetilde v\rangle}.\]
	Retain the sign of $2i\partial_t\log d$. For the remaining square coefficient apply \eqref{eq:transport-weight-log-bounds} and \eqref{eq:cubic-weight-bound} to obtain
	\begin{equation}\label{eq:mixed-multiplier-coefficient}
		\begin{aligned}
			&\sum_{\mathfrak A_K}\iint\Omega_0^*(b_N^{\rm vel}+k_N^{\rm vel}\langle\widetilde v\rangle^{2p_N})\left|\Lambda_{\alpha,\beta}-2i\partial_t\log d\right|\frac{|Z|^2}{W^2}\dd\widehat v\dd x\\
			&\qquad\leq C\sup_{t,x}\left\{\eps\theta^{-\delta_{\rm gap}}(1+m_\eps^{(3-\gamma)/2})(1+J_H)\right\}\frac{\mathcal D_K^{\rm vel}}\eps.
		\end{aligned}
	\end{equation}
	By \eqref{eq:single-weight-pointwise}, $J_H\leq C\sqrt{e_{\rm vel}}$. Thus the supremum in \eqref{eq:mixed-multiplier-coefficient} over $0\leq t<T$ and $x\in\R^3$ tends to zero by \eqref{eq:adaptive-absorption}.

	In the commutator of \eqref{eq:mixed-exact-equation}, the two factors are $Z_{\alpha,\beta}$ and $Z_{\alpha+e_l,\beta-e_l}$, both of order $N$ with power $p_N$. Since
	\[\frac\theta d\leq C\theta^p,\qquad \theta^{p-3-\gamma}=\theta^{-\delta_{\rm gap}}\leq C,\]
	the Cauchy-Schwarz inequality bounds its pairing by
	\begin{equation}\label{eq:mixed-transport-commutator-bound}
		\begin{aligned}
			&\left|\sum_{(\alpha,\beta)\in\mathfrak A_K}\ \sum_{l:\beta_l>0}
			 \iint\Omega_0^*(b_N^{\rm vel}+k_N^{\rm vel}\langle\widetilde v\rangle^{2p_N})
			 \frac{\beta_l\theta Z_{\alpha,\beta}Z_{\alpha+e_l,\beta-e_l}}{d W^2}\dd\widehat v\dd x\right|\\
			&\quad\leq C(\mathcal D_K^{\rm vel})^{1/2}(\mathcal D_0+\mathcal D_K^{\rm vel})^{1/2}
			 \leq C(\mathcal D_0+\mathcal D_K^{\rm vel}) .
		\end{aligned}
	\end{equation}
	For $\beta-e_l=0$ and $N=K$, we apply Lemma~\ref{lem:weighted-elliptic-recovery} with $p_K=0$.
	
	We choose
	\[0<e_{\rm vel}\leq\min\left\{1,e_0,\left(\frac{c_1}{8C_1}\right)^2\right\},\qquad\zeta_0:=\frac{c_1}{8}.\]
	Let $C_2$ be a bound for the coefficients in \eqref{eq:mixed-multiplier-coefficient} and \eqref{eq:mixed-transport-commutator-bound}. By \eqref{eq:adaptive-absorption}, we choose $0<\eps_{\rm vel}\leq\eps_0$ such that
	\begin{equation}\label{eq:mixed-absorption-choice}
	 \begin{aligned}
	  C_2\sup_{t,x}\big\{\eps\theta^{-\delta_{\rm gap}}(1+m_\eps^{(3-\gamma)/2})(1+C\sqrt{e_{\rm vel}})\big\}&\leq c_1/8,\\
	  C_2\eps_{\rm vel}&\leq c_1/8,\qquad 0<\eps\leq\eps_{\rm vel}.
	 \end{aligned}
	\end{equation}
	The supremum in \eqref{eq:mixed-absorption-choice} is over $0\leq t<T$ and $x\in\R^3$. By \eqref{eq:mixed-collision-lower-bound}-\eqref{eq:mixed-transport-commutator-bound}, \eqref{eq:mixed-absorption-choice} and $2i\partial_t\log d\leq0$, with $\eps,\mathcal B\leq1$ we get
	\[\mathcal V_K'+\frac{c_1-4(c_1/8)}\eps\mathcal D_K^{\rm vel}\leq\frac{C_{\rm vel}}\eps\mathcal D_0+C_{\rm vel}\eps(T-t)^{-1+\eta}(1+\mathcal B).\]
	This is \eqref{eq:mixed-velocity-differential} with $c_{\rm vel}=c_1/2$. If $f_0=M(0)$, then
	\[H(0,x,v)=0\quad\Longrightarrow\quad Z_{\alpha,\beta}(0,x,v)=0\quad\Longrightarrow\quad \mathcal V_K(0)=0.\qedhere\]
\end{proof}

\subsection{Augmented global energy and the choice of constants}\label{subsec:augmented-global}

	Take $0<\kappa_{\rm reg}\leq\min\{1,c_0/(2C_{\rm vel})\}$ and $c_{\rm reg}=\min\{c_0/2,c_{\rm vel}\}$. In the sum of \eqref{eq:new-closed-global} and $\kappa_{\rm reg}$ times \eqref{eq:mixed-velocity-differential},
	\[(c_0-\kappa_{\rm reg}C_{\rm vel})\mathcal D_0+\kappa_{\rm reg}c_{\rm vel}\mathcal D_K^{\rm vel}\geq c_{\rm reg}\mathcal D_{\rm reg}.\]
	Use $\mathcal B\leq\mathcal B_{\rm reg}$ and $1+\mathcal B\leq1+e_{\rm vel}$ in the sum to obtain \eqref{eq:augmented-global-inequality}. By \eqref{eq:roadmap-continuation-threshold} and \eqref{eq:augmented-bootstrap-comparison},
	\[\mathcal B+\mathcal V_K\leq\kappa_{\rm reg}^{-1}\mathcal B_{\rm reg}\leq\kappa_{\rm reg}^{-1}e_{\rm reg}\leq e_{\rm vel}/2\quad\text{if }\mathcal B_{\rm reg}\leq e_{\rm reg}.\]
	If $f_0=M(0)$, then $\mathcal B_{\rm reg}(0)=\mathcal B(0)$.

	Fix $\eps>0$ and $S<T$. On $[0,S]$, by \eqref{eq:fixed-tail-weight}, \eqref{eq:new-global-weights} and the compact parameter bounds,
	\[W^{-2}\asymp_{S,\eps}\langle v\rangle^\ell,\qquad d\asymp_S\langle x\rangle,\qquad\theta\asymp_S1,\qquad\Omega_0^*\asymp_S\langle x\rangle^{2M_x-3}.\]
	For
	\[\|g\|_{\mathcal Y_{\rm mix}}^2:=\sum_{|\alpha|+|\beta|\leq K}\iint\langle x\rangle^{2M_x-3}\langle v\rangle^\ell|\partial_x^\alpha\partial_v^\beta g|^2\dd v\dd x,\]
	use $\dd\widehat v=\theta^{-3}\dd v$ in \eqref{eq:physical-spatial-recovery} and \eqref{eq:mixed-velocity-energy} to get
	\[\|H\|_{\mathcal Y_{\rm mix}}^2\leq C_{S,\eps}(I_0^{\rm mic}+I_0^{\rm low}+\mathcal V_K)\leq C_{S,\eps}\mathcal B_{\rm reg}.\]
	Differentiate the Maxwellians and use \eqref{eq:macro-equivalence-new} and \eqref{eq:new-norm-comparison}. Then
	\begin{align*}
		\|M-\bar M\|_{\mathcal Y_{\rm mix}}^2&\leq C_{S,\eps}\mathcal Z_K^2\leq C_{S,\eps}\mathcal B,\\
		\|\bar M-M_\infty\|_{\mathcal Y_{\rm mix}}^2&\leq C_S,\qquad
			\supp_x(\bar M-M_\infty)\subset\overline B_{2R_0}.
	\end{align*}
	Use $f-M_\infty=H+(M-\bar M)+(\bar M-M_\infty)$ to obtain
	\begin{equation}\label{eq:fixed-mixed-control}
		\sup_{0\leq t\leq S}\|f(t)-M_\infty\|_{\mathcal Y_{\rm mix}}^2\leq C_{S,\eps}\left(1+\sup_{0\leq t\leq S}\mathcal B_{\rm reg}(t)\right).
	\end{equation}
	Here $\mathcal Y_{\rm mix}\hookrightarrow C^1(\R^6)$ by $K\geq14$. Higher differentiability follows from Propositions~\ref{thm:B-local} and \ref{prop:B-time-regularity}.

	The choices below are independent of $S<T$.
\begin{enumerate}
	\item Fix $\gamma$, the profile, $r$ and the exterior radius. Set $a=r/2$, $M_x>3$ and take $K\geq\max\{14,2k_0+2\}$ even.
	\item Take $0<d_*<\min\{\eta,1/4\}$ and the projection on $X^K_{\bar\varpi-a}$ from Proposition~\ref{prop:euler-spectrum} with \eqref{eq:two-order-projection-agreement}. Fix $P$ when $E_{\rm u}\ne\{0\}$. Take $\max\{0,2\kappa\}<\delta<d_*$ and $\kappa,\delta,h_*,\lambda,\lambda_E$ satisfying \eqref{eq:new-spectral-choices}. Then fix $\mathcal R_*,t_{\rm obs},R_c$ and the radii in \eqref{eq:fixed-spectral-cutoff}.
	\item Set $(q_0,q_1)=(8,2K+8)$ and choose $\ell$ for the collision estimates and Appendix~\ref{app:polynomial-local}. With $\ell$ fixed, use \eqref{eq:spatial-collision-hierarchy}, \eqref{eq:mixed-weight-hierarchy} uniformly for $m\geq m_0$. Increase $m_0\geq\ell$ for the Gaussian-tail comparisons and \eqref{eq:closed-inner-damping-comparison}. Fix $\kappa_{\rm w}$ by the tail bounds of Section~\ref{sec:expanded-critical}.
	\item Take $e_0$ from Proposition~\ref{thm:new-closed-pair} and $0<e_{\rm vel}\leq e_0$ from Proposition~\ref{prop:mixed-velocity-estimate}. Choose $\kappa_{\rm reg}$ in \eqref{eq:augmented-global}, and $C_{\rm rec}$ in \eqref{eq:roadmap-continuation-macroscopic}, \eqref{eq:roadmap-continuation-pointwise} and \eqref{eq:selection-local-macro-bound}. At fixed $\delta_*,b_0$ take $e_{\rm reg},e_{\rm sh}$ by \eqref{eq:roadmap-continuation-threshold}.
	\item Take $0<T\leq T_{\rm kin}$ with $R_2T^{1/r}<R_{\rm e}/2$ and $CT\leq2h_*-\lambda_E$. By \eqref{eq:adaptive-absorption}, choose $\eps$ for \eqref{eq:transport-absorption-choice} and \eqref{eq:mixed-absorption-choice}. Apply Theorem~\ref{thm:coupled-shooting} to \eqref{eq:positive-initial-Maxwellian} and take $0<\eps_0\leq\eps_{\rm sh}$.
\end{enumerate}
	Equations \eqref{eq:fitted-metric-ABS} and \eqref{eq:fitted-metric-W} are uniform on the initial parameter ball for $0<\eps\leq\eps_0$.

\section{Correction of the unstable modes}\label{sec:selection}

\subsection{Initial data, energy bounds and continuation}\label{subsec:selection-local}

	\begin{proof}[Proof of Proposition~\ref{prop:positive-realization}]
	By \eqref{eq:positive-initial-field} and the support of the functionals in Proposition~\ref{prop:projection-local},
	\begin{equation}\label{eq:initial-profile-support}
		\begin{aligned}
			\Pu(\chi_{\rm init}\Phi_j)&=\Phi_j,\qquad \supp\mathbf z_0^\sharp\subset\overline B_{R_1},\\
			\|\mathbf z_0^\sharp\|_{H^m_y}&\leq C_m|\mathbf a_0|_P\quad(m\geq0),\qquad
			\|\mathbf z_0^\sharp\|_{L^\infty_y}\leq C_0|\mathbf a_0|_P .
		\end{aligned}
	\end{equation}
	Use $x=T^{1/r}y$ in \eqref{eq:positive-initial-field} and set $(\rho_0,u_0,\mathcal T_0)=(\theta(0)^3n_0,u(0)+\theta(0)V_0,\theta(0)^2\Theta_0)$ as in \eqref{eq:positive-initial-Maxwellian}. By continuity at $\mathbf z=0$, choose $\delta_{\rm in}>0$ with
	\begin{equation}\label{eq:initial-parameter-range}
	 |\mathbf z|\leq\delta_{\rm in}
	 \ \Longrightarrow\quad B>-1,\quad |V|\leq1/8,\quad
	 \big((1+B)^3\mathrm e^{-3S/5},(1+B)^2\mathrm e^{3S/5}\big)\in[7/8,9/8]^2 .
	\end{equation}
	Choose $R_{\rm in}$ with $2C_0R_{\rm in}\leq\delta_{\rm in}$ using \eqref{eq:initial-profile-support}. For $|\mathbf a_0|_P<2R_{\rm in}$, \eqref{eq:initial-parameter-range} gives $\rho_0,\mathcal T_0>0$. Taking the moments of $f_0=\mathcal M(\rho_0,u_0,\mathcal T_0)$ and using $\chi_{\rm sp}\chi_{\rm init}=\chi_{\rm init}$, we get
	\begin{equation}\label{eq:initial-fitted-data}
		\begin{aligned}
			M(0)&=f_0,\qquad H(0)=0,\\
			\Pu(\chi_{\rm sp}\mathsf T\mathbf z_0^\sharp)&=\Pu\mathsf T\mathbf z_0^\sharp=\sum_{j=1}^J(\mathbf a_0)_j\Phi_j,\qquad \mathbf a(0)=\mathbf a_0 .
		\end{aligned}
	\end{equation}
	For $O\in SO(3)$ the radial fields in \eqref{eq:positive-initial-field} give
	\[
	 V_0(Ox)=OV_0(x),\qquad (B_0,S_0^{\rm ent})(Ox)=(B_0,S_0^{\rm ent})(x),\qquad
	 f_0(Ox,Ov)=f_0(x,v).
	\]
	Moreover, we have
	\[f_0=\bar M(0)\quad\text{for }|x|\geq R_1T^{1/r},\qquad (\rho_0,u_0,\mathcal T_0)=(\theta_\infty^3,0,\theta_\infty^2)\quad\text{for }|x|\geq2R_0.\]
	By the Gaussian integration, we have
	\[\begin{aligned}
		\HH(f_0\mid M_\infty)&=\int_{B_{2R_0}}\left(\rho_0\log\frac{\rho_0}{\theta_\infty^3}-\rho_0+\theta_\infty^3+\frac{\rho_0|u_0|^2}{2\theta_\infty^2}\right)\dd x\\
		&\quad+\frac32\int_{B_{2R_0}}\rho_0\left(\frac{\mathcal T_0}{\theta_\infty^2}-1-\log\frac{\mathcal T_0}{\theta_\infty^2}\right)\dd x<\infty.
	\end{aligned}\]
	Equations \eqref{eq:background} and \eqref{eq:initial-parameter-range} give
	\[\rho_0\geq\tfrac78\theta(0,x)^3\geq\tfrac78\theta_*^3,\qquad \mathcal T_0\geq\tfrac78\theta(0,x)^2\geq\tfrac78\theta_*^2.\]
	At fixed $T$, \eqref{eq:positive-initial-field}, \eqref{eq:positive-initial-Maxwellian} and these bounds give $f_0\in C^\infty(\mathbf a_0,x,v)$. At $J=0$, \eqref{eq:positive-initial-field} gives $\mathbf z_0=0$, so $f_0=\bar M(0)$.
\end{proof}

\begin{proof}[Proof of Lemma~\ref{lem:selection-initial-energy}]
	By \eqref{eq:initial-fitted-data},
	\[G_i(0)=P_i(0)\ (i=0,K),\qquad I_b^{\rm low}(0)=\mathcal V_K(0)=0,\qquad \mathcal B_{\rm reg}(0)=\mathcal B(0).\]
	Since $J_x=1$ on $\supp\mathbf z_0$ by \eqref{eq:initial-profile-support} and \eqref{eq:fixed-spectral-cutoff}, \eqref{eq:weighted-macro-variable} gives
	\[\mathsf V(0,y)=\sum_{j=1}^J(\mathbf a_0)_j\chi_{\rm init}(y)\Phi_j(y),\qquad \|\mathsf V(0)\|_{X^K_{\bar\varpi-b}}\leq C|\mathbf a_0|_P\quad(b=0,a).\]
	By $W^{-2}\leq M^{-1}$, \eqref{eq:new-Gaussian-match} and \eqref{eq:observation-zero-order-bound},
	\begin{equation}\label{eq:initial-energy-pairings}
		\begin{aligned}
			I_b(0)&\leq\tfrac12\iint\Omega_b^*\frac{|P_K(0)|^2+\kappa_b^{\rm E}|P_0(0)|^2}{M(0)}\dd\widehat v\dd x
			=A_b^{\rm E}(0)\leq C|\mathbf a_0|_P^2,\\
			\mathcal R_*(\mathsf V(0),\mathsf V(0))&\leq C\|\mathsf V(0)\|_{L^2(B_{R_c})}^2\leq C|\mathbf a_0|_P^2 .
		\end{aligned}
	\end{equation}
	By \eqref{eq:spatial-master-energies} and $\mathcal B_{\rm reg}(0)=\mathcal B(0)$, \eqref{eq:initial-energy-pairings} gives
	\[\mathcal E(0;\mathbf a_0)+\mathcal B_{\rm reg}(0;\mathbf a_0)\leq C_I|\mathbf a_0|_P^2,\qquad \mathcal E(0;\mathbf 0)=\mathcal B_{\rm reg}(0;\mathbf 0)=0.\]
	For $J>0$ take the derivative in direction $h\in\R^J$ at fixed $(y,\widehat v)$. The measure $\Omega_b^*\dd x\dd\widehat v$ and $m_\eps$ are independent of $\mathbf a_0$. In \eqref{eq:matched-change-of-variables} the inputs to $D_{\mathbf z}M$ are linear in $\mathbf a_0$. By \eqref{eq:explicit-Maxwellian-derivative}, \eqref{eq:initial-profile-support} and $0\leq\partial_a\Phi\leq1$,
	\begin{equation}\label{eq:initial-gaussian-derivatives}
		\begin{aligned}
			|P_i|&\leq C|\mathbf a_0|_P\langle\widetilde v\rangle^2M,\\
			|D_{\mathbf a_0}P_i[h]|&\leq C(1+|\mathbf a_0|_P)|h|_P\langle\widetilde v\rangle^4M,\\
			|D_{\mathbf a_0}\log W[h]|&\leq C|h|_P\langle\widetilde v\rangle^2,\qquad m_\eps\geq m_0,\quad |\mathbf a_0|_P\leq R_{\rm in}.
		\end{aligned}
	\end{equation}
	On $|y|<R_1$, $J_x=1$. By \eqref{eq:initial-gaussian-derivatives} and $W^{-2}\leq M^{-1}$, for $i=0,K$,
	\begin{align*}
		\frac{|P_iD_{\mathbf a_0}P_i[h]|+|P_i|^2|D_{\mathbf a_0}\log W[h]|}{W^2}
		&\leq C(|\mathbf a_0|_P+|\mathbf a_0|_P^2)|h|_P\langle\widetilde v\rangle^6M\\
		&\leq C(|\mathbf a_0|_P+|\mathbf a_0|_P^2)|h|_P\mathrm e^{-c|\widehat v|^2}.
	\end{align*}
	The bound holds on a neighborhood of the initial ball, by compactness of the parameter range. Use \eqref{eq:matched-change-of-variables} and $G_i(0)=P_i(0)$ to differentiate under the integral.
	\begin{equation}\label{eq:rw10-initial-energy-derivative}
		\begin{aligned}
			D_{\mathbf a_0}I_b(0)[h]=\iint\frac{\Omega_b^*}{W^2}\bigl(&P_KD_{\mathbf a_0}P_K[h]+\kappa_b^{\rm E}P_0D_{\mathbf a_0}P_0[h]\\
				&-(|P_K|^2+\kappa_b^{\rm E}|P_0|^2)D_{\mathbf a_0}\log W[h]\bigr)\dd\widehat v\dd x .
		\end{aligned}
	\end{equation}
	Integrate on $B_{R_1}\times\R^3_{\widehat v}$ and use \eqref{eq:observation-zero-order-bound}. Then
	\begin{align*}
		|D_{\mathbf a_0}I_b(0)[h]|&\leq C(|\mathbf a_0|_P+|\mathbf a_0|_P^2)|h|_P,\\
		|D_{\mathbf a_0}\mathcal R_*(\mathsf V(0),\mathsf V(0))[h]|
			&=2|\mathcal R_*(\mathsf V(0),D_{\mathbf a_0}\mathsf V(0)[h])|
			\leq C|\mathbf a_0|_P|h|_P .
	\end{align*}
	Take $\sup_{|h|_P=1}$ in the derivative of \eqref{eq:spatial-master-energies} to get \eqref{eq:shooting-initial}.
\end{proof}

\begin{proof}[Proof of Lemma~\ref{lem:selection-energy-validity}]
	By \eqref{eq:shooting-bootstrap}, \eqref{eq:augmented-bootstrap-comparison} and \eqref{eq:roadmap-continuation-threshold},
	\[\mathcal B+\mathcal V_K\leq\kappa_{\rm reg}^{-1}\mathcal B_{\rm reg}\leq\kappa_{\rm reg}^{-1}e_{\rm sh}\leq e_{\rm vel}/4.\]
	The initial parameters satisfy \eqref{eq:initial-parameter-range}. While $(n,V,\Theta)\in\mathcal U$, \eqref{eq:roadmap-continuation-macroscopic} gives
	\[\sup_x\big(|n-1|+|V|+|\Theta-1|\big)\leq C_{\rm rec}e_{\rm sh}^{1/2}\leq1/8 .\]
	Hence $7/8\leq n,\Theta\leq9/8$ and $|V|\leq1/8$ up to $S$ by $L^\infty_x$ continuity. From \eqref{eq:roadmap-continuation-pointwise} and \eqref{eq:selection-local-macro-bound} we get
	\begin{equation}\label{eq:selection-recovered-smallness}
		\sup_{x,v}\langle\widetilde v\rangle^8|H/W|\leq C_{\rm rec}\sqrt{e_{\rm sh}}\leq\delta_*/8,\qquad
		Z_{K-1,R_2}\leq C_{\rm rec}\sqrt{e_{\rm sh}}\leq b_0/2 .
	\end{equation}
	Apply Propositions~\ref{thm:new-closed-pair}, \ref{prop:mixed-velocity-estimate} and \ref{thm:actual-spectral-forcing}.
\end{proof}

\begin{lemma}\label{lem:selection-local-family}
	Fix $0<T\leq T_{\rm kin}$. There is $R_{\rm loc}\in(0,R_{\rm in})$ independent of $\eps$ such that for $0<\eps\leq\eps_{\rm vel}$ all data \eqref{eq:positive-initial-Maxwellian} with $|\mathbf a_0|_P\leq R_{\rm loc}$ have unique smooth, strictly positive solutions on $[0,h_{\eps,T}]$. For $O\in SO(3)$, $f(t,Ox,Ov)=f(t,x,v)$, and
	\begin{equation}\label{eq:selection-local-initial-energy}
		\mathcal E(0)+\mathcal B_{\rm reg}(0)+|\mathbf a_0|_P^2\leq e_{\rm sh}/8 .
	\end{equation}
	The integrals in \eqref{eq:new-kinetic-energies}, \eqref{eq:spatial-master-energies} and \eqref{eq:mixed-velocity-energy} are finite and continuous on compact continuation intervals, with
	\[(n,V,\Theta)\in C_tL_x^\infty,\qquad \mathcal E,\mathcal B,\mathcal V_K,\mathbf a\in C_t^1,\qquad F_{\rm u}\in C_t .\]
	Suppose the solution for $\mathbf a_0$ satisfies \eqref{eq:shooting-bootstrap} on $[0,S]$, $S<T$. For $|\mathbf b_0|_P\leq R_{\rm loc}$ sufficiently close to $\mathbf a_0$, the solution for $\mathbf b_0$ exists on $[0,S]$ and
	\begin{equation}\label{eq:selection-continuous-dependence}
		\begin{aligned}
			\sup_{0\leq t\leq S}\big(&|\mathcal E(t;\mathbf b_0)-\mathcal E(t;\mathbf a_0)|
			+|\mathcal B_{\rm reg}(t;\mathbf b_0)-\mathcal B_{\rm reg}(t;\mathbf a_0)|\\
			&+|\mathbf a(t;\mathbf b_0)-\mathbf a(t;\mathbf a_0)|_P\big)\longrightarrow0 \qquad(\mathbf b_0\to\mathbf a_0).
		\end{aligned}
	\end{equation}
\end{lemma}
\begin{proof}
	Take $M_{\rm ref}=\bar M(0)$ and $(q_0,q_1)=(8,2K+8)$ in \eqref{eq:B-reference-norms}. At $s\mathbf a_0$, $0\leq s\leq1$, $|\mathbf a_0|_P\leq R_{\rm in}$, \eqref{eq:initial-parameter-range} gives $\Theta_0\leq9/8$ and $|V_0|\leq1/8$. In particular,
	\[\frac{f_0(s\mathbf a_0)^2}{M_{\rm ref}}\leq C\exp\left(-\frac{|\widehat v-V_0|^2}{\Theta_0}+\frac{|\widehat v|^2}{2}\right)\leq C\mathrm e^{-c|\widehat v|^2}.\]
	For each $(\alpha,\beta)$, differentiating \eqref{eq:physical-Maxwellian} gives $\partial_x^\alpha\partial_v^\beta f_0$ and $\partial_x^\alpha\partial_v^\beta D_{\mathbf a_0}f_0[h]$ as $f_0$ times polynomials in $\widehat v$ with bounded coefficients at fixed $T$. The coefficients of the second polynomial are $O(|h|_P)$. Thus $W_q\geq\sqrt{M_{\rm ref}}\langle\widehat v\rangle^{-q}$ and
	\[f_0(\mathbf a_0)-M_{\rm ref}=\int_0^1D_{\mathbf a_0}f_0(s\mathbf a_0)[\mathbf a_0]\dd s\qquad\Longrightarrow\qquad |f_0(\mathbf a_0)-M_{\rm ref}|_{q_0}\leq C_{\rm in}|\mathbf a_0|_P,\]
	with $C_{\rm in}$ independent of $\eps$. Set
	\[R_{\rm loc}:=\min\left\{\frac{R_{\rm in}}2,\frac{\delta_*}{8C_{\rm in}},\sqrt{\frac{e_{\rm sh}}{8(C_I+1)}}\right\}.\]
	By the Gaussian estimate and \eqref{eq:initial-profile-support},
	\begin{equation}\label{eq:selection-local-data-norms}
		\begin{aligned}
			\sup_{|\mathbf a_0|_P\leq2R_{\rm loc}}|f_0(\mathbf a_0)-M_{\rm ref}|_{q_0}
			&\leq 2C_{\rm in}R_{\rm loc}\leq\delta_*/4,\\
			\sup_{|\mathbf a_0|_P\leq2R_{\rm loc}}|f_0(\mathbf a_0)-M_{\rm ref}|_{q_1,N,m_x}
			&<\infty\qquad(N\in\mathbb N_0,\ m_x\geq0).
		\end{aligned}
	\end{equation}
	Propositions~\ref{thm:B-local} and \ref{prop:B-time-regularity} give the common interval $[0,h_{\eps,T}]$. For \eqref{eq:selection-local-initial-energy}, use \eqref{eq:shooting-initial} and $(C_I+1)R_{\rm loc}^2\leq e_{\rm sh}/8$. Rotating $(v_*,\omega)$ gives, for $f_O(t,x,v)=f(t,Ox,Ov)$,
	\[Df_O=\eps^{-1}Q(f_O,f_O),\qquad f_O(0)=f_0,\qquad |f_O-M_{\rm ref}|_{q_0}=|f-M_{\rm ref}|_{q_0}.\]
	Uniqueness gives $f_O=f$.

	Take $m_x>M_x+K+4$ in \eqref{eq:B-spatial-decay-bound} and \eqref{eq:B-spatial-decay-time}. Use the majorant \eqref{eq:B-energy-integrability} for the integrands and their first time derivatives, and \eqref{eq:roadmap-unstable-coordinates} for $\mathbf a,F_{\rm u}$.

	At fixed $\eps,T$, the Gaussian derivative bounds give $|f_0(\mathbf b_0)-f_0(\mathbf a_0)|_{q_1,N,m_x}\to0$. By \eqref{eq:local-reference-difference},
	\[\sup_{0\leq t\leq h_{\eps,T}}|f(t;\mathbf b_0)-f(t;\mathbf a_0)|_{q_1,N}\leq C_N|f_0(\mathbf b_0)-f_0(\mathbf a_0)|_{q_1,N}\longrightarrow0.\]
	At a restart $t_0\leq S$, freeze $M_{\rm ref}=M(t_0;\mathbf a_0)$. If the solutions converge at $t_0$ as $\mathbf b_0\to\mathbf a_0$, \eqref{eq:selection-recovered-smallness} gives
	\[\max_{\mathbf c_0\in\{\mathbf a_0,\mathbf b_0\}}|f(t_0;\mathbf c_0)-M_{\rm ref}|_{q_0}\leq\delta_*/8+|f(t_0;\mathbf b_0)-f(t_0;\mathbf a_0)|_{q_0}<\delta_*/4.\]
	By \eqref{eq:roadmap-continuation-macroscopic} and Proposition~\ref{thm:B-restart}, the reference bounds and restart length are uniform for $t_0\leq S$. Iterate \eqref{eq:local-reference-difference} and apply \eqref{eq:B-energy-integrability} and \eqref{eq:roadmap-unstable-coordinates} to obtain \eqref{eq:selection-continuous-dependence}.
\end{proof}

\begin{lemma}\label{lem:selection-continuation}
	Let $0<T\leq T_{\rm kin}$, $0<\eps\leq\eps_{\rm vel}$, $|\mathbf a_0|_P\leq R_{\rm loc}$, and let the solution of Lemma~\ref{lem:selection-local-family} exist on $[0,S_*)$, $0<S_*<T$. If
	\[\sup_{0\leq t<S_*}\big(\mathcal E(t)+\mathcal B_{\rm reg}(t)+|\mathbf a(t)|_P^2\big)<e_{\rm sh},\]
	then $f$ extends smoothly to $[0,S_*+h]$ for some $0<h<T-S_*$.
\end{lemma}
\begin{proof}
	On the compact subintervals of $[0,S_*)$, Lemma~\ref{lem:selection-energy-validity} and \eqref{eq:selection-recovered-smallness} give
	\[(n,V,\Theta)\in\operatorname{int}\mathcal U,\qquad \sup_{t<S_*}\mathcal B_{\rm reg}(t)<e_{\rm sh}\leq e_{\rm reg}/2,\qquad \sup_{t<S_*,\,x,v}\langle\widetilde v\rangle^8|H/W|\leq\delta_*/8.\]
	Fix $S_*<S<T$, and take the common restart length $h_0\in(0,T-S_*)$ from Proposition~\ref{thm:B-restart}. Restart at $S_*-h_0<t_0<S_*$. By uniqueness the solutions agree on the overlap.
\end{proof}

\subsection{Selection by Brouwer's fixed point theorem}\label{subsec:selection-brouwer}

\begin{proof}[Proof of Theorem~\ref{thm:coupled-shooting}]
	By Lemma~\ref{lem:selection-energy-validity}, \eqref{eq:new-closed-inner}, \eqref{eq:augmented-global-inequality} and \eqref{eq:new-actual-mode-bound} hold under \eqref{eq:shooting-bootstrap}. For $k_*>0$ put
	\begin{equation}\label{eq:selection-exit-functional}
	 \mathcal Q_{\rm sh}(t)=(T-t)^{-d_*}\big(|\mathbf a(t)|_P^2-k_*\mathcal E(t)\big),
	 \qquad g_*:=d_*-2\kappa>0 .
	\end{equation}
	For $J>0$, differentiating the coordinate term in \eqref{eq:selection-exit-functional} and using \eqref{eq:new-actual-mode-bound}, \eqref{eq:unstable-coordinate-lower-bound} and $2\kappa\geq d_*-\|P\|^{-1}$, we obtain
	\begin{align*}
		\frac\dd{\dd t}\bigl((T-t)^{-d_*}|\mathbf a|_P^2\bigr)
		&=\frac{d_*}{(T-t)^{1+d_*}}|\mathbf a|_P^2+\frac{2}{(T-t)^{1+d_*}}\langle A_{\rm u}\mathbf a,\mathbf a\rangle_P+2(T-t)^{-d_*}\langle F_{\rm u},\mathbf a\rangle_P\\
		&\geq\frac{g_*}{(T-t)^{1+d_*}}|\mathbf a|_P^2-2(T-t)^{-d_*}|F_{\rm u}|_P|\mathbf a|_P.
	\end{align*}
	By \eqref{eq:new-actual-mode-bound} and Young's inequality, we have
	\begin{align*}
		C(T-t)^{(-1+\eta)/2-d_*}|\mathbf a|_P\mathcal D_a^{1/2}&\leq\frac{g_*}{4(T-t)^{1+d_*}}|\mathbf a|_P^2+C_g(T-t)^{\eta-d_*}\mathcal D_a,\\
		C(T-t)^{-1-d_*}|\mathbf a|_P\mathcal E&\leq\frac{g_*}{4(T-t)^{1+d_*}}|\mathbf a|_P^2+C_g(T-t)^{-1-d_*}\mathcal E^2.
	\end{align*}
	As a consequence, we have
	\begin{equation}\label{eq:mode-barrier-derivative}
		\frac\dd{\dd t}\bigl((T-t)^{-d_*}|\mathbf a|_P^2\bigr)\geq\frac{g_*}{2(T-t)^{1+d_*}}|\mathbf a|_P^2-C_g(T-t)^{\eta-d_*}\mathcal D_a-C_g(T-t)^{-1-d_*}\mathcal E^2.
	\end{equation}
	For $J=0$, the left-hand side of \eqref{eq:mode-barrier-derivative} is zero and its right-hand side is nonpositive, and hence this inequality remains valid. By the energy inequality \eqref{eq:new-closed-inner}, we have
	\begin{equation}\label{eq:weighted-E-barrier}
		\begin{aligned}
			\frac\dd{\dd t}\bigl((T-t)^{-d_*}\mathcal E\bigr)
			\leq{}&-\frac{\lambda-d_*}{(T-t)^{1+d_*}}\mathcal E-\frac c\eps(T-t)^{-d_*}\mathcal D_a\\
			&+\frac C{(T-t)^{1+d_*}}|\mathbf a|_P^2+C\eps(T-t)^{-1+\eta-d_*}.
		\end{aligned}
	\end{equation}
	Choose $k_*,e_1,\eps_{\rm sh}>0$ in this order so that
	\begin{align*}
		k_*C&\leq\min\{g_*/4,(\lambda-d_*)/2\},\qquad k_*C_I\leq1/4,\\
		e_1&\leq\min\{e_{\rm sh}/2,k_*(\lambda-d_*)/(2C_g)\},\\
		\eps_{\rm sh}&\leq\eps_{\rm vel},\qquad C_g\eps_{\rm sh}T^\eta\leq k_*c/2 .
	\end{align*}
	If $\eps\leq\eps_{\rm sh}$ and $\mathcal E+\mathcal B_{\rm reg}+|\mathbf a|_P^2\leq e_1$, then
	\begin{equation}\label{eq:selection-barrier-absorption}
	 \begin{aligned}
	  C_g\mathcal E^2\leq\frac{k_*(\lambda-d_*)}{2}\mathcal E,\quad
	  \frac{k_*c}{\eps}-C_g(T-t)^\eta\geq\frac{k_*c}{2\eps}.
	 \end{aligned}
	\end{equation}
	Subtract $k_*$ times \eqref{eq:weighted-E-barrier} from \eqref{eq:mode-barrier-derivative}. Applying \eqref{eq:selection-barrier-absorption} and $k_*C\leq g_*/4$, we obtain
	\begin{equation}\label{eq:outward-barrier}
		\begin{aligned}
			\mathcal Q_{\rm sh}'&\geq\frac{g_*}{4(T-t)^{1+d_*}}|\mathbf a|_P^2
				+\frac{k_*(\lambda-d_*)}{2(T-t)^{1+d_*}}\mathcal E\\
			&\quad+\frac{k_*c}{2\eps}(T-t)^{-d_*}\mathcal D_a-k_*C\eps(T-t)^{-1+\eta-d_*}.
		\end{aligned}
	\end{equation}
	Take $R_{\rm sh}>4k_*C/g_*+1$ independently of $S$ and set
	\begin{equation}\label{eq:selection-exit-level}
		L_{\rm sh}=R_{\rm sh}\eps T^{\eta-d_*}.
	\end{equation}
	At $\mathcal Q_{\rm sh}=L_{\rm sh}$ we have $(T-t)^{-d_*}|\mathbf a|_P^2\geq L_{\rm sh}$. Hence \eqref{eq:outward-barrier} gives
	\begin{equation}\label{eq:strict-outward}
		(T-t)\mathcal Q_{\rm sh}'\geq\tfrac14g_*L_{\rm sh}-k_*C\eps(T-t)^{\eta-d_*}
		\geq(\tfrac14 g_*R_{\rm sh}-k_*C)\eps T^{\eta-d_*}>0 .
	\end{equation}

	For $J=0$, take $\mathcal A_\eps=\{0\}$. For $J\geq1$, \eqref{eq:shooting-initial} and $k_*C_I\leq1/4$ give
	\begin{equation}\label{eq:initial-barrier-size}
	 \frac34T^{-d_*}|\mathbf a_0|_P^2
	 \leq\mathcal Q_{\rm sh}(0;\mathbf a_0)=T^{-d_*}\bigl(|\mathbf a_0|_P^2-k_*\mathcal E(0;\mathbf a_0)\bigr)
	 \leq T^{-d_*}|\mathbf a_0|_P^2 .
	\end{equation}
	For $|\mathbf e|_P=1$ and $0<\varrho\leq R_{\rm loc}$,
	\begin{equation}\label{eq:shooting-initial-radial-derivative}
	 \frac\dd{\dd\varrho}\mathcal Q_{\rm sh}(0;\varrho\mathbf e)
	 =T^{-d_*}\left(2\varrho-k_*\langle\nabla_{\mathbf a_0}\mathcal E(0;\varrho\mathbf e),\mathbf e\rangle_P\right)
	 \geq\frac74T^{-d_*}\varrho>0.
	\end{equation}
	Decrease $\eps_{\rm sh}$ so that
	\[R_{\rm sh}\eps_{\rm sh}T^\eta<\tfrac34R_{\rm loc}^2,\qquad \tfrac43(C_I+1)R_{\rm sh}\eps_{\rm sh}T^\eta\leq e_1/2.\]
	Define
	\begin{equation}\label{eq:selection-parameter-ball}
		\mathcal A_\eps:=\{\mathbf a_0:|\mathbf a_0|_P\leq R_{\rm loc},\ \mathcal Q_{\rm sh}(0;\mathbf a_0)\leq L_{\rm sh}\}.
	\end{equation}
	For $|\mathbf e|_P=1$, \eqref{eq:initial-barrier-size} gives
	\[\mathcal Q_{\rm sh}(0;\mathbf 0)=0<L_{\rm sh}<\tfrac34T^{-d_*}R_{\rm loc}^2\leq\mathcal Q_{\rm sh}(0;R_{\rm loc}\mathbf e).\]
	By \eqref{eq:shooting-initial-radial-derivative} there is a unique root $\rho_\eps(\mathbf e)\in(0,R_{\rm loc})$ continuous in $\mathbf e$, with
	\[\mathcal Q_{\rm sh}(0;\rho_\eps(\mathbf e)\mathbf e)=L_{\rm sh},\qquad
		L_{\rm sh}T^{d_*}\leq\rho_\eps(\mathbf e)^2\leq\tfrac43L_{\rm sh}T^{d_*}<R_{\rm loc}^2.\]
	Let $\mathbb B_P=\{\xi:|\xi|_P\leq1\}$. By \eqref{eq:selection-parameter-ball}, the radial map from $\mathbb B_P$ onto $\mathcal A_\eps$ is
	\begin{equation}\label{eq:selection-radial-homeomorphism}
		\Psi_\eps(0)=0,\qquad \Psi_\eps(s\mathbf e)=s\rho_\eps(\mathbf e)\mathbf e,\qquad 0<s\leq1,\quad|\mathbf e|_P=1 .
	\end{equation}
	The inverse $\Psi_\eps^{-1}(\varrho\mathbf e)=\varrho\mathbf e/\rho_\eps(\mathbf e)$ extends continuously by zero since $\rho_\eps^2\geq L_{\rm sh}T^{d_*}$. Thus $\Psi_\eps$ is a homeomorphism with $\Psi_\eps(\partial\mathbb B_P)=\partial\mathcal A_\eps$. On $\mathcal A_\eps$, \eqref{eq:initial-barrier-size} and \eqref{eq:shooting-initial} give
	\begin{equation}\label{eq:selection-initial-ball-bounds}
		|\mathbf a_0|_P^2\leq\tfrac43 L_{\rm sh}T^{d_*},\qquad
		T^{-d_*}\mathcal E(0;\mathbf a_0)\leq\tfrac43 C_IL_{\rm sh},\qquad
		\mathcal B_{\rm reg}(0;\mathbf a_0)\leq\tfrac43 C_IL_{\rm sh}T^{d_*}.
	\end{equation}
	These bounds include $\mathbf a_0=0$ when $J=0$. Also,
	\[\mathcal E(0;\mathbf a_0)+\mathcal B_{\rm reg}(0;\mathbf a_0)+|\mathbf a_0|_P^2\leq(C_I+1)|\mathbf a_0|_P^2\leq\tfrac43(C_I+1)R_{\rm sh}\eps T^\eta\leq e_1/2.\]

	For $\mathbf a_0\in\mathcal A_\eps$, use $(T-t)^{-d_*}|\mathbf a|_P^2\leq L_{\rm sh}+k_*(T-t)^{-d_*}\mathcal E$ in \eqref{eq:weighted-E-barrier} while \eqref{eq:shooting-bootstrap} and $\mathcal Q_{\rm sh}\leq L_{\rm sh}$ hold. Since $k_*C\leq(\lambda-d_*)/2$,
	\[\frac\dd{\dd t}\bigl((T-t)^{-d_*}\mathcal E\bigr)+\frac{\lambda-d_*}{2(T-t)}(T-t)^{-d_*}\mathcal E
	\leq\frac{CL_{\rm sh}+C\eps T^{\eta-d_*}}{T-t}.\]
	Multiplication by $(T-t)^{-(\lambda-d_*)/2}$ gives
	\begin{align*}
		(T-t)^{-d_*}\mathcal E(t)
		&\leq \left(\frac{T-t}{T}\right)^{(\lambda-d_*)/2}T^{-d_*}\mathcal E(0)
		+\frac{2C}{\lambda-d_*}\bigl(L_{\rm sh}+\eps T^{\eta-d_*}\bigr)
		\leq CL_{\rm sh}.
	\end{align*}
	The last bound uses \eqref{eq:selection-exit-level} and \eqref{eq:selection-initial-ball-bounds}. With \eqref{eq:selection-exit-functional} this gives
	\begin{equation}\label{eq:barrier-improves-bootstrap}
		\mathcal E(t)+|\mathbf a(t)|_P^2\leq C_3L_{\rm sh}(T-t)^{d_*}.
	\end{equation}
	Using \eqref{eq:barrier-improves-bootstrap}, we integrate \eqref{eq:augmented-global-inequality} against $\mathrm e^{-Ct}$. The bounds
	\[\int_0^t\frac{\mathcal E(q)}{T-q}\dd q\leq C_3L_{\rm sh}\frac{T^{d_*}-(T-t)^{d_*}}{d_*},\qquad \int_0^t\eps(T-q)^{\eta-1}\dd q=\frac\eps\eta\bigl(T^\eta-(T-t)^\eta\bigr)\]
	give
	\begin{equation}\label{eq:selection-global-bound}
		\sup_{0\leq q\leq t}\mathcal B_{\rm reg}(q)+\frac{c_{\rm reg}}\eps\int_0^t\mathcal D_{\rm reg}(q)\dd q\leq C_T\left(\mathcal B_{\rm reg}(0)+\frac{L_{\rm sh}T^{d_*}}{d_*}+\frac{\eps T^\eta}{\eta}\right)\leq C_4\eps T^\eta.
	\end{equation}
	Reduce $\eps_{\rm sh}$ so that $(C_3R_{\rm sh}+C_4)\eps_{\rm sh}T^\eta\leq e_1/2$. By \eqref{eq:selection-exit-level}, \eqref{eq:barrier-improves-bootstrap} and \eqref{eq:selection-global-bound},
	\begin{equation}\label{eq:selection-strict-energy}
		\mathcal E(t)+\mathcal B_{\rm reg}(t)+|\mathbf a(t)|_P^2\leq(C_3R_{\rm sh}+C_4)\eps T^\eta\leq e_1/2 .
	\end{equation}
	While $\mathcal Q_{\rm sh}\leq L_{\rm sh}$, \eqref{eq:selection-strict-energy} gives continuation below $T$ by Lemma~\ref{lem:selection-continuation}.

	If $J=0$, then $\mathcal Q_{\rm sh}=-k_*(T-t)^{-d_*}\mathcal E\leq0<L_{\rm sh}$. Thus \eqref{eq:selection-strict-energy} continues the solution up to $T$. For $J\geq1$ fix $0<S<T$. Suppose all $\mathbf a_0\in\mathcal A_\eps$ hit the exit level before $S$. Put
	\begin{equation}\label{eq:selection-first-exit-time}
		t_{\rm exit}(\mathbf a_0):=\inf\{t\geq0:\mathcal Q_{\rm sh}(t;\mathbf a_0)=L_{\rm sh}\}<S,\qquad t_{\rm exit}|_{\partial\mathcal A_\eps}=0 .
	\end{equation}
	At $t_0=t_{\rm exit}(\mathbf a_0)$ the sum in \eqref{eq:selection-strict-energy} is at most $e_1/2$, so the solution continues past $t_0$ with sum below $e_1$. If $t_0>0$, \eqref{eq:strict-outward} and the first hitting property give, for small $h>0$,
	\[\max_{0\leq t\leq t_0-h}\mathcal Q_{\rm sh}(t;\mathbf a_0)<L_{\rm sh}<\mathcal Q_{\rm sh}(t_0+h;\mathbf a_0).\]
	By \eqref{eq:selection-continuous-dependence} the same inequalities hold for nearby $\mathbf b_0$, so $|t_{\rm exit}(\mathbf b_0)-t_0|<h$. If $t_0=0$, then $\mathcal Q_{\rm sh}(h;\mathbf b_0)>L_{\rm sh}$ for small $h>0$ and nearby $\mathbf b_0$, giving $0\leq t_{\rm exit}(\mathbf b_0)<h$. Hence $t_{\rm exit}$ is continuous.
	
	By \eqref{eq:selection-exit-functional} at $t=t_{\rm exit}(\mathbf a_0)$,
	\begin{equation}\label{eq:selection-nonzero-exit}
		|\mathbf a(t_{\rm exit};\mathbf a_0)|_P^2
		=(T-t_{\rm exit})^{d_*}L_{\rm sh}+k_*\mathcal E(t_{\rm exit};\mathbf a_0)
		\geq (T-t_{\rm exit})^{d_*}L_{\rm sh}>0 .
	\end{equation}
	The lower bound in \eqref{eq:selection-nonzero-exit} and the parametrization \eqref{eq:selection-radial-homeomorphism} allow us to set
	\begin{equation}\label{eq:selection-brouwer-map}
		\mathcal F_S(\xi):=-\frac{\mathbf a\bigl(t_{\rm exit}(\Psi_\eps(\xi));\Psi_\eps(\xi)\bigr)}{\bigl|\mathbf a\bigl(t_{\rm exit}(\Psi_\eps(\xi));\Psi_\eps(\xi)\bigr)\bigr|_P},\qquad \xi\in\mathbb B_P.
	\end{equation}
	By \eqref{eq:selection-continuous-dependence}, continuity of $t_{\rm exit}$ and \eqref{eq:selection-nonzero-exit}, this is a continuous map $\mathbb B_P\to\partial\mathbb B_P$. Equations \eqref{eq:selection-first-exit-time}, \eqref{eq:selection-radial-homeomorphism} and \eqref{eq:initial-fitted-data} give
	\[\mathcal F_S(\xi)=-\frac{\Psi_\eps(\xi)}{|\Psi_\eps(\xi)|_P}=-\xi,\qquad |\xi|_P=1.\]
	A fixed point given by Brouwer would have $|\xi|_P=1$ and $\xi=-\xi$. Consequently some $\mathbf a_0(S)\in\mathcal A_\eps$ satisfies $\mathcal Q_{\rm sh}(t;\mathbf a_0(S))\leq L_{\rm sh}$ on $[0,S]$.

	Take $S_j\uparrow T$ with $\mathbf a_0(S_j)\to\mathbf a_0^*\in\mathcal A_\eps$ along a subsequence, using $\mathbf a_0^*=0$ if $J=0$. By \eqref{eq:selection-initial-ball-bounds}, the initial sum is below $e_1$. On $[0,S]$ with the sum below $e_1$, \eqref{eq:selection-continuous-dependence} and \eqref{eq:selection-strict-energy} for $S_j>S$ give
	\[\sup_{0\leq t\leq S}\big(\mathcal E(t;\mathbf a_0^*)+\mathcal B_{\rm reg}(t;\mathbf a_0^*)+|\mathbf a(t;\mathbf a_0^*)|_P^2\big)\leq e_1/2.\]
	Apply Lemma~\ref{lem:selection-continuation} to reach every $S<T$. The limit in \eqref{eq:barrier-improves-bootstrap} gives \eqref{eq:shooting-pointwise-result}. By \eqref{eq:local-reference-difference} the dissipation integrands converge pointwise. Fatou's lemma in \eqref{eq:selection-global-bound} gives
	\[\sup_{0\leq t\leq S}\mathcal B_{\rm reg}(t;\mathbf a_0^*)+\frac{c_{\rm reg}}\eps\int_0^S\mathcal D_{\rm reg}(t;\mathbf a_0^*)\dd t\leq C_4\eps T^\eta.\]
	Let $S\uparrow T$ to obtain \eqref{eq:shooting-pointwise-result-2}. Use \eqref{eq:selection-initial-ball-bounds} for the initial size and Lemma~\ref{lem:selection-energy-validity} for the parameter range.

	For the selected solution, multiply \eqref{eq:new-closed-inner} by $(T-t)^{-\delta}$ and use \eqref{eq:barrier-improves-bootstrap}.
	\begin{equation}\label{eq:selection-delta-energy}
	 \begin{aligned}
	  \frac{\dd}{\dd t}\big((T-t)^{-\delta}\mathcal E\big)
	    +\frac{\lambda-\delta}{(T-t)^{1+\delta}}\mathcal E
	    +\frac c\eps(T-t)^{-\delta}\mathcal D_a
	   &\leq C L_{\rm sh}(T-t)^{d_*-\delta-1}+C\eps(T-t)^{\eta-\delta-1}.
	 \end{aligned}
	\end{equation}
	By \eqref{eq:shooting-exponents}, $\delta<d_*<\min\{\eta,\lambda\}$. Integrating \eqref{eq:selection-delta-energy} up to $S<T$ and letting $S\uparrow T$ yields
	\[\begin{aligned}
		&\sup_{0\leq t<T}(T-t)^{-\delta}\mathcal E(t)+\frac c\eps\int_0^T(T-t)^{-\delta}\mathcal D_a(t)\dd t+\int_0^T(T-t)^{-1-\delta}\mathcal E(t)\dd t\\
		&\qquad\leq C\left(T^{-\delta}\mathcal E(0)+\frac{L_{\rm sh}T^{d_*-\delta}}{d_*-\delta}+\frac{\eps T^{\eta-\delta}}{\eta-\delta}\right)\leq C\eps T^{\eta-\delta}.
	\end{aligned}\]
	This proves \eqref{eq:shooting-weighted-result}.
	
	If $J=0$, the proof is complete. For $J>0$, by \eqref{eq:unstable-coordinate-lower-bound}, we have $2\langle A_{\rm u}\xi,\xi\rangle_P\geq-2\kappa|\xi|_P^2$. Hence, we have
	\[\frac\dd{\dd q}|\mathcal U_{\rm u}(t,q)\xi|_P^2\leq\frac{2\kappa}{T-q}|\mathcal U_{\rm u}(t,q)\xi|_P^2,\qquad \|\mathcal U_{\rm u}(t,q)\|_P\leq\left(\frac{T-t}{T-q}\right)^\kappa\quad(t\leq q<T).\]
	By the variation of constants on $[t,S]$, we obtain
	\[\mathbf a(t)=\mathcal U_{\rm u}(t,S)\mathbf a(S)-\int_t^S\mathcal U_{\rm u}(t,q)F_{\rm u}(q)\dd q.\]
	The bound \eqref{eq:barrier-improves-bootstrap} implies
	\[|\mathcal U_{\rm u}(t,S)\mathbf a(S)|_P\leq C L_{\rm sh}^{1/2}(T-t)^\kappa(T-S)^{d_*/2-\kappa}\longrightarrow0\quad(S\uparrow T),\]
	since $d_*>2\kappa$. Moreover, by \eqref{eq:new-actual-mode-bound} and \eqref{eq:barrier-improves-bootstrap}, we have
	\[\begin{aligned}
		\int_0^T(T-q)^{1-\delta}|F_{\rm u}(q)|_P^2\dd q
		&\leq CT^\eta\int_0^T(T-q)^{-\delta}\mathcal D_a(q)\dd q+CL_{\rm sh}^2\int_0^T(T-q)^{2d_*-\delta-1}\dd q\\
		&\leq C\eps^2T^{2\eta-\delta}+\frac{CL_{\rm sh}^2T^{2d_*-\delta}}{2d_*-\delta}<\infty.
	\end{aligned}\]
	Since $\delta>2\kappa$, by the Cauchy--Schwarz inequality, we have
	\[\begin{aligned}
		\int_t^T|\mathcal U_{\rm u}(t,q)F_{\rm u}(q)|_P\dd q
		&\leq(T-t)^\kappa\left(\int_t^T(T-q)^{\delta-2\kappa-1}\dd q\right)^{1/2}\left(\int_t^T(T-q)^{1-\delta}|F_{\rm u}(q)|_P^2\dd q\right)^{1/2}\\
		&=\frac{(T-t)^{\delta/2}}{\sqrt{\delta-2\kappa}}\left(\int_t^T(T-q)^{1-\delta}|F_{\rm u}(q)|_P^2\dd q\right)^{1/2}<\infty.
	\end{aligned}\]
	Thus the backward integral is absolutely convergent. Let $S\uparrow T$ in the variation of constants to obtain \eqref{eq:selected-terminal-integral}.

\end{proof}

\appendix
\section{A polynomially weighted local construction near a local Maxwellian}\label{app:polynomial-local}

\subsection{Frozen reference and pointwise collision bounds}
Fix smooth reference parameters $A_{\rm ref},u_{\rm ref},\theta_{\rm ref},m$ depending on $x$ alone. Use physical velocity $v$ in
\begin{align}\label{eq:B-reference}
	M_{\rm ref}(x,v)&:=A_{\rm ref}(x)\mu(\widetilde v_{\rm ref}),\qquad \widetilde v_{\rm ref}:=(v-u_{\rm ref}(x))/\theta_{\rm ref}(x),\notag\\
	W(x,v)&:=A_{\rm ref}(x)^{1/2}w_{m(x),\ell}(\widetilde v_{\rm ref}),\qquad W_q(x,v):=W(x,v)\langle \widetilde v_{\rm ref}\rangle^{-q}.
\end{align}
The mass is $\int M_{\rm ref}(x,v)\dd v=A_{\rm ref}(x)\theta_{\rm ref}(x)^3$. Assume
\[
 \begin{gathered}
  0<A_{{\rm ref},-}\leq A_{\rm ref}\leq A_{{\rm ref},+},\qquad
  0<\theta_{{\rm ref},-}\leq\theta_{\rm ref}\leq\theta_{{\rm ref},+},\\
  |u_{\rm ref}|\leq u_{{\rm ref},+},\qquad m_0\leq m\leq m_+ .
 \end{gathered}
\]
Assume bounded spatial derivatives of $A_{\rm ref},u_{\rm ref},\theta_{\rm ref}$ at all orders and $\|\nabla_xm\|_\infty<\infty$. The local time uses the displayed bounds and $\|\nabla_x(A_{\rm ref},u_{\rm ref},\theta_{\rm ref},m)\|_\infty$.

For the remainder $h=f-M_{\rm ref}$, set
\begin{equation}\label{eq:B-reference-norms}
 |h|_q:=\sup_{x,v}\frac{|h(x,v)|}{W_q(x,v)},\qquad
 |h|_{q,N}:=\sum_{|\alpha|+|\beta|\leq N}|\partial_x^\alpha\partial_v^\beta h|_q .
\end{equation}
Use $(q_0,q_1)=(8,2K+8)$ at all orders $N$. Choose $\ell,m_0$ for the collision estimates over $q\in[0,q_1]$. The norm with spatial decay is
\begin{equation}\label{eq:B-spatial-reference-norms}
	|h|_{q,N,m_x}:=\sum_{|\alpha|+|\beta|\leq N} \sup_{x,v}\langle x\rangle^{m_x} |\partial_x^\alpha\partial_v^\beta h(x,v)|/W_q(x,v), \qquad m_x\geq0.
\end{equation}

For $\nu[h](x,v):=\int|v-v_*|^\gamma h(x,v_*)\dd v_*$, set
\begin{align*}
	\nu_{\rm ref}&:=\nu[M_{\rm ref}],\\
	\mathcal K_{\rm ref}h&:=Q^+(M_{\rm ref},h)+Q^+(h,M_{\rm ref})-M_{\rm ref}\nu[h],\\
	Q(M_{\rm ref},h)+Q(h,M_{\rm ref})&=-\nu_{\rm ref}h+\mathcal K_{\rm ref}h .
\end{align*}

\begin{lemma}\label{lem:B-rows}
	For $\delta_{\rm loc}>0$ small and $q\in\{q_0,q_1\}$, take $\ell,m_0$ sufficiently large and fixed. The kernel $k_q$ of $W_q^{-1}\mathcal K_{\rm ref}(W_q\cdot)$ satisfies
	\begin{equation}\label{eq:B-linear-row}
		\int|k_q(x,v,v')|\dd v'\leq\delta_{\rm loc}\nu_{\rm ref}(x,v)+C_{\rm ref}.
	\end{equation}
	Equation \eqref{eq:B-linear-row} also holds for the kernel of $W_q^{-1}[Q^+(M_{\rm ref},W_q\cdot)+Q^+(W_q\cdot,M_{\rm ref})]$. The constant $C_{\rm ref}$ uses only the bounds after \eqref{eq:B-reference}. For $Q$ and $Q^+$, we have
	\begin{align}\label{eq:B-nonlinear-row}
		\frac{|Q(g,h)(x,v)|}{W_q(x,v)}&\leq C\nu_{\rm ref}(x,v)|g|_q|h|_q,\notag\\
		\frac{|Q(g,h)(x,v)|}{W_{q_1}(x,v)}&\leq C_{q_1}\nu_{\rm ref}(x,v)(|g|_{q_0}|h|_{q_1}+|g|_{q_1}|h|_{q_0}).
	\end{align}
	The constants in \eqref{eq:B-nonlinear-row} depend only on $A_{{\rm ref},-},A_{{\rm ref},+},\ell,q_1,\gamma$. We use $C_{q_1}\geq C$ in the two bounds.
\end{lemma}
\begin{proof}
	At fixed $x$ use normalized $v$, $w_q=w_{m,\ell}\langle v\rangle^{-q}$ and $\nu=\nu[\mu]$. For the conjugated gain $k_q^+$ and the starred loss, \eqref{eq:polynomial-absolute-row} and \eqref{eq:new-moments} give
	\begin{equation}\label{eq:rwA-normalized-rows}
		\begin{aligned}
			\nu(v)^{-1}\int k_q^+(v,v')\dd v'
			&\leq C_{q_1}(\ell^{-1}+|v|^{-2}+\ell^{-1/2}|v|^{-1})\qquad(|v|\geq\sqrt\ell),\\
			\frac{\mu(v)}{w_q(v)\nu(v)}\nu[w_q](v)&\leq C\sqrt\mu(v)\langle v\rangle^q.
		\end{aligned}
	\end{equation}
	Choose $C_{q_1}\ell^{-1}\leq\delta_{\rm loc}/4$, then $R_{\rm row}$ large so the sum in \eqref{eq:rwA-normalized-rows} is at most $\delta_{\rm loc}$ for $|v|\geq R_{\rm row}$. Inside the ball, use $|a'|\geq\rho-R_{\rm row}$ in \eqref{eq:new-grad-kernel} and \eqref{eq:polynomial-gain-ratio}. Then
	\[\sup_{|v|\leq R_{\rm row}}\int k_q^+(v,v')\dd v'\leq C_{\ell,q_1,R_{\rm row}}\left(\int_0^1\rho\dd\rho+\int_1^\infty(1+\rho)^\gamma(1+\rho^2)^{-\ell/4}\dd\rho\right)<\infty .\]
	For the starred loss on this ball use \eqref{eq:rwA-normalized-rows}.

	Take $L>\gamma+5$ and $\ell\geq2L$ in \eqref{eq:two-channel-new}. The first channel satisfies
	\[\int |v-v_*|^\gamma\langle v_*\rangle^{-L}\dd v_*\leq C\langle v\rangle^\gamma\int\langle v_*\rangle^{\gamma-L}\dd v_*\leq C_L\nu(v).\]
	For the other channel write $v=v_\perp+ae$, $v_*=y+se$ with $v_\perp,y\in e^\perp$. By \eqref{eq:gain-collision-coordinates}, $\rho=a-s>0$ and
	\begin{align*}
		v'&=v_\perp+se,\qquad v_*'=y+ae,\qquad \dd v_*=\dd y\dd s,\\
		B_\gamma\dd\omega&=\pi^{-1}|v-v_*|^{\gamma-1}(a-s)\dd e\leq\pi^{-1}|v-v_*|^\gamma\dd e .
	\end{align*}
	With $A=1+|v_\perp|^2$ and $B=1+a^2$, $\langle v'\rangle^2=A+s^2$ and $\langle v_*'\rangle^2=B+|y|^2$. Use $|v-v_*|^2\leq2(|v'|^2+|v_*'|^2)$ to get
	\[|v-v_*|^\gamma\leq C_\gamma\sum_{\zeta\in\{0,\gamma\}}(A+s^2)^{\zeta/2}(B+|y|^2)^{(\gamma-\zeta)/2}.\]
	Extend to $s\in\R$ and put $s=\sqrt A\,s_1$, $y=\sqrt B\,y_1$. Then
	\begin{equation}\label{eq:pointwise-gain-plane-integral}
		\begin{aligned}
			&\int_{\R}\int_{e^\perp}
				\frac{\langle v\rangle^L|v-v_*|^\gamma}{(A+s^2)^{L/2}(B+|y|^2)^{L/2}}\dd y\dd s\\
			&\qquad\leq C_L\sum_{\zeta\in\{0,\gamma\}}\langle v\rangle^L A^{(1+\zeta-L)/2}B^{(2+\gamma-\zeta-L)/2}.
		\end{aligned}
	\end{equation}
	For $R=|v|\geq1$, $\tau=|a|/R$ has angular measure $4\pi\dd\tau$ on $[0,1]$. Use the two substitutions
	\begin{align*}
		0\leq\tau\leq\tfrac12&:\quad A\asymp R^2,\quad u=R\tau,\quad \dd\tau=R^{-1}\dd u,\\
		\tfrac12\leq\tau\leq1&:\quad B\asymp R^2,\quad u=R\sqrt{1-\tau^2},\quad |\dd\tau|\leq2uR^{-2}\dd u.
	\end{align*}
	Then
	\begin{align*}
		&\int_{\mathbb S^2}\langle v\rangle^L A^{(1+\zeta-L)/2}B^{(2+\gamma-\zeta-L)/2}\dd e\\
		\leq&\, C_L R^\zeta\int_0^\infty(1+u^2)^{(2+\gamma-\zeta-L)/2}\dd u+C_L R^{\gamma-\zeta}\int_0^\infty u(1+u^2)^{(1+\zeta-L)/2}\dd u\\
		\leq&\, C_L\nu(v).
	\end{align*}
	For $R\leq1$, use $1\leq A,B\leq2$ in \eqref{eq:pointwise-gain-plane-integral}. Together with the integral of $|v-v_*|^\gamma\langle v_*\rangle^{-L}$ this proves
	\begin{equation}\label{eq:appendix-bilinear-row}
		\iint B_\gamma(v-v_*)R_w\dd\omega\dd v_*\leq C_L\nu(v).
	\end{equation}
	
	The inequality $\langle v\rangle\leq\langle v'\rangle\langle v_*'\rangle$ gives $w_q(v')w_q(v_*')/w_q(v)\leq R_w$. By \eqref{eq:appendix-bilinear-row} and \eqref{eq:new-moments},
	\begin{equation}\label{eq:rwA-bilinear-rows}
		\begin{aligned}
			\iint B_\gamma \frac{w_q(v')w_q(v_*')}{w_q(v)}\dd\omega\dd v_*&\leq C_\ell\nu(v),\\
			\nu[w_q](v)&\leq C_\ell\nu(v).
		\end{aligned}
	\end{equation}
	For $r_q=q_1-q_0$, we have
	\begin{align*}
		\frac{\langle v\rangle^{r_q}|g_*'h'|}{w_{q_0}(v)}
		\leq C_{q_1}\frac{|\langle v_*'\rangle^{r_q}g_*'h'|+|g_*'\langle v'\rangle^{r_q}h'|}{w_{q_0}(v)},\quad
		\langle v\rangle^{r_q}Q^-(g,h)=Q^-(g,\langle v\rangle^{r_q}h).
	\end{align*}
	Apply \eqref{eq:rwA-bilinear-rows} at $q=q_0$ to the right sides. For physical $v$ and $z=\widetilde v_{\rm ref}$, set $\widetilde g(z)=A_{\rm ref}^{-1/2}g(u_{\rm ref}+\theta_{\rm ref}z)$ and define $\widetilde h$ in the same way. By \eqref{Eq.collision-operator},
	\begin{align}\label{eq:reference-frequency-scaling}
		\nu_{\rm ref}(v)&=A_{\rm ref}\theta_{\rm ref}^{3+\gamma}\nu(z),\\
		\frac{|Q(g,h)(v)|}{W_q(v)\nu_{\rm ref}(v)}
		&=A_{\rm ref}^{-1/2}\frac{|Q(\widetilde g,\widetilde h)(z)|}{w_q(z)\nu(z)}.\notag
	\end{align}
	The normalized linear row is multiplied by $A_{\rm ref}\theta_{\rm ref}^{3+\gamma}$, so \eqref{eq:rwA-normalized-rows} gives \eqref{eq:B-linear-row}. The nonlinear estimate \eqref{eq:B-nonlinear-row} follows from the second line of \eqref{eq:reference-frequency-scaling}.
\end{proof}

\subsection{Damped characteristics and local construction}

	Set $a_{\rm col}=\eps^{-1}\nu_{\rm ref}$, $b_q=v\cdot\nabla_x\log W_q$ and $E=|\widetilde v_{\rm ref}|^2$. At fixed $v$, \eqref{eq:B-reference} gives, with $\Phi=\Phi_{m,\ell}(E)$,
	\begin{align*}
		\nabla_x\log W_q
		&=\tfrac12\nabla_x\log A_{\rm ref}-\tfrac14\partial_m\Phi\,\nabla_xm
		-\left(\tfrac14\partial_E\Phi+\tfrac q{2(1+E)}\right)\nabla_x E,\\
		\nabla_x E&=-\frac{2}{\theta_{\rm ref}}\left(\sum_i\widetilde v_{{\rm ref},i}\nabla_xu_{{\rm ref},i}+E\nabla_x\theta_{\rm ref}\right).
	\end{align*}
	By \eqref{eq:Phi-derivatives}, $|\partial_m\Phi|\leq C$ and $(1+E)\partial_E\Phi\leq1+m_+$. Therefore $|\nabla_x\log W_q|\leq C_{\rm ref}$. By \eqref{eq:reference-frequency-scaling} and $\gamma>1$,
	\begin{equation}\label{eq:B-transport-split}
		|b_q|\leq C_{\rm ref}\langle v\rangle\leq\delta_{\rm loc}a_{\rm col}+C_{M_{\rm ref},\eps,\delta_{\rm loc}},\qquad q=q_0,q_1 .
	\end{equation}
	For $0\leq s\leq t\leq T_{\rm loc}$ let $X(s)=x-(t-s)v$ and set
	\[E_q(t,s;x,v):=\exp\left[-\int_s^t(a_{\rm col}+b_q)(X(\tau),v)\dd\tau\right].\]
	By \eqref{eq:B-transport-split},
	\begin{align}\label{eq:B-characteristic}
		E_q(t,s)&\leq \mathrm e^{CT_{\rm loc}}\exp\left[-(1-\delta_{\rm loc})\int_s^t a_{\rm col}(X(\tau),v)\dd\tau\right],\notag\\
		\int_0^t E_q(t,s)a_{\rm col}(X(s),v)\dd s
		&\leq\frac{\mathrm e^{CT_{\rm loc}}}{1-\delta_{\rm loc}}
			\left(1-\exp\left[-(1-\delta_{\rm loc})\int_0^t a_{\rm col}(X(\tau),v)\dd\tau\right]\right),\notag\\
		\int_0^t E_q(t,s)\dd s&\leq \mathrm e^{CT_{\rm loc}}T_{\rm loc}.
	\end{align}
	The second line is the integral of the $s$ derivative of the exponential in the first line. For $(\partial_t+v\cdot\nabla_x+a_{\rm col})G=F$, integrate along $X(s)$ to get
	\begin{equation}\label{eq:local-damped-characteristic-formula}
		\frac{G(t,x,v)}{W_q(x,v)}=E_q(t,0)\frac{G(0,X(0),v)}{W_q(X(0),v)}
			+\int_0^t E_q(t,s)\frac{F(s,X(s),v)}{W_q(X(s),v)}\dd s .
	\end{equation}

\begin{proposition}\label{thm:B-local}
	Let $\delta_*>0$ be sufficiently small in terms of the fixed amplitude bounds and $\ell,q_1$. For smooth $h_0$ with $f_0=M_{\rm ref}+h_0\geq0$, assume
	\[|h_0|_{q_0}\leq\delta_*,\qquad |h_0|_{q_1,N}<\infty\quad(N\geq0).\]
	For some $T_{\rm loc}>0$, \eqref{Eq.Boltzmann} has a unique classical solution $f=M_{\rm ref}+h$ in the class
	\[f(0)=f_0,\qquad f\geq0,\qquad \sup_{t\leq T_{\rm loc}}|h(t)|_{q_0}\leq2\delta_*,\qquad \sup_{t\leq T_{\rm loc}}|h(t)|_{q_1,N}<\infty\quad(N\geq0).\]
	The time $T_{\rm loc}$ uses $\eps,\delta_*$ and the $C^1_x$ bounds of the reference parameters. The solution is smooth in $(x,v)$ with $f,\partial_tf$ continuous. If $f_0>0$, then $f>0$. Differences for bounded data in $|\cdot|_{q_1,N}$ with a common reference satisfy \eqref{eq:local-reference-difference}.
\end{proposition}
\begin{proof}
	Starting at $h^0=0$, solve along \eqref{eq:local-damped-characteristic-formula}
	\begin{align}\label{eq:B-iteration}
		(D+a_{\rm col})h^{n+1}
		=\eps^{-1}(\mathcal K_{\rm ref}h^n+Q(h^n,h^n))-v\cdot\nabla_xM_{\rm ref},\quad
		h^{n+1}(0)=h_0 .
	\end{align}
	Set $Z_n=\sup_{t\leq T_{\rm loc}}|h^n(t)|_{q_0}$ and $C_{\rm src}=|v\cdot\nabla_x M_{\rm ref}|_{q_0}$. The latter is bounded by the Gaussian moments and the first spatial derivatives of the reference. The right-hand side of \eqref{eq:B-iteration} divided by $W_{q_0}$ is bounded in absolute value by
	\[a_{\rm col}(\delta_{\rm loc}Z_n+CZ_n^2)+\eps^{-1}C_{\rm ref}Z_n+C_{\rm src}\]
	by \eqref{eq:B-linear-row} and \eqref{eq:B-nonlinear-row}. Integrate this in \eqref{eq:local-damped-characteristic-formula} by \eqref{eq:B-characteristic}.
	\begin{equation}\label{eq:rwA-low-iteration}
		\begin{aligned}
			Z_{n+1}&\leq\mathrm e^{CT_{\rm loc}}\left[|h_0|_{q_0}+\frac{\delta_{\rm loc}Z_n+CZ_n^2}{1-\delta_{\rm loc}}
			+T_{\rm loc}(\eps^{-1}C_{\rm ref}Z_n+C_{\rm src})\right].
		\end{aligned}
	\end{equation}
	With $\delta_{\rm loc}$ small in Lemma~\ref{lem:B-rows}, choose $\delta_*$ and then $T_{\rm loc}$ by
	\begin{equation}\label{eq:local-row-smallness}
		\begin{gathered}
			\mathrm e^{CT_{\rm loc}}\left(\frac{\delta_{\rm loc}}{1-\delta_{\rm loc}}+\eps^{-1}C_{\rm ref}T_{\rm loc}\right)\leq\tfrac1{16},\qquad
			\frac{4C\delta_*\mathrm e^{CT_{\rm loc}}}{1-\delta_{\rm loc}}\leq\tfrac18,\\
			\mathrm e^{CT_{\rm loc}}\leq\tfrac98,\qquad \mathrm e^{CT_{\rm loc}}T_{\rm loc}C_{\rm src}\leq\tfrac18\delta_*,\qquad
			\frac{9(\delta_{\rm loc}+4C_{q_1}\delta_*)}{8(1-\delta_{\rm loc})}<\tfrac12,\\
			\delta_{\rm loc}+C\delta_*\leq\tfrac14,\qquad
			\mathrm e^{CT_{\rm loc}}\left(\frac{\delta_{\rm loc}+C\delta_*}{1-\delta_{\rm loc}-C\delta_*}
				+\eps^{-1}C_{\rm ref}T_{\rm loc}\right)\leq\tfrac12 .
		\end{gathered}
	\end{equation}
	Equations \eqref{eq:rwA-low-iteration} and \eqref{eq:local-row-smallness} give
	\[Z_n\leq2\delta_*\ \Longrightarrow\ Z_{n+1}\leq(9/8+1/8+1/8+1/8)\delta_*=\tfrac32\delta_*.\]
	Subtracting successive equations in \eqref{eq:B-iteration} gives, for $g^n=h^n-h^{n-1}$,
	\[(\partial_t+v\cdot\nabla_x+a_{\rm col})g^{n+1}=\eps^{-1}\big(\mathcal K_{\rm ref}g^n+Q(g^n,h^n)+Q(h^{n-1},g^n)\big),\qquad g^{n+1}(0)=0.\]
	Use $Z_n+Z_{n-1}\leq4\delta_*$ and the coefficient $1/16+1/8$ in \eqref{eq:local-row-smallness}. For $g=h-\widetilde h$ the initial term is $h_0-\widetilde h_0$. Thus
	\begin{equation}\label{eq:rwA-low-difference}
		\begin{aligned}
			\sup_{t\leq T_{\rm loc}}|h^{n+1}-h^n|_{q_0}
			&\leq\tfrac3{16}\sup_{t\leq T_{\rm loc}}|h^n-h^{n-1}|_{q_0}\qquad(n\geq1),\\
			\sup_{t\leq T_{\rm loc}}|h-\widetilde h|_{q_0}
			&\leq\tfrac98|h_0-\widetilde h_0|_{q_0}+\tfrac3{16}\sup_{t\leq T_{\rm loc}}|h-\widetilde h|_{q_0}.
		\end{aligned}
	\end{equation}
	The first line gives convergence in $|\cdot|_{q_0}$. The second gives uniqueness and $\sup_{t\leq T_{\rm loc}}|h-\widetilde h|_{q_0}\leq(18/13)|h_0-\widetilde h_0|_{q_0}$.

	For $M\geq1$, put $Z_{q_1,M}(t)=\max_{0\leq n\leq M}\sup_{s\leq t}|h^n(s)|_{q_1}$. Since $|h^n|_{q_0}\leq2\delta_*$, the collision terms in \eqref{eq:B-iteration} divided by $W_{q_1}$ are bounded by
	\[a_{\rm col}(\delta_{\rm loc}+4C_{q_1}\delta_*)|h^n|_{q_1}+\eps^{-1}C_{\rm ref}|h^n|_{q_1}.\]
	By \eqref{eq:B-characteristic},
	\begin{equation}\label{eq:local-high-norm-iteration}
		Z_{q_1,M}(t)\leq C(1+|h_0|_{q_1})+\kappa_{\rm it}Z_{q_1,M}(t)+C\int_0^tZ_{q_1,M}(s)\dd s ,
	\end{equation}
	where
	\begin{equation}\label{eq:local-iteration-coefficient}
		\kappa_{\rm it}:=\frac{\mathrm e^{CT_{\rm loc}}(\delta_{\rm loc}+4C_{q_1}\delta_*)}{1-\delta_{\rm loc}}<\tfrac12
	\end{equation}
	by \eqref{eq:local-row-smallness}. Hence $Z_{q_1,M}(t)\leq 2C(1+|h_0|_{q_1})\mathrm e^{2Ct}$, uniformly in $M$.

	Write $\mathfrak a=(\alpha,\beta)$, $D^{\mathfrak a}=\partial_x^\alpha\partial_v^\beta$ and $N=|\alpha|+|\beta|$. In this notation \eqref{eq:rw9-collision-derivatives} is
	\[
	 D^{\mathfrak a}Q(g,h)=\sum_{\mathfrak c\leq\mathfrak a}\binom{\mathfrak a}{\mathfrak c}
	                   Q(D^{\mathfrak c}g,D^{\mathfrak a-\mathfrak c}h).
	\]
	For $N\geq1$, apply $D^{\mathfrak a}$ to the equation for $h$ and use $Q(M_{\rm ref},M_{\rm ref})=0$. We get
	\begin{align}\label{eq:B-differentiated}
		(\partial_t+v\cdot\nabla_x+a_{\rm col})D^{\mathfrak a}h
		&=\eps^{-1}\bigl(\mathcal K_{\rm ref}D^{\mathfrak a}h+Q(h,D^{\mathfrak a}h)+Q(D^{\mathfrak a}h,h)\bigr)\notag\\
		&\quad+\eps^{-1}\mathcal F_{\mathfrak a}-D^{\mathfrak a}(v\cdot\nabla_xM_{\rm ref})-\sum_{\beta_l>0}\beta_lD^{\alpha+e_l,\beta-e_l}h,
	\end{align}
	where
	\begin{align*}
		\mathcal F_{\mathfrak a}
		&=\sum_{0<\mathfrak c\leq\mathfrak a}\binom{\mathfrak a}{\mathfrak c}\bigl(Q(D^{\mathfrak c}M_{\rm ref},D^{\mathfrak a-\mathfrak c}h)+Q(D^{\mathfrak a-\mathfrak c}h,D^{\mathfrak c}M_{\rm ref})\bigr)+\sum_{\substack{\mathfrak c\leq\mathfrak a\\\mathfrak c\neq0,\mathfrak a}}\binom{\mathfrak a}{\mathfrak c}Q(D^{\mathfrak c}h,D^{\mathfrak a-\mathfrak c}h).
	\end{align*}
	By \eqref{eq:B-nonlinear-row} and the Gaussian bounds on $D^{\mathfrak c}M_{\rm ref}$, for $q=q_0,q_1$,
	\begin{equation}\label{eq:rwA-derivative-sources}
		\begin{aligned}
			&\sup_{x,v}\frac{|\mathcal F_{\mathfrak a}|}{\nu_{\rm ref}W_q}
			\leq C_N\left(|h|_{q,N-1}+\sum_{l=1}^{N-1}|h|_{q,l}|h|_{q,N-l}\right),\\
			&\qquad\qquad\left|\sum_{\beta_l>0}\beta_lD^{\alpha+e_l,\beta-e_l}h\right|_q\leq N|h|_{q,N}.
		\end{aligned}
	\end{equation}
	For any $g$ and $|h|_{q_0}\leq2\delta_*$, \eqref{eq:B-nonlinear-row} also gives
	\begin{align}\label{eq:local-highest-rows}
		|Q(h,g)|+|Q(g,h)|&\leq 4C\delta_*\nu_{\rm ref}W_{q_0}|g|_{q_0},\notag\\
		|Q(h,g)|+|Q(g,h)|&\leq 2C_{q_1}\nu_{\rm ref}W_{q_1}(2\delta_*|g|_{q_1}+|h|_{q_1}|g|_{q_0}).
	\end{align}
	For the iterates in \eqref{eq:B-iteration}, the differentiated equation has the additional loss term
	\[-\eps^{-1}\sum_{0<\mathfrak c\leq\mathfrak a}\binom{\mathfrak a}{\mathfrak c}(D^{\mathfrak c}\nu_{\rm ref})D^{\mathfrak a-\mathfrak c}(h^{n+1}-h^n).\]
	Changing variables $v_*=v+z$ in $\nu_{\rm ref}$ and using the Gaussian derivative bounds gives
	\begin{equation}\label{eq:local-frequency-derivatives}
		D^{\alpha,\beta}\nu_{\rm ref}(x,v)=\int|v-v_*|^\gamma\partial_x^\alpha\partial_{v_*}^\beta M_{\rm ref}(x,v_*)\dd v_*,\qquad |D^{\alpha,\beta}\nu_{\rm ref}|\leq C_N\nu_{\rm ref}.
	\end{equation}
	By \eqref{eq:local-frequency-derivatives}, the loss commutator divided by $W_q$ is bounded by $C_Na_{\rm col}(|h^{n+1}|_{q,N-1}+|h^n|_{q,N-1})$. Set
	\[Z_{N,q,M}(t):=\sum_{|\mathfrak a|\leq N}\max_{0\leq n\leq M}\sup_{s\leq t}|D^{\mathfrak a}h^n(s)|_q .\]
	At $q_0$, apply \eqref{eq:local-highest-rows} and \eqref{eq:local-row-smallness} for each $\mathfrak a$ before summing. In \eqref{eq:B-characteristic}, integrate collision terms against $a_{\rm col}E_{q_0}$ and the transport commutator in \eqref{eq:rwA-derivative-sources} against $E_{q_0}$. Thus
	\begin{equation}\label{eq:B-derivative-gronwall}
		Z_{N,q,M}(t)\leq C_N+\tfrac12Z_{N,q,M}(t)+C_N\int_0^tZ_{N,q,M}(s)\dd s .
	\end{equation}
	Hence $Z_{N,q_0,M}(t)\leq2C_N\exp(2C_Nt)$. At $q_1$, bound $|h^n|_{q_1}|D^{\mathfrak a}h^n|_{q_0}$ by \eqref{eq:local-high-norm-iteration} and this $q_0$ estimate. The remaining coefficient is $\kappa_{\rm it}<1/2$. Thus \eqref{eq:B-derivative-gronwall} holds at $q_1$, with $C_N$ depending on $|h_0|_{q_1,N}$, the bounds through order $N-1$ and reference derivatives through $N+1$.

	For $K_1\Subset K_2\subset\R^6$ and $l<N$, interpolate the difference on $K_1$ using \eqref{eq:B-derivative-gronwall} on $K_2$ and \eqref{eq:rwA-low-difference} to get
	\[\sup_{t\leq T_{\rm loc}}\|h^n(t)-h^{n'}(t)\|_{C^l(K_1)}\leq C\left(\sup_{t\leq T_{\rm loc}}|h^n(t)-h^{n'}(t)|_{q_0}\right)^{1-l/N}\longrightarrow0.\]
	Pass to the differentiated equation using \eqref{eq:B-nonlinear-row}. Since $(|v|+\nu_{\rm ref})W_{q_1}/W_{q_0}\leq C\langle v\rangle^{\gamma-2K}\leq C$, the equation gives
	\begin{equation}\label{eq:local-time-low-norm}
		|\partial_t h(t)|_{q_0}\leq C_{M_{\rm ref},\eps}\bigl(1+|h(t)|_{q_1,1}+|h(t)|_{q_1}^2\bigr).
	\end{equation}

	For $g=h-\widetilde h$, set $\Delta_{N,q}(t)=\sum_{|\mathfrak a|\leq N}\sup_{s\leq t}|D^{\mathfrak a}g(s)|_q$ and $\Delta_{-1,q}=0$. With $G=D^{\mathfrak a}g$ at $|\mathfrak a|=N$, the leading terms satisfy
	\begin{align*}
		\frac{|\mathcal K_{\rm ref}G+Q(G,h)+Q(\widetilde h,G)|}{W_q}
		&\leq\big((\delta_{\rm loc}+4C_{q_1}\delta_*)\nu_{\rm ref}+C_{\rm ref}\big)|G|_q\\
		&\quad+\mathbf1_{\{q=q_1\}}C_{q_1}\nu_{\rm ref}(|h|_{q_1}+|\widetilde h|_{q_1})|G|_{q_0},
	\end{align*}
	by \eqref{eq:B-linear-row} and \eqref{eq:B-nonlinear-row}. Upon division by $W_q$, the other collision terms have bound $C_N\nu_{\rm ref}\Delta_{N-1,q_1}$. By \eqref{eq:B-characteristic}, \eqref{eq:local-iteration-coefficient} and the transport commutator in \eqref{eq:rwA-derivative-sources},
	\begin{align}\label{eq:B-difference-order-norms}
		\tfrac12\Delta_{N,q_0}(t)&\leq C_N\bigl(\Delta_{N,q_1}(0)+\Delta_{N-1,q_1}(t)\bigr)+C_N\int_0^t\Delta_{N,q_0}(s)\dd s,\notag\\
		\tfrac12\Delta_{N,q_1}(t)&\leq C_N\bigl(\Delta_{N,q_1}(0)+\Delta_{N-1,q_1}(t)+\Delta_{N,q_0}(t)\bigr)+C_N\int_0^t\Delta_{N,q_1}(s)\dd s .
	\end{align}
	Starting with $\Delta_{-1,q_1}=0$, apply Gronwall's inequality in \eqref{eq:B-difference-order-norms} in the order
	\[\Delta_{N-1,q_1}(t)\leq C_{N-1}\Delta_{N,q_1}(0)\ \Longrightarrow\
	\Delta_{N,q_0}(t)\leq C_N\Delta_{N,q_1}(0)\ \Longrightarrow\ \Delta_{N,q_1}(t)\leq C_N\Delta_{N,q_1}(0).\]
	This gives
	\begin{equation}\label{eq:local-reference-difference}
		\sup_{0\leq t\leq T_{\rm loc}}|h(t)-\widetilde h(t)|_{q_1,N}\leq C_N|h_0-\widetilde h_0|_{q_1,N}
	\end{equation}
	The constant is uniform for data in bounded subsets of $|\cdot|_{q_1,N}$.

	For $f_\pm=\max\{\pm f,0\}$, approximate the negative part by smooth convex functions in \eqref{Eq.Boltzmann}. The gain terms with two factors of the same sign are nonnegative, so
	\[(\partial_t+v\cdot\nabla_x+\eps^{-1}\nu[f])f_-\leq\eps^{-1}\big(Q^+(f_+,f_-)+Q^+(f_-,f_+)\big).\]
	Use $f_-\leq |h|$ and $f_+\leq M_{\rm ref}+|h|$. By \eqref{eq:rwA-bilinear-rows} and \eqref{eq:reference-frequency-scaling}, $|\nu[h]|\leq C\delta_*\nu_{\rm ref}$. By Lemma~\ref{lem:B-rows} and \eqref{eq:B-transport-split},
	\begin{align*}
		\eps^{-1}\nu[f]+b_{q_0}&\geq (1-\delta_{\rm loc}-C\delta_*)a_{\rm col}-C,\\
		\frac{Q^+(f_+,f_-)+Q^+(f_-,f_+)}{\eps W_{q_0}}
		&\leq\big((\delta_{\rm loc}+C\delta_*)a_{\rm col}+\eps^{-1}C_{\rm ref}\big)|f_-|_{q_0}.
	\end{align*}
	Use \eqref{eq:B-characteristic} with damping $1-\delta_{\rm loc}-C\delta_*>0$. By \eqref{eq:local-row-smallness} and $(f_0)_-=0$,
	\[\sup_{t\leq T_{\rm loc}}|f_-(t)|_{q_0}\leq\tfrac12\sup_{t\leq T_{\rm loc}}|f_-(t)|_{q_0}.\]
	Hence $f_-=0$. If $f_0>0$, integrate $Df+\eps^{-1}\nu[f]f=\eps^{-1}Q^+(f,f)\geq0$ along the same characteristics, using $\nu[f]\leq C\langle v\rangle^\gamma$.
	\[f(t,x,v)\geq f_0(x-tv,v)\exp\left[-\eps^{-1}\int_0^t\nu[f](s,x-(t-s)v,v)\dd s\right]>0.\qedhere\]
\end{proof}

\subsection{Time regularity and spatial decay}

	By differentiating \eqref{Eq.Boltzmann} at $t=0$, we prescribe the initial time derivatives by
\begin{equation}\label{eq:B-time-initial}
	f_{0,0}:=f_0,\qquad f_{0,j+1}:=-v\cdot\nabla_x f_{0,j} +\eps^{-1}\sum_{i=0}^j\binom ji Q(f_{0,i},f_{0,j-i}),\qquad j\geq0.
\end{equation}

\begin{proposition}\label{prop:B-time-regularity}
	Let $f$ be the solution of Proposition~\ref{thm:B-local} and assume $|f_{0,j}|_{q_1,N}<\infty$ for $j\geq1$, $N\geq0$ in \eqref{eq:B-time-initial}. Then $f$ is smooth on $[0,T_{\rm loc}]\times\R^3_x\times\R^3_v$ with one-sided time derivatives at $0,T_{\rm loc}$, and
	\[\sup_{t\leq T_{\rm loc}}|\partial_t^j f(t)|_{q_1,N}<\infty,\qquad \partial_t^j f(0)=f_{0,j}\qquad(j\geq1,\ N\geq0).\]
	Use $T_{\rm loc}$ from Proposition~\ref{thm:B-local}. The initial norms are finite if
	\[\sup_{x,v}\langle v\rangle^L|D^{\mathfrak a}h_0|<\infty\qquad(L\geq0,\ \mathfrak a\in\mathbb N_0^6).\]
\end{proposition}
\begin{proof}
	Suppose $|g(x,v)|+|h(x,v)|\leq C_m\langle v\rangle^{-m}$. Conservation of energy gives
	\begin{align*}
		|v-v_*|^\gamma\big(|g_*'h'|+|g_*h|\big)
		&\leq C_m|v-v_*|^\gamma(1+|v|^2+|v_*|^2)^{-m/2}\\
		&\leq C_m\langle v\rangle^{\gamma-m/2}\langle v_*\rangle^{\gamma-m/2}.
	\end{align*}
	Integrate with $m/2>\gamma+3+L$ and apply \eqref{eq:rw9-collision-derivatives} to get $\sup_{x,v}\langle v\rangle^L|D^{\mathfrak a}Q(g,h)|<\infty$ for rapidly decreasing input derivatives. Induct in \eqref{eq:B-time-initial}, using $W_{q_1}\asymp\langle v\rangle^{-\ell/2-q_1}$.

	Fix $N_t\geq1$. Choose $\chi_{\rm time}\in C^\infty([0,\infty);[0,1])$ equal to one near zero and zero on $[2,\infty)$, and $0<t_{\rm c}<T_{\rm loc}/4$ so that
	\begin{equation}\label{eq:B-time-seed}
	 h^0(t)=h_0+\chi_{\rm time}(t/t_{\rm c})\sum_{j=1}^{N_t}\frac{t^j}{j!}f_{0,j},
	 \qquad \sum_{j=1}^{N_t}\frac{(2t_{\rm c})^j}{j!}|f_{0,j}|_{q_0}\leq\delta_*.
	\end{equation}
	By \eqref{eq:B-time-seed}, $\sup_{t\leq T_{\rm loc}}|h^0|_{q_0}\leq2\delta_*$. Thus \eqref{eq:rwA-low-iteration} and \eqref{eq:rwA-low-difference} still give $h^n\to h$. Differentiating \eqref{eq:B-iteration} for $g_j^n=\partial_t^jh^n$ gives
	\begin{equation}\label{eq:B-time-iterate}
		\begin{aligned}
			(D+a_{\rm col})g_j^{n+1}
			&=\eps^{-1}\mathcal K_{\rm ref}g_j^n+\eps^{-1}\sum_{i=0}^j\binom ji Q(g_i^n,g_{j-i}^n)\\
			&\quad-\mathbf1_{\{j=0\}}v\cdot\nabla_xM_{\rm ref},\qquad 0\leq j\leq N_t .
		\end{aligned}
	\end{equation}
	For $f^n=M_{\rm ref}+h^n$, \eqref{eq:B-iteration} gives
	\begin{equation}\label{eq:physical-local-iteration}
		Df^{n+1}=\eps^{-1}Q(f^n,f^n)-a_{\rm col}(f^{n+1}-f^n).
	\end{equation}
	The seed \eqref{eq:B-time-seed} gives $\partial_t^j f^0(0)=f_{0,j}$ for $0\leq j\leq N_t$. Induct on $n$ and then on $j$, starting with $f^{n+1}(0)=f_0$. When $\partial_t^j(f^{n+1}-f^n)(0)=0$ and $j<N_t$, \eqref{eq:physical-local-iteration} gives
	\[\partial_t^{j+1}f^{n+1}(0)=-v\cdot\nabla_xf_{0,j}+\eps^{-1}\sum_{i=0}^j\binom ji Q(f_{0,i},f_{0,j-i})=f_{0,j+1}.\]
	Thus
	\begin{equation}\label{eq:B-time-iterate-initial}
		g_0^n(0)=h_0,\qquad g_j^n(0)=f_{0,j}\qquad(n\geq0,\ 1\leq j\leq N_t).
	\end{equation}

	For $j=0$ use \eqref{eq:local-high-norm-iteration} and \eqref{eq:B-derivative-gronwall} with the seed norms. Induct lexicographically on $(j,N)$, $N=|\mathfrak a|$. For $j\geq1$, the terms of order $(j,N)$ from $D^{\mathfrak a}$ on the right of \eqref{eq:B-time-iterate} are
	\[\eps^{-1}\big(\mathcal K_{\rm ref}D^{\mathfrak a}g_j^n+Q(h^n,D^{\mathfrak a}g_j^n)+Q(D^{\mathfrak a}g_j^n,h^n)\big).\]
	For $(i,\mathfrak c)\notin\{(0,0),(j,\mathfrak a)\}$ both factors in the Leibniz sum have lower order. In the loss commutator, $|\mathfrak a-\mathfrak c|\leq N-1$. After division by $W_q$, the lower-order collision terms are bounded by $C_{j,N,q}a_{\rm col}$ by \eqref{eq:B-nonlinear-row} and \eqref{eq:local-frequency-derivatives}. The transport commutator is bounded by $N|g_j^{n+1}|_{q,N}$.

	Define $Z_{j,N,q,M}$ as $Z_{N,q,M}$ with $h^n$ replaced by $g_j^n$. At $q_0$, use \eqref{eq:local-highest-rows} and \eqref{eq:B-characteristic} for each $\mathfrak a$ before summing. By \eqref{eq:local-iteration-coefficient},
	\begin{equation}\label{eq:B-mixed-iterate-gronwall}
		Z_{j,N,q,M}(t)\leq C_{j,N,q}+\kappa_{\rm it}Z_{j,N,q,M}(t)+C_{j,N,q}\int_0^t Z_{j,N,q,M}(s)\dd s .
	\end{equation}
	At $q_1$ put $|h^n|_{q_1}|D^{\mathfrak a}g_j^n|_{q_0}$ in $C_{j,N,q_1}$ using \eqref{eq:local-high-norm-iteration} and \eqref{eq:B-mixed-iterate-gronwall} at $q_0$. In both norms $\kappa_{\rm it}<1/2$. Gronwall's inequality and $M\to\infty$ give
	\begin{equation}\label{eq:B-time-uniform-iterates}
		\sup_{n\geq0}\sup_{0\leq t\leq T_{\rm loc}}|\partial_t^j h^n(t)|_{q_1,N}<\infty\qquad(0\leq j\leq N_t,\ N\geq0).
	\end{equation}

	For $N_t\geq j+2$, use \eqref{eq:B-time-uniform-iterates} at $(j+1,N)$ and $(j,N+1)$ to take a diagonal subsequence through $(j,N)$ by Arzel\`a-Ascoli. The zeroth limit is $h$ by \eqref{eq:rwA-low-difference}. Pass to the limit in the spatial integral identities and in
	\[\partial_t^{i-1}D^{\mathfrak a}h^n(t)-\partial_t^{i-1}D^{\mathfrak a}h^n(s)=\int_s^t\partial_t^iD^{\mathfrak a}h^n(\sigma)\dd\sigma\qquad(1\leq i\leq j,\ |\mathfrak a|\leq N).\]
	The limits give the derivatives of $h$, with initial values \eqref{eq:B-time-iterate-initial} and bounds \eqref{eq:B-time-uniform-iterates}.
\end{proof}

	Fix $m_x\geq0$ with $|M_{\rm ref}-M_\infty|_{q_1,N,m_x}+|h_0|_{q_1,N,m_x}<\infty$ for $N\geq0$. Since collisions act at fixed $x$, $L_x=\langle x\rangle^{m_x}$ satisfies
	\[L_x\mathcal K_{\rm ref}g=\mathcal K_{\rm ref}(L_xg),\qquad L_xQ(g,h)=Q(L_xg,h)=Q(g,L_xh).\]
	For $L_xg/W_q$ the transport coefficient is $b_q-v\cdot\nabla_x\log L_x=b_q-m_x(v\cdot x)/\langle x\rangle^2$. Use \eqref{eq:B-transport-split} and $\gamma>1$ to obtain
	\begin{equation}\label{eq:spatial-decay-transport}
		|b_q-v\cdot\nabla_x\log L_x|\leq C_{m_x,{\rm ref}}\langle v\rangle\leq\delta_{\rm loc}a_{\rm col}+C_{m_x,M_{\rm ref},\eps,\delta_{\rm loc}} .
	\end{equation}
	Taking moments of $M_{\rm ref}-M_\infty$ and differentiating \eqref{eq:B-reference} give
	\[L_x|\partial_x^\alpha(A_{\rm ref},u_{\rm ref},\theta_{\rm ref})|\leq C_\alpha\quad(|\alpha|\geq1),\qquad |L_xD^{\mathfrak a}(v\cdot\nabla_x M_{\rm ref})|_q<\infty .\]
	On $[s_0,s_0+h]$ with $\exp(C_{m_x,M_{\rm ref},\eps,\delta_{\rm loc}}h)\leq9/8$, let $Z_{N,q,M}^{(m_x)}$ denote $Z_{N,q,M}$ with $L_xD^{\mathfrak a}h^n$ and $s\in[s_0,t]$. Put $L_x$ on the full derivative in \eqref{eq:local-highest-rows}. By \eqref{eq:spatial-decay-transport} in \eqref{eq:B-characteristic}, first for $q_0$ and then for $q_1$,
	\begin{equation}\label{eq:spatial-decay-derivative-gronwall}
		\tfrac12 Z_{N,q,M}^{(m_x)}(t)\leq C_{N,m_x}\big(1+Z_{N,q,M}^{(m_x)}(s_0)\big)
			+C_{N,m_x}\int_{s_0}^t Z_{N,q,M}^{(m_x)}(s)\dd s .
	\end{equation}
	Iterate the bound on the subintervals and let $M\to\infty$ to get
	\begin{equation}\label{eq:B-spatial-decay-bound}
		\sup_{0\leq t\leq T_{\rm loc}}|h(t)|_{q_1,N,m_x}<\infty,\qquad N\geq0 .
	\end{equation}
	Apply \eqref{eq:B-difference-order-norms} with $L_x(h-\widetilde h)$ to get \eqref{eq:local-reference-difference} in $|\cdot|_{q_1,N,m_x}$ on bounded sets.

	If $|f_{0,j}|_{q_1,N,m_x}<\infty$ for $j\geq1$, $N\geq0$, use \eqref{eq:B-time-seed} and the same estimate \eqref{eq:spatial-decay-derivative-gronwall} in \eqref{eq:B-time-iterate}. Place $L_x$ on the full derivative or a lower-order factor. Then
	\begin{equation}\label{eq:B-spatial-decay-time}
		\sup_{0\leq t\leq T_{\rm loc}}|\partial_t^j f(t)|_{q_1,N,m_x}<\infty\qquad(j\geq1,\ N\geq0).
	\end{equation}
	For \eqref{eq:positive-initial-Maxwellian}, $(f_0,Q(f_0,f_0))=(M_\infty,0)$ on $|x|\geq2R_0$. Thus \eqref{eq:B-time-initial} gives $\supp_x f_{0,j}\subset\overline B_{2R_0}$, $j\geq1$, and $|f_{0,j}|_{q_1,N,m_x}<\infty$.

	Take $m_x>M_x+K+4$. By \eqref{eq:B-spatial-decay-bound}, \eqref{eq:B-spatial-decay-time} and $f-M_\infty=h+(M_{\rm ref}-M_\infty)$,
	\begin{equation}\label{eq:local-spatial-moment-decay}
		\int\langle v\rangle^2\bigl(|\partial_x^\alpha(f-M_\infty)|+|\partial_t\partial_x^\alpha f|\bigr)\dd v\leq C_\alpha\langle x\rangle^{-m_x}\qquad(|\alpha|\geq0).
	\end{equation}
	On compact intervals with $\inf\rho_f,\inf T_f>0$, \eqref{eq:physical-parameters} and \eqref{eq:local-spatial-moment-decay} give
	\[|\partial_t^j\partial_x^\alpha(\rho_f-\theta_\infty^3,u_f,T_f-\theta_\infty^2)|\leq C_\alpha\langle x\rangle^{-m_x}\qquad(j=0,1).\]
	In \eqref{eq:physical-Maxwellian} this gives, for $|\mathfrak a|\leq K$,
	\[|D^{\mathfrak a}(M-M_\infty)|+|\partial_tD^{\mathfrak a}M|\leq C\langle x\rangle^{-m_x}\mathrm e^{-c|v|^2}.\]
	By \eqref{eq:B-spatial-decay-bound}, \eqref{eq:B-spatial-decay-time} and $H=(f-M_\infty)-(M-M_\infty)$,
	\[|D^{\mathfrak a}H|+|\partial_tD^{\mathfrak a}H|\leq C\langle x\rangle^{-m_x}W_{q_1}.\]
	On the compact interval,
	\[\frac{W_{q_1}^2}{W^2}\widehat\nu_M\langle\widetilde v\rangle^{2p_{|\mathfrak a|}}\leq C\langle\widetilde v_{\rm ref}\rangle^{\gamma-2(q_1-2K)}.\]
	Thus the energy and dissipation integrals are dominated by
	\begin{equation}\label{eq:B-energy-integrability}
		C\int_{\R^3}\langle x\rangle^{2M_x-3+2K-2m_x}\dd x\int_{\R^3}\langle\widetilde v_{\rm ref}\rangle^{\gamma-2(q_1-2K)}\dd\widetilde v_{\rm ref}<\infty .
	\end{equation}
	Indeed, the exponents are $2M_x-3+2K-2m_x<-11$ and $\gamma-2(q_1-2K)=\gamma-16<-3$. On $[0,S]$, \eqref{eq:Phi-differentiated} and \eqref{eq:fitted-log-weight} give
	\[(1+|\widetilde v|^2)\partial_a\Phi\leq C_{S,\eps},\quad |\partial_t(|\widetilde v|^2)|\leq C_{S,\eps}\langle\widetilde v\rangle^2\qquad\Longrightarrow\qquad |\partial_t\log W|\leq C_{S,\eps}.\]
	Thus the majorant in \eqref{eq:B-energy-integrability} also controls the first time derivatives of the energy integrands. For nearby data use \eqref{eq:local-reference-difference} with a common reference.

\subsection{Low-norm recovery and continuation}

\begin{lemma}\label{lem:B-pointwise}
	Let $H=f-\mathcal M(P_f)$, $K\geq14$ and $(n,V,\Theta)\in\mathcal U$. If $\mathcal B_{\rm reg}<\infty$ and $\mathcal B\leq e_{\rm vel}$ from Proposition~\ref{prop:mixed-velocity-estimate}, then
	\begin{equation}\label{eq:B-energy-sup}
		\sup_{x,v}\langle\widetilde v\rangle^8|H/W|\leq C\mathcal B_{\rm reg}^{1/2},
	\end{equation}
	with $C$ uniform in $m_\eps$.
\end{lemma}
\begin{proof}
	Fix $(t,x_0,v_0)$ with $t<T$. Put $(d_0,\theta_0)=(d(t,x_0),\theta(t,x_0))$, $Z=\langle\widetilde v(t,x_0,v_0)\rangle$, and
	\[\mathcal O=B_{cd_0Z^{-2}}(x_0)\times B_{c\theta_0 Z^{-1}}(v_0).\]
	From \eqref{eq:coefficient-derivative-assumptions}, $d|\nabla_x\widetilde v|\leq C\langle\widetilde v\rangle$ and $\theta|\nabla_v\widetilde v|\leq C$. For $c$ small,
	\[(d,\theta,J_x,\langle\widetilde v\rangle)\asymp(d_0,\theta_0,J_x(x_0),Z)\qquad\text{on }\mathcal O.\]
	Integrate $d_0|\nabla_x\log W|\leq CZ^2$ and $\theta_0|\nabla_v\log W|\leq CZ$ from \eqref{eq:rw56-pointwise-multiplier} and \eqref{eq:Phi-derivatives}. Then
	\begin{equation}\label{eq:restart-box-weight}
		\left|\log\frac{W(x,v)}{W(x_0,v_0)}\right|
			\leq C\left(Z^2\frac{|x-x_0|}{d_0}+Z\frac{|v-v_0|}{\theta_0}\right)\leq 2Cc,\qquad (x,v)\in\mathcal O .
	\end{equation}
	Using $\Omega_0^*\asymp J_x^2d^{-3}$, $\dd\widehat v=\theta^{-3}\dd v$, we have
	\begin{equation}\label{eq:restart-box-measure}
		\iint_{\mathcal O}\Omega_0^*\theta^{-3}\dd x\dd v\asymp J_x(x_0)^2d_0^{-3}\theta_0^{-3}(d_0Z^{-2})^3(\theta_0 Z^{-1})^3=J_x(x_0)^2Z^{-9}.
	\end{equation}
	For $X=(x-x_0)/(cd_0Z^{-2})$, $Y=(v-v_0)/(c\theta_0 Z^{-1})$, we have for $|\alpha|=i\leq2$, $|\beta|=j\leq2$,
	\[\partial_X^\alpha\partial_Y^\beta H=(cd_0Z^{-2})^i(c\theta_0 Z^{-1})^j\partial_x^\alpha\partial_v^\beta H.\]
	Use the $H^2$ pointwise estimate in each variable. By \eqref{eq:restart-box-weight}, \eqref{eq:restart-box-measure}, \eqref{eq:new-kinetic-energies} and \eqref{eq:mixed-velocity-energy},
	\begin{align*}
		Z^8\frac{|H(x_0,v_0)|}{W(x_0,v_0)}
		&\leq C J_x(x_0)^{-1}\mathcal B_{\rm reg}^{1/2}
			\sum_{0\leq i,j\leq2}Z^{8+9/2-2i-j-p_{i+j}}\\
		&=C J_x(x_0)^{-1}\mathcal B_{\rm reg}^{1/2}
			\sum_{0\leq i,j\leq2}Z^{25/2-2K+j}
		\leq C\mathcal B_{\rm reg}^{1/2}.
	\end{align*}
	Here $25/2-2K+j\leq29/2-2K\leq0$ and $J_x(x_0)^{-1}\leq C$ by \eqref{eq:spatial-weight-J}. Take the supremum over $(x_0,v_0)$.
\end{proof}

\begin{proposition}\label{thm:B-restart}
	Fix $\eps>0$ and $0<S_*\leq S<T$. Let $f$ be the classical solution with datum from Proposition~\ref{prop:positive-realization} on $[0,S_*)$. Assume
	\[
	 (n,V,\Theta)(t,x)\in\mathcal U,\qquad
	 \sup_{t<S_*}\mathcal B_{\rm reg}(t)\leq e_{\rm reg}.
	\]
	In \eqref{eq:roadmap-continuation-threshold} use $C_{\rm rec}$ at least the constant in \eqref{eq:B-energy-sup}, and use $\delta_*$ from Proposition~\ref{thm:B-local} with $(q_0,q_1)=(8,2K+8)$. Then $f$ extends smoothly across $S_*$ with a common restart length for fixed $(S,\eps)$.
\end{proposition}
\begin{proof}
	At $0\leq t_0<S_*$, freeze $M_{\rm ref}=\mathcal M(P_f(t_0))$ and the weight in \eqref{eq:fitted-velocity-weight}, so that
	\[(A_{\rm ref},u_{\rm ref},\theta_{\rm ref},m)=(n\Theta^{-3/2},u+\theta V,\theta\sqrt\Theta,m_\eps)(t_0,x).\]
	By $d\geq(T-S)^{1/r}$, $\theta\geq\theta_*$ and \eqref{eq:background-bounds}, \eqref{eq:roadmap-continuation-macroscopic}, \eqref{eq:adaptive-m}, \eqref{eq:Dm}, we have
	\begin{equation}\label{eq:rwA-restart-coefficients}
		\begin{gathered}
			0<A_{{\rm ref},-}\leq A_{\rm ref}\leq A_{{\rm ref},+},\qquad c_S\leq\theta_{\rm ref}\leq C_{S,\eps},\qquad m_0\leq m\leq C_{S,\eps},\\
			|u_{\rm ref}|+|\nabla_x A_{\rm ref}|+|\nabla_xu_{\rm ref}|+|\nabla_x\theta_{\rm ref}|+|\nabla_xm|\leq C_{S,\eps}.
		\end{gathered}
	\end{equation}
	The amplitude bounds use only $\mathcal U$, so $\delta_*$ is independent of $S$. Equations \eqref{eq:augmented-bootstrap-comparison} and \eqref{eq:roadmap-continuation-threshold} give
	\[\mathcal B(t_0)+\mathcal V_K(t_0)\leq\kappa_{\rm reg}^{-1}\mathcal B_{\rm reg}(t_0)\leq\kappa_{\rm reg}^{-1}e_{\rm reg}\leq e_{\rm vel}/2.\]
	Now $W_{q_0}=W(t_0)\langle\widetilde v(t_0)\rangle^{-8}$. By \eqref{eq:B-energy-sup},
	\begin{equation}\label{eq:rwA-restart-low}
		|f(t_0)-M_{\rm ref}|_{q_0}\leq C_{\rm rec}\mathcal B_{\rm reg}(t_0)^{1/2}\leq C_{\rm rec}\sqrt{e_{\rm reg}}\leq\delta_*/8.
	\end{equation}

	For both the preceding and the new reference, \eqref{eq:Phi-fixed-tail} gives
	\begin{equation}\label{eq:rwA-restart-weight}
		c_{t_0}\langle v\rangle^{-\ell/2-q}\leq W_q(x,v)\leq C_{t_0}\langle v\rangle^{-\ell/2-q}\qquad(q=q_0,q_1).
	\end{equation}
	At $t_0$, \eqref{eq:local-spatial-moment-decay} for the preceding reference bounds the derivatives of the three moments in \eqref{eq:physical-parameters}. In particular,
	\[|\partial_x^\alpha(\rho_f-\theta_\infty^3,u_f,T_f-\theta_\infty^2)(t_0,x)|\leq C_{t_0,\alpha}\langle x\rangle^{-m_x}.\]
	By \eqref{eq:physical-Maxwellian}, $|D^{\mathfrak a}(M_{\rm ref}-M_\infty)|\leq C_{t_0,\mathfrak a}\langle x\rangle^{-m_x}\mathrm e^{-c_{t_0}|v|^2}$. With \eqref{eq:B-spatial-decay-bound} for the preceding reference and \eqref{eq:rwA-restart-weight},
	\[|M_{\rm ref}-M_\infty|_{q_1,N,m_x}+|f(t_0)-M_{\rm ref}|_{q_1,N,m_x}<\infty\qquad(N\geq0).\]
	For $t_0>0$, $f_{0,j}=\partial_t^jf(t_0)$ in \eqref{eq:B-time-initial}. By \eqref{eq:B-spatial-decay-time} and \eqref{eq:rwA-restart-weight},
	\[|f_{0,j}|_{q_1,N,m_x}<\infty\qquad(j\geq1,\ N\geq0).\]
	At $t_0=0$ the same norms follow from \eqref{eq:positive-initial-Maxwellian} and \eqref{eq:B-time-initial}. The bounds \eqref{eq:rwA-restart-coefficients} and \eqref{eq:rwA-restart-low} give a common $h_S>0$ by Proposition~\ref{thm:B-local}. Apply Proposition~\ref{prop:B-time-regularity} on $[t_0,t_0+h_S]$.

	By \eqref{eq:local-time-low-norm} and \eqref{eq:rwA-restart-weight}, the original solution is continuous at $t_0$ in the new frozen $q_0$ norm. By \eqref{eq:rwA-restart-low}, both solutions are in the $2\delta_*$ ball just after $t_0$ and agree by \eqref{eq:rwA-low-difference}. The bound $\frac32\delta_*$ from \eqref{eq:rwA-low-iteration} continues equality through the overlap. For $\max\{0,S_*-h_S\}<t_0<S_*$ we have $t_0+h_S>S_*$.

\end{proof}

\section*{Acknowledgments}
F. Shao is supported by the China National Postdoctoral Program for Innovative Talents under Grant BX2026005. Z. Zhang is partially supported by the National Natural Science Foundation of China under Grant 12288101.

\section*{Disclosure of AI assistance}
The authors began investigating the problem studied in this paper in November 2025. The weighted coercivity estimates for the linearized Boltzmann collision operator were established in February 2026, and the full energy framework was developed in August 2026. In preparing the manuscript, the authors used the AI tools (GPT 5.6-sol and GPT 6-Astra) to assist with language polishing and checking proofs. The authors take full responsibility for the accuracy, originality, and integrity of the entire paper.

\end{document}